\documentclass[letterpaper]{amsart} 
 
\usepackage{amssymb,latexsym} 
\usepackage{psfrag}
\usepackage{epsfig}
 
\numberwithin{equation}{section}

\newtheorem{theorem}{Theorem}[section] 
\newtheorem{proposition}[theorem]{Proposition} 
 
\newtheorem{corollary}[theorem]{Corollary} 
\newtheorem{lemma}[theorem]{Lemma} 
 
\theoremstyle{definition} 
 
\newtheorem{remark}[theorem]{Remark} 
 
\newtheorem{definition}[theorem]{Definition} 
\newtheorem{problem}[theorem]{Problem}

\newcommand{\credit}[1]{\smallskip\noindent {\textbf{#1.\ }}} 
 
\def\proof{\smallskip\noindent {\bf Proof.\ }} 
\def\endproof{\hfill$\square$\medskip} 

\def\FF{\mathbb{F}}

\begin{document}

\title{Recognizing Cluster Algebras of Finite type}

\author{Ahmet I. Seven}

\address{Northeastern University, Boston, MA 02115, USA and Middle East Technical University, Ankara, Turkey}
\email{aseven@math.metu.edu.tr}

\thanks{The author's research was supported in part
by Andrei Zelevinsky's NSF grant \#DMS-0200299.}

\date{November 27, 2005}

\maketitle

\section{Introduction}

\label{sec:introduction}

Cluster algebras were introduced in \cite{CAI} by Fomin and Zelevinsky
to provide an algebraic framework for the study of canonical bases in quantum groups. Since their introduction, it has also been observed that
cluster algebras are closely related with different areas in mathematics. For example;
they provide a natural algebraic set-up to study recursively defined rational functions in combinatorics and number theory \cite{LP}. In geometry, they introduce natural Poisson transformations \cite{GSV1}. 
In representation theory, they form a natural algebraic framework to study positivity \cite{CAIII}.

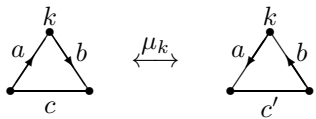
\begin{figure}[ht] 
\setlength{\unitlength}{1.5pt} 
\begin{picture}(30,17)(-5,0) 
\put(0,0){\line(1,0){20}} 
\put(0,0){\line(2,3){10}} 
\put(0,0){\vector(2,3){6}} 
\put(10,15){\line(2,-3){10}} 
\put(10,15){\vector(2,-3){6}} 
\put(0,0){\circle*{2}} 
\put(20,0){\circle*{2}} 
\put(10,15){\circle*{2}} 
\put(2,10){\makebox(0,0){$a$}} 
\put(18,10){\makebox(0,0){$b$}} 
\put(10,-4){\makebox(0,0){$c$}} 
\put(10,19){\makebox(0,0){$k$}} 
\end{picture} 
$ 
\begin{array}{c} 
\stackrel{\textstyle\mu_k}{\longleftrightarrow} 
\\[.3in] 
\end{array} 
$ 
\setlength{\unitlength}{1.5pt} 
\begin{picture}(30,17)(-5,0) 
\put(0,0){\line(1,0){20}} 
\put(0,0){\line(2,3){10}} 
\put(10,15){\vector(-2,-3){6}} 
\put(10,15){\line(2,-3){10}} 
\put(20,0){\vector(-2,3){6}} 
\put(0,0){\circle*{2}} 
\put(20,0){\circle*{2}} 
\put(10,15){\circle*{2}} 
\put(2,10){\makebox(0,0){$a$}} 
\put(18,10){\makebox(0,0){$b$}} 
\put(10,-4){\makebox(0,0){$c'$}} 
\put(10,19){\makebox(0,0){$k$}} 
\end{picture} 
 
\vspace{-.2in} 
\caption{Diagram mutation} 
\label{fig:diagram-mutation-general} 
\end{figure}

One of the most striking results in the theory of cluster algebras
due to S. Fomin and A. Zelevinsky is the classification of
cluster algebras of finite type, which turns out to be identical
to the Cartan-Killing classification \cite{CAI}. This result can be stated purely combinatorially in terms of certain transformations, called  mutations, on certain graphs. To be more precise, let us assume that $\Gamma$ is a finite directed graph whose edges are weighted with positive integers. We call $\Gamma$ a \emph{diagram} if it has the following property: the product of weights along any cycle is a perfect square. 
For any vertex $k$ in $\Gamma$, the \emph{mutation $\mu_k$ in the direction~$k$} is the transformation that changes $\Gamma$ as follows:

\begin{itemize} 
\item The orientations of all edges incident to~$k$ are reversed, 
their weights intact. 
\item 
For any vertices $i$ and $j$ which are connected in 
$\Gamma$ via a two-edge oriented path going through~$k$ (refer to 
Figure~\ref{fig:diagram-mutation-general} for the rest of notation), 
the direction of the edge $(i,j)$ in $\mu_k(\Gamma)$ and its weight $c'$ are uniquely determined by the rule 
\begin{equation} 
\label{eq:weight-relation-general} 
\pm\sqrt {c} \pm\sqrt {c'} = \sqrt {ab} \,, 
\end{equation} 
where the sign before $\sqrt {c}$ 
(resp., before $\sqrt {c'}$) 
is ``$+$'' if $i,j,k$ form an oriented cycle 
in~$\Gamma$ (resp., in~$\mu_k(\Gamma)$), and is ``$-$'' otherwise. 
Here either $c$ or $c'$ can be equal to~$0$, which means that the corresponding edge is absent. 
 
\item 
The rest of the edges and their weights in $\Gamma$ 
remain unchanged. 
\end{itemize} 

\noindent 
It is not hard to show that the resulting weighted graph is a diagram; in particular, its edge-weights are positive integers.
It is also easy to check that $\mu_k$ is involutive, i.e. $\mu_k^{2}(\Gamma)=\Gamma$. 
Two diagrams $\Gamma$ and $\Gamma'$ related by a sequence of diagram mutations are called \emph{mutation equivalent}. A diagram is called $2-$finite if every mutation equivalent diagram has all edge weights equal to 1,2 or 3.
The combinatorial part of the classification theorem in \cite{CAII} is the following: a diagram is 2-finite if and only if it is mutation equivalent to a Dynkin diagram, i.e. a diagram whose underlying undirected graph is a Dynkin graph. However, in \cite{CAII}, an algorithm for checking whether a given diagram is mutation equivalent to a Dynkin diagram is not given. In particular, we do not know how many mutations one needs to perform to show that a given diagram is mutation equivalent to a particular Dynkin diagram, say, $E_8$. This makes the following recognition problem natural:

\begin{problem}{\emph{Recognition Problem for 2-finite diagrams:}}
\label{pr:pb}
How to recognize whether a given diagram $\Gamma$ is $2$-finite without having to perform an unspecified number of mutations.
\end{problem}

In this paper, we solve Problem~\ref{pr:pb} completely by providing the list of all \emph{minimal} 2-infinite diagrams (Section~\ref{subsec:list}).
The list contains all extended Dynkin diagrams but
also has 6 more infinite series, and a substantial number of exceptional diagrams with at most 9 vertices. 
For the proof of this fact, we first show that any diagram in our list is minimal 2-infinite. To prove that any minimal 2-infinite diagram is indeed in our list, we use an inductive argument. 
The basis of the induction is the following fact: our list contains any two-vertex diagram with the edge weight greater than or equal to $4$.
The inductive step is the following statement: if a diagram $\Gamma$ contains a subdiagram that belongs to our list, then, for any vertex $k$ in $\Gamma$, the diagram $\mu_k(\Gamma)$ also contains a subdiagram from our list (Lemmas~\ref{lem:star} and \ref{lem:xxx}). 
Those two properties imply that our list contains all minimal 2-infinite diagrams. To be more precise, let us assume
that $\Gamma$ is a minimal 2-infinite diagram. Then, by definition, there is a sequence of mutations $\mu_r,...,\mu_1$ such that the diagram $\Gamma'=\mu_r\circ...\circ\mu_1(\Gamma)$ contains an edge whose weight is greater than or equal to 4. Here we note that $\Gamma=\mu_1\circ...\circ\mu_r(\Gamma')$ because mutations are involutive. 
Thus, by induction on $k$, the diagram $\Gamma$ contains a subdiagram $\Gamma'$ from our list. Since $\Gamma$ is \emph{minimal} 2-infinite we have $\Gamma=\Gamma'$, showing that our list is a complete list 
of minimal 2-infinite diagrams. We give a more detailed outline of this proof in Section~\ref{sec:main-th}.

We have used some computer assistance to produce our list of diagrams and 
prove that it really is the list of minimal 2-infinite diagrams. 
More specifically, we use 
computer to obtain exceptional minimal 2-infinite diagrams, which are the minimal 2-infinite diagrams
that do not appear in series. We first used a theoretical argument to show that those
exceptional diagrams can have at most 9 vertices, even so it presented a challenge for us to compute them
explicitly because we needed, in one form or another, a fast method to check using a computer 
if a given diagram is 2-finite. In Section 5, we develop such a method for simply-laced diagrams,
here a diagram is called simply-laced if all of its edges have weight equal to $1$.
The basic idea of our method is to 
view the underlying graph of a diagram as an alternating bilinear form on a vector space
over the 2-element field, and describe an arbitrary simply-laced 
2-finite diagram using  
algebraic invariants of the corresponding bilinear form
\footnote{After the first version of this paper appeared, M. Barot, C. Geiss and A. Zelevinsky 
had the paper "Cluster algebras of finite type and positive symmetrizable matrices" 
(to appear in J. London Math. Soc.), where they obtained 
a description of 2-finite diagrams using bilinear forms over integers.}. A nice combinatorial set-up
to carry out this idea is provided by a class of (undirected) graph transformations
called basic moves, which were introduced and studied in \cite{BHII,S}.
A basic move is a simpler operation than a mutation; there is 
also a classification of graphs under basic moves using algebraic
and combinatorial invariants which can be easily implemented \cite{J,S}.
We take advantage of this classification thanks to our following 
characterization: a simply-laced diagram that does not 
contain any non-oriented cycle is 2-finite if and only if its
underlying graph can be obtained
from a Dynkin graph using basic moves (Theorem~\ref{th:mut=BM}).
Using this description, we design and implement the algorithm in 
Section~\ref{subsec:alg}, computing the exceptional minimal
2-infinite diagrams. Our computer program is available at \cite{W}.




In addition to giving an explicit description of minimal 2-infinite diagrams,
we determine representatives for their mutation classes. In particular, we prove that any minimal 2-infinite diagram with at least 5 vertices is mutation equivalent to an extended Dynkin diagram (Theorem~\ref{cor:main}).
We also remark that one can enlarge the set of extended Dynkin diagrams by including some other representatives giving the following "intermediate" recognition criterion: a diagram is 2-infinite if and only, using at most $9$  mutations, it can be transformed into a diagram which contains one of the distinguished representatives (Remark~\ref{th:ps-main}). 

Another long list of directed graphs (quivers) was obtained by Happel and Vossieck
in \cite{HV} to classify finite dimensional algebras
which are of minimal infinite representation type. 
We observed that our list of simply-laced minimal 2-infinite diagrams 
is the same as the list of Happel and Vossieck upto a natural operation 
of replacing the dotted edges indicating relations of quivers in \cite{HV} 
by arrows in the reverse direction. 
This remarkable coincidence of such long lists suggests a close relation between the associated finite dimensional algebras
and cluster algebras, 
which we will explore in a seperate publication. Let us also note that our list 
contains non-simply-laced diagrams while the diagrams in \cite{HV} are 
simply-laced.

The paper is organized as follows. In Section~\ref{sec:def}, we give basic definitions. In Section~\ref{sec:main-th}, we state our main results and outline their proofs. 
In Section~\ref{sec:series}, we prove some statements that allow us to compute series of minimal 2-infinite diagrams. In Section~\ref{sec:min<=9}, we compute exceptional simply-laced minimal 2-infinite diagrams using basic moves. In Sections~\ref{sec:pf-main} and  \ref{subsec:lem-xxx}, we prove our main results. Section~\ref{subsec:list} is our list of minimal 2-infinite diagrams.












\section{Basic Definitions}

\label{sec:def}

In this section, we recall some definitions and statements from \cite{CAI,CAII}. 
We start with the skew-symmetrizability property of an 
integer matrix \cite[Definition~4.4]{CAI}.

\begin{definition} 
\label{def:skew-symmetrizable} 
Let $B$ be a $n \times n$ matrix whose entries are integers. 
The matrix $B$ is called 
skew-symmetrizable if there exists a diagonal matrix $D$ 
with positive diagonal entries such that $DB$ is skew-symmetric. 
\end{definition}

For any skew-symetrizable matrix $B$, Fomin and Zelevinsky introduced a weighted directed graph as follows (\cite[Definition~7.3]{CAII}).

\begin{definition} 
\label{def:diagram-of-B} 
Let $n$ be a positive integer and let $I=\{1,2,...,n\}$. 
The \emph{diagram} of a skew-symmetrizable integer 
matrix~$B=(b_{ij})_{i,j\in I}$ is the weighted directed graph $\Gamma (B)$ 
with the vertex set $I$ such that there is a directed edge from $i$ to $j$ 
if and only if $b_{ij} > 0$, and this edge is assigned the weight 
$|b_{ij}b_{ji}|\,$. 
\end{definition} 
According to \cite[Lemma~7.4]{CAII}; 
if $B$ is a skew-symmetrizable matrix, then, for all $k \geq 3$ and all 
$i_1, \dots, i_k\,$, it satisfies
\begin{equation} 
\label{eq:cycle=cycle}
b_{i_1 i_2} b_{i_2 i_3} \cdots b_{i_k i_1} = 
(-1)^k b_{i_2 i_1} b_{i_3 i_2} \cdots b_{i_1 i_k}\,. 
\end{equation}
In particular, if the edges $e_1,e_2,...,e_r$ with weights
$w_1,w_2,...,w_r$ form an induced cycle (which is not necessarily oriented) in $\Gamma(B)$, 
then the product $w_1w_2...w_r$ is a perfect square. 
Thus we can naturally define a \emph{diagram} as follows:

\begin{definition}
\label{lem:def-diagram}
A diagram $\Gamma$ is a finite directed graph whose edges are weighted with positive integers such that the product of weights along any cycle is a perfect square. 
\end{definition}
\noindent
By some abuse of notation, we denote by the same symbol $\Gamma$ 
the underlying undirected graph of a diagram. If an edge $e=[i,j]$ has weight equal to 1, then we call $e$ \emph{weightless} and do not specify its weight 
in the picture. If all the edges are weightless, then we call $\Gamma$ 
\emph{simply-laced}. 
By a \emph{subdiagram} of $\Gamma$, we always mean a diagram $\Gamma'$ 
obtained from $\Gamma$ by taking an induced directed subgraph on a subset of
vertices and keeping all its edge weights the same as in $\Gamma$ 
\cite[Definition~9.1]{CAII}. We will denote this by $\Gamma' \subset \Gamma$.

For any vertex $k$ in a diagram $\Gamma$, there is the associated mutation $\mu_k$ which changes $\Gamma$ as described in Fig.~\ref{fig:diagram-mutation-general}. This operation naturally defines an equivalence relation on the set of all diagrams. More precisely, two diagrams are called mutation equivalent if they can be obtained from each other by applying a sequence of mutations. 
An important type of diagrams that behave very nicely under mutations are 2-finite diagrams:

\begin{definition} 
\label{def:2-finite} 
A diagram $\Gamma$ is called \emph{$2$-finite} if any diagram $\Gamma'$ 
which is mutation equivalent to $\Gamma$ 
has all edge weights equal to $1,2$ or~$3$. A diagram is called \emph{$2$-infinite}
if it is not 2-finite.
\end{definition}
\noindent
Let us note that a subdiagram of  a $2$-finite diagram is $2$-finite.
We also note that there are only finitely many diagrams which are mutation equivalent to a given 2-finite diagram.

\begin{figure}[ht] 
\vspace{-.2in} 
\[ 
\begin{array}{ccl} 
A_n && 
\setlength{\unitlength}{1.5pt} 
\begin{picture}(140,17)(0,-2) 
\put(0,0){\line(1,0){140}} 
\multiput(0,0)(20,0){8}{\circle*{2}} 
\end{picture}\\
B_n
&& 
\setlength{\unitlength}{1.5pt} 
\begin{picture}(140,17)(0,-2) 
\put(0,0){\line(1,0){140}} 
\multiput(0,0)(20,0){8}{\circle*{2}} 
\put(10,4){\makebox(0,0){$2$}} 
\end{picture} 
\\[.2in] 
D_n 
&& 
\setlength{\unitlength}{1.5pt} 
\begin{picture}(140,17)(0,-2) 
\put(20,0){\line(1,0){120}} 
\put(0,10){\line(2,-1){20}} 
\put(0,-10){\line(2,1){20}} 
\multiput(20,0)(20,0){7}{\circle*{2}} 
\put(0,10){\circle*{2}} 
\put(0,-10){\circle*{2}} 
\end{picture} 
\\[.2in] 
E_6 
&& 
\setlength{\unitlength}{1.5pt} 
\begin{picture}(140,17)(0,-2) 
\put(0,0){\line(1,0){80}} 
\put(40,0){\line(0,-1){20}} 
\put(40,-20){\circle*{2}} 
\multiput(0,0)(20,0){5}{\circle*{2}} 
\end{picture} 
\\[.4in] 
E_7 
&& 
\setlength{\unitlength}{1.5pt} 
\begin{picture}(140,17)(0,-2) 
\put(0,0){\line(1,0){100}} 
\put(40,0){\line(0,-1){20}} 
\put(40,-20){\circle*{2}} 
\multiput(0,0)(20,0){6}{\circle*{2}} 
\end{picture} 
\\[.4in] 
E_8 
&& 
\setlength{\unitlength}{1.5pt} 
\begin{picture}(140,17)(0,-2) 
\put(0,0){\line(1,0){120}} 
\put(40,0){\line(0,-1){20}} 
\put(40,-20){\circle*{2}} 
\multiput(0,0)(20,0){7}{\circle*{2}} 
\end{picture} 
\\[.45in] 
F_4 
&& 
\setlength{\unitlength}{1.5pt} 
\begin{picture}(140,17)(0,-2) 
\put(0,0){\line(1,0){60}} 
\multiput(0,0)(20,0){4}{\circle*{2}} 
\put(30,4){\makebox(0,0){$2$}} 
\end{picture} 
\\[.1in] 
G_2 
&& 
\setlength{\unitlength}{1.5pt} 
\begin{picture}(140,17)(0,-2) 
\put(0,0){\line(1,0){20}} 
\multiput(0,0)(20,0){2}{\circle*{2}} 
\put(10,4){\makebox(0,0){$3$}} 
\end{picture} 
\end{array} 
\] 
\vspace{-.1in} 
\caption{Dynkin diagrams} 
\label{fig:dynkin-diagrams} 
\end{figure}
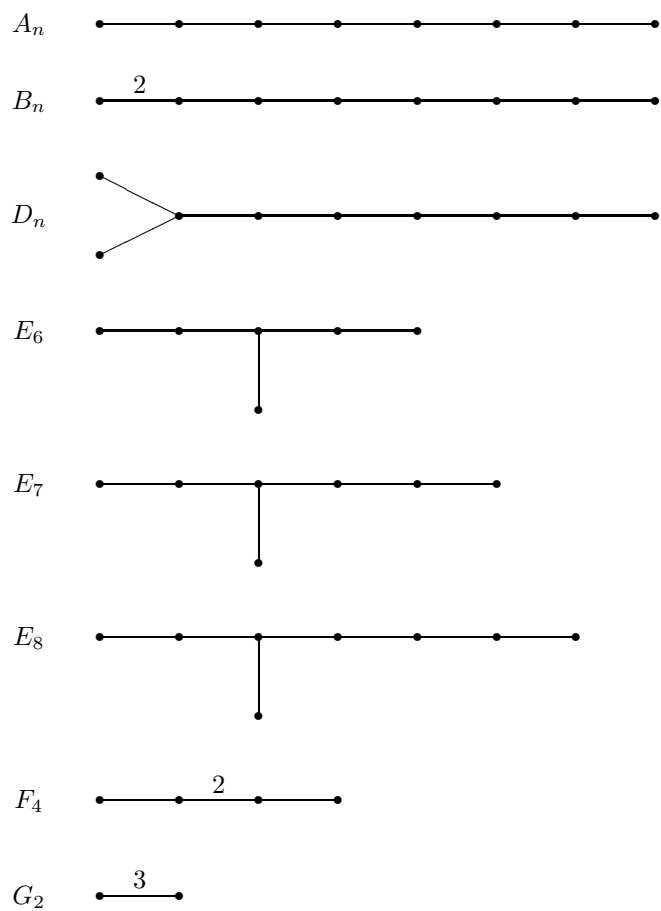

2-finite diagrams were classified by Fomin and Zelevinsky in \cite{CAII}.
Their classification is identical to the Cartan-Killing classification. More precisely:

\begin{theorem}
\label{th:2-finite-class}
A diagram is 2-finite if and only if it is mutation equivalent to an arbitrarily oriented Dynkin diagram (Fig.~\ref{fig:dynkin-diagrams}). 
\end{theorem}

It is a natural problem to give an explicit description of 2-finite diagrams (Problem~\ref{pr:pb}). A conceptual solution to this problem could be obtained by finding  the list of all \emph{minimal} 2-infinite diagrams. More precisely:

\begin{definition} 
\label{def:min-2-infinite} 
A diagram $\Gamma$ is called \emph{minimal $2$-infinite} if 
it is $2$-infinite and any proper subdiagram of 
$\Gamma$ is $2$-finite.
\end{definition}

Clearly one has the following:
\begin{align}
\label{eq:2-infinite} 
&\text{a diagram $\Gamma$ is $2$-infinite if and only if it contains a subdiagram which}
\\
\nonumber
&\text{is minimal 2-infinite.}
\end{align} 

In Section~\ref{subsec:list}, we give a complete list of minimal 2-infinite diagrams, solving Problem~\ref{pr:pb}. 

\pagebreak

\pagebreak

\begin{figure}[ht] 
\[ 
\begin{array}{ccl} 
A_n^{(1)}
&& 
\setlength{\unitlength}{1.5pt} 
\begin{picture}(140,60)(0,-2) 
\put(60,0){\circle*{2.0}} 
\put(60,20){\circle*{2.0}}
\put(40,40){\circle*{2.0}}
\put(20,40){\circle*{2.0}}
\put(0,20){\circle*{2.0}}
\put(20,-20){\circle*{2.0}}
\put(0,0){\circle*{2.0}}
\put(40,-20){\circle*{2.0}}

\put(60,0){\line(-1,-1){20}}
\put(60,0){\line(0,1){20}}
\put(60,20){\line(-1,1){20}}
\put(40,40){\line(-1,0){20}}
\put(20,40){\line(-1,-1){20}}
\put(0,20){\line(0,-1){20}}
\put(0,0){\line(1,-1){20}}
\put(20,-20){\line(1,0){20}}

\put(104,0){\makebox(0,0){non-oriented}}
\end{picture} 
\\[.3in] 
B_n^{(1)} 
&& 
\setlength{\unitlength}{1.5pt} 
\begin{picture}(140,15)(0,-2) 
\put(20,0){\line(1,0){120}} 
\put(0,10){\line(2,-1){20}} 
\put(0,-10){\line(2,1){20}} 
\multiput(20,0)(20,0){7}{\circle*{2}} 
\put(0,10){\circle*{2}} 
\put(0,-10){\circle*{2}} 
\put(130,4){\makebox(0,0){$2$}} 
\end{picture} 
\\[.2in] 
C_n^{(1)}
&& 
\setlength{\unitlength}{1.5pt} 
\begin{picture}(140,17)(0,-2) 
\put(0,0){\line(1,0){140}} 
\multiput(0,0)(20,0){8}{\circle*{2}} 
\put(10,4){\makebox(0,0){$2$}} 
\put(130,4){\makebox(0,0){$2$}} 
\end{picture} 
\\[.1in] 
D_n^{(1)} 
&& 
\setlength{\unitlength}{1.5pt} 
\begin{picture}(140,17)(0,-2) 
\put(20,0){\line(1,0){100}} 
\put(0,10){\line(2,-1){20}} 
\put(0,-10){\line(2,1){20}} 
\put(120,0){\line(2,-1){20}} 
\put(120,0){\line(2,1){20}} 
\multiput(20,0)(20,0){6}{\circle*{2}} 
\put(0,10){\circle*{2}} 
\put(0,-10){\circle*{2}} 
\put(140,10){\circle*{2}} 
\put(140,-10){\circle*{2}} 
\end{picture} 
\\[.2in] 
E_6^{(1)} 
&& 
\setlength{\unitlength}{1.5pt} 
\begin{picture}(140,17)(0,-2) 
\put(0,0){\line(1,0){80}} 
\put(40,0){\line(0,-1){40}} 
\put(40,-20){\circle*{2}} 
\put(40,-40){\circle*{2}} 
\multiput(0,0)(20,0){5}{\circle*{2}} 
\end{picture} 
\\[.7in] 
E_7^{(1)} 
&& 
\setlength{\unitlength}{1.5pt} 
\begin{picture}(140,17)(0,-2) 
\put(0,0){\line(1,0){120}} 
\put(60,0){\line(0,-1){20}} 
\put(60,-20){\circle*{2}} 
\multiput(0,0)(20,0){7}{\circle*{2}} 
\end{picture} 
\\[.25in] 
E_8^{(1)} 
&& 
\setlength{\unitlength}{1.5pt} 
\begin{picture}(140,17)(0,-2) 
\put(0,0){\line(1,0){140}} 
\put(40,0){\line(0,-1){20}} 
\put(40,-20){\circle*{2}} 
\multiput(0,0)(20,0){8}{\circle*{2}} 
\end{picture} 
\\[.3in] 
F_4^{(1)} 
&& 
\setlength{\unitlength}{1.5pt} 
\begin{picture}(140,17)(0,-2) 
\put(0,0){\line(1,0){80}} 
\multiput(0,0)(20,0){5}{\circle*{2}} 
\put(30,4){\makebox(0,0){$2$}} 
\end{picture} 
\\[.1in] 
G_2^{(1)}(a) 
&& 
\setlength{\unitlength}{1.5pt} 
\begin{picture}(140,17)(0,-2) 
\put(0,0){\line(1,0){40}} 
\multiput(0,0)(20,0){3}{\circle*{2}} 
\put(10,4){\makebox(0,0){$3$}} 
\put(30,4){\makebox(0,0){$a$}} 
\put(70,0){\makebox(0,0){$a=1,2,3$}} 
\end{picture}\\[.1in] 
I_2(a)
&& 
\setlength{\unitlength}{1.5pt} 
\begin{picture}(65,17)(20,-2) 
\put(20,0){\line(1,0){30}} 
\put(20,0){\circle*{2}} 
\put(50,0){\circle*{2}} 
\put(35,4){\makebox(0,0){$a$}} 
\put(80,0){\makebox(0,0){$a\geq 4$}} 
\end{picture} 
\end{array} 
\] 
\caption{Extended Dynkin diagrams} 
\label{fig:extended-dynkin-diagrams} 
\end{figure}
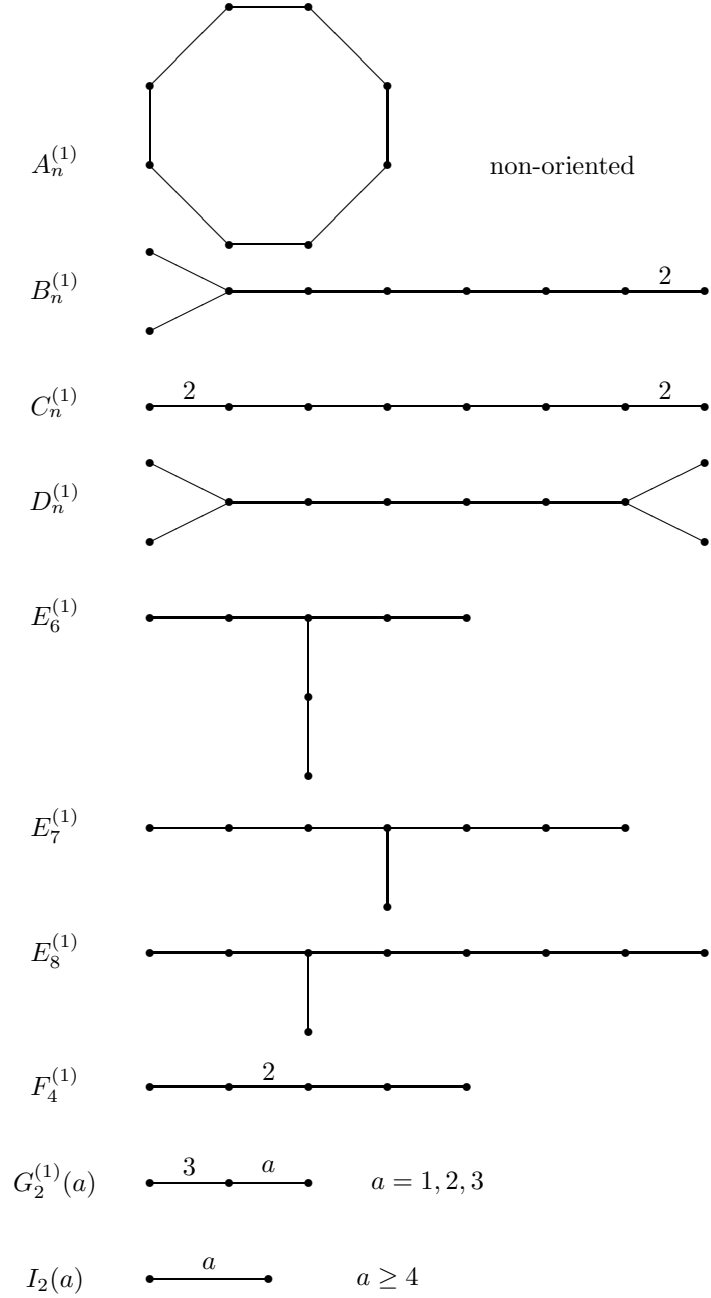 



\section{Main Results}
\label{sec:main-th}

\emph{
Throughout the paper, we assume that all diagrams are connected. We also assume, unless otherwise stated, that any diagram has an arbitrary orientation which does not contain any non-oriented cycle.}

\vspace{.1in}
Our main result is the following statement:

\begin{theorem}
\label{th:main}

The list of minimal 2-infinite diagrams consists precisely of the diagrams given in Section~\ref{subsec:list}.


\end{theorem}

We also determine representatives for mutation classes of minimal 2-infinite diagrams as follows:

\begin{theorem}
\label{cor:main}

Any minimal 2-infinite diagram is either one of the diagrams in 
in Table 2 (Section~\ref{subsec:list}) or it is
mutation equivalent to an extended Dynkin diagram (Fig.~\ref{fig:extended-dynkin-diagrams}). 

\end{theorem}

We prove Theorem~\ref{th:main} as follows. We first show that any diagram in our list is minimal
2-infinite:

\begin{lemma}
\label{lem:List-min-2inf}
Any diagram $\Gamma$ in Section~\ref{subsec:list} is minimal 2-infinite.
\end{lemma}
\noindent
Next, to complete the proof of Theorem~\ref{th:main}, we show that any minimal
2-infinite diagram is indeed one of the diagrams in Section~\ref{subsec:list}. For this,
we recall that any 2-infinite diagram, in particular any minimal one, is mutation 
equivalent to a diagram which contains a subdiagram of the form $I_2(a), a\geq4$: 

\setlength{\unitlength}{2pt}   
\begin{picture}(30,10)(0,-2) 
\setlength{\unitlength}{1.5pt}
\put(80,0){\circle*{2.0}}
\put(110,0){\circle*{2.0}}
\put(80,0){\line(1,0){30}} 
\put(95,4){\makebox(0,0){$a$}} 
\put(125,0){\makebox(0,0){$a\geq 4$}} 
\end{picture}

\vspace{.1in} 
\noindent
To be more precise, let us assume that $X'$ is a 2-infinite diagram and 
$\mu_r,...,\mu_1$ a sequence of mutations such that the diagram 
$X=\mu_r\circ...\circ\mu_1(X')$ 
contains a subdiagram of the form $I_2(a), a\geq4$. Here we note that 
$X'=\mu_1\circ...\circ\mu_r(X)$ because mutations are involutive. 
We prove, by induction on $r$, that $X'$ contains a subdiagram from our list,
so if $X'$ is \emph{minimal} 2-infinite then this subdiagram must be 
$X'$ itself (because any diagram in our list is 2-infinite), 
proving Theorem~\ref{th:main}. The basis of the induction is the fact 
that any diagram of the form $I_2(a), a\geq4$ is included in 
Section~\ref{subsec:list} (Table 1). The inductive step is the following
statement:
\begin{align} 
\label{eq:star2}
&\text{
If a diagram $X$ contains a subdiagram $\Gamma$ from Section~\ref{subsec:list}, 
then, for any vertex $k$ in $X$,}
\\
\nonumber
&\text{
the diagram $\mu_k(X)$ also contains a subdiagram from Section~\ref{subsec:list}.}
\end{align} 

\noindent
To establish \eqref{eq:star2}, we consider it in two possible cases: 
the vertex $k$ being contained in $\Gamma$ or not. If $k$ is a vertex of
$\Gamma$, we show that $\mu_k(\Gamma)$
contains a subdiagram from our list. If $k$ is not in $\Gamma$, we denote
by $\Gamma k$ the minimal subdiagram of $X$ that contains $\Gamma$ and $k$,
and show that $\mu_k(\Gamma k)$ contains a subdiagram from our list.
For this we will assume, without loss of generality, that $k$ is connected
to (at least two vertices in) $\Gamma$ because otherwise $\mu_k$ does not effect
$\Gamma$ (Fig~\ref{fig:diagram-mutation-general}). 
In short, we prove the following two statements to show that \eqref{eq:star2} is satisfied:

\begin{lemma}
\label{lem:star}

Let $\Gamma$ be an arbitrary diagram in Section~\ref{subsec:list}.
If $k$ is a vertex in $\Gamma$, then $\mu_k(\Gamma)$ contains a subdiagram 
$\Gamma'$ which is in Section~\ref{subsec:list}. 
Furthermore, if $\Gamma$ is in Table 1, then $\Gamma'$ can be chosen from Table 1. 
\end{lemma}

\begin{lemma}
\label{lem:xxx}

Suppose that $\Gamma$ is an arbitrary diagram in Section~\ref{subsec:list}.
Let $\Gamma k$ be a diagram obtained from $\Gamma$ by adjoining a vertex $k$. 
Then $\mu_k(\Gamma k)$ contains a subdiagram 
$\Gamma'$ which is in Section~\ref{subsec:list}. 
\end{lemma}


Let us note that Lemmas~\ref{lem:List-min-2inf}, \ref{lem:star}, \ref{lem:xxx} 
prove Theorem~\ref{th:main}. We prove those lemmas
in Sections~\ref{subsec:List-in-M},~\ref{subsec:lem-star} and ~\ref{subsec:lem-xxx}.
Here we discuss our proof of Lemma~\ref{lem:xxx}, which is more involved. To prove this 
lemma we assume, without loss of generality, that $\Gamma k$ does not have 
any subdiagram which contains
$k$ and belongs to our list (otherwise Lemma~\ref{lem:star} applies). This assumption 
greatly restricts possible non-simply-laced $\Gamma k$ and 
we manage to obtain the lemma for such $\Gamma k$ using a case-by-case 
analysis. To treat simply-laced $\Gamma k$, it turns out to be convenient for us
to consider them in two classes: those that do not contain any subdiagram which is mutation
equivalent to the Dynkin diagram $E_6$ and those that do. If a simply-laced $\Gamma k$
belongs to the first class, then $\Gamma$ is in Table 1 (because any simply-laced diagram
in other tables contains a subdiagram which is mutation equivalent to $E_6$). For such
$\Gamma k$, we obtain the lemma from the following stronger statement:

\begin{proposition}
\label{pr:main-add-seq}
Suppose that $\Gamma$ is a simply-laced diagram in Table 1 (Section~\ref{subsec:list})
, i.e. $\Gamma$ is one of the following diagrams: $A_n^{(1)}$, $D_n^{(1)}$, $D_n^{(1)}(m,r)$, $D_n^{(1)}(r)$, $D_n^{(1)}(m,r,s)$.
Let $\Gamma k$ be a simply-laced diagram obtained from $\Gamma$ by adjoining a vertex $k$. 
Suppose that $k$ is connected to at least two vertices in $\Gamma$. 
Suppose also that
\begin{align} 
\label{eq:knotinE6}&\text{the vertex $k$ is not contained in any subdiagram $E\subset \Gamma k$
such that $E$ is}
\\
\nonumber
&\text{mutation equivalent to $E_6$}.
\end{align} 
Then (precisely) one of the following holds:
\begin{align} 
\label{eq:knotinM-seq}&\text{$k$ is contained in a diagram $\Gamma''\subset \Gamma k$ such that $\Gamma''$ is in Table 1,}
\\
\label{eq:kct2-seq}&\text{the diagram $\mu_k(\Gamma k)$ is in Table 1.}
\\
\nonumber
\end{align} 
\end{proposition}

\noindent
Let us note that if \eqref{eq:knotinM-seq} is satisfied, then Lemma~\ref{lem:star} 
applies, giving the same conclusion as Lemma~\ref{lem:xxx}. Now
to complete the proof of Lemma~\ref{lem:xxx}, we need to establish it
for $\Gamma k$ that contains a subdiagram, say $E$, which is mutation
equivalent to $E_6$. 
For this we first show that $\mu_k(\Gamma k)$
contains a minimal 2-infinite diagram which has at most $9$ vertices, then
we show that any such minimal 2-infinite diagram is contained in our
list. The first part is obvious if $\Gamma k$, thus $\mu_k(\Gamma k)$,
has at most $9$ vertices (recall that $\Gamma k$ is 2-infinite, 
so $\mu_k(\Gamma k)$ is also 2-infinite thus contains
a minimal 2-infinite diagram). For larger $\Gamma k$, we first observe 
the following fact in Corollary~\ref{lem:E6-in-2-infinite}: any (simply-laced) diagram which 
has at least 
$9$ vertices and contains a subdiagram mutation equivalent to $E_6$ is 2-infinite.
To use this fact in our set-up, we also observe that if $\Gamma k$ has
at least $10$ vertices then there exists a connected subdiagram $Xk\subset \Gamma k$ of
9 vertices which contains both $E$ and $k$ (Section~\ref{subsec:lem-xxx-1}). Then $Xk$ must be 2-infinite by the 
mentioned fact, so it contains a minimal 2-infinite subdiagram, say $M$, which
has at most $9$ vertices. We note
that $M$ contains $k$ because any subdiagram of $Xk$ that does not
contain $k$ is a proper subdiagram of $\Gamma$, 
and any proper subdiagram of $\Gamma$ 
is 2-finite because $\Gamma$ is minimal 2-infinite. Let us also note that 
$\mu_k(M)$, which is a subdiagram of $\mu_k(\Gamma k)$, is also 2-infinite,
so it contains a minimal 2-infinite diagram which has at most 9 vertices.
Thus, to complete the proof of 
Lemma~\ref{lem:xxx}, it is enough to show that any minimal 2-infinite diagram $M$ with 
at most 9 vertices is contained in Section~\ref{subsec:list}:

\begin{proposition}
\label{pr:simply-laced-9}
Any simply-laced minimal 2-infinite diagram which has at most 9 vertices is contained in Section~\ref{subsec:list}.
\end{proposition}

We prove the proposition using some theory that we develop in Section~\ref{sec:min<=9} along with some
computer assistance. To motivate for our proof, let us first note that we may
obtain all simply-laced minimal 2-infinite diagrams (with at most 9 vertices) 
as follows: first we compute all simply-laced 2-finite diagrams (with at most 8 vertices) 
mutating the corresponding Dynkin diagrams,
then extend any 2-finite diagram by connecting, in all possible ways, 
an additional vertex such that
the resulting diagram is 2-infinite and any proper subdiagram is 2-finite. 
To implement the second step of
this algorithm, we need an efficient method to check, possibly using a computer, 
if a given simply-laced diagram is 2-finite. Our basic idea to develop such a method is 
to view the underlying graph of a diagram as an alternating bilinear form on a vector space
over the 2-element field, and characterize an arbitrary (simply-laced) 
2-finite diagram using  
algebraic invariants of the corresponding bilinear form. A nice combinatorial set-up
to carry out this idea is provided by a class of (undirected) graph transformations
called basic moves, which were introduced and studied in \cite{BHII,S}.
A basic move
changes a graph in a way similar to a mutation does modulo 2: 
it introduces or deletes edges containing a fixed vertex connected to
a given vertex (thus a basic move is assigned to a pair of vertices
connected to each other, for a precise description see Definition~\ref{def:BM}). 
We note in Proposition~\ref{pr:mut-BM} that the underlying graphs of mutation-equivalent
simply-laced diagrams can be obtained from each other 
by a sequence of basic moves. 
We prove the converse of this statement
for 2-finite diagrams: any simply-laced diagram that does not 
contain any non-oriented cycle is 2-finite if and only if its
underlying graph can be obtained
from a Dynkin graph using basic moves (Theorem~\ref{th:mut=BM}).
The advantage of characterizing 2-finite diagrams using basic moves
is that a basic move is a simpler operation than a mutation; there is 
also a classification of graphs under basic moves using algebraic
and combinatorial invariants which can be easily implemented \cite{J,S}.
In Proposition~\ref{pr:BM-Arf}, we give such a characterization for graphs that can be 
obtained from Dynkin graphs with 6,7 or 8 vertices using basic moves.
Using this description, we design and implement the algorithm in 
Section~\ref{subsec:alg},  
obtaining all simply-laced minimal 2-infinite diagrams that contain a subdiagram
which is mutation equivalent to $E_6$ (our computer program is
available at \cite{W}). For the remaining simply-laced 
minimal 2-infinite  diagrams, we prove 
that they must belong to Table 1 
(Corollary~\ref{cor:BM-D}), completing the proof of Proposition~\ref{pr:simply-laced-9}.

We will prove Theorem~\ref{th:main} in Section~\ref{subsec:lem-xxx}. We prove Theorem~\ref{cor:main} in Section~\ref{subsec:cor-main}.

\section{Series of minimal 2-infinite diagrams: Proof of Proposition~\ref{pr:main-add-seq}}

\label{sec:series}

To prove Proposition~\ref{pr:main-add-seq}, it will be convenient for us to prove first a slightly stronger
statement for the diagram $A_n^{(1)}$:

\begin{proposition}
\label{pr:main-add-seq-A}
In the situation of Proposition~\ref{pr:main-add-seq}, if $\Gamma$ is of type $A_n^{(1)}$,
then $\mu_k(\Gamma k)$ is 
one of the following diagrams: $A_n^{(1)}$, $D_n^{(1)}$, $D_n^{(1)}(m,r)$, $D_n^{(1)}(r)$, $D_n^{(1)}(m,r,s)$.

\end{proposition}

\subsection{Proof of Proposition~\ref{pr:main-add-seq-A}}

\label{subsec:series-1}

Let us index the vertices in $\Gamma$ by $\{1,...,n\}$. 
Let us also write 
$$\{i\in \Gamma: k \mathrm{\:is \: connected \: to \:} i\}=\{i_1,...,i_r\}$$
where $1 \leq i_1  < i_2 <...< i_r \leq n$ and  $r\geq 2$.
Since $k$ is not contained in any non-oriented cycle in $\Gamma k$, the number $r$ is even. 

\begin{figure}[ht]

\setlength{\unitlength}{1.8pt}

\begin{center}

\begin{picture}(120,80)(-110,-30)

\put(-50,10){\line(-3,1){30}}
\put(-50,10){\line(-1,3){10}}
\put(-50,10){\line(1,3){10}}
\put(-50,10){\line(3,1){30}}
\put(-50,10){\line(3,-1){30}}
\put(-50,10){\line(1,-3){10}}
\put(-50,10){\line(-1,-3){10}}
\put(-50,10){\line(-3,-1){30}}
\put(-50,10){\circle*{2.0}}
\put(-47,13){\makebox(0,0){$k$}}

\thicklines



\put(-20,20){\circle*{2.0}}
\put(-20,0){\circle*{2.0}} 
\put(-20,20){\circle*{2.0}}
\put(-40,40){\circle*{2.0}}
\put(-60,40){\circle*{2.0}}
\put(-80,20){\circle*{2.0}}
\put(-60,-20){\circle*{2.0}}
\put(-80,0){\circle*{2.0}}
\put(-40,-20){\circle*{2.0}}

\put(-20,0){\line(-1,-1){20}}
\put(-20,0){\line(0,1){20}}
\put(-20,20){\line(-1,1){20}}
\put(-40,40){\line(-1,0){20}}
\put(-60,40){\line(-1,-1){20}}
\put(-80,20){\line(0,-1){20}}
\put(-80,0){\line(1,-1){20}}
\put(-60,-20){\line(1,0){20}}

\put(-19,-4){\makebox(0,0){$i_6$}}
\put(-60,43){\makebox(0,0){$i_3$}}
\put(-40,43){\makebox(0,0){$i_4$}}
\put(-83,0){\makebox(0,0){$i_1$}}
\put(-83,20){\makebox(0,0){$i_2$}}
\put(-40,-23){\makebox(0,0){$i_7$}}
\put(-60,-23){\makebox(0,0){$i_8$}}
\put(-19,23){\makebox(0,0){$i_5$}}

\end{picture}
\end{center}

\caption{The subdiagram on vertices $\{i_1,i_2,i_4,i_5,i_7,k\}$ is mutation equivalent to $E_6$} 

\label{fig:Ahat-0}

\end{figure}
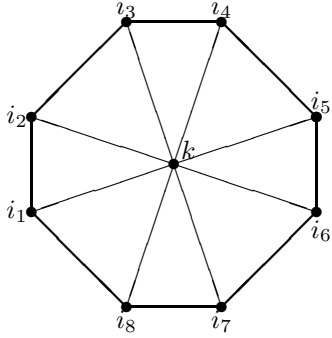

We prove the lemma using a case by case analysis as follows:

\credit{Case 1} 
\emph{$r\geq 8$.}
In this case the subdiagram with the vertices $\{i_1,i_1+1,i_4,i_4+1,i_7,k\}$ is always mutation equivalent to $E_6$ (Fig.~\ref{fig:Ahat-0})
, contradicting \eqref{eq:knotinE6}.

\credit{Case 2} 
\emph{r=6}.

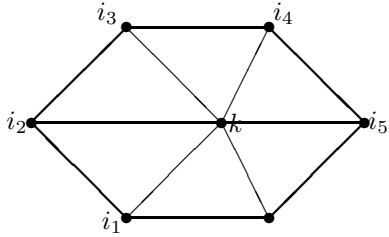
\begin{figure}[ht]

\setlength{\unitlength}{1.8pt}

\begin{center}


\begin{picture}(70,50)(5,-25)

\put(30,20){\line(1,-1){20}}
\put(50,0){\line(-1,-1){20}}
\put(53,0){\makebox(0,0){$k$}}
\put(50,0){\circle*{2.0}}
\put(50,0){\line(-1,0){40}}
\put(50,0){\line(1,0){30}}
\put(50,0){\line(1,2){10}}
\put(50,0){\line(1,-2){10}}

\thicklines

\put(60,20){\circle*{2.0}}
\put(30,20){\circle*{2.0}}
\put(30,-20){\circle*{2.0}}
\put(60,-20){\circle*{2.0}}
\put(80,0){\circle*{2.0}}
\put(10,0){\circle*{2.0}}

\put(27,-21){\makebox(0,0){$i_1$}}
\put(26,23){\makebox(0,0){$i_3$}}
\put(7,0){\makebox(0,0){$i_2$}}
\put(63,23){\makebox(0,0){$i_4$}}
\put(83,0){\makebox(0,0){$i_5$}}

\put(10,0){\line(1,1){20}}
\put(30,-20){\line(-1,1){20}}
\put(30,20){\line(1,0){30}}
\put(60,20){\line(1,-1){20}}
\put(80,0){\line(-1,-1){20}}
\put(60,-20){\line(-1,0){30}}

\end{picture}

\end{center}

\caption{The subdiagram on vertices $\{i_1,i_2,...,i_5,k\}$ is mutation equivalent to $E_6$} 

\label{fig:Ahat-11}

\end{figure}

\credit{Subcase 2.1} 
\emph{$\Gamma$ has length 6.}
Then the subdiagram with the vertices $\{i_1,i_2,...,i_5,k\}$ is 
mutation equivalent to $E_6$ (Fig.~\ref{fig:Ahat-11})
, contradicting \eqref{eq:knotinE6}.

\credit{Subcase 2.2} 
\emph{$\Gamma$ has length greater than 6.}
Let us note that $k$ is contained in a cycle $C\subset \Gamma k$ of length greater than 3.

\begin{figure}[ht]

\setlength{\unitlength}{1.8pt}

\begin{center}

\begin{picture}(70,50)(5,-25)

\put(53,0){\makebox(0,0){$k$}}
\put(50,0){\circle*{2.0}}
\put(50,0){\line(-1,0){40}}
\put(50,0){\line(1,2){10}}
\put(50,0){\line(1,-2){10}}
\put(80,0){\line(-1,-1){10}}
\put(60,-20){\line(-1,0){10}}
\put(63,-23){\makebox(0,0){$i_5$}}

\thicklines

\put(60,20){\circle*{2.0}}
\put(30,20){\circle*{2.0}}
\put(30,-20){\circle*{2.0}}
\put(60,-20){\circle*{2.0}}
\put(80,0){\circle*{2.0}}
\put(10,0){\circle*{2.0}}

\put(27,-21){\makebox(0,0){$n$}}
\put(26,23){\makebox(0,0){$2$}}
\put(7,0){\makebox(0,0){$1$}}
\put(63,23){\makebox(0,0){$3$}}
\put(83,0){\makebox(0,0){$4$}}

\put(10,0){\line(1,1){20}}
\put(30,-20){\line(-1,1){20}}
\put(30,20){\line(1,0){30}}
\put(60,20){\line(1,-1){20}}

\end{picture}

\end{center}

\caption{The subdiagram on vertices $\{k,1,2,3,n,i_5\}$ is mutation equivalent to $E_6$} 

\label{fig:Ahat-121}

\end{figure}
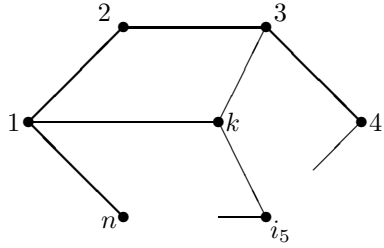

\credit{Subsubcase 2.2.1}
\emph{$k$ is contained in a cycle $C\subset \Gamma k$ of length equal to 4}. 
Let us assume, without loss of generality, that $C=[k,i_1=1,2,i_2=3]$
and $k$ is not connected to the vertex $2$. Then the subdiagram 
with the vertices $C \cup \{n,i_5\}$ in $\Gamma k$ is mutation equivalent to $E_6$ (Fig.~\ref{fig:Ahat-121}), 
contradicting \eqref{eq:knotinE6}.

\begin{figure}[ht]

\setlength{\unitlength}{1.8pt}

\begin{center}

\begin{picture}(70,50)(5,-25)
\put(47,-3){\makebox(0,0){$k$}}
\put(50,0){\circle*{2.0}}
\put(50,0){\line(-1,0){40}}
\put(50,0){\line(1,0){30}}
\put(50,0){\line(1,2){10}}
\put(50,0){\line(1,-2){10}}
\put(80,0){\line(-1,-1){10}}
\put(60,-20){\line(-1,0){10}}
\put(63,-23){\makebox(0,0){$i_4$}}
\put(60,20){\line(1,-1){10}}

\thicklines

\put(60,20){\circle*{2.0}}
\put(30,20){\circle*{2.0}}
\put(30,-20){\circle*{2.0}}
\put(60,-20){\circle*{2.0}}
\put(80,0){\circle*{2.0}}
\put(10,0){\circle*{2.0}}
\put(20,10){\circle*{2.0}}
\put(45,20){\circle*{2.0}}

\put(27,-21){\makebox(0,0){$n$}}
\put(7,0){\makebox(0,0){$1$}}
\put(48,23){\makebox(0,0){$p$}}
\put(63,23){\makebox(0,0){$i_2$}}
\put(83,0){\makebox(0,0){$i_3$}}

\put(10,0){\line(1,1){20}}
\put(30,-20){\line(-1,1){20}}
\put(30,20){\line(1,0){30}}

\end{picture}

\end{center}

\caption{The diagram on vertices $\{k,1,2,p,i_2,i_4\}$, $p \geq 3$, contains (a subdiagram mutation equivalent to) $E_6$} 

\label{fig:Ahat-122}

\end{figure}
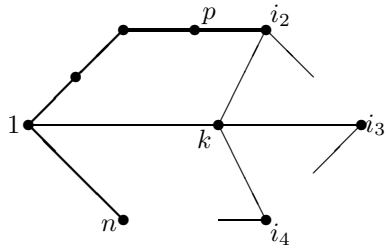

\credit{Subsubcase 2.2.2}
\emph{$k$ is contained in a cycle $C\subset \Gamma k$ of length greater than 4}.
Let us assume, without loss of generality, that $C=[k,i_1=1,2,...,r,i_2]$
where $p \geq 3$ and $k$ is not connected to any vertex in $\{2,...,p\}$.
Then the subdiagram with the vertices $\{k,1,2,p,i_2,i_4\}$is mutation equivalent to $E_6$ 
(Fig.~\ref{fig:Ahat-122}), 
contradicting \eqref{eq:knotinE6}.

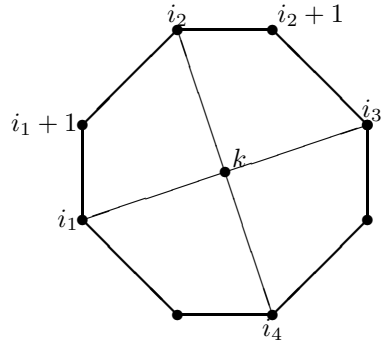
\begin{figure}[ht]

\setlength{\unitlength}{1.8pt}

\begin{center}

\begin{picture}(120,80)(-110,-30)

\put(-50,10){\line(-1,3){10}}
\put(-50,10){\line(3,1){30}}
\put(-50,10){\line(1,-3){10}}
\put(-50,10){\line(-3,-1){30}}
\put(-50,10){\circle*{2.0}}
\put(-47,13){\makebox(0,0){$k$}}

\thicklines



\put(-20,20){\circle*{2.0}}
\put(-20,0){\circle*{2.0}} 
\put(-20,20){\circle*{2.0}}
\put(-40,40){\circle*{2.0}}
\put(-60,40){\circle*{2.0}}
\put(-80,20){\circle*{2.0}}
\put(-60,-20){\circle*{2.0}}
\put(-80,0){\circle*{2.0}}
\put(-40,-20){\circle*{2.0}}

\put(-20,0){\line(-1,-1){20}}
\put(-20,0){\line(0,1){20}}
\put(-20,20){\line(-1,1){20}}
\put(-40,40){\line(-1,0){20}}
\put(-60,40){\line(-1,-1){20}}
\put(-80,20){\line(0,-1){20}}
\put(-80,0){\line(1,-1){20}}
\put(-60,-20){\line(1,0){20}}

\put(-60,43){\makebox(0,0){$i_2$}}
\put(-32,43){\makebox(0,0){$i_2+1$}}
\put(-83,0){\makebox(0,0){$i_1$}}
\put(-88,20){\makebox(0,0){$i_1+1$}}
\put(-40,-23){\makebox(0,0){$i_4$}}
\put(-19,23){\makebox(0,0){$i_3$}}

\end{picture}
\end{center}

\caption{The subdiagram on vertices $\{k,i_1,i_1+1,i_2,i_2+1,i_4\}$ is mutation equivalent to $E_6$} 

\label{fig:Ahat-2}

\end{figure}

\credit{Case 3} 
\emph{r=4}.
Let us first assume that the diagram $\Gamma k$ does not contain a triangle.
Then, e.g., the subdiagram with the vertices $\{k,i_1,i_1+1,i_2,i_2+1,i_4\}$  is always mutation equivalent to $E_6$ (Fig.~\ref{fig:Ahat-2}), contradicting  \eqref{eq:knotinE6}. Let us now assume that $\Gamma k$ contains at least one 
triangle and consider the subcases.

\credit{Subcase 3.1} \emph{The diagram $\Gamma k$ contains precisely one triangle}.

\credit{Subsubcase 3.1.1} \emph{The vertex $k$ is contained in a cycle $C'\subset \Gamma k$ of length greater than 4.}
Let us assume without loss of generality that $C'=[k,i_1=1,2,...,i_2=r]$. Since $\Gamma k$
contains precisely one triangle, we may assume that $i_3$ is not connected to $i_2$.
Then the subdiagram with the vertices $C'\cup \{i_3\}$ is mutation equivalent to $E_6$, contradicting \eqref{eq:knotinE6}.

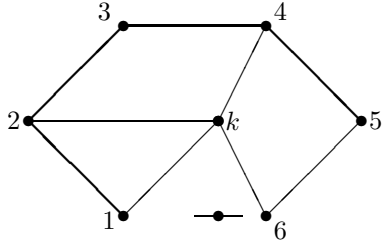
\begin{figure}[ht]

\setlength{\unitlength}{1.8pt}

\begin{center}

\begin{picture}(70,50)(5,-25)
\put(50,0){\line(-1,-1){20}}
\put(53,0){\makebox(0,0){$k$}}
\put(50,0){\circle*{2.0}}
\put(50,0){\circle*{2.0}}
\put(50,-20){\circle*{2.0}}
\put(50,0){\line(-1,0){40}}
\put(50,0){\line(1,2){10}}
\put(50,0){\line(1,-2){10}}
\put(80,0){\line(-1,-1){20}}
\put(55,-20){\line(-1,0){10}}
\put(63,-23){\makebox(0,0){$6$}}

\thicklines

\put(60,20){\circle*{2.0}}
\put(30,20){\circle*{2.0}}
\put(30,-20){\circle*{2.0}}
\put(60,-20){\circle*{2.0}}
\put(80,0){\circle*{2.0}}
\put(10,0){\circle*{2.0}}

\put(27,-21){\makebox(0,0){$1$}}
\put(26,23){\makebox(0,0){$3$}}
\put(7,0){\makebox(0,0){$2$}}
\put(63,23){\makebox(0,0){$4$}}
\put(83,0){\makebox(0,0){$5$}}

\put(10,0){\line(1,1){20}}
\put(30,-20){\line(-1,1){20}}
\put(30,20){\line(1,0){30}}
\put(60,20){\line(1,-1){20}}

\end{picture}

\end{center}

\caption{The subdiagram on vertices $\{k,1,2,3,4,6\}$ is mutation equivalent to $E_6$} 

\label{fig:Ahat-212}

\end{figure}

\credit{Subsubcase 3.1.2} \emph{The vertex $k$ is not contained in any cycle $C'\subset \Gamma k$ 
such that $C'$ has length greater than 4.}
Let us first assume, without loss of generality, that $T=[k,i_1=1,i_2=2]$
is the unique triangle in $\Gamma k$. Since $k$ is not contained in a cycle of length greater
than 4, we have $i_3=4$ and $i_5=6$. Then the subdiagram with the vertices $\{k,1,2,3,4,6\}$ is mutation equivalent to $E_6$ (Fig.~\ref{fig:Ahat-212}), contradicting \eqref{eq:knotinE6}.

\credit{Subcase 3.2} \emph{The diagram $\Gamma k$ contains precisely two triangles, say 
$T_1$ and $T_2$.}

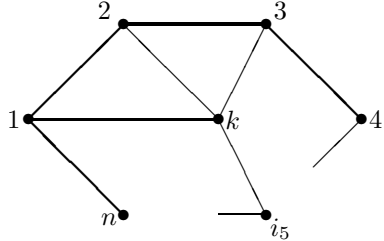
\begin{figure}[ht]

\setlength{\unitlength}{1.8pt}

\begin{center}

\begin{picture}(70,50)(5,-25)
\put(53,0){\makebox(0,0){$k$}}
\put(50,0){\circle*{2.0}}
\put(50,0){\line(-1,0){40}}
\put(50,0){\line(-1,1){20}}
\put(50,0){\line(1,2){10}}
\put(50,0){\line(1,-2){10}}
\put(80,0){\line(-1,-1){10}}
\put(60,-20){\line(-1,0){10}}
\put(63,-23){\makebox(0,0){$i_5$}}

\thicklines

\put(60,20){\circle*{2.0}}
\put(30,20){\circle*{2.0}}
\put(30,-20){\circle*{2.0}}
\put(60,-20){\circle*{2.0}}
\put(80,0){\circle*{2.0}}
\put(10,0){\circle*{2.0}}

\put(27,-21){\makebox(0,0){$n$}}
\put(26,23){\makebox(0,0){$2$}}
\put(7,0){\makebox(0,0){$1$}}
\put(63,23){\makebox(0,0){$3$}}
\put(83,0){\makebox(0,0){$4$}}

\put(10,0){\line(1,1){20}}
\put(30,-20){\line(-1,1){20}}
\put(30,20){\line(1,0){30}}
\put(60,20){\line(1,-1){20}}

\end{picture}

\end{center}

\caption{The subdiagram on vertices $\{n,1,2,3,4,k\}$ is mutation equivalent to $E_6$} 

\label{fig:Ahat-221}

\end{figure}

\credit{Subsubcase 3.2.1} \emph{$T_1$ and $T_2$ share a common edge.}
Let us assume, without loss of generality, that $T_1=[k,i_1=1,i_2=2]$
and $T_2=[k,i_2=2,i_3=3]$. Then $k$ is not
connected to any vertex in $\{4,n\}$, thus the subdiagram with the vertices  
$T_1 \cup T_2 \cup \{k,4,n\}$ (Fig.~\ref{fig:Ahat-221})
is mutation equivalent to $E_6$, contradicting \eqref{eq:knotinE6}.  

\credit{Subsubcase 3.2.2} \emph{$T_1$ and $T_2$ do not share a common edge.}
Then $\mu_k(\Gamma k)$ is of type $D_n^{(1)}(m,r,s)$. 

\credit{Subcase 3.3} \emph{The diagram $\Gamma k$ contains precisely three triangles, say 
$T_1,T_2$ and $T_3$.}
Then the diagram $\mu_k(\Gamma k)$ is of type $D_n^{(1)}(m,r)$. 

There is no other subcases in Case 3.

\credit{Case 4} 
\emph{r=2}.

\credit{Subcase 4.1} \emph{The vertices $i_1$ and $i_2$ are connected to each other.}
Then $\mu_k(\Gamma k)$ is of type $A_n^{(1)}$.

\credit{Subcase 4.2}\emph{The vertices $i_1$ and $i_2$ are not connected to each other.} 
Then the vertex $k$ is  contained in precisely two cycles, say $C_1$ and $C_2$ in $\Gamma k$. 
We may assume, without loss of generality, that the length of $C_1$
is less than or equal to the length of $C_2$.

\credit{Subsubcase 4.2.1} \emph{The cycle $C_1$ has length 4.} 
Then $\mu_k(\Gamma k)$ is of type $D_n^{(1)}(r)$. 

\credit{Subsubcase 4.2.2} \emph{The cycle $C_1$ has length greater than 4.} 
Let us assume, without loss of generality that, $i_1=1$
and $i_2 \geq 4$. Then the subdiagram with the vertices $\{k,i_2-2,i_2-1,i_2,i_2+1,i_2+2\}$ 
is an $E_6$, contradicting \eqref{eq:knotinE6}.





\subsection{Proof of Proposition~\ref{pr:main-add-seq}}
\label{subsec:series-2}
If $\Gamma$ is of type $A_n^{(1)}$, then the proposition follows from 
Proposition~\ref{pr:main-add-seq-A}. We will show that the proposition holds for each of the remaining diagrams $D_n^{(1)}$, $D_n^{(1)}(m,r)$, $D_n^{(1)}(r)$, $D_n^{(1)}(m,r,s)$. To do so, we will consider all possible cases in which the vertex $k$ may be 
connected to the diagram $\Gamma$ and show that \eqref{eq:knotinM-seq} or 
\eqref{eq:kct2-seq} is satisfied.

\begin{lemma}
\label{lem:caseDhat}

If $\Gamma$ is of type $D_n^{(1)}$, then 
Proposition~\ref{pr:main-add-seq} holds.

\end{lemma}

\proof
We assume that $\Gamma$ is indexed as in Table 1.

\credit{Case 1} \emph{$k$ is not connected to any vertex in the set $\{b_1,b_2,c_1,c_2\}$}.
Let us write $\{a_i: k \mathrm{\: is \: connected \: to\: }a_i\}=\{a_{i_1},...,a_{i_r}\}$ where
$1\leq i_1<...<i_r\leq m$ and $r \geq 2$. 

\credit{Subcase 1.1} \emph{$r \geq 3$}.
Let us first assume that $i_r < m$. Then the subdiagram with the vertices 
$$\{b_1,b_2,a_1, a_2,...,a_{i_1},k,a_{i_r}, a_{i_{r+1}},...,a_m,c_1,c_2\}$$ is of type
$D_n^{(1)}$ (Fig.~\ref{fig:Dhat-11}).
Let us now assume that $i_r = m$. Then the subdiagram with the vertices 
$$\{b_1,b_2,a_1, a_2,...,a_{i_1},k,a_{i_r},c_1,c_2\}$$ is of type
$D_n^{(1)}$.

\begin{figure}[ht]

\setlength{\unitlength}{1.8pt}

\begin{center}

\begin{picture}(120,60)(-20,-30)

\put(40,20){\line(1,-1){20}}
\put(40,20){\line(-1,-1){20}}
\put(40,20){\line(0,-1){10}}
\put(22,-3){\makebox(0,0){$a_{i_1}$}}
\put(62,-3){\makebox(0,0){$a_{i_r}$}}
\put(41,23){\makebox(0,0){$k$}}
\put(30,0){\circle*{2.0}}
\put(50,0){\circle*{2.0}}

\thicklines

\put(0,0){\circle*{2.0}}

\put(-20,20){\circle*{2.0}}  \put(-14,20){\makebox(0,0){$b_2$}}
\put(-20,-20){\circle*{2.0}} \put(-14,-20){\makebox(0,0){$b_1$}}
\put(20,0){\circle*{2.0}}
\put(40,0){\circle*{2.0}}
\put(60,0){\circle*{2.0}}
\put(80,0){\circle*{2.0}}
\put(100,20){\circle*{2.0}}
\put(100,-20){\circle*{2.0}}
\put(40,20){\circle*{2.0}}
\put(104,20){\makebox(0,0){$c_2$}}
\put(104,-20){\makebox(0,0){$c_1$}}

\put(0,0){\line(1,0){80}}
\put(0,0){\line(-1,-1){20}}
\put(0,0){\line(-1,1){20}}
\put(80,0){\line(1,1){20}}
\put(80,0){\line(1,-1){20}}

\put(78,-3){\makebox(0,0){$a_m$}}
\put(3,-3){\makebox(0,0){$a_1$}}

\end{picture}
\end{center}

\caption{} 

\label{fig:Dhat-11}

\end{figure}

\credit{Subcase 1.2} \emph{$r=2$}.
Let us first assume that $a_{i_1}$ and $a_{i_2}$ are not connected to each other.
If $i_2 < m$, then the subdiagram with the vertices 
$$\{b_1,b_2,a_1, a_2,...,a_{i_1},k,a_{i_2}, a_{i_2+1},...,a_m,c_1,c_2\}$$ is of type
$D_n^{(1)}$.
If $i_2=m$, then the subdiagram with the vertices 
$$\{b_1,b_2,a_1, a_2,...,a_{i_1},k,a_{i_2},c_1,c_2\}$$ is of type
$D_n^{(1)}$. If $a_{i_1}$ and $a_{i_2}$ are connected to each other (Fig.~\ref{fig:Dhat-12}), then $\mu_k(\Gamma k)$ is of type $D_n^{(1)}$.

\begin{figure}[ht]

\setlength{\unitlength}{1.8pt}

\begin{center}

\begin{picture}(120,60)(-20,-30)

\put(30,20){\line(1,-2){10}}
\put(30,20){\line(-1,-2){10}}
\put(22,-3){\makebox(0,0){$a_{i_1}$}}
\put(42,-3){\makebox(0,0){$a_{i_2}$}}
\put(31,23){\makebox(0,0){$k$}}
\put(30,20){\circle*{2.0}}

\thicklines

\put(0,0){\circle*{2.0}}

\put(-20,20){\circle*{2.0}}  \put(-14,20){\makebox(0,0){$b_2$}}
\put(-20,-20){\circle*{2.0}} \put(-14,-20){\makebox(0,0){$b_1$}}
\put(20,0){\circle*{2.0}}
\put(40,0){\circle*{2.0}}
\put(60,0){\circle*{2.0}}
\put(80,0){\circle*{2.0}}
\put(100,20){\circle*{2.0}}
\put(100,-20){\circle*{2.0}}
\put(104,20){\makebox(0,0){$c_2$}}
\put(104,-20){\makebox(0,0){$c_1$}}

\put(0,0){\line(1,0){80}}
\put(0,0){\line(-1,-1){20}}
\put(0,0){\line(-1,1){20}}
\put(80,0){\line(1,1){20}}
\put(80,0){\line(1,-1){20}}

\put(78,-3){\makebox(0,0){$a_m$}}
\put(3,-3){\makebox(0,0){$a_1$}}

\end{picture}
\end{center}

\caption{} 

\label{fig:Dhat-12}

\end{figure}

\credit{Case 2} \emph{$k$ is connected to a vertex, say $b_1$ in the set $\{b_1,b_2,c_1,c_2\}$}



\credit{Subcase 2.1} \emph{$k$ is not connected to $b_2$}. 

\credit{Subsubcase 2.1.1} \emph{$m \geq 2$}. Then we have the following:
\begin{itemize}
\item
If $m=2$, then either $k$ is contained in a non-oriented cycle or the subdiagram with the vertices 
the vertices $\{k,b_1,b_2,a_1,a_2,c_1\}$ is mutation equivalent to $E_6$ contradicting \eqref{eq:knotinE6}.
\item 
If $m> 2$, then either $k$ is contained in a non-oriented cycle or the subdiagram with the vertices $\{k,b_1,b_2,a_1,a_2,a_3\}$ is mutation 
equivalent to $E_6$, contradicting \eqref{eq:knotinE6}.
\end{itemize}

\credit{Subsubcase 2.1.2} \emph{$m=1$}. Then we have the following:
\begin{itemize}
\item
Suppose that $k$ is not connected to $a_1$. 
If $k$ is connected to both $c_1$ and $c_2$ then $k$ is contained in a 
non-oriented cycle, otherwise $\mu_k(\Gamma k)$ is of type $D_n^{(1)}(m,r)$.
\item
Suppose that $k$ is connected to $a_1$. Then $k$ is connected to a vertex in $\{c_1,c_2\}$. Let us first assume that $k$ is not
connected to any vertex in $C=\{c_1,c_2\}$. Then the subdiagram with the vertices 
$\{k,a_1,b_2,c_1,c_2\}$ is of type $D_n^{(1)}$. Let us now assume that
$k$ is connected to a vertex in $C$. If $k$ is connected to both $c_1$
and $c_2$ then then the subdiagram with the vertices $\{k,b_1,a_1,c_1,c_2\}$ is of
type $D_n^{(1)}(r)$; otherwise the diagram $\mu_k(\Gamma k)$ is of type $D_n^{(1)}$.   
\end{itemize}

\credit{Subcase 2.2} \emph{$k$ is connected to $b_2$}.

\credit{Subsubcase 2.2.1} \emph{$k$ is not connected to $a_1$}. Then we have the following:
\begin{itemize}
\item[(i)]
Suppose that $k$ is connected to a vertex $v$ in $\Gamma$ such that $v\ne b_1,b_2$.
Then $k$ is contained in a non-oriented cycle. 
\item[(ii)] 
If (i) does not hold, then $\mu_k(\Gamma k)$ is of type $D_n^{(1)}(m,r)$.
\end{itemize}


\credit{Subsubcase 2.2.2} \emph{k is connected to $a_1$.}Then we have the following: 
\begin{itemize}
\item[(i)]
Suppose that $k$ is connected to a vertex $v$ in $\Gamma$ such that $v\ne b_1,b_2,a_1$.
For such a vertex $v$, let us denote by $P_v$ the shortest path $[a_1,....,v]$ that
connects $a_1$ to $v$. We may assume, without loss of generality, that for any vertex
$v'$ in $P_v$ such that $v'\ne a_1,v$, the vertex $k$ is not connected to $v'$.
Then the subdiagram with the vertices $\{k,b_1,b_2\}\cup P_v$ is of type $D_n^{(1)}(r)$.  
\item[(ii)] If (i) does not hold, then $\mu_k(\Gamma k)$ is of type $D_n^{(1)}$.
\end{itemize}


\endproof

\begin{lemma}
\label{lem:caseDhat2}
If $\Gamma$ is of type $D_n^{(1)}(m,r)$, then Proposition~\ref{pr:main-add-seq} holds.
\end{lemma}

\proof 
We assume that $\Gamma$ is indexed as in Table 1.

\credit{Case 1} \emph{
$k$ is connected to a vertex,say $b_1$, in the set $\{b_1,b_2\}$}.

\credit{Subcase 1.1} \emph{$k$ is not connected to $b_2$}. Then we have the following:
\begin{itemize}
\item
If $m=1$, then either $k$ is contained in a non-oriented cycle or the subdiagram with the vertices $\{k,b_1,b_2,a_1,c_2,c_3\}$ is mutation equivalent to $E_6$, contradicting \eqref{eq:knotinE6}.
\item
If $m=2$, then the subdiagram with the vertices 
the vertices $\{k,b_1,b_2,a_1,a_2,c_1\}$ is mutation equivalent to $E_6$,
contradicting \eqref{eq:knotinE6}.
\item
If $m> 2$, then the subdiagram with the vertices $\{k,b_1,b_2,a_1,a_2,a_3\}$ is mutation 
equivalent to $E_6$ contradicting, \eqref{eq:knotinE6}.
\end{itemize}

\credit{Subcase 1.2} \emph{k is connected to $b_2$}. 

\credit{Subsubcase 1.2.1} \emph{$k$ is not connected to $a_1$}. 
\begin{itemize}
\item[(i)]
Suppose that $k$ is connected to a vertex $v$ in $\Gamma$ such that $v\ne b_1,b_2$. Then $k$ is necessarily contained in a non-oriented cycle. 
\item[(ii)]
If (i) does not hold, then $\mu_k(\Gamma k)$ is of type $D_n^{(1)}(m,r,s)$.
\end{itemize}

\credit{Subsubcase 1.2.2} \emph{k is connected to $a_1$}. Then we have the following:
\begin{itemize}
\item[(i)]
Suppose that $k$ is connected to a vertex $v$ in $\Gamma$ such that $v\ne b_1,b_2,a_1$.
For such a vertex $v$, let us denote by $P_v$ the shortest path $[a_1,....,v]$ that
connects $a_1$ to $v$. We may assume, without loss of generality, that for any vertex
$v'$ in $P_v$ such that $v'\ne a_1,v$, the vertex $k$ is not connected to $v'$.
Then the subdiagram with the vertices $\{k,b_1,b_2\}\cup P_v$ is of type $D_n^{(1)}(r)$.  
\item[(ii)]
If (i) does not hold, then $\mu_k(\Gamma k)$ is of type $D_n^{(1)}(m,r)$.
\end{itemize}

\credit{Case 2} \emph{$k$ is not connected to any vertex in the set $\{b_1,b_2\}$}.

\credit{Subcase 2.1} \emph{$k$ is not connected to any vertex in $\{c_1,...,c_r\}$}.
Let us write 
$$\{a_i: k \mathrm{\:\:is\:\: connected\:\: to\:\:} a_i\}=\{a_{i_1},...,a_{i_s}\}$$
where $i_1<...<i_s$, $s \geq 2$. 

\credit{Subcase 2.1.1} \emph{$s \geq 3$}.
Then the subdiagram with the vertices 
$$\{b_1,b_2,a_1, a_2,...,a_{i_1},k,a_{i_s}, a_{i_{s+1}},...,a_m,c_1,c_2,...,c_r\}$$ is of type
$D_n^{(1)}(m,r)$ (Fig.~\ref{fig:Dhat-211}).
 
\begin{figure}[ht]

\setlength{\unitlength}{1.8pt}

\begin{center}

\begin{picture}(120,80)(5,-30)

\put(40,20){\line(1,-1){20}}
\put(40,20){\line(-1,-1){20}}
\put(40,20){\line(0,-1){10}}
\put(22,-3){\makebox(0,0){$a_{i_1}$}}
\put(62,-3){\makebox(0,0){$a_{i_r}$}}
\put(41,23){\makebox(0,0){$k$}}
\put(30,0){\circle*{2.0}}
\put(50,0){\circle*{2.0}}

\thicklines

\put(0,0){\circle*{2.0}}

\put(-20,20){\circle*{2.0}}  \put(-14,20){\makebox(0,0){$b_2$}}
\put(-20,-20){\circle*{2.0}} \put(-14,-20){\makebox(0,0){$b_1$}}
\put(20,0){\circle*{2.0}}
\put(40,0){\circle*{2.0}}
\put(60,0){\circle*{2.0}}
\put(80,0){\circle*{2.0}}
\put(100,20){\circle*{2.0}}
\put(100,0){\circle*{2.0}}
\put(120,40){\circle*{2.0}}
\put(140,40){\circle*{2.0}}
\put(160,20){\circle*{2.0}}
\put(160,0){\circle*{2.0}}
\put(120,-20){\circle*{2.0}}
\put(140,-20){\circle*{2.0}}

\put(40,20){\circle*{2.0}}
\put(98,23){\makebox(0,0){$c_2$}}
\put(98,-3){\makebox(0,0){$c_1$}}
\put(118,43){\makebox(0,0){$c_3$}}
\put(118,-23){\makebox(0,0){$c_r$}}

\put(0,0){\line(1,0){100}}
\put(0,0){\line(-1,-1){20}}
\put(0,0){\line(-1,1){20}}
\put(80,0){\line(1,1){20}}
\put(100,0){\line(0,1){20}}
\put(100,0){\line(1,-1){20}}
\put(100,20){\line(1,1){20}}
\put(120,40){\line(1,0){20}}
\put(140,40){\line(1,-1){20}}
\put(160,20){\line(0,-1){20}}
\put(160,0){\line(-1,-1){20}}
\put(140,-20){\line(-1,0){20}}

\put(78,-3){\makebox(0,0){$a_m$}}
\put(3,-3){\makebox(0,0){$a_1$}}

\end{picture}
\end{center}

\caption{} 

\label{fig:Dhat-211}

\end{figure}

\credit{Subcase 2.1.2} \emph{
$s=2$}.
If $a_{i_1}$ and $a_{i_2}$ are not connected to each other,
then the subdiagram with the vertices
$$\{b_1,b_2,a_1, a_2,...,a_{i_1},k,a_{i_2}, a_{i_2+1},...,a_m,c_1,c_2,...,c_r\}$$ is of type
$D_n^{(1)}(m,r)$. 
If $a_{i_1}$ and $a_{i_2}$ are connected to each other, then $\mu_k(\Gamma k)$ is of type $D_n^{(1)}(m,r)$.

\credit{Subcase 2.2} \emph{ $k$ is
connected to a vertex in $C=\{c_1,...,c_r\}$}.
Let us write 
$$\{c_i: k \mathrm{\:is \:connected \:to}\: c_i\}=\{c_{i_1},...,c_{i_s}\}$$
where $1\leq i_1<...<i_s\leq r$. (Note that $C$ is an oriented cycle by our assumption the beginning of Section~\ref{sec:main-th}).

\credit{Subsubcase 2.2.1} \emph{$s \geq 3$}. 
Then $k$ is necessarily contained in 
a non-oriented cycle. 

\credit{Subsubcase 2.2.2} \emph{s=2}. 
We note that if $c_{i_1}$ and $c_{i_2}$ are not connected to each other
then $k$ is contained in a non-oriented cycle. Let us now
assume that $c_{i_1}$ and $c_{i_2}$ are connected to each other
and consider subcases.

\credit{Subsubsubcase 2.2.2.1} \emph{
$k$ is not connected to any vertex in 
$\{a_1,...,a_m\}$}. 
\begin{itemize}
\item
If $\{c_{i_1},c_{i_2}\} \ne \{c_1,c_2\}$, then 
$\mu_k(\Gamma k)$ is of type $D_n^{(1)}(m,r)$. 

\item
If $\{c_{i_1},c_{i_2}\} = \{c_1,c_2\}$, then the subdiagram formed
by $C \cup \{k,a_m\}$ is 
of type $D_n^{(1)}(r)$ 
(Fig.~\ref{fig:Dhat-2221}).
\end{itemize}

\begin{figure}[ht]

\setlength{\unitlength}{1.8pt}

\begin{center}

\begin{picture}(120,80)(5,-30)

\put(120,0){\line(-1,0){20}}
\put(120,0){\line(-1,1){20}}
\put(22,-3){\makebox(0,0){$a_{i_1}$}}
\put(62,-3){\makebox(0,0){$a_{i_r}$}}
\put(123,0){\makebox(0,0){$k$}}
\put(120,0){\circle*{2.0}}

\thicklines

\put(0,0){\circle*{2.0}}

\put(-20,20){\circle*{2.0}}  \put(-14,20){\makebox(0,0){$b_2$}}
\put(-20,-20){\circle*{2.0}} \put(-14,-20){\makebox(0,0){$b_1$}}
\put(20,0){\circle*{2.0}}
\put(40,0){\circle*{2.0}}
\put(60,0){\circle*{2.0}}
\put(80,0){\circle*{2.0}}
\put(100,20){\circle*{2.0}}
\put(100,0){\circle*{2.0}}
\put(120,40){\circle*{2.0}}
\put(140,40){\circle*{2.0}}
\put(160,20){\circle*{2.0}}
\put(160,0){\circle*{2.0}}
\put(120,-20){\circle*{2.0}}
\put(140,-20){\circle*{2.0}}

\put(98,23){\makebox(0,0){$c_2$}}
\put(98,-3){\makebox(0,0){$c_1$}}
\put(118,43){\makebox(0,0){$c_3$}}
\put(118,-23){\makebox(0,0){$c_r$}}

\put(0,0){\line(1,0){100}}
\put(0,0){\line(-1,-1){20}}
\put(0,0){\line(-1,1){20}}
\put(80,0){\line(1,1){20}}
\put(100,0){\line(0,1){20}}
\put(100,0){\line(1,-1){20}}
\put(100,20){\line(1,1){20}}
\put(120,40){\line(1,0){20}}
\put(140,40){\line(1,-1){20}}
\put(160,20){\line(0,-1){20}}
\put(160,0){\line(-1,-1){20}}
\put(140,-20){\line(-1,0){20}}

\put(78,-3){\makebox(0,0){$a_m$}}
\put(3,-3){\makebox(0,0){$a_1$}}

\end{picture}
\end{center}

\caption{} 

\label{fig:Dhat-2221}

\end{figure}

\credit{Subsubsubcase 2.2.2.2} \emph{
 $k$ is connected to a vertex in $\{a_1,...,a_m\}$}.
Then $k$ is contained in a non-oriented 
cycle. 






\credit{Subsubcase 2.2.3} \emph{s=1}. 
Then $k$ is either contained in a non-oriented cycle or it is contained in a diagram which is mutation equivalent to $E_6$, contradicting \eqref{eq:knotinE6}.





\endproof

\begin{lemma}
\label{lem:caseDhat4}
If $\Gamma$ is of type $D_n^{(1)}(m,r,s)$, then Proposition~\ref{pr:main-add-seq} holds.
\end{lemma}

\proof 
We assume that $\Gamma$ is indexed as in Table 1. 
Let us denote $C_1=\{b_1,...,b_s\}$ and $C_2=\{c_1,...,c_r\}$.



\credit{Case 1} \emph{
$k$ is connected to a vertex in $C_1$ or $C_2$}.
We assume, without loss of generality, that $k$ is connected to 
a vertex in $C_1$.
We also write 
$$\{b_i: k \mathrm{\:is \:connected \:to}\: b_i\}=\{b_{i_1},...,b_{i_t}\}$$
where $1\leq i_1<...<i_t\leq r$

\credit{Subcase 1.1} \emph{$t>3$}. 
Then $k$ is contained in 
a non-oriented cycle. 

\credit{Subcase 1.2} \emph{$t=2$}. 
Let us first note the following:
\begin{itemize}
\item[(i)]
If $b_{i_1}$ and $b_{i_2}$ are not connected to each other
then $k$ is contained in a non-oriented cycle. 
\item[(ii)]
If $k$ is connected to a vertex $v$ in $\Gamma$ such that
$v\ne b_{i_1},b_{i_2}$, then $k$ is contained in a non-oriented cycle.
\end{itemize}
Let us now
assume that (i) and (ii) do not hold. Then we have the following: 
\begin{itemize}
\item[(a)]
If $\{c_{i_1},c_{i_2}\} \ne \{c_1,c_2\}$, then 
$\mu_k(\Gamma k)$ is of type $D_n^{(1)}(m,r,s)$.
\item[(b)]
If $\{c_{i_1},c_{i_2}\} = \{c_1,c_2\}$, then the subdiagram induced
by $C \cup \{k,a_1\}$ is of type $D_n^{(1)}(r)$.
\end{itemize}







\credit{Subcase 1.3} \emph{$t=1$}. 
Then $k$ is contained in a diagram which is mutation equivalent to $E_6$, 
contradicting \eqref{eq:knotinE6}.

\credit{Case 2} \emph{
$k$ is not connected to any vertex in 
$C_1$ and not connected to any vertex in $C_2$}.
Let us write $\{a_i: k \mathrm{\: is\: connected \:to\:} a_i\: \}=\{a_{i_1},...,a_{i_l}\}$
where $i_1<...<i_l$. If $l\geq 3$, then the subdiagram with the vertices
$C_1 \cup \{a_1, a_2,...,a_{i_1},k,a_{i_l}, a_{i_{l+1}},...,a_m\} \cup C_2$ is of type
$D^4$. 
Let us now assume that $l=2$. If      
the vertices $a_{i_1}$ and $a_{i_2}$ are not connected to each other,
then the subdiagram with the vertices
$$C_1 \cup \{a_1, a_2,...,a_{i_1},k,a_{i_2}, a_{i_2+1},...,a_m,c_1,c_2\}$$ is of type
$D_n^{(1)}(m,r,s)$. 
If $a_{i_1}$ and $a_{i_2}$ are connected to each other, then $\mu_k(\Gamma k)$ is of type $D_n^{(1)}(m,r)$.
\endproof

\begin{lemma}
\label{lem:caseDhat3}

If $\Gamma$ is of type $D_n^{(1)}(r)$, then Proposition~\ref{pr:main-add-seq} holds.

\end{lemma}
\proof 
We assume that $\Gamma$ is indexed as in Table 1.

\credit{Case 1} \emph{
$k$ is not connected to any vertex in $\{a_1,c_1\}$}.  
Let us write 
$$\{b_i: k \mathrm{\:is \:connected \:to}\: b_i\}=\{b_{i_1},...,b_{i_t}\}$$
where $1\leq i_1<...<i_s\leq r, s\geq 2$.

\credit{Subcase 1.1} \emph{$s>3$}. 
Then $k$ is contained in 
a non-oriented cycle. 

\credit{Subcase 1.2} \emph{$s=2$}. 
We note that if $b_{i_1}$ and $b_{i_2}$ are not connected to each other,
then $k$ is contained in a non-oriented cycle. 
Let us now
assume that $b_{i_1}$ and $b_{i_2}$ are not connected to each other, then we have the following: 
\begin{itemize}
\item
If $\{c_{i_1},c_{i_2}\} \ne \{c_1,c_2\}$, then 
$\mu_k(\Gamma k)$ is of type $D_n^{(1)}(r)$.
\item
If $\{c_{i_1},c_{i_2}\} = \{c_1,c_2\}$, then the sub-diagram induced
by $C \cup \{k,a_1\}$ is of type $D_n^{(1)}(r)$.
\end{itemize}



 

\credit{Case 2} \emph{
 $k$ is connected to precisely one vertex, say $a_1$, in $\{a_1,c_1\}$}. 
If $k$ is connected to more than
one vertex in $B=\{b_1,...,b_r\}$, then it is necessarily contained in a
non-oriented cycle, so we assume that $k$ is connected to
precisely one vertex in $B$. 

\credit{Subcase 2.1}\emph{
$k$ is not connected to any vertex in $\{b_1,b_2\}$}.
Then $k$ is connected to a vertex in $\{b_3,...,b_r\}$ and hence
contained in a non-oriented cycle. 

\credit{Subcase 2.2}\emph{
$k$ is connected to a vertex, say $b_2$, in $\{b_1,b_2\}$}.
We note that if $k$ is connected to $b_1$, then $k$ is contained in
a non-oriented cycle. Let us now assume that $k$ is not connected to $b_1$.
Then we have the following:
\begin{itemize}
\item
If $r>5$, then the subdiagram with the vertices $\{b_r,b_1,b_2,b_3,b_4,k\}$ is
the tree $E_6$. 
\item
If $r=5$, then the subdiagram with the vertices $B \cup {k}$ is
mutation equivalent to $E_6$, contradicting \eqref{eq:knotinE6}.
\item
If $r=4$, then the subdiagram with the vertices $B \cup \{k,a_1\}$ is
mutation equivalent to $E_6$, contradicting \eqref{eq:knotinE6}.
\item
If $r=3$, then $\mu_k(\Gamma k)$ is of type $D_n^{(1)}(r)$.
\end{itemize}

\credit{Case 3} \emph{$k$ is connected to $a_1$ and $c_1$.} 
In this case, the vertex $k$ is always contained in a non-oriented cycle. 

\endproof

\section{Simply-laced minimal 2-infinite diagrams with at most 9 vertices: Proof of Proposition~\ref{pr:simply-laced-9} and another characterization of simply-laced 2-finite diagrams}
\label{sec:min<=9}

In this section we will characterize simply-laced 2-finite diagrams using basic moves, a class of graph transformations. Using this characterization we will prove Proposition~\ref{pr:simply-laced-9} in 
Section~\ref{subsec:simply-laced-9}.

\begin{definition}
\label{def:BM}

Suppose that $\bar{\Gamma}$ is a weightless undirected graph and that $a,c$ are two vertices which are connected to each other. The basic move $\phi_{c,a}$ is the transformation that changes $\bar{\Gamma}$ as follows: it connects $c$ to vertices that are connected to $a$ but not connected to $c$; at the same time it disconnects vertices from $c$ if they are connected to $a$. We call two graphs $\bar{\Gamma}$ and $\bar{\Gamma}'$ \emph{BM(basic move)-equivalent} if they can be obtained from each other by a sequence of basic moves.

For a diagram $\Gamma$, we denote by $\bar{\Gamma}$ the undirected graph which is defined as follows: the vertex set of $\bar{\Gamma}$ is the same as that of $\Gamma$ and two vertices $i,j$ in $\bar{\Gamma}$ are connected to each other if and only if they are connected in $\Gamma$ by an edge whose weight is an odd integer.
By a \emph{subgraph} of $\bar{\Gamma}$, we mean a graph 
obtained from $\bar{\Gamma}$ by taking an induced subgraph on a subset of
vertices.
\end{definition}

Basic moves were introduced in \cite{BHII}. Here we note that they
are related to mutations as follows:

\begin{proposition}
\label{pr:mut-BM}

Suppose that a diagram $\Gamma$ is mutation equivalent to $\Gamma'$. Then $\bar{\Gamma}$ is BM-equivalent to $\bar{\Gamma'}$.

\end{proposition}

\proof
It is enough to establish the proposition for $\Gamma'=\mu_k(\Gamma)$ where $k$ is an arbitrary vertex in $\Gamma$.
Let us assume that $c_1,c_2,...,c_r$ are the vertices which are connected to $k$ by a directed edge pointing towards $k$ with an odd weight. Then  
$\bar{\Gamma'}=\phi_{c_r,k}\circ...\phi_{c_2,k}\circ\phi_{c_1,k}(\bar{\Gamma})$.

\endproof

We observe that the converse of Proposition~\ref{pr:mut-BM} holds for 2-finite diagrams. More precisely:

\begin{theorem}
\label{th:mut=BM}
Suppose that $\Gamma$ is a simply-laced connected diagram that does not contain any non-oriented cycle. Let $\bar{\Gamma}$ be the underlying undirected graph of $\Gamma$. 
Then the diagram $\Gamma$ is 2-finite if and only if $\bar{\Gamma}$ is 
BM-equivalent to a Dynkin graph.
\end{theorem}

\noindent
We prove the theorem in Section~\ref{subsec:pf:mut=BM} after some preparation.

\subsection{Basic moves as linear transformations}
\label{subsec:BM-lin}
In this subsection we give some definitions related to basic moves; the details can be found in \cite{S}.
Let $\bar{\Gamma}$ be a connected graph with vertex set $I$. We denote by $V$ the vector space with the basis $I$. The graph $\bar{\Gamma}$ naturally defines an alternating $\FF_2$-valued bilinear form $\Omega$ on $V$ as follows: for any $i,j \in I$, $\Omega(i,j)=1$ if and only if $i$ and $j$ are connected to each other in $\bar{\Gamma}$. Then a basic move $\phi_{c,a}$ corresponds to the following change of basis for $\Omega$: 
$$\phi_{c,a}(I)=I-\{c\} \cup \{c+a\}.$$
For any vector subspace $U$ of $V$, we denote by $U_0$ the kernel of the restriction of $\Omega$ on $U$, i.e. $U_0=\{u\in U: \Omega(u,u')=0 \mathrm{\:for\:all\:}u'\in U\}$. 
We denote by $Q$ the $\FF_2$-valued quadratic form which is defined as follows: $Q(u+v)=Q(u)+Q(v)+\Omega(u,v), (u,v \in V)$ and $Q(i)=1$ for all $i \in I$.
We define $U_{00}=U_0 \cap Q^{-1}(0)$. Clearly $U_{00}$ is a vector subspace of $U_{0}$. If $V_0=V_{00}$, then the Arf invariant of $Q$ is defined as follows: $Arf(Q)=\sum Q(e_i)Q(f_i)$ where $\{e_1,f_1,...,e_r,f_r,h_1,...,h_p\}$ is a symplectic basis, i.e. a basis such that $\Omega(e_i,f_i)=1$ and the rest of the values of $\Omega$ are $0$.

We will use the following simple fact in our proof of 
Theorem~\ref{th:mut=BM}:

\begin{proposition}
\label{pr:BM-ker}
Suppose that $U$ is a vector subspace of codimension one in $V$. Then 
dim($V_{00}$) $\geq$ dim($U_{00}$)- 1. 
\end{proposition}
\noindent
Let $v\in V$ be a non-zero vector which is not in $U$. 
Let us assume that $K=\{x_1,...,x_p\}$ is a basis of $U_{00}$ and write $r(K)=\#\{x_i\in K:\Omega(v,x_i)=1\}$. If $r(K)=0$, then $K$ is also a basis for $V_{00}$, thus dim($V_{00}$)$=$dim($U_{00}$). Suppose that $r(K)=1$ and $\Omega(v,x_i)=1$. Then $K-\{x_i\}$ is a basis for $V_{00}$. Let us now assume that $r(K)>1$ and, assume without loss of generality, that $\Omega(v,x_i)=1$ for $i=1,2$. Then $K_1=\{x_1,x_1+x_2,x_3,...,x_p\}$ is also a basis for $U_{00}$ while $r(K_1)=r(K)-1$ (because $\Omega(v,x_1+x_2)=0$), thus the proposition follows by induction.

\subsection{Basic moves and Dynkin graphs}

In this subsection we give some properties of graphs which are BM-equivalent to 
Dynkin graphs and discuss their implications on mutations. Throughout this subsection,
we assume that $\bar{\Gamma}$ is a connected graph in the set-up of (sub)section~\ref{subsec:BM-lin}.

\begin{proposition}
\label{pr:BM-E6}
Suppose that $\bar{\Gamma}$ contains $E_6$. If the number of vertices in $\bar{\Gamma}$ is greater than or equal to $9$, then it is not BM-equivalent to any Dynkin graph.
\end{proposition}
\noindent
Since $\bar{\Gamma}$ contains $E_6$, any tree which is BM-equivalent to $\bar{\Gamma}$ also contains $E_6$ 
\cite[Theorem~2.7]{S}. Since no Dynkin graph with 9 or more vertices contains $E_6$, the graph $\bar{\Gamma}$ can not be equivalent to any Dynkin graph.

\begin{corollary}
\label{lem:E6-in-2-infinite}

Suppose that $\Gamma$ is a simply-laced diagram that contains a subdiagram which is mutation equivalent to $E_6$.
If the number of vertices in $\Gamma$ is greater than or equal to 9, then $\Gamma$ is $2$-infinite.

\end{corollary}

\begin{proposition}
\label{pr:BM-Arf}
In the set-up of Section~\ref{subsec:BM-lin}, we have the following description of the BM-equivalence classes of the simply-laced Dynkin graphs with $6$, $7$, or $8$ vertices:

\begin{itemize}
\item
Suppose that $\bar{\Gamma}$ has precisely 6 vertices. Then $\bar{\Gamma}$ is BM-equivalent to a Dynkin graph if and only if one of the following holds:
\begin{itemize}
\item[(a)]
$V_0=\{0\}$,
\item[(b)]
dim($V_0$)$=2$ and $V_0\ne V_{00}$.
\end{itemize}
\item
Suppose that $\bar{\Gamma}$ has precisely 7 vertices. Then $\bar{\Gamma}$ is BM-equivalent to a Dynkin graph if and only if dim($V_0$)$=1$ and 
one of the following holds:
\begin{itemize}
\item[(a)]
$V_0\ne V_{00}$,
\item[(b)]
$V_0=V_{00}$ and $Arf(Q)=0$
\end{itemize}
\item
Suppose that $\bar{\Gamma}$ has precisely 8 vertices. 
Then $\bar{\Gamma}$ is BM-equivalent to a Dynkin graph if and only if 
one of the following holds:
\begin{itemize}
\item[(a)]
$V_0=\{0\}$ and $Arf(Q)=0$.
\item[(b)]
dim($V_0$)$=1$ and $\bar{\Gamma}$ is BM-equivalent to the Dynkin graph $D_8$. (Note that an explicit description of graphs which are BM-equivalent to $D_8$ is given in \cite[Theorems~2.7,~2.10]{S}).
\end{itemize}
\end{itemize}

\end{proposition}
\noindent
This is a special case of \cite[Theorem~4.1]{BHII} and \cite[Theorem~3.8]{J} (see also \cite[Theorem~2.7]{S}).

\begin{proposition}
\label{pr:BM-V00}
Suppose that $\bar{\Gamma}$ contains a subgraph which is BM-equivalent to $E_6$. If dim($V_{00}$)$\geq 1$, then $\bar{\Gamma}$ is not BM-equivalent to any Dynkin graph.
\end{proposition}
\noindent
This statement is also a special case of \cite[Theorem~3.8]{J} and \cite[Theorem~2.7]{S}.

\begin{figure}[ht]

\setlength{\unitlength}{1.8pt}

\begin{center}

\begin{picture}(100,45)(-10,0)

\thicklines

\put(10,0){\circle*{2.0}$^{a_1}$}
\put(30,0){\circle*{2.0}$^{a_2}$}
\put(50,0){\circle*{2.0}$^{a_3}$}
\put(70,0){\circle*{2.0}$^{a_4}$}
\put(90,0){\circle*{2.0}$^{a_5}$}
\put(50,20){\circle*{2.0}$^{a_6}$}
\put(50,40){\circle*{2.0}$^{a_7}$}


\put(10,0){\line(1,0){80}}
\put(50,0){\line(0,1){40}}

\end{picture}

\end{center}

\caption{The extended Dynkin graph $E_6^{(1)} \,$.} 
\label{fig:BM-E7Fhat}
\end{figure}
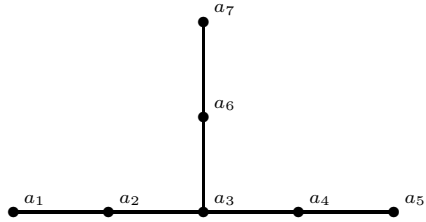

\begin{proposition}
\label{pr:BM-E6hat}
Suppose that $\bar{\Gamma}$ contains a subgraph $X$ which is $E_6^{(1)}$ or $E_7^{(1)}$. Then $\bar{\Gamma}$ is not BM-equivalent to any Dynkin graph.
\end{proposition}
\noindent
If the number of vertices in $\bar{\Gamma}$ is greater than or equal to 9, then the statement follows from Proposition~\ref{pr:BM-E6} because $X$ contains $E_6$ as a subgraph. Let us now assume that $\Gamma$ has at most 8 vertices. 
If $X$ is $E_7^{(1)}$, then $X=\bar{\Gamma}$ which is not BM-equivalent to any Dynkin graph by \cite{J}.
Let us now assume that $X$ is $E_6^{(1)}$ indexed as in 
Fig~\ref{fig:BM-E7Fhat}. Then $\bar{\Gamma}$ can not be BM-equivalent to any Dynkin graph with seven vertices \cite[Theorem~2.3]{S}, so we may assume that $\bar{\Gamma}$ has precisely $8$ vertices. We may also assume that $V_0=0$, otherwise $\bar{\Gamma}$ is not equivalent to any Dynkin graph by Proposition~\ref{pr:BM-V00}. Under these assumptions, by Proposition~\ref{pr:BM-Arf}, the graph 
$\bar{\Gamma}$ is BM-equivalent to a Dynkin graph if and only if $Arf(Q)=0$.
We will establish a symplectic basis and show that $Arf(Q)=1$, which will prove the proposition. We note that the set $$\{e_1=a_1,f_1=a_2,e_2=a_4,f_2=a_5,e_3=a_6,f_3=a_7,e_4=a_1+a_3+a_5+a_7\}$$ could be completed to a symplectic basis, i.e. there is a non-zero vector $f_4\in V$ such that $\Omega(e_4,f_4)=1$ and $\Omega(e_i,f_4)=\Omega(f_i,f_4)=0$ for $i=1,2,3$, which gives $Arf(Q)=1$ because $Q(a_1+a_3+a_5+a_7)=0$ and $Q(a_i)=1$ for $i=1,...,7$. 
 
\begin{proposition}
\label{pr:BM-Dhat}
If $\bar{\Gamma}$ contains a subgraph $X$ which is of type $D_n^{(1)}$, then it is not BM-equivalent to any Dynkin graph.
\end{proposition}
\noindent
Let us denote by 
$U\subset V$ the vector subspace which is spanned by the vertices in $X$
(recall that we work in the set-up of Section~\ref{subsec:BM-lin}). By \cite[Proposition~3.1]{S}, we have the following: if $X$ is a $D_4^{(1)}$ (resp. $D_n^{(1)}, n\geq 5$), then dim($U_{00}$)$=3$ (resp. dim($U_{00}$)$\geq 2$). 
If $\bar{\Gamma}$ does not contain any subgraph which is BM-equivalent to $E_6$, then dim($V_{00}$)$\geq 2$ 
by \cite[Theorem~2.9]{S}, implying the conclusion of the proposition by \cite[Theorem~2.3,2.7]{S}.
Let us now assume that $\bar{\Gamma}$ contains a subgraph which is BM-equivalent to $E_6$. If $\bar{\Gamma}$ has at least $9$ vertices, then the statement follows from Proposition~\ref{pr:BM-E6}.
If $\bar{\Gamma}$ has $6$ or $7$ vertices, then dim($U_{00}$)$\geq 1$ by Proposition~\ref{pr:BM-ker}, so the proposition follows from Proposition~\ref{pr:BM-V00}.
If $\bar{\Gamma}$ has $8$ vertices, then it contains a subgraph which is BM-equivalent to $E_6^{(1)}$ by
the classification of graphs under basic moves \cite[Theorems~2.3,~2.7]{S} and \cite[Theorem~3.8]{J}, thus the proposition follows from Proposition~\ref{pr:BM-E6hat}.

Finally, we will need the following statement on mutations.

\begin{figure}[ht]

\setlength{\unitlength}{1.8pt}

\begin{center}
\includegraphics[width=3.6in]{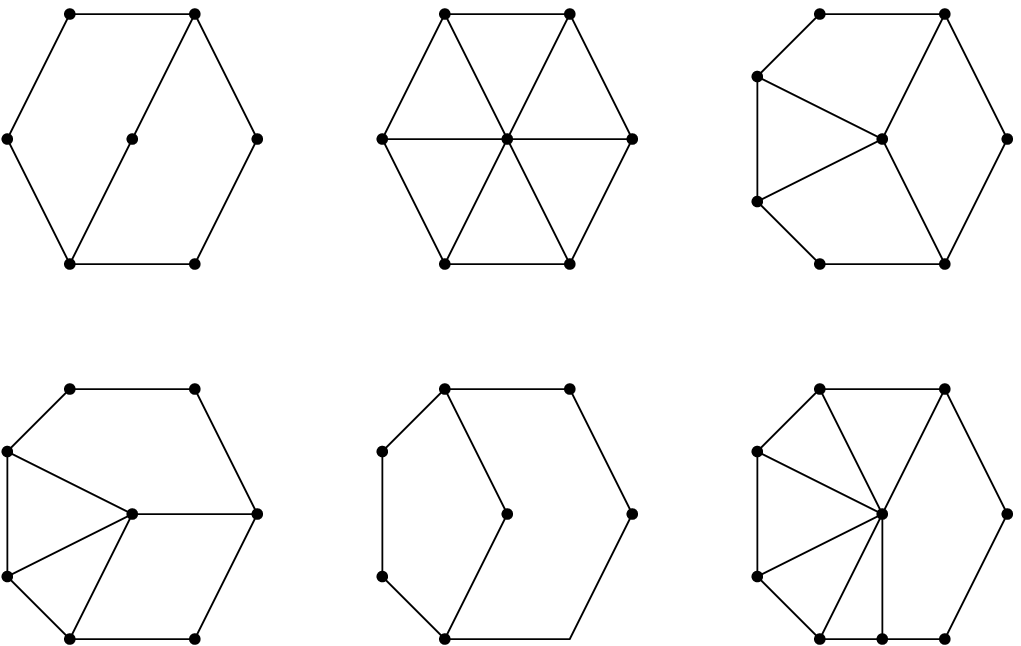}

\end{center}

\caption{} 
\label{fig:Ck}
\end{figure}

\begin{proposition}
\label{pr:BM-Ck}
Let $C$ be a diagram which is a non-oriented cycle whose length is less than or equal to 7 and let $Ck$ be a diagram obtained by adjoining a vertex $k$ to $C$. Suppose that $k$ is a vertex connected to at least two vertices of $C$. Suppose also that $k$ is not contained in any non-oriented cycle in $Ck$, in particular $k$ is connected to an even number of vertices in $C$. Then $Ck$ is mutation equivalent to a simply-laced extended Dynkin diagram. 
\end{proposition}
\noindent
If $Ck$ does not contain any subgraph which is mutation equivalent to $E_6$, then the statement is the same as Proposition~\ref{pr:main-add-seq-A}. 
If $Ck$ contains a subgraph which is mutation equivalent to $E_6$, then $Ck$ is one of the diagrams in Fig.\ref{fig:Ck}. 
It follows from a direct check that $Ck$ is mutation equivalent to a simply-laced extended Dynkin diagram. 
(In fact, the diagram $\mu_k(Ck)$ is in Table 4 or 5). We will prove in Lemma~\ref{lem:min-eDynkin} that each diagram in Table 4 (resp. Table 5) is mutation equivalent to $E_6^{(1)}$ (resp. $E_7^{(1)}$).

\subsection{Proof of Theorem~\ref{th:mut=BM}}
\label{subsec:pf:mut=BM}

In view of Proposition~\ref{pr:mut-BM}, it is enough to prove that 

(*) if $\Gamma$ is not 2-finite, then $\bar{\Gamma}$ is not BM-equivalent to any Dynkin graph.

\noindent
To establish (*), let us first notice that 
there is a sequence
of mutations $\mu_1,...,\mu_k$, $k\geq 1$ and a simply-laced diagram $\Gamma'$ such that 
\begin{itemize} 
\item[{\rm(i)}] 
$\Gamma'=\mu_k \circ...\circ \mu_1(\Gamma)$,
\item[{\rm(ii)}]
the diagram $\Gamma'$ contains a non-oriented cycle $C$,
\item[{\rm(iii)}]
for any $i:1,...,k-1$, the diagram $\Gamma_i=\mu_i \circ ...\circ \mu_1(\Gamma)$ is simply-laced 
and it does not contain any non-oriented cycle.
\end{itemize}
\noindent
Let us note that $k$ is not contained in $C$ and it is connected to an even number of vertices in $C$ (because $\mu_k(\Gamma')=\Gamma_{k-1}$ is a simply-laced diagram which does not contain any non-oriented cycle). We denote by $Ck$ the subdiagram with the vertices $C \cup {k}$. Let us first assume that $C$ has length less than or equal to $7$. Then $\mu_k(Ck)\subset \Gamma_{k-1}$ is mutation equivalent to an extended Dynkin diagram, say $X$, (Proposition~\ref{pr:BM-Ck}). Note that $\bar{Ck}$ is BM-equivalent to $\bar{X}$ by Proposition~\ref{pr:mut-BM}.   
Since $\Gamma_{k-1}$ does not contain any non-oriented cycle, the diagram $X$ is one of the following: $D_n^{(1)}$, $E_6^{(1)}$, $E_7^{(1)}$. Thus (*) follows from Proposition~\ref{pr:BM-Dhat} and Proposition~\ref{pr:BM-E6hat}. Let us now assume that the length of $C$ is greater than or equal to $8$. If $Ck$ does not contain any subdiagram which is mutation equivalent to $E_6$, then $Ck$ is mutation equivalent to $D_n^{(1)}$ by Proposition~\ref{pr:main-add-seq}, so 
(*) follows from Proposition~\ref{pr:BM-Dhat}; otherwise it follows from Proposition~\ref{pr:BM-E6}.

\noindent
\begin{corollary}
\label{cor:BM-E6-inv}
Suppose that $\Gamma$ is a diagram that contains (a subdiagram which is mutation equivalent to) $E_6$. Let $\Gamma'$ be a simply-laced diagram which is mutation equivalent to $\Gamma$. Suppose also that $\Gamma'$ does not contain any non-oriented cycle. Then $\Gamma$ contains a subdiagram which is mutation equivalent to $E_6$.
\end{corollary}
\noindent
By \cite[Theorem~2.7]{S}, the diagram $\bar{\Gamma'}$ contains a subgraph, say $\bar{E}$, which is BM-equivalent to $E_6$. Since $\bar{\Gamma'}$ does not contain any non-oriented cycle, the subdiagram $E$ is mutation equivalent to (an orientation of) $E_6$ by the theorem.

\begin{corollary}
\label{cor:BM-E6-fin}

Suppose that $\Gamma$ is a simply-laced 2-finite diagram. Suppose also that $\Gamma$ contains a subdiagram which is mutation equivalent to $E_6$. Then $\Gamma$ is mutation equivalent to one of the following diagrams: 
$E_6$, $E_7$, $E_8$.
\end{corollary}
\noindent
The statement immediately follows from Corollary~\ref{cor:BM-E6-inv} (and the classification of 2-finite diagrams). 


\begin{corollary}
\label{cor:BM-D}

Suppose that $\Gamma$ is a simply-laced minimal 2-infinite diagram that does not contain any subdiagram which is mutation equivalent to $E_6$. Then $\Gamma$ is in Table 1, i.e. $\Gamma$ is one of the following diagrams: $A_n^{(1)}$, $D_n^{(1)}$, $D_n^{(1)}(m,r)$, $D_n^{(1)}(r)$, $D_n^{(1)}(m,r,s)$.

\end{corollary}
\noindent
If $\Gamma$ contains a non-oriented cycle $C$, then $\Gamma=C$ (i.e. of type $A_n^{(1)}$) because $\Gamma$ is minimal 2-infinite. Let us now assume that $\Gamma$ does not contain any non-oriented cycle. Then there is a sequence
of mutations $\mu_1,...,\mu_k$, $k\geq 1$ and a simply-laced diagram $\Gamma'$ such that 
\begin{itemize} 
\item[{\rm(i)}] 
$\Gamma'=\mu_k \circ...\circ \mu_1(\Gamma)$,
\item[{\rm(ii)}]
the diagram $\Gamma'$ contains a non-oriented cycle $C$,
\item[{\rm(iii)}]
for any $i:1,...,k-1$, the diagram $\Gamma_i=\mu_i \circ ...\circ \mu_1(\Gamma)$ is simply-laced
and it does not contain any non-oriented cycle.
\end{itemize}
We note that $k$ is not contained in $C$ and it is connected to an even number of vertices in $C$ (because $\mu_k(\Gamma')=\Gamma_{k-1}$ is a simply-laced diagram which does not contain any non-oriented cycle). We also note that, for any $i:1,...,k-1$, the diagram $\Gamma_i$ does not contain any subdiagram which is mutation equivalent to $E_6$ (Corollary~\ref{cor:BM-E6-inv}). This implies that the vertex $k$ is not contained in any subdiagram which is mutation equal to $E_6$ in $Ck$ (otherwise $\mu_k(Ck)\subset \mu_k(\Gamma')=\Gamma_{k-1}$ contains one). By Proposition~\ref{pr:main-add-seq-A}, the diagram $\mu_k(Ck)$ is one of the diagrams $D_n^{(1)}$, $D_n^{(1)}(m,r)$, $D_n^{(1)}(r)$, $D_n^{(1)}(m,r,s)$. By Proposition~\ref{pr:main-add-seq} and Lemma~\ref{lem:star}, for any $j=1,...,k-1$, the diagram 
$\Gamma^j=\mu_{k-j} \circ...\mu_{k-1}\circ \mu_k(\Gamma')$, in particular $\Gamma^{k-1}=\Gamma$, contains one of the diagrams $D_n^{(1)}$, $D_n^{(1)}(m,r)$, $D_n^{(1)}(r)$, $D_n^{(1)}(m,r,s)$. Since $\Gamma$ is minimal 2-infinite, the diagram $\Gamma$ must be one of $D_n^{(1)}$, $D_n^{(1)}(m,r)$, $D_n^{(1)}(r)$, $D_n^{(1)}(m,r,s)$.

\begin{corollary}
\label{cor:BM-A}
A simply-laced diagram $\Gamma$ is mutation equivalent to the Dynkin diagram of type $A$ if and only if it does not have subdiagrams of the following form:
\begin{align} 
\label{eq:A} 
&\text{a non-oriented cycle,}
\\
\nonumber
&\text{a tree of type $D_4$,}
\\
\nonumber 
&\text{two triangles sharing a common edge,} 
\\
\nonumber 
&\text{any (oriented) cycle whose length is greater than or equal to 4.} 
\end{align} 
\end{corollary}
\noindent
In view of Theorem~\ref{th:mut=BM}, the corollary follows from \cite[Theorem~2.7]{S}.

\begin{corollary}
\label{cor:BM-E6-min}

Suppose that $\Gamma$ is a simply-laced minimal 2-infinite diagram. Suppose also that $\Gamma$ contains a subdiagram which is mutation equivalent to $E_6$. Then $\Gamma$ is one of the diagrams in Tables 4-6 in Section~\ref{subsec:list}.
\end{corollary}
\noindent
We obtained this corollary by implementing the following algorithm.

\subsection{Algorithm to compute minimal 2-infinite diagrams that contains a subdiagram which is mutation equivalent to $E_6$}
\label{subsec:alg}
Let us denote the set of those diagrams by $\bold{E}$. Let $\Gamma\in \bold{E}$ be a diagram with $n$ vertices. Then $7\leq n \leq9$ by Corollary~\ref{lem:E6-in-2-infinite}. If $n=8$ or $n=9$, then $\Gamma$ contains a subdiagram which is mutation equivalent to $E_{n-1}^{(1)}$ (Corollary~\ref{cor:BM-E6-fin}). Thus we can compute $\bold{E}$ as follows:

\vspace{.1in}
\noindent
For $n=7,8,9$: 
\begin{itemize}
\item
compute the mutation class of $E_{n-1}$,
\item
for each diagram $\Gamma$ in the mutation class of $E_{n-1}$, and each non-empty subset $X$ of vertices in $\Gamma$,  
\begin{itemize}
\item
form the graph $\Gamma k$ by adjoining a vertex $k$ to $\Gamma$ as follows: $k$ is connected to all of the vertices in $X$ but not connected any other vertex in $\Gamma$,
\item
check if $\Gamma k$ is minimal 2-infinite using Theorem~\ref{th:mut=BM} and Proposition~\ref{pr:BM-Arf},
\item
(if $\Gamma k$ is minimal 2-infinite, check whether it is in the mutation class of $E_{n-1}^{(1)}$).
\end{itemize}

\end{itemize}
\noindent
We implemented this algorithm using the computer algebra package Maple (v.8). Our code is available at \cite{W}.   

\subsection{Proof of Proposition~\ref{pr:simply-laced-9}}
\label{subsec:simply-laced-9}
If $\Gamma$ contains a subdiagram which is mutation equivalent to $E_6$, then the theorem follows from Corollary~\ref{cor:BM-E6-min}. 
If $\Gamma$ does not contain any subdiagram which is mutation equivalent to $E_6$, then the theorem follows from Corollary~\ref{cor:BM-D}.

\section{Proof of Theorem~\ref{th:main} }

\label{sec:pf-main}


We prove the theorem following the outline given in Section~\ref{sec:main-th}, i.e. by proving 
Lemmas~\ref{lem:List-min-2inf}, \ref{lem:star}, \ref{lem:xxx}. Since Lemma~\ref{lem:xxx} is 
somewhat more involved, we will prove it in a seperate section, Section~\ref{subsec:lem-xxx}, 
for convenience.

\subsection{Proof of Lemma~\ref{lem:List-min-2inf}}
\label{subsec:List-in-M}









We first prove that $\Gamma$ is 2-infinite. According to \cite[Propositions~9.3,~9.7]{CAII}), any non-oriented oriented cycle is 2-infinite. In particular any diagram in Table 2 is 2-infinite. By \cite[Propositions~9.3,~9.7]{CAII}), any extended Dynkin diagram is also 2-infinite. Thus for $\Gamma$ which is not in Table 2 the fact 
that it is 2-infinite follows from the following stronger statement:

\begin{lemma}
\label{lem:min-eDynkin}

Suppose that $\Gamma$ in Section~\ref{subsec:list} is not one of the diagrams in Table 2. 
Then $\Gamma$ is mutation equivalent to an extended Dynkin diagram (Fig.~\ref{fig:extended-dynkin-diagrams}). 

\end{lemma}

We prove the lemma using a case by case analysis. Throughout the proof, we  assume that $\Gamma$ is indexed as in the associated table  (Section~\ref{subsec:list}).

\credit{Case 1}\emph{$\Gamma$ is in Table 1.}  

\credit{Subcase 1.1}\emph{$\Gamma$ is of type $B_n^{(1)}(m,r)$ (resp. $B_n^{(1)}(r)$).} 
Then we have the following:
\begin{itemize}
\item
If $r=3$, then $\mu_{b_1}(\Gamma)$ is of type $B_n^{(1)}$. 
\item
If $r>3$, then $\mu_{b_1}(\Gamma)$ is of type $B_n^{(1)}(m+1,r-1)$ 
(resp. $B^{(1)}(1,r-1)$), so it is mutation equivalent to $B_n^{(1)}$ by induction on $r$.
\end{itemize}

\credit{Subcase 1.2}\emph{$\Gamma$ is of type $D_n^{(1)}(m,r)$.} 
Then we have the following:
\begin{itemize}
\item
If $r=3$, then $\mu_{c_1}(\Gamma)$ is of type $D_n^{(1)}$. 
\item
If $r>3$, then $\mu_{c_1}(\Gamma)$ is of type $D_n^{(1)}(m+1,r-1)$, so it is mutation equivalent to $D_n^{(1)}$ by induction on $r$. 
\end{itemize}

\credit{Subcase 1.3}\emph{$\Gamma$ is of type $D_n^{(1)}(r)$.} 
Then we have the following:
\begin{itemize}
\item
If $r=3$, then $\mu_{b_1}(\Gamma)$ is of type $D_4^{(1)}$.
\item
If $r>3$, then $\mu_{b_1}(\Gamma)$ is of type $D_n^{(1)}(1,r-1)$ which is mutation equivalent to $D_n^{(1)}$ by Subcase 1.1.
\end{itemize}

\credit{Subcase 1.4}\emph{$\Gamma$ is of type $D_n^{(1)}(m,r,s)$.} 
Then we have the following:
\begin{itemize}
\item
If $r=3$, then $\mu_{c_1}(\Gamma)$ is of type $D_n^{(1)}(m+1,s)$, which is mutation equivalent to $D_n^{(1)}$ by Subcase 1.2. 
\item
If $r>3$, then $\mu_{c_1}(\Gamma)$ is of type $D_n^{(1)}(m+1,r-1,s)$, so it is mutation equivalent to $D_n^{(1)}$ by induction on $r$. 
\end{itemize}

\credit{Case 2}\emph{$\Gamma$ is in Table 3.}  

\credit{Subcase 2.1}\emph{$\Gamma$ contains a triangle, say $T$, which is not adjacent to any cycle.} Then $T$ contains a vertex, say $k$, which is adjacent to precisely two edges one of them being weightless. Then the mutation $\mu_k$ eliminates the triangle $T$. Continuing this process one has   $F_4^{(1)}$ or $G_2^{(1)}$. 

\credit{Subcase 2.2}\emph{Subcase 2.1 does not hold}. In this case $\Gamma$ has two cycles and there is a vertex $k$ which is contained in both cycles.  If the cycles are triangles, then $\mu_k(\Gamma)$ is as in Subcase 2.1; otherwise $\mu_k(\Gamma)$ has two ajacent triangles, and applying this case once more will reduce it to the Subcase 2.1.

\credit{Case 3}\emph{$\Gamma$ is in Tables 4-6.}  
Then, it follows from a computer check that $\Gamma$ is mutation equivalent to $E_6^{(1)}$, $E_7^{(1)}$ or $E_8^{(1)}$ \cite{W}. 

We are done with the proof of Lemma~\ref{lem:min-eDynkin}, thus have shown that any diagram $\Gamma$ in Section~\ref{subsec:list} is 2-infinite. Now we will prove that $\Gamma$ is \emph{minimal}:

\begin{lemma}
\label{lem:List-min}

Suppose that $\Gamma$ is an arbitrary diagram in Section~\ref{subsec:list}. Then any connected subdiagram $\Gamma'$ obtained from $\Gamma$ by deleting a vertex $b$ is mutation equivalent to a Dynkin diagram.

\end{lemma}

We prove the lemma using a case by case analysis. Throughout the proof, we  assume that $\Gamma$ is indexed as in the associated table (Section~\ref{subsec:list}).

\credit{Case 1}\emph{$\Gamma$ is in Table 1.}  
The lemma obviously holds for the extended Dynkin diagrams and the diagrams $I_2(a), a\geq 4$.

\credit{Subcase 1.1}\emph{$\Gamma$ is of type $D_n^{(1)}(m,r)$.}

\credit{Subsubcase 1.1.1}\emph{$b=b_1$ or $b=b_2$.}
Then the diagram $\mu_{c_{r-1}}...\mu_{c_1}(\Gamma')$ is of type $D$.

\credit{Subsubcase 1.1.2}\emph{$b=c_i$ for some $i: 1<i<r$.}
Then the diagram $\mu_{c_{i-1}}...\mu_{c_1}(\Gamma')$ is of type $D$.

The remaining (sub)subcases are trivial.

\credit{Subcase 1.2}\emph{$\Gamma$ is of type $D_n^{(1)}(r)$.} 

\credit{Subsubcase 1.2.1}\emph{$b=a_1$ or $b=c_1$}
Let us assume, without loss of generality, that $b=a_1$. Then
$\mu_{c_1}(\Gamma')$ is an (oriented) cycle which is mutation equivalent to a Dynkin graph of type $D$ \cite{CAII}.

\credit{Subsubcase 1.2.2}\emph{$b=b_i$ for some $i: 1<i<r$.}
Then the diagram $\mu_{c_{i-1}}...\mu_{c_1}(\Gamma')$ is of type $D$.

The remaining (sub)subcases are trivial.

\credit{Subcase 1.3}\emph{$\Gamma$ is of type $D_n^{(1)}(m,r,s)$.} 
By the same arguments as in Subcase 1.1, the diagram $\Gamma'$ is mutation equivalent to a tree of type $A$ or $D$.

\credit{Subcase 1.4}\emph{$\Gamma$ is of type $B_n^{(1)}(m,r)$ or $B_n^{(1)}(r)$.} 
In this case $\Gamma'$ is mutation equivalent to a tree of type $A,B$ or $D$.

\credit{Case 2}\emph{$\Gamma$ is in Tables 2-3.}  
In this case the lemma is almost obvious.

\credit{Case 3}\emph{$\Gamma$ is in Tables 4-6.}  
For this case we obtained the lemma 
by a computer check (see Section~\ref{subsec:alg} and \cite{W}).


\subsection{Proof of Lemma~\ref{lem:star}}
\label{subsec:lem-star}
Let us also show this lemma by considering possible cases:

\credit{Case 1}\emph{$\Gamma$ is in Table 1.}  
In this case we have one of the following:
\begin{itemize}
\item
$\mu_k(\Gamma k)$ is in Table 1.
\item
The subdiagram $\Gamma'\subset \mu_k(\Gamma k)$ formed by all of the vertices which are different from $k$ is in Table 1.
\end{itemize}
For most $k$ we observed those two statements in our proof of Lemma~\ref{lem:min-eDynkin}. For the remaining vertices they follow from a direct check.

\credit{Case 2}\emph{$\Gamma$ is in Tables 2,3.}  
In this case we have one of the following:
\begin{itemize}
\item
$\mu_k(\Gamma k)$ is in Table 2,3.
\item
The subdiagram $\Gamma'\subset \mu_k(\Gamma k)$ formed by all of the vertices which are different from $k$ is in 
Section~\ref{subsec:list}.
\end{itemize}
Those two statements also follow from a direct check.

\credit{Case 3}\emph{$\Gamma$ is in Tables 4-6.}  
In this case the lemma 
follows from Proposition~\ref{pr:simply-laced-9} 
(because $\Gamma k$ is 2-infinite so is $\mu_k(\Gamma k)$)).


\section{Proof of Lemma~\ref{lem:xxx}}
\label{subsec:lem-xxx}


We have divided the proof into several subsections. In each
subsection, we assume that the diagrams $\Gamma$ and $\Gamma k$
are as in the title (of the subsection).
We will also assume, without loss of generality, that $k$ is
connected to at least two vertices in $\Gamma$, (otherwise
$\Gamma \subset \mu_k(\Gamma k)$, so we may take $\Gamma'=\Gamma$). 
We will show that
one of the following two statements is satisfied:
\begin{align} 
\label{eq:knotinM}&\text{$k$ is contained in a diagram $\Gamma''\subset \Gamma k$ such that 
$\Gamma''$ is in Section~\ref{subsec:list},}
\\
\label{eq:kct2}&\text{the diagram $\mu_k(\Gamma k)$ contains a subdiagram $\Gamma'$ which is
in Section~\ref{subsec:list}}
\\
\nonumber
\end{align} 
Let us note that \eqref{eq:knotinM} implies \eqref{eq:kct2}
by Lemma~\ref{lem:star}, so Lemma~\ref{lem:xxx} holds in any case. 

\subsection{The diagram $\Gamma k$ is simply laced}
\label{subsec:lem-xxx-1}
If $\Gamma$ has at most 8 vertices, then obviously $\mu_k(\Gamma k)$ has at most 9 vertices and it is 2-infinite
, so we have \eqref{eq:kct2} by Proposition~\ref{pr:simply-laced-9}. Let us now assume that 
$\Gamma$ has at least 9 vertices and that
\begin{align} 
\label{eq:kinE6}&\text{$\Gamma k$ contains in a subdiagram $E\subset \Gamma k$
which is mutation}
\\
\nonumber
&\text{equivalent to the Dynkin diagram $E_6$}.
\end{align}
We claim that there exists a connected subdiagram $Xk\subset \Gamma k$ of
9 vertices which contains both $E$ and $k$. To prove this claim we argue as follows: 
if $\Gamma$ is in Table 1, then it follows by inspection that $\Gamma$ does not
contain any subdiagram which is mutation equivalent to $E_6$, so the subdiagram $E$ 
must contain the vertex $k$ and we may take $Xk$ to be any connected subdiagram 
which has 9 vertices and contains $E$. If $\Gamma$ is not in Table 1, then it 
must belong to Table 6 with precisely 9 vertices.
We note, by inspection, that $\Gamma$ contains
a subdiagram  which is mutation equivalent to $E_6$, thus in our set-up of 
\eqref{eq:kinE6} we may assume that $E \subset \Gamma$.
Let $X\subset \Gamma$ be any (connected) subdiagram with precisely $8$ vertices such that $E\subset X$. Note that $X$ is 2-finite because it is a proper subdiagram of $\Gamma$, which is minimal 2-infinite 
(Lemma~\ref{lem:List-min-2inf}). Since $k$ is connected to at least two vertices in $\Gamma$, it is connected to at least one vertex in $\Gamma$, so we may take $Xk$ to be the subdiagram formed by $X$ and $k$.
Then $Xk$ is 2-infinite by Corollary~\ref{lem:E6-in-2-infinite}, thus it contains a subdiagram $X'$
from Section~\ref{subsec:list} by Proposition~\ref{pr:simply-laced-9}. 
Since $X$ is 2-finite, the vertex $k$ is necessarily contained in $X'$, thus \eqref{eq:knotinM} is satisfied.

If \eqref{eq:kinE6} is not satisfied, then Lemma~\ref{lem:xxx} 
follows from Proposition~\ref{pr:main-add-seq}.

\subsection{The diagram $\Gamma k$ is non-simply-laced while $\Gamma$ is simply-laced}
\label{sec:pflemma6}
A glance at Table 2,3 shows that 
\begin{align}
\label{eq:edge-3}
~\text{if $\Gamma k$ contains an edge $e$ whose weight is greater than or equal to 3, then $e$ is} 
\\
\nonumber
~\text{contained in a subdiagram $\Gamma'\subset \Gamma k$ which is one of the diagrams in Table 2 or Table 3,}
\end{align}
therefore \eqref{eq:knotinM} is satisfied. Let us now assume that 
the condition of \eqref{eq:edge-3} is not satisfied. 
Then, by Definition~\ref{lem:def-diagram}, any edge which contains $k$ has weight equal to $2$. 

\credit{Case 1}\emph{$k$ is connected to two vertices, say $i,j$, which are not connected to each other.} 
Then the subdiagram with the vertices $\{k,i,j\}$ is of type $C_n^{(1)}$. 

\credit{Case 2}\emph{Case 1 does not hold.} 
Then $k$ is connected to precisely two vertices, say $i,j$, in $\Gamma$. Furthermore, the vertices $i$ and $j$ are connected to each other.

\credit{Subcase 2.1}\emph{$\Gamma$ is of type $A_n^{(1)}$, i.e. a non-oriented cycle.}
Then the diagram $\mu_k(\Gamma k)$ is of type $B_n^{(1)}(r)$.

\credit{Subcase 2.2}\emph{$\Gamma$ is not of type $A_n^{(1)}$.}

\credit{Subsubcase 2.2.1}\emph{The edge $[i,j]$ is contained in a cycle, say $C$, in $\Gamma$.}
Then the subdiagram with $C\cup\{k\}$ is of type $B_n^{(1)}(r)$.

\credit{Subsubcase 2.2.2}\emph{The edge $[i,j]$ is not contained in any cycle in $\Gamma$.}
We first note, by inspection, that $\Gamma$ contains a subdiagram $X$ which is one of the following: a tree of type $D_4$, two triangles sharing a common edge, a (oriented) cycle whose length is greater than or equal to 4. 
Let $P$ be the shortest path that connects $k$ to $X$. Then the subdiagram with $P$ and $X$ contains a subdiagram which is of type $B_n^{(1)}$ or $B_n^{(1)}(m,r)$.

\subsection{The diagram $\Gamma$ is a non-simply-laced diagram in Table 1}
\label{sec:pflemma7}


\begin{lemma}
\label{lem:caseChat}

If $\Gamma$ is of type $C_n^{(1)}$, then 
Lemma~\ref{lem:xxx} holds.
\end{lemma}

\proof 
In view of \eqref{eq:edge-3}, we may assume that the maximum edge-weight in $\Gamma k$ is $2$.

\credit{Case 1} \emph{
$k$ is not connected to any vertex in the set $\{b_1,c_1\}$.}

\credit{Subcase 1.1} \emph{
$k$ is not contained in any edge whose weight is equal to 2 in $\Gamma k$.}
Then we have the following:
\begin{itemize}
\item[(i)]
Suppose that $k$ is connected to two distinct vertices $a_i,a_j$ 
in $\{a_1,...,a_m\}$ and $a_i,a_j$ are not connected 
to each other by an edge. Let $r=\min\{a_i:k \mathrm{\:is\: connected\: to\:} a_i\}$
and $s=\max\{i:k \mathrm{\:is \:connected \:to\:} a_i\}$.
Then the subdiagram with the vertices $$\{b_1,a_1,...,a_r,k,a_s,....,a_m,c_1\}$$
is of type $C_n^{(1)}$ (Fig.~\ref{fig:Chat-11}).
\item[(ii)]
If (i) does not hold, then $k$ is connected to precisely
two vertices $a_i,a_j$ in $\{a_1,...,a_m\}$ such that 
$a_i,a_j$ are connected to each other by an edge.
Then $\mu_k(\Gamma k)$ is of type $C_n^{(1)}$. 
\end{itemize}

\begin{figure}[ht]

\setlength{\unitlength}{1.8pt}

\begin{center}

\begin{picture}(70,30)(-10,-10)

\put(10,0){\line(1,1){10}}
\put(40,0){\line(-2,1){20}}
\put(20,10){\line(0,-1){6}}
\put(20,10){\line(1,-1){6}}

\thicklines

\put(0,0){\circle*{2.0}}
\put(10,0){\circle*{2.0}}
\put(-10,0){\circle*{2.0}}
\put(20,0){\circle*{2.0}}
\put(30,0){\circle*{2.0}}
\put(40,0){\circle*{2.0}}
\put(50,0){\circle*{2.0}}
\put(60,0){\circle*{2.0}}
\put(20,10){\circle*{2.0}}

\put(-10,0){\line(1,0){70}}

\put(-5,3){\makebox(0,0){$2$}}
\put(55,3){\makebox(0,0){$2$}}
\put(20,14){\makebox(0,0){$k$}}
\put(10,-3){\makebox(0,0){$r$}}
\put(40,-3){\makebox(0,0){$s$}}
\put(60,-3){\makebox(0,0){$c_1$}}
\put(-10,-3){\makebox(0,0){$b_1$}}

\end{picture}
\end{center}

\caption{} 

\label{fig:Chat-11}

\end{figure}

\credit{Subcase 1.2} \emph{$k$ is contained in an edge whose weight is equal to 2.}
Then for any $a_i$ connected to $k$, the edge $[k,a_i]$ is weighted 2.
Let $P$ be the shortest path that connects $k$ to $c_1$.
Then the subdiagram with the vertices $P$ is of type $C_n^{(1)}$   
(Fig.~\ref{fig:Chat-12}).

\begin{figure}[ht]

\setlength{\unitlength}{1.8pt}

\begin{center}

\begin{picture}(70,30)(-10,-10)

\put(20,10){\line(-1,-1){6}}
\put(40,0){\line(-2,1){20}}
\put(20,10){\line(0,-1){6}}
\put(20,10){\line(1,-1){6}}

\thicklines

\put(0,0){\circle*{2.0}}
\put(10,0){\circle*{2.0}}
\put(-10,0){\circle*{2.0}}
\put(20,0){\circle*{2.0}}
\put(30,0){\circle*{2.0}}
\put(40,0){\circle*{2.0}}
\put(50,0){\circle*{2.0}}
\put(60,0){\circle*{2.0}}
\put(20,10){\circle*{2.0}}

\put(-10,0){\line(1,0){70}}

\put(-5,3){\makebox(0,0){$2$}}
\put(55,3){\makebox(0,0){$2$}}
\put(20,14){\makebox(0,0){$k$}}
\put(30,8){\makebox(0,0){$2$}}
\put(60,-3){\makebox(0,0){$c_1$}}
\put(-10,-3){\makebox(0,0){$b_1$}}

\end{picture}
\end{center}

\caption{} 

\label{fig:Chat-12}

\end{figure}


\credit{Case 2} \emph{
$k$ is connected to precisely one vertex, say $b_1$, in the set $\{b_1,c_1\}$.} 
We note that the vertex $k$ 
is connected to at least one vertex in  $\{a_1,...,a_m\}$.


\credit{Subcase 2.1} \emph{The edge $[k,b_1]$ is weightless}.
We note that
for any vertex $a$ in $\{a_1,...,a_m\}$ such that 
$k$ is connected to $a$, the edge $[a,k]$ is weighted 2.
Let us write $s=\max\{i:k\:\mathrm{is \:connected\: to\:} a_i\}$. 
Then the diagram $\Gamma'$ induced by $\{k,a_s,...,a_m,c_1\}$
is of type $C_n^{(1)}$. 

\credit{Subcase 2.2} \emph{The  edge $[k,b_1]$ is weighted}.
We note that, for any vertex $a$ in $\{a_1,...,a_m\}$ such that 
$k$ is connected to $a$, the edge $[a,k]$ is weightless.
Then we have the following:
\begin{itemize}
\item
If $k$ is not connected to $a_1$, then
the subdiagram with $\{k,b_1,a_1\}$ is of type $C_n^{(1)}$ 
\item
Suppose that $k$ is connected to $a_1$. If $k$ is not connected to 
any vertex in $A=\{a_2,....,a_m\}$, then $\mu_k(\Gamma k)$ is a $C_n^{(1)}$.
Let us now assume that $k$ is connected to a vertex in $A$ and let
$s=\min\{i:k\:\mathrm{is \: connected \:to\:} a_i \in A\}$. Then the diagram induced
by $\{k,b_1,a_1,...,a_s\}$ is a $B_n^{(1)}(r)$. 
\end{itemize}

\credit{Case 3} \emph{k is connected to $b_1$ and $c_1$.} 
Let us first assume that the edge $[k,b_1]$ has weight equal to 2.
Then the edge $[k,c_1]$ has also weight equal to 2,
so the subdiagram with the vertices $\{b_1,k,c_1\}$ is of type $C_n^{(1)}$.
Let us now assume that the edge $[k,b_1]$, hence $[k,c_1]$, is weightless
and consider subcases.

\credit{Subcase 3.1} \emph{$k$ is connected to a vertex in $A=\{a_1,...,a_m\}$.}
We note that for any vertex $a\in A$ such that $k$ is connected
to $a$, the edge $[k,a]$ is weighted $2$.
Let us write
$$\{a_i:k \mathrm{\:is\: connected\: to\:} a_i\}=\{a_{i_1},...,a_{i_r}\}, r\geq 1,$$
where $1\leq i_1 <i_2<...<i_r\leq m$. 

\credit{Subsubcase 3.1.1} \emph{$r\geq 3$.}
Then the subdiagram with the vertices $\{a_{i_1},k,a_{i_r}\}$ is of type $C_n^{(1)}$.

\credit{Subsubcase 3.1.2} \emph{$r=2$.}
If the vertices $a_{i_1}$ and $a_{i_2}$ are not connected to each other,
then the subdiagram with the vertices $\{a_{i_1},k,a_{i_2}\}$ is of type $C_n^{(1)}$.
Let us now assume that $a_{i_1}$ and $a_{i_2}$ are connected to each other.
Then we have the following:

\begin{itemize}
\item[(i)] If $i_1\ne 1$, then the subdiagram with the vertices $\{a_1,b_1,k,a_{i_2}\}$
is of type $C_n^{(1)}$ (Fig~\ref{fig:Chat-312i}).
\item[(ii)]
If $i_2\ne m$ then the subdiagram with the vertices $\{a_m,c_1,k,a_{i_1}\}$
is of type $C_n^{(1)}$.
\item[(iii)]
If neither of the conditions in (i,ii) holds,
then $i_1=1$ and $i_2=2=m$, so $\mu_k(\Gamma k)$ is of type $F_4^{(1)}(3^{2};0)$.
\end{itemize}

\begin{figure}[ht]

\setlength{\unitlength}{1.8pt}

\begin{center}

\begin{picture}(70,30)(-10,-10)

\put(-12,1){\line(4,1){30}}
\put(20,10){\line(0,-1){10}}
\put(20,10){\line(1,-1){10}}
\put(60,0){\line(-4,1){40}}

\thicklines

\put(0,0){\circle*{2.0}}
\put(10,0){\circle*{2.0}}
\put(-10,0){\circle*{2.0}}
\put(20,0){\circle*{2.0}}
\put(30,0){\circle*{2.0}}
\put(40,0){\circle*{2.0}}
\put(50,0){\circle*{2.0}}
\put(60,0){\circle*{2.0}}
\put(20,10){\circle*{2.0}}

\put(-10,0){\line(1,0){70}}

\put(-5,-3){\makebox(0,0){$2$}}
\put(55,-3){\makebox(0,0){$2$}}
\put(20,14){\makebox(0,0){$k$}}
\put(30,5){\makebox(0,0){$2$}}
\put(18,5){\makebox(0,0){$2$}}
\put(60,-3){\makebox(0,0){$c_1$}}
\put(-10,-3){\makebox(0,0){$b_1$}}
\put(0,-3){\makebox(0,0){$a_1$}}
\put(20,-3){\makebox(0,0){$a_{i_1}$}}
\put(30,-3){\makebox(0,0){$a_{i_2}$}}

\end{picture}
\end{center}

\caption{} 

\label{fig:Chat-312i}

\end{figure}


\credit{Subsubcase 3.1.3} \emph{$r=1$.} Then we have the following:
\begin{itemize}
\item[(i)] If $1<i_1<m$, then the subdiagram with the vertices $\{k,a_{i_1}, a_{i_1-1},a_{i_1+1}\}$ is of type $B_n^{(1)}$ (Fig.~\ref{fig:Chat-313i}).
\item[(ii)]
Suppose that the condition of (i) does not hold.
Then, without loss of generality, we may assume that
$i_1=1$. If $m>2$, then the subdiagram with the vertices $\{a_1,k,c_1,a_m\}$is of type
$C_n^{(1)}$. If $m=2$, then the diagram $\mu_k(\Gamma k\}$ is of type
$F_4^{(1)}(3^{2};0;1^1)$. If $m=1$, then $\mu_k(\Gamma k\}$ is of type
$B_n^{(1)}$.
\end{itemize}

\begin{figure}[ht]

\setlength{\unitlength}{1.8pt}

\begin{center}

\begin{picture}(70,35)(-10,-10)

\put(-10,0){\line(3,2){30}}
\put(20,20){\line(1,-2){10}}
\put(60,0){\line(-2,1){40}}

\thicklines

\put(0,0){\circle*{2.0}}
\put(10,0){\circle*{2.0}}
\put(-10,0){\circle*{2.0}}
\put(20,0){\circle*{2.0}}
\put(30,0){\circle*{2.0}}
\put(40,0){\circle*{2.0}}
\put(50,0){\circle*{2.0}}
\put(60,0){\circle*{2.0}}
\put(20,20){\circle*{2.0}}

\put(-10,0){\line(1,0){70}}

\put(-5,-3){\makebox(0,0){$2$}}
\put(55,-3){\makebox(0,0){$2$}}
\put(20,24){\makebox(0,0){$k$}}
\put(23,8){\makebox(0,0){$2$}}
\put(60,-3){\makebox(0,0){$c_1$}}
\put(-10,-3){\makebox(0,0){$b_1$}}
\put(0,-3){\makebox(0,0){$a_1$}}
\put(30,-3){\makebox(0,0){$a_{i_1}$}}
\put(50,-3){\makebox(0,0){$a_{m}$}}

\end{picture}
\end{center}

\caption{} 

\label{fig:Chat-313i}

\end{figure}

\credit{Subcase 3.2}\emph{
$k$ is not connected to any vertex in $\{a_1,...,a_m\}$.}
Then we have the following:
\begin{itemize}
\item If $m>2$ , then the subdiagram with the vertices $\{a_1,b_1,k,c_1,a_m\}$
is of type $C_n^{(1)}$ (Fig.~\ref{fig:Chat-32i}).
\item
If $m=2$ then $\mu_k(\Gamma k)$ is of type $F_4^{(1)}(4^{1};3^1)$.
\item
If $m=1$ then then $\mu_k(\Gamma k)$ is of type $B_n^{(1)}(r)$. 
\end{itemize}

\begin{figure}[ht]

\setlength{\unitlength}{1.8pt}

\begin{center}

\begin{picture}(70,30)(-10,-10)

\put(-12,1){\line(4,1){30}}
\put(60,0){\line(-4,1){40}}

\thicklines

\put(0,0){\circle*{2.0}}
\put(10,0){\circle*{2.0}}
\put(-10,0){\circle*{2.0}}
\put(20,0){\circle*{2.0}}
\put(30,0){\circle*{2.0}}
\put(40,0){\circle*{2.0}}
\put(50,0){\circle*{2.0}}
\put(60,0){\circle*{2.0}}
\put(20,10){\circle*{2.0}}

\put(-10,0){\line(1,0){70}}

\put(-5,-3){\makebox(0,0){$2$}}
\put(55,-3){\makebox(0,0){$2$}}
\put(20,14){\makebox(0,0){$k$}}
\put(60,-3){\makebox(0,0){$c_1$}}
\put(-10,-3){\makebox(0,0){$b_1$}}
\put(0,-3){\makebox(0,0){$a_1$}}

\end{picture}
\end{center}

\caption{} 

\label{fig:Chat-32i}

\end{figure}








\endproof

\begin{lemma}
\label{lem:caseBhat}
If $\Gamma$ is of type $B_n^{(1)}$, then Lemma~\ref{lem:xxx} holds.
\end{lemma}

\proof
In view of \eqref{eq:edge-3}, we may assume that the maximum edge-weight in $\Gamma k$ is $2$.

\credit{Case 1} \emph{$k$ is not connected to $c_1$.}

\credit{Subcase 1.1}\emph{ 
$k$ is contained in an edge whose weight is equal to 2.} 
Let us first note that, for any
vertex $v$ such that $k$ is connected to $v$, the edge $[k,v]$ is weighted 2.
Then we have the following:

\begin{itemize}
\item
Suppose that $k$ is connected to a vertex in $\{a_1,...,a_m\}$ and 
let $$s=\max\{i: k\mathrm{ \:is \:connected \:to\:} a_i\}.$$ Then the subdiagram with the vertices
$\{k,a_s,...,a_m,c_1\}$ is of type $C_n^{(1)}$ (Fig.~\ref{fig:Bhat-11i}).
\item
Suppose that $k$ is not connected to any vertex in $\{a_1,...,a_m\}$.
Then $k$ is connected to $b_1$ and $b_2$,
 thus the subdiagram with the vertices $\{k,b_1,a_1,...,a_m,c_1\}$ is of type $C_n^{(1)}$.
\end{itemize}

\begin{figure}[ht]

\setlength{\unitlength}{1.8pt}

\begin{center}

\begin{picture}(120,60)(-20,-30)

\put(40,20){\line(1,-1){20}}
\put(40,20){\line(-1,0){10}}
\put(40,20){\line(-1,-1){10}}
\put(40,20){\line(0,-1){10}}

\thicklines

\put(0,0){\circle*{2.0}$^{a_1}$}

\put(-20,20){\circle*{2.0}$^{b_2}$}
\put(-20,-20){\circle*{2.0}} \put(-14,-20){\makebox(0,0){$b_1$}}
\put(20,0){\circle*{2.0}$^{a_2}$}
\put(40,0){\circle*{2.0}}
\put(60,0){\circle*{2.0}$\: ^{a_s}$}
\put(80,0){\circle*{2.0}$\: ^{a_m}$}
\put(100,0){\circle*{2.0}$^{c_1}$}
\put(40,20){\circle*{2.0}}
\put(43,23){\makebox(0,0){$k$}}

\put(0,0){\line(1,0){100}}
\put(0,0){\line(-1,-1){20}}
\put(0,0){\line(-1,1){20}}

\put(90,-3){\makebox(0,0){$2$}}
\put(53,10){\makebox(0,0){$2$}}

\end{picture}
\end{center}

\caption{} 

\label{fig:Bhat-11i}

\end{figure}

\credit{Subcase 1.2}\emph{ 
$k$ is not contained in any edge weighted 2}. 

\credit{Subsubcase 1.2.1}\emph{ 
 $k$ is connected to a vertex in $\{a_1,....,a_m\}$}.
Let $s=\max\{i: k \mathrm{\:is \:connected \:to\:} a_i\}$. 
\begin{itemize}
\item
Suppose that $s>1$. If $k$ is not connected to $a_{s-1}$,
then the subdiagram with the vertices $\{a_{s-1},a_s,...,a_m,c_1,k\}$ is of type $B_n^{(1)}$ 
(Fig.~\ref{fig:Bhat121ia})
Let us now assume that $k$ is connected to $a_{s-1}$. Then we have the following:
\begin{itemize}
\item
If $k$ is not connected to any vertex $v$ such that $v\ne a_s,a_{s-1}$, 
then the diagram $\mu_k(\Gamma k)$ is of type $B_n^{(1)}$. 
\item
Suppose that $k$ is connected to vertex $v$ such that $v\ne a_s,a_{s-1}$. 
Let $v'$  be such a vertex that is closest to $a_{s-1}$ 
and let $P$ be the shortest undirected path that connects $a_s$ to $v'$. 
Then the subdiagram with
$P \cup \{a_s,a_{s+1},...,a_m,c_1,k\}$ is of type $B_n^{(1)}(m,r)$ 
(Fig.~\ref{fig:Bhat-121ib}).
\end{itemize}
\item
Suppose that $s=1$. 
If there is a vertex, say $b_1$, in $\{b_1,b_2\}$ such that $k$ is not
connected to $b_1$, then
the subdiagram with the vertices $\{k,b_1,a_1,....,a_m,c_1\}$ is of type $B_n^{(1)}$ 
.
 If k is connected to $b_1$ and $b_2$ then $\mu_k(\Gamma k)$ is of type 
$B_n^{(1)}$. 
\end{itemize}

\begin{figure}[ht]

\setlength{\unitlength}{1.8pt}

\begin{center}

\begin{picture}(120,60)(-20,-30)

\put(40,20){\line(1,-1){20}}
\put(40,20){\line(-1,0){10}}
\put(40,20){\line(-1,-1){20}}

\thicklines

\put(0,0){\circle*{2.0}$^{a_1}$}

\put(-20,20){\circle*{2.0}$^{b_2}$}
\put(-20,-20){\circle*{2.0}} \put(-14,-20){\makebox(0,0){$b_1$}}
\put(20,0){\circle*{2.0}$^{a_2}$}
\put(40,0){\circle*{2.0}}
\put(60,0){\circle*{2.0}}
\put(80,0){\circle*{2.0}$\: ^{a_m}$}
\put(100,0){\circle*{2.0}$^{c_1}$}
\put(40,20){\circle*{2.0}}
\put(43,23){\makebox(0,0){$k$}}

\put(0,0){\line(1,0){100}}
\put(0,0){\line(-1,-1){20}}
\put(0,0){\line(-1,1){20}}

\put(90,-3){\makebox(0,0){$2$}}
\put(40,-3){\makebox(0,0){$a_{s-1}$}}
\put(60,-3){\makebox(0,0){$a_{s-1}$}}

\end{picture}
\end{center}

\caption{} 

\label{fig:Bhat121ia}

\end{figure}

\begin{figure}[ht]

\setlength{\unitlength}{1.8pt}

\begin{center}

\begin{picture}(120,60)(-20,-30)

\put(40,20){\line(1,-1){20}}
\put(40,20){\line(-1,0){60}}
\put(40,20){\line(-1,-1){20}}
\put(40,20){\line(0,-1){20}}

\thicklines

\put(0,0){\circle*{2.0}$^{a_1}$}

\put(-20,20){\circle*{2.0}} 
\put(-20,-20){\circle*{2.0}} 
\put(20,0){\circle*{2.0}}
\put(40,0){\circle*{2.0}}
\put(60,0){\circle*{2.0}}
\put(80,0){\circle*{2.0}$\: ^{a_m}$}
\put(100,0){\circle*{2.0}$^{c_1}$}
\put(40,20){\circle*{2.0}}
\put(43,23){\makebox(0,0){$k$}}

\put(0,0){\line(1,0){100}}
\put(0,0){\line(-1,-1){20}}
\put(0,0){\line(-1,1){20}}

\put(90,-3){\makebox(0,0){$2$}}
\put(60,-3){\makebox(0,0){$a_s$}}
\put(40,-3){\makebox(0,0){$a_{s-1}$}}
\put(20,-3){\makebox(0,0){$v'$}}

\end{picture}
\end{center}

\caption{} 

\label{fig:Bhat-121ib}

\end{figure}

\credit{Subsubcase 1.2.1}\emph{ 
 $k$ is not connected to any vertex in $\{a_1,....,a_m\}$}. 
Then $\mu_k(\Gamma k)$ is of type $B_n^{(1)}(m,r)$. 


\credit{Case 2} \emph{k is connected to $c_1$.}

\credit{Subcase 2.1}\emph{The edge $[k,c_1]$ is weighted 2}.
We first note that the subdiagram with the vertices $\{b_1,b_2,a_1,...,a_m,k\}$ is weightless.
\begin{itemize}
\item
If $k$ is not connected to $a_m$, then
the subdiagram with the vertices $\{k,c_1,a_m\}$ is of type $C_n^{(1)}$. 
\item Suppose that $k$ is connected
to $a_m$. If $k$ is connected to a vertex $v$ such that $v\ne a_m,c_1$,  then $k$ is contained in a diagram
of type $B_n^{(1)}(r)$ (Fig.~\ref{fig:Bhat-21ii}),
otherwise $\mu_k(\Gamma k)$ is of type $B_n^{(1)}(r)$. 
\end{itemize}

\begin{figure}[ht]

\setlength{\unitlength}{1.8pt}

\begin{center}

\begin{picture}(120,60)(-20,-30)

\put(80,20){\line(1,-1){20}}
\put(80,20){\line(-1,0){10}}
\put(80,20){\line(-2,-1){40}}
\put(80,20){\line(0,-1){20}}

\thicklines

\put(0,0){\circle*{2.0}$^{a_1}$}

\put(-20,20){\circle*{2.0}} 
\put(-20,-20){\circle*{2.0}} 
\put(20,0){\circle*{2.0}}
\put(40,0){\circle*{2.0}}
\put(60,0){\circle*{2.0}}
\put(80,0){\circle*{2.0}$\: ^{a_m}$}
\put(100,0){\circle*{2.0}$^{c_1}$}
\put(80,20){\circle*{2.0}}
\put(83,23){\makebox(0,0){$k$}}

\put(0,0){\line(1,0){100}}
\put(0,0){\line(-1,-1){20}}
\put(0,0){\line(-1,1){20}}

\put(90,-3){\makebox(0,0){$2$}}
\put(91,12){\makebox(0,0){$2$}}
\put(40,-3){\makebox(0,0){$v$}}

\end{picture}
\end{center}

\caption{} 

\label{fig:Bhat-21ii}

\end{figure}

\credit{Subcase 2.2} \emph{The edge $[k,c_1]$ is weightless}.
We note that, for any vertex $v\ne c_1$ such that $k$ is connected to $v$, the edge $[k,v]$ is weighted 2. 
Let us first assume that $k$ is not connected any vertex in $\{b_1,b_2\}$ and write
$$s=\min\{i: k \mathrm{\:is \:connected \:to\:} a_i\}.$$ Then the subdiagram with the vertices $\{b_1,b_2,a_1,...,a_s,k\}$
is of type $B_n^{(1)}$ (Fig.~\ref{fig:Bhat-22}).
Let us now assume that 
$k$ is connected to a vertex in $\{b_1,b_2\}$.
If $k$ is connected to both $b_1$ and $b_2$, then 
the subdiagram with the vertices $\{k,b_1,b_2\}$ is of type $C_n^{(1)}$. 

To consider the remaining (sub)subcases, we assume, without loss of generality,
that $k$ is connected to $b_2$ and not connected to $b_1$.

\begin{figure}[ht]

\setlength{\unitlength}{1.8pt}

\begin{center}

\begin{picture}(120,60)(-20,-30)

\put(60,20){\line(2,-1){40}}
\put(60,20){\line(-1,-1){20}}
\put(60,20){\line(0,-1){20}}

\thicklines

\put(0,0){\circle*{2.0}$^{a_1}$}

\put(-20,20){\circle*{2.0}} 
\put(-20,-20){\circle*{2.0}} 
\put(20,0){\circle*{2.0}}
\put(40,0){\circle*{2.0}}
\put(60,0){\circle*{2.0}}
\put(80,0){\circle*{2.0}$\: ^{a_m}$}
\put(100,0){\circle*{2.0}$^{c_1}$}
\put(60,20){\circle*{2.0}}

\put(0,0){\line(1,0){100}}
\put(0,0){\line(-1,-1){20}}
\put(0,0){\line(-1,1){20}}

\put(90,-3){\makebox(0,0){$2$}}
\put(48,11){\makebox(0,0){$2$}}
\put(60,24){\makebox(0,0){$k$}}
\put(40,-3){\makebox(0,0){$s$}}
\put(62,11){\makebox(0,0){$2$}}

\end{picture}
\end{center}

\caption{} 

\label{fig:Bhat-22}

\end{figure}

\credit{Subsubcase 2.2.1} \emph{$m>1$.}
If $k$ is not connected to any vertex in $\{a_1,...,a_m\}$, then the subdiagram with the vertices
$\{b_2,k,c_1,a_m\}$ is of type $C_n^{(1)}$. 
Let us now assume that $k$ is connected to a vertex
in $\{a_1,...,a_m\}$ and write
$s=\max\{i:k \mathrm{\:is \:connected \:to\:} a_i\}$. If $s>1$, then the subdiagram with the vertices 
$\{b_2,k,a_s\}$ is of type $C_n^{(1)}$ (Fig.~\ref{fig:Bhat-221}).
If $s=1$, then the subdiagram with the vertices
$\{b_2,k,c_1,a_m\}$ is of type $C_n^{(1)}$. 

\begin{figure}[ht]

\setlength{\unitlength}{1.8pt}

\begin{center}

\begin{picture}(120,60)(-20,-30)

\put(20,20){\line(-1,0){40}}
\put(20,20){\line(4,-1){80}}
\put(20,20){\line(-1,-1){20}}
\put(20,20){\line(0,-1){20}}

\thicklines

\put(0,0){\circle*{2.0}$^{a_1}$}

\put(-20,20){\circle*{2.0}} 
\put(-20,-20){\circle*{2.0}} 
\put(20,0){\circle*{2.0}}
\put(40,0){\circle*{2.0}}
\put(60,0){\circle*{2.0}}
\put(80,0){\circle*{2.0}$\: ^{a_m}$}
\put(100,0){\circle*{2.0}$^{c_1}$}
\put(20,20){\circle*{2.0}}

\put(0,0){\line(1,0){100}}
\put(0,0){\line(-1,-1){20}}
\put(0,0){\line(-1,1){20}}

\put(90,-3){\makebox(0,0){$2$}}
\put(0,23){\makebox(0,0){$2$}}
\put(20,24){\makebox(0,0){$k$}}
\put(20,-3){\makebox(0,0){$a_s$}}
\put(8,10){\makebox(0,0){$2$}}
\put(-20,24){\makebox(0,0){$b_2$}}
\put(-20,-23){\makebox(0,0){$b_1$}}
\put(23,10){\makebox(0,0){$2$}}

\end{picture}
\end{center}

\caption{} 

\label{fig:Bhat-221}

\end{figure}


\credit{Subsubcase 2.2.2} \emph{$m=1$.}
If $k$ is not connected to $a_1$, then $\mu_k(\Gamma k)$ is of type $F_4^{(1)}(3^{2};0;1^1)$. 
Similarly if $k$ is connected to $a_1$, then $\mu_k(\Gamma k)$ is $F_4^{(1)}(3^{1};1;1^{2})$. 

\endproof

\begin{lemma}
\label{lem:caseBhat2}

If $\Gamma$ is of type $B_n^{(1)}(m,r)$, then Lemma~\ref{lem:xxx} holds.

\end{lemma}

\proof
In view of \eqref{eq:edge-3}, we may assume that the maximum edge-weight in $\Gamma k$ is $2$.

\credit{Case 1} \emph{$k$ is not connected to $c_1$.}

\credit{Subcase 1.1}\emph{ 
$k$ is contained in an edge whose weight is equal to 2}. 
Let us first note that, for any
vertex $v$ such that $k$ is connected to $v$, the edge $[k,v]$ is weighted 2.
Let $P$ be the shortest (undirected) path that connects $k$ to $c_1$ in $\Gamma k$.
Then $P$ induces a diagram which is of type $C_n^{(1)}$. 

\credit{Subcase 1.2}\emph{ 
$k$ is not contained in any edge whose weight is equal to 2}. 

\credit{Subsubcase 1.2.1}\emph{$k$ is connected to a vertex in $\{a_1,....,a_m\}$}.
Let $s=\max\{i: k  \mathrm{\:is \:connected \:to\:} a_i\}$. 
\begin{itemize}
\item
Suppose that $s>1$. We note that if $k$ is not connected to $a_{s-1}$,
then the subdiagram with the vertices in $\{a_{s-1},a_s,...,a_m,c_1,k\}$ is of type $B_n^{(1)}$. 
Let us now assume that $k$ is connected to $a_{s-1}$. 
If $k$ is not connected to any vertex $v$ such that $v\ne a_s,a_{s-1}$, 
then the diagram $\mu_k(\Gamma k)$ is of type $B_n^{(1)}(m,r)$. 
Let us now assume
that $k$ is connected to a vertex $v$ such that $v\ne a_s,a_{s-1}$, and
assume that $v'$ is such a vertex that is closest to $a_{s-1}$. 
Let $P$ be the shortest undirected path that connects $a_{s-1}$ to $v'$. 
Then the subdiagram with
$P \cup \{a_s,a_{s+1},...,a_m,c_1,k\}$ is of type $B_n^{(1)}(m,r)$. 
\item 
Suppose that $s=1$. 
If $k$ is not connected to $b_1$ (resp. $b_2$), then
the subdiagram with $\{k,b_1 (\mathrm{resp. \:b_2}),a_1,....,a_m,c_1\}$ is of type $B_n^{(1)}$. 
If k is connected to $b_1$ and $b_2$ then $k$ is contained in a non-oriented
cycle.  
\end{itemize}

\credit{Subsubcase 1.2.2}\emph{  k is not connected to any vertex in $\{a_1,....,a_m\}$}.
Let us write 
$$D=\{b_i: k \:\mathrm{is\: connected\: to} \:v\}=\{b_{i_1},...,b_{i_s}\}$$
where $i_1<...<i_s,s \geq 2$.
\begin{itemize}
\item[(i)]
Suppose that there exist two vertices
$b_i,b_j \in D$ such that $b_i$ and $b_j$ are not connected
to each other. (Note that this assumption
holds when $s>2$). Then the vertex $k$ is necessarily contained in a 
non-oriented cycle (Fig.~\ref{fig:Bhat2-121i}).
\item[(ii)]
Suppose that the condition of (i) does not hold. Then we have
$D=\{b_{i_1},b_{i_2}\}$ and the vertices $b_{i_1},b_{i_2}$ are 
connected to each other. If $D=\{b_1,b_2\}$, then 
the subdiagram with $\{k,b_1,b_2,a_1,....,a_m\}$ is of type $B_n^{(1)}(m,r)$ (Fig.~\ref{fig:Bhat2-121ii});
otherwise $\mu_k(\Gamma k)$ is of type $B_n^{(1)}(m,r)$.
\end{itemize}

\begin{figure}[ht]

\setlength{\unitlength}{1.8pt}

\begin{center}

\begin{picture}(120,80)(-80,-30)

\put(-60,0){\line(-1,0){10}}
\put(-60,0){\line(-1,1){20}}
\put(-60,0){\line(0,-1){20}}

\thicklines

\put(0,0){\circle*{2.0}$^{a_1}$}

\put(-20,20){\circle*{2.0}}
\put(-20,0){\circle*{2.0}} 
\put(20,0){\circle*{2.0}$^{a_2}$}
\put(40,0){\circle*{2.0}}
\put(60,0){\circle*{2.0}}
\put(80,0){\circle*{2.0}$\: ^{a_m}$}
\put(100,0){\circle*{2.0}$^{c_1}$}
\put(-20,20){\circle*{2.0}}
\put(-40,40){\circle*{2.0}}
\put(-60,40){\circle*{2.0}}
\put(-80,20){\circle*{2.0}}
\put(-60,-20){\circle*{2.0}}
\put(-80,0){\circle*{2.0}}
\put(-40,-20){\circle*{2.0}}

\put(-20,0){\line(1,0){120}}
\put(-20,0){\line(-1,-1){20}}
\put(-20,0){\line(0,1){20}}
\put(0,0){\line(-1,1){20}}
\put(-20,20){\line(-1,1){20}}
\put(-40,40){\line(-1,0){20}}
\put(-60,40){\line(-1,-1){20}}
\put(-80,20){\line(0,-1){20}}
\put(-80,0){\line(1,-1){20}}
\put(-60,-20){\line(1,0){20}}

\put(90,-3){\makebox(0,0){$2$}}
\put(-19,-4){\makebox(0,0){$b_1$}}
\put(-83,20){\makebox(0,0){$b_i$}}
\put(-19,23){\makebox(0,0){$b_2$}}
\put(-60,4){\makebox(0,0){$k$}}
\put(-63,-20){\makebox(0,0){$b_j$}}

\end{picture}
\end{center}

\caption{} 

\label{fig:Bhat2-121i}

\end{figure}

\begin{figure}[ht]

\setlength{\unitlength}{1.8pt}

\begin{center}

\begin{picture}(120,80)(-80,-30)

\put(-40,0){\line(1,0){20}}
\put(-40,0){\line(1,1){20}}

\thicklines

\put(0,0){\circle*{2.0}$^{a_1}$}

\put(-20,20){\circle*{2.0}}
\put(-20,0){\circle*{2.0}} 
\put(20,0){\circle*{2.0}$^{a_2}$}
\put(40,0){\circle*{2.0}}
\put(60,0){\circle*{2.0}}
\put(80,0){\circle*{2.0}$\: ^{a_m}$}
\put(100,0){\circle*{2.0}$^{c_1}$}
\put(-20,20){\circle*{2.0}}
\put(-40,40){\circle*{2.0}}
\put(-60,40){\circle*{2.0}}
\put(-80,20){\circle*{2.0}}
\put(-60,-20){\circle*{2.0}}
\put(-80,0){\circle*{2.0}}
\put(-40,-20){\circle*{2.0}}
\put(-40,0){\circle*{2.0}}

\put(-20,0){\line(1,0){120}}
\put(-20,0){\line(-1,-1){20}}
\put(-20,0){\line(0,1){20}}
\put(0,0){\line(-1,1){20}}
\put(-20,20){\line(-1,1){20}}
\put(-40,40){\line(-1,0){20}}
\put(-60,40){\line(-1,-1){20}}
\put(-80,20){\line(0,-1){20}}
\put(-80,0){\line(1,-1){20}}
\put(-60,-20){\line(1,0){20}}

\put(90,-3){\makebox(0,0){$2$}}
\put(-19,-4){\makebox(0,0){$b_1$}}
\put(-19,23){\makebox(0,0){$b_2$}}
\put(-42,3){\makebox(0,0){$k$}}

\end{picture}
\end{center}

\caption{} 

\label{fig:Bhat2-121ii}

\end{figure}


\credit{Case 2} \emph{k is connected to $c_1$.}

\credit{Subcase 2.1}\emph{The edge $[k,c_1]$ is weighted 2}.
We first note that the subdiagram with the vertices $\{b_1,b_2,a_1,...,a_m,k\}$ is weightless.
If $k$ is not connected to $a_m$ then
the subdiagram with the vertices $\{k,c_1,a_m\}$ is of $C_n^{(1)}$. 
Let us now assume that $k$ is connected
to $a_m$. If k is connected to any  vertex $v$ such that $v\ne a_m,c_1$,  then $k$ is contained in a diagram
of type $B_n^{(1)}(r)$
, otherwise $\mu_k(\Gamma k)$ is of type $B_n^{(1)}(m,r)$. 

\credit{Subcase 2.2} \emph{The edge $[k,c_1]$ is weightless}.
We note that, for any vertex $v\ne c_1$ such that $k$ is connected to $v$, the edge $[k,v]$ is weighted 2. 

\noindent
Let us first assume that k is not connected any vertex in $\{b_1,b_2,...,b_r\}$ and write
$s=\min\{i: k \mathrm{\:is \:connected \:to\:} a_i\}$. Then the subdiagram with the vertices $$\{k,a_s,...,a_1,b_1,b_2,...,b_r\}$$
is of type $B_n^{(1)}(m,r)$. 

\noindent
Let us now assume that 
$k$ is connected to a vertex in $\{b_1,b_2,...,b_r\}$ and write 
$D=\{b_i:k \mathrm{\:is \: connected\: to\:} b_i\}=\{b_{i_1},...,b_{i_s}\}$
where $1\leq i_1<i_2<...<i_s\leq r$.

\credit{Subsubcase 2.2.1} \emph{$s\geq 2$}. 
Then there exist two vertices
$b_i,b_j \in D$ such that $b_i$ and $b_j$ are not connected
to each other. Then the subdiagram with the vertices $\{k,b_i,b_j\}$ is of
type $C_n^{(1)}$.

\credit{Subsubcase 2.2.2} \emph{$s=2$.}
If $b_{i_1}$ and $b_{i_2}$ are not connected to each other, then
the subdiagram with the vertices $\{k,b_{i_1},b_{i_2}\}$ is of
type $C_n^{(1)}$; otherwise 
the subdiagram with the vertices $\{k,b_1,b_2,...,b_r\}$ is of type $B_n^{(1)}(r)$ 
(Fig.~\ref{fig:Bhat2-222}).

\credit{Subsubcase 2.2.3} \emph{$s=1$}
\begin{itemize}
\item
Suppose that $k$ is not connected to any vertex in $\{a_1,...,a_m\}$.
If $b_i$ is not connected to $a_m$, then the subdiagram with the vertices $\{b_i,k,c_1,a_m\}$
of type $C_n^{(1)}$. Let us now assume that $b_i$ is connected to $a_m$.
Then we have $m=1$ and the subdiagram with the vertices $\{k,b_1,b_2,a_1,c_1\}$ is of type $F_4^{(1)}(4^{1};3^1)$.
\item
Suppose that $k$ is connected to a vertex in $\{a_1,....,a_m\}$. Let us write $a_s=\min\{a_j:k \mathrm{\:is \:connected \:to\:} a_j\}$. If $b_i$ and $a_j$ are not
connected, then the subdiagram with the vertices $\{b_i,k,a_j\}$ is of type $C_n^{(1)}$,
otherwise the subdiagram with the vertices $\{b_1,b_2,k,a_m\}$ is of type $B_n^{(1)}(r)$
(Fig.~\ref{fig:Bhat2-223}). 
\end{itemize}

\begin{figure}[ht]

\setlength{\unitlength}{1.8pt}

\begin{center}

\begin{picture}(120,80)(-80,-30)

\put(-20,40){\line(-1,0){20}}
\put(-20,40){\line(0,-1){20}}
\put(-20,40){\line(3,-1){120}}
\put(-20,40){\line(2,-1){20}}
\put(-20,40){\line(1,-1){10}}

\thicklines

\put(0,0){\circle*{2.0}$^{a_1}$}

\put(-20,20){\circle*{2.0}}
\put(-20,0){\circle*{2.0}} 
\put(20,0){\circle*{2.0}$^{a_2}$}
\put(40,0){\circle*{2.0}}
\put(60,0){\circle*{2.0}}
\put(80,0){\circle*{2.0}$\: ^{a_m}$}
\put(100,0){\circle*{2.0}$^{c_1}$}
\put(-20,20){\circle*{2.0}}
\put(-40,40){\circle*{2.0}}
\put(-60,40){\circle*{2.0}}
\put(-80,20){\circle*{2.0}}
\put(-60,-20){\circle*{2.0}}
\put(-80,0){\circle*{2.0}}
\put(-40,-20){\circle*{2.0}}
\put(-20,40){\circle*{2.0}}

\put(-20,0){\line(1,0){120}}
\put(-20,0){\line(-1,-1){20}}
\put(-20,0){\line(0,1){20}}
\put(0,0){\line(-1,1){20}}
\put(-20,20){\line(-1,1){20}}
\put(-40,40){\line(-1,0){20}}
\put(-60,40){\line(-1,-1){20}}
\put(-80,20){\line(0,-1){20}}
\put(-80,0){\line(1,-1){20}}
\put(-60,-20){\line(1,0){20}}

\put(90,-3){\makebox(0,0){$2$}}
\put(-19,-4){\makebox(0,0){$b_1$}}
\put(-17,23){\makebox(0,0){$b_2$}}
\put(-20,44){\makebox(0,0){$k$}}
\put(-22,30){\makebox(0,0){$2$}}
\put(-30,43){\makebox(0,0){$2$}}

\end{picture}
\end{center}

\caption{} 

\label{fig:Bhat2-222}

\end{figure}

\begin{figure}[ht]

\setlength{\unitlength}{1.8pt}

\begin{center}

\begin{picture}(120,80)(-80,-30)

\put(0,20){\line(-1,0){20}}
\put(0,20){\line(0,-1){20}}
\put(0,20){\line(5,-1){100}}

\thicklines

\put(0,0){\circle*{2.0}$^{a_1}$}

\put(-20,20){\circle*{2.0}}
\put(-20,0){\circle*{2.0}} 
\put(20,0){\circle*{2.0}$^{a_2}$}
\put(40,0){\circle*{2.0}}
\put(60,0){\circle*{2.0}}
\put(80,0){\circle*{2.0}$\: ^{a_m}$}
\put(100,0){\circle*{2.0}$^{c_1}$}
\put(-20,20){\circle*{2.0}}
\put(-40,40){\circle*{2.0}}
\put(-60,40){\circle*{2.0}}
\put(-80,20){\circle*{2.0}}
\put(-60,-20){\circle*{2.0}}
\put(-80,0){\circle*{2.0}}
\put(-40,-20){\circle*{2.0}}
\put(-20,20){\circle*{2.0}}

\put(-20,0){\line(1,0){120}}
\put(-20,0){\line(-1,-1){20}}
\put(-20,0){\line(0,1){20}}
\put(0,0){\line(-1,1){20}}
\put(-20,20){\line(-1,1){20}}
\put(-40,40){\line(-1,0){20}}
\put(-60,40){\line(-1,-1){20}}
\put(-80,20){\line(0,-1){20}}
\put(-80,0){\line(1,-1){20}}
\put(-60,-20){\line(1,0){20}}

\put(90,-3){\makebox(0,0){$2$}}
\put(-19,-4){\makebox(0,0){$b_1$}}
\put(-17,23){\makebox(0,0){$b_2$}}
\put(0,24){\makebox(0,0){$k$}}
\put(1,10){\makebox(0,0){$2$}}
\put(-10,23){\makebox(0,0){$2$}}

\end{picture}
\end{center}

\caption{} 

\label{fig:Bhat2-223}

\end{figure}

\endproof

\begin{lemma}
\label{lem:caseBhat22}
If $\Gamma$ is of type $B_n^{(1)}(r)$, then Lemma~\ref{lem:xxx} holds.
\end{lemma}

\proof
In view of \eqref{eq:edge-3}, we may assume that the maximum edge-weight in $\Gamma k$ is $2$.
\credit{Case 1} \emph{$k$ is not connected to $c_1$.}

\credit{Subcase 1.1}\emph{ 
$k$ is contained in an edge whose weight is equal to 2}. 
In this case, for any
vertex $v$ connected to $k$, the edge $[k,v]$ is weighted 2.
Let $P$ be the shortest (undirected) path that connects $k$ to $c_1$ in $\Gamma k$.
Then $P$ induces a diagram which is of type $C_n^{(1)}$. 

\credit{Subcase 1.2}\emph{ 
$k$ is not contained in any edge whose weight is equal to 2}. 
Let us write 
$$D=\{b_i: k \mathrm{\:is \:connected \:to\:} v\}=\{b_{i_1},...,b_{i_s}\}$$
where $i_1<...<i_s$, $s\geq 2$.
\begin{itemize}
\item[(i)]
Suppose that there exist there exists two vertices
$b_i,b_j \in D$ such that $b_i$ and $b_j$ are not connected
to each other. (Note that this assumption holds when $s>2$). 
Then the vertex $k$ is contained in a non-oriented cycle.

\item[(ii)]
Suppose that (i) does not hold. Then, we have
$D=\{b_{i_1},b_{i_2}\}$ and the vertices $b_{i_1},b_{i_2}$ are 
connected to each other. If $D=\{b_1,b_2\}$ then 
the subdiagram with the vertices $\{k,b_1,b_2,c_1\}$ is of type $B_n^{(1)}(r)$,
otherwise $\mu_k(\Gamma k)$ is of type $B_n^{(1)}(r)$.
\end{itemize}
 

\credit{Case 2} \emph{$k$ is connected to $c_1$.}

\credit{Subcase 2.1}\emph{The edge $[k,c_1]$ is weighted 2}.
We first note that the subdiagram with the vertices $\{b_1,b_2,...,b_r,k\}$ is weightless. We also note the following:
\begin{itemize} 
\item[(i)]
If $k$ is connected to both $b_1$ and $b_2$, then $k$ is contained in a non-oriented cycle.
\item[(ii)]
Suppose that (i) does not hold. We may assume, without loss of generality, that $k$ is not connected to $b_1$.
We note that the subdiagram with the vertices 
$\{k,b_1,c_1\}$ is of type $C_n^{(1)}$.
\end{itemize}

\credit{Subcase 2.2} \emph{The edge $[k,c_1]$ is weightless}.
We note that, for any vertex $v\ne c_1$ such that k is connected to $v$, the
edge $[k,v]$ is weighted 2. 
Let us write 
$D=\{b_i:k \mathrm{\:is \:connected \:to\:} b_i\}=\{b_{i_1},...,b_{i_s}\}$
where $1\leq i_1<i_2<...<i_s\leq r$, $s \geq 1$.
We first assume that $k$ is connected to a vertex in $\{b_1,b_2\}$.
Then we have the following: if $k$ is connected to both $b_1$ and $b_2$, then $k$ is contained in a non-oriented cycle.
To consider the remaining possibilities,
let us now assume, without loss of generality, that
$k$ is connected to $b_2$ and not connected to $b_1$.
Then we have the following: if $r>3$, then the subdiagram with the vertices $\{b_1,b_2,b_3,k\}$ is of type 
$B_n^{(1)}$ (Fig.~\ref{fig:Bhat22-22}); otherwise $\mu_k(\Gamma k)$ is of type $F_4^{(1)}(4^{1};3^1)$.

\endproof

\begin{figure}[ht]

\setlength{\unitlength}{1.8pt}

\begin{center}

\begin{picture}(120,80)(-130,-30)

\put(0,20){\line(0,-1){20}}
\put(0,20){\line(-1,0){20}}

\thicklines

\put(0,0){\circle*{2.0}}
\put(0,20){\circle*{2.0}}
\put(-20,20){\circle*{2.0}}
\put(-20,0){\circle*{2.0}} 
\put(-20,20){\circle*{2.0}}
\put(-40,40){\circle*{2.0}}
\put(-60,40){\circle*{2.0}}
\put(-80,20){\circle*{2.0}}
\put(-60,-20){\circle*{2.0}}
\put(-80,0){\circle*{2.0}}
\put(-40,-20){\circle*{2.0}}

\put(-20,0){\line(1,0){20}}
\put(-20,0){\line(-1,-1){20}}
\put(-20,0){\line(0,1){20}}
\put(0,0){\line(-1,1){20}}
\put(-20,20){\line(-1,1){20}}
\put(-40,40){\line(-1,0){20}}
\put(-60,40){\line(-1,-1){20}}
\put(-80,20){\line(0,-1){20}}
\put(-80,0){\line(1,-1){20}}
\put(-60,-20){\line(1,0){20}}

\put(-10,-3){\makebox(0,0){$2$}}
\put(-9,12){\makebox(0,0){$2$}}
\put(-19,-4){\makebox(0,0){$b_1$}}
\put(-40,43){\makebox(0,0){$b_3$}}
\put(-83,0){\makebox(0,0){$b_i$}}
\put(-40,-23){\makebox(0,0){$b_r$}}
\put(-19,23){\makebox(0,0){$b_2$}}
\put(-10,23){\makebox(0,0){$2$}}
\put(0,23){\makebox(0,0){$k$}}
\put(4,0){\makebox(0,0){$c_1$}}

\end{picture}
\end{center}

\caption{} 

\label{fig:Bhat22-22}

\end{figure}

\begin{lemma}
\label{lem:caseI2}

If $\Gamma$ is of type $I_2(a), a\geq 4$, then Lemma~\ref{lem:xxx} holds.

\end{lemma}

We leave the proof the lemma to the reader as an easy exercise.

\subsection{The diagram $\Gamma$ is a (non-simply-laced) diagram in Table 2}
\label{sec:pflemma8}
For such $\Gamma$ and $\Gamma k$, one could easily check that Lemma~\ref{lem:xxx} holds.

\subsection{The diagram $\Gamma$ is a (non-simply-laced) diagram in Table 3}
\label{sec:pflemma9}


In view of \eqref{eq:edge-3}, we only need to consider the diagrams which do not contain any edge whose weight is equal to $3$.

\begin{lemma}
\label{lem:caseFhat}

If $\Gamma$ is of type $F_4^{(1)}$, then 
Lemma~\ref{lem:xxx} holds.

\end{lemma}

\begin{figure}[ht]

\setlength{\unitlength}{1.8pt}

\begin{center}

\begin{picture}(105,30)(-10,-10)

\thicklines

\put(10,0){\circle*{2.0}$^{a_1}$}
\put(30,0){\circle*{2.0}$^{a_2}$}
\put(50,0){\circle*{2.0}$^{b_1}$}
\put(70,0){\circle*{2.0}$^{b_2}$}
\put(90,0){\circle*{2.0}$^{c_1}$}

\put(60,-3){\makebox(0,0){$2$}}

\put(10,0){\line(1,0){80}}

\end{picture}

\end{center}

\caption{The extended Dynkin diagram $F_4^{(1)} \,$.} 

\label{fig:Fhat}

\end{figure}
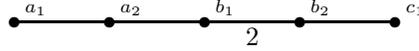

\proof
Throughout the proof, we assume that $\Gamma$ is indexed as in Fig.~\ref{fig:Fhat}.

\credit{Case 1} \emph{$k$ is not connected to any vertex in $\{b_1,b_2\}$.}

\credit{Subcase 1.1}\emph{$k$ is not connected to $c_1$.}

Let us note that the vertex $k$ is connected to both $a_1$ and $a_2$.
If the edge $[k,a_2]$ (hence the edge $[k,a_1]$) is weighted 2
, then the subdiagram with the vertices $\{k,a_2,b_1,b_2\}$ is of type $C_n^{(1)}$ 
(Fig.~\ref{fig:Fhat-11})
; otherwise the subdiagram with the vertices $\{k,a_2,b_1,b_2,c_1\}$ is of type $F_4^{(1)}$. 

\begin{figure}[ht]

\setlength{\unitlength}{1.8pt}

\begin{center}

\begin{picture}(70,40)(-10,-20)

\put(20,-20){\line(-1,2){10}}
\put(20,-20){\line(1,2){10}}

\thicklines

\put(20,-20){\circle*{2.0}}
\put(10,0){\circle*{2.0}$^{a_1}$}
\put(30,0){\circle*{2.0}$^{a_2}$}
\put(50,0){\circle*{2.0}$^{b_1}$}
\put(70,0){\circle*{2.0}$^{b_2}$}
\put(90,0){\circle*{2.0}$^{c_1}$}

\put(60,-3){\makebox(0,0){$2$}}
\put(20,-23){\makebox(0,0){$k$}}
\put(27,-10){\makebox(0,0){$2$}}
\put(13,-10){\makebox(0,0){$2$}}

\put(10,0){\line(1,0){80}}

\end{picture}

\end{center}

\caption{} 

\label{fig:Fhat-11}

\end{figure}

\credit{Subcase 1.2}\emph{$k$ is connected to $c_1$.}

\begin{itemize}
\item 
If the edge $[k,c_1]$ is weighted 2,
then the subdiagram with the vertices $\{k,c_1,b_2,b_1\}$ is of type $C_n^{(1)}$ (Fig.~\ref{fig:Fhat-12i}).
 
\item
Suppose that the edge $[k,c_1]$ is weightless. 
If $k$ is connected to $a_2$, then
the subdiagram with the vertices $\{k,a_2,b_1,b_2\}$ is of type $C_n^{(1)}$;
otherwise the subdiagram with the vertices $\{k,a_1,a_2,b_1,b_2\}$ is of type $C_n^{(1)}$ (Fig.~\ref{fig:Fhat-12ii}). 
\end{itemize}

\begin{figure}[ht]
\setlength{\unitlength}{1.8pt}

\begin{center}

\begin{picture}(70,40)(-10,-20)

\put(70,-20){\line(-3,1){10}}
\put(70,-20){\line(1,1){20}}

\thicklines

\put(70,-20){\circle*{2.0}}
\put(10,0){\circle*{2.0}$^{a_1}$}
\put(30,0){\circle*{2.0}$^{a_2}$}
\put(50,0){\circle*{2.0}$^{b_1}$}
\put(70,0){\circle*{2.0}$^{b_2}$}
\put(90,0){\circle*{2.0}$^{c_1}$}

\put(60,-2){\makebox(0,0){$2$}}
\put(70,-23){\makebox(0,0){$k$}}
\put(83,-10){\makebox(0,0){$2$}}

\put(10,0){\line(1,0){80}}

\end{picture}

\end{center}

\caption{} 

\label{fig:Fhat-12i}

\end{figure}

\begin{figure}[ht]

\setlength{\unitlength}{1.8pt}

\begin{center}

\begin{picture}(70,40)(-10,-20)

\put(70,-20){\line(-3,1){60}}
\put(70,-20){\line(1,1){20}}

\thicklines

\put(70,-20){\circle*{2.0}}
\put(10,0){\circle*{2.0}$^{a_1}$}
\put(30,0){\circle*{2.0}$^{a_2}$}
\put(50,0){\circle*{2.0}$^{b_1}$}
\put(70,0){\circle*{2.0}$^{b_2}$}
\put(90,0){\circle*{2.0}$^{c_1}$}

\put(60,-2){\makebox(0,0){$2$}}
\put(70,-23){\makebox(0,0){$k$}}
\put(33,-10){\makebox(0,0){$2$}}

\put(10,0){\line(1,0){80}}

\end{picture}

\end{center}

\caption{} 

\label{fig:Fhat-12ii}

\end{figure}

\credit{Case 2} \emph{$k$ is connected to both $b_1$ and $b_2$.}

\credit{Subcase 2.1}\emph{The edge $[k,b_1]$ is weighted 2.}

We first note that the edge $[k,b_2]$ is weightless.
We also note that if $k$ is connected to a vertex $v$ in $\{a_1,a_2\}$, 
then the edge $[k,v]$ is weighted 2.

\begin{itemize}
\item
Suppose that $k$ is connected to $a_1$. Then the subdiagram with the vertices
$\{a_1,k,b_1\}$ is of type $C_n^{(1)}$ 
(Fig.~\ref{fig:Fhat-21i}). 

\item
Suppose that $k$ is connected to $c_1$. Then the subdiagram with the vertices
$\{k,b_1,b_2,c_1\}$ is of type $B_n^{(1)}(r)$ 
(Fig.~\ref{fig:Fhat-21ii}). 

\item
Suppose that $k$ is not connected to any vertex in $\{a_1,c_1\}$.
If $k$ is connected to $a_2$, then the subdiagram with the vertices 
$\{a_1,a_2,k,b_2,c_1\}$ is of type $F_4^{(1)}$ 
; otherwise
the subdiagram with the vertices $\{a_2,b_1,b_2,k,c_1\}$ is of type $F_4^{(1)}(3^{1};1;1^2)$.  
\end{itemize}

\begin{figure}[ht]

\setlength{\unitlength}{1.8pt}

\begin{center}

\begin{picture}(70,40)(-10,-20)

\put(60,-20){\line(-5,2){50}}
\put(60,-20){\line(1,2){10}}
\put(60,-20){\line(-1,2){10}}
\put(60,-20){\line(3,2){10}}

\thicklines

\put(60,-20){\circle*{2.0}}
\put(10,0){\circle*{2.0}$^{a_1}$}
\put(30,0){\circle*{2.0}$^{a_2}$}
\put(50,0){\circle*{2.0}$^{b_1}$}
\put(70,0){\circle*{2.0}$^{b_2}$}
\put(90,0){\circle*{2.0}$^{c_1}$}

\put(60,-2){\makebox(0,0){$2$}}
\put(60,-23){\makebox(0,0){$k$}}
\put(53,-10){\makebox(0,0){$2$}}
\put(30,-10){\makebox(0,0){$2$}}

\put(10,0){\line(1,0){80}}

\end{picture}

\end{center}

\caption{} 

\label{fig:Fhat-21i}

\end{figure}

\begin{figure}[ht]

\setlength{\unitlength}{1.8pt}

\begin{center}

\begin{picture}(70,40)(-10,-20)

\put(60,-20){\line(-5,2){10}}
\put(60,-20){\line(1,2){10}}
\put(60,-20){\line(-1,2){10}}
\put(60,-20){\line(3,2){30}}

\thicklines

\put(60,-20){\circle*{2.0}}
\put(10,0){\circle*{2.0}$^{a_1}$}
\put(30,0){\circle*{2.0}$^{a_2}$}
\put(50,0){\circle*{2.0}$^{b_1}$}
\put(70,0){\circle*{2.0}$^{b_2}$}
\put(90,0){\circle*{2.0}$^{c_1}$}

\put(60,-2){\makebox(0,0){$2$}}
\put(60,-23){\makebox(0,0){$k$}}
\put(53,-10){\makebox(0,0){$2$}}
\put(30,-10){\makebox(0,0){$2$}}

\put(10,0){\line(1,0){80}}

\end{picture}

\end{center}

\caption{} 

\label{fig:Fhat-21ii}

\end{figure}

\credit{Subcase 2.2}\emph{The edge $[k,b_2]$ is weighted 2.}

We note that the edge $[k,b_1]$ is weightless. We also note that if 
$k$ is connected to $c_1$, the edge $[k,c_1]$ is weighted 2.
Let us now consider the subcases.

\begin{itemize}
\item
Suppose that $k$ is connected to $a_1$ or $a_2$.
If $k$ is connected to $a_2$, then the edge $[k,a_2]$ is weightless and the
subdiagram with the vertices $\{a_2,k,b_1,b_2\}$
is of type $B_n^{(1)}(r)$, otherwise the subdiagram with the vertices $\{a_1,a_2,b_1,b_2,k\}$
is of type $B_n^{(1)}(r)$ (Fig.~\ref{fig:Fhat-22i}).
\item
Suppose that $k$ is not connected to any vertex in $\{a_1,a_2\}$.
If $k$ is connected to $c_1$, then the subdiagram with the vertices 
$\{a_2,b_1,k,b_2,c_1\}$ of type $F_4^{(1)}(3^{2};0;1^1)$. 
If $k$ is not connected to $c_1$, then the subdiagram with the vertices 
$\{a_2,b_1,k,b_2,c_1\}$ of type $F_4^{(1)}(3^{1};1;1^2)$. 
\end{itemize}

\begin{figure}[ht]

\setlength{\unitlength}{1.8pt}

\begin{center}

\begin{picture}(70,40)(-10,-20)

\put(60,-20){\line(-5,2){50}}
\put(60,-20){\line(1,2){10}}
\put(60,-20){\line(-1,2){10}}
\put(60,-20){\line(3,2){10}}

\thicklines

\put(60,-20){\circle*{2.0}}
\put(10,0){\circle*{2.0}$^{a_1}$}
\put(30,0){\circle*{2.0}$^{a_2}$}
\put(50,0){\circle*{2.0}$^{b_1}$}
\put(70,0){\circle*{2.0}$^{b_2}$}
\put(90,0){\circle*{2.0}$^{c_1}$}

\put(60,-2){\makebox(0,0){$2$}}
\put(60,-23){\makebox(0,0){$k$}}
\put(69,-10){\makebox(0,0){$2$}}

\put(10,0){\line(1,0){80}}

\end{picture}

\end{center}

\caption{} 

\label{fig:Fhat-22i}

\end{figure}

\credit{Case 3} \emph{$k$ is connected to $b_1$ not $b_2$.}

We note that if the edge $[k,b_1]$ is weighted 2,
then the subdiagram with the vertices $\{k,b_1,b_2\}$ is of type $C_n^{(1)}$
(Fig.~\ref{fig:Fhat-3}).
To consider the remaining subcases we assume that the edge $[k,b_1]$
is weightless.

\credit{Subcase 3.1}\emph{$k$ is not connected to $a_2$.}

The subdiagram with the vertices 
$\{k,b_1,a_2,b_2\}$ is of type $B_n^{(1)}$ (Fig.~\ref{fig:Fhat-31}).

\credit{Subcase 3.2}\emph{k is connected to $a_2$.}

Let $X$ denote the subdiagram with the vertices 
$\{k,a_2,b_1,b_2,c_1\}$.
If $k$ is not connected to $c_1$, then $X$ is of type $F_4^{(1)}(3^{1};1;2^{1})_1$.
If $k$ is connected to $c_1$, then $X$ is of type $F_4^{(1)}(4^{1};3^1)$.

\begin{figure}[ht]

\setlength{\unitlength}{1.8pt}

\begin{center}

\begin{picture}(70,40)(-10,-20)

\put(60,-20){\line(-1,2){10}}
\put(60,-20){\line(2,1){10}}
\put(60,-20){\line(-2,1){10}}

\thicklines

\put(60,-20){\circle*{2.0}}
\put(10,0){\circle*{2.0}$^{a_1}$}
\put(30,0){\circle*{2.0}$^{a_2}$}
\put(50,0){\circle*{2.0}$^{b_1}$}
\put(70,0){\circle*{2.0}$^{b_2}$}
\put(90,0){\circle*{2.0}$^{c_1}$}

\put(60,-2){\makebox(0,0){$2$}}
\put(60,-23){\makebox(0,0){$k$}}
\put(53,-10){\makebox(0,0){$2$}}

\put(10,0){\line(1,0){80}}

\end{picture}

\end{center}

\caption{} 

\label{fig:Fhat-3}

\end{figure}

\begin{figure}[ht]

\setlength{\unitlength}{1.8pt}

\begin{center}

\begin{picture}(70,40)(-10,-20)

\put(60,-20){\line(-1,2){10}}
\put(60,-20){\line(2,1){10}}
\put(60,-20){\line(-5,1){20}}

\thicklines

\put(60,-20){\circle*{2.0}}
\put(10,0){\circle*{2.0}$^{a_1}$}
\put(30,0){\circle*{2.0}$^{a_2}$}
\put(50,0){\circle*{2.0}$^{b_1}$}
\put(70,0){\circle*{2.0}$^{b_2}$}
\put(90,0){\circle*{2.0}$^{c_1}$}

\put(60,-2){\makebox(0,0){$2$}}
\put(60,-23){\makebox(0,0){$k$}}

\put(10,0){\line(1,0){80}}

\end{picture}

\end{center}

\caption{} 

\label{fig:Fhat-31}

\end{figure}

\credit{Case 4} \emph{$k$ is connected to $b_2$ not $b_1$.}

We note that if the edge $[k,b_2]$ is weighted 2,
then the subdiagram with the vertices $\{k,b_1,b_2\}$ is of type $C_n^{(1)}$.
To consider the remaining subcases we assume that the edge $[k,b_2]$
is weightless.

\credit{Subcase 4.1}\emph{$k$ is not connected to $a_2$.}

\begin{itemize}
\item
If $k$ is connected to $a_1$ then the subdiagram with the vertices 
$\{a_1,k,b_2,b_1\}$ is of type $C_n^{(1)}$
(Fig~\ref{fig:Fhat-41i}).

\item
If $k$ is not connected to $a_1$,
then $k$ is connected to $c_1$ and
the subdiagram with the vertices 
$\{a_2,b_1,b_2,c_1,k\}$ is of type $F_4^{(1)}(3^{1};1;2^{1})_1$ 
(Fig.~\ref{fig:Fhat-41ii}).
\end{itemize}

\credit{Subcase 4.2}\emph{$k$ is connected to $a_2$}

We note that the edge $[k,a_2]$ is weighted 2.

\begin{itemize}
\item
Supose that $k$ is connected to $a_1$.
Then the subdiagram with the vertices 
$\{a_1,k,b_2,b_1\}$ is of type $B_n^{(1)}$. 

\item
Suppose that $k$ is not connected to $a_1$.
Then the subdiagram with the vertices 
$\{a_1,a_2,k,b_1\}$ is of type $C_n^{(1)}$
(Fig.~\ref{fig:Fhat-42ii}).
\end{itemize}

\endproof

\begin{figure}[ht]

\setlength{\unitlength}{1.8pt}

\begin{center}

\begin{picture}(70,40)(-10,-20)

\put(60,-20){\line(-5,2){50}}
\put(60,-20){\line(1,2){10}}
\put(60,-20){\line(1,1){10}}

\thicklines

\put(60,-20){\circle*{2.0}}
\put(10,0){\circle*{2.0}$^{a_1}$}
\put(30,0){\circle*{2.0}$^{a_2}$}
\put(50,0){\circle*{2.0}$^{b_1}$}
\put(70,0){\circle*{2.0}$^{b_2}$}
\put(90,0){\circle*{2.0}$^{c_1}$}

\put(60,-2){\makebox(0,0){$2$}}
\put(60,-23){\makebox(0,0){$k$}}
\put(30,-10){\makebox(0,0){$2$}}

\put(10,0){\line(1,0){80}}

\end{picture}

\end{center}

\caption{} 

\label{fig:Fhat-41i}

\end{figure}

\begin{figure}[ht]

\setlength{\unitlength}{1.8pt}

\begin{center}

\begin{picture}(70,40)(-10,-20)

\put(80,-20){\line(-1,2){10}}
\put(80,-20){\line(1,2){10}}

\thicklines

\put(80,-20){\circle*{2.0}}
\put(10,0){\circle*{2.0}$^{a_1}$}
\put(30,0){\circle*{2.0}$^{a_2}$}
\put(50,0){\circle*{2.0}$^{b_1}$}
\put(70,0){\circle*{2.0}$^{b_2}$}
\put(90,0){\circle*{2.0}$^{c_1}$}

\put(60,-2){\makebox(0,0){$2$}}
\put(80,-23){\makebox(0,0){$k$}}

\put(10,0){\line(1,0){80}}

\end{picture}

\end{center}

\caption{} 

\label{fig:Fhat-41ii}

\end{figure}

\begin{figure}[ht]

\setlength{\unitlength}{1.8pt}

\begin{center}

\begin{picture}(70,40)(-10,-20)

\put(50,-20){\line(-1,1){20}}
\put(50,-20){\line(1,1){20}}
\put(50,-20){\line(2,1){10}}

\thicklines

\put(50,-20){\circle*{2.0}}
\put(10,0){\circle*{2.0}$^{a_1}$}
\put(30,0){\circle*{2.0}$^{a_2}$}
\put(50,0){\circle*{2.0}$^{b_1}$}
\put(70,0){\circle*{2.0}$^{b_2}$}
\put(90,0){\circle*{2.0}$^{c_1}$}

\put(60,-2){\makebox(0,0){$2$}}
\put(50,-23){\makebox(0,0){$k$}}
\put(38,-10){\makebox(0,0){$2$}}

\put(10,0){\line(1,0){80}}

\end{picture}

\end{center}

\caption{} 

\label{fig:Fhat-42ii}

\end{figure}

\begin{lemma}
\label{lem:caseFhat11}

If $\Gamma$ is of type $F_4^{(1)}(3^{1};1;2^{1})_1$, then Lemma~\ref{lem:xxx} holds.

\end{lemma}

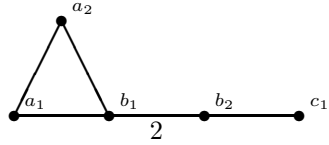
\begin{figure}[ht]

\setlength{\unitlength}{1.8pt}

\begin{center}

\begin{picture}(85,40)(-10,-10)

\thicklines

\put(10,0){\circle*{2.0}$^{a_1}$}
\put(30,0){\circle*{2.0}$^{b_1}$}
\put(50,0){\circle*{2.0}$^{b_2}$}
\put(70,0){\circle*{2.0}$^{c_1}$}
\put(20,20){\circle*{2.0}$^{a_2}$}

\put(40,-3){\makebox(0,0){$2$}}

\put(10,0){\line(1,0){60}}
\put(20,20){\line(-1,-2){10}}
\put(20,20){\line(1,-2){10}}

\end{picture}

\end{center}

\caption{The diagram $F_4^{(1)}(3^{1};1;2^{1})_1$}

\label{fig:Fhat11}

\end{figure}

\proof
Let us assume that $\Gamma$ is indexed as in Fig.~\ref{fig:Fhat11}.

\credit{Case 1} \emph{$k$ is not connected to any vertex in $\{b_1,b_2\}$.}

\credit{Subcase 1.1}\emph{$k$ is not connected to $c_1$.}

In this case the vertex $k$ is connected to both $a_1$ and $a_2$.
If the edge $[k,a_1]$ is weighted 2, then the 
subdiagram with the vertices $\{k,a_1,b_1,b_2\}$ is of type $C_n^{(1)}$; 
otherwise the 
subdiagram with the vertices $\{k,a_1,b_1,b_2,c_1\}$ is of type $F_4^{(1)}$ 
(Fig.~\ref{fig:Fhat11-11}).

\begin{figure}[ht]

\setlength{\unitlength}{1.8pt}

\begin{center}

\begin{picture}(70,40)(-10,-10)

\put(0,20){\line(1,-2){10}}
\put(0,20){\line(1,0){20}}

\thicklines

\put(0,20){\circle*{2.0}}
\put(10,0){\circle*{2.0}$^{a_1}$}
\put(30,0){\circle*{2.0}$^{b_1}$}
\put(50,0){\circle*{2.0}$^{b_2}$}
\put(70,0){\circle*{2.0}$^{c_1}$}
\put(20,20){\circle*{2.0}$^{a_2}$}

\put(40,-3){\makebox(0,0){$2$}}
\put(0,23){\makebox(0,0){$k$}}
\put(10,23){\makebox(0,0){($2$)}}
\put(2,10){\makebox(0,0){($2$)}}

\put(10,0){\line(1,0){60}}
\put(20,20){\line(-1,-2){10}}
\put(20,20){\line(1,-2){10}}

\end{picture}

\end{center}

\caption{} 

\label{fig:Fhat11-11}

\end{figure}

\credit{Subcase 1.2}\emph{$k$ is connected to $c_1$.}

\begin{itemize}
\item
If the edge $[k,c_1]$ is weighted 2, then the subdiagram with the vertices 
$\{k,c_1,b_2,b_1\}$ is of type $C_n^{(1)}$ (Fig.\ref{fig:Fhat11-12i}).
 
\item
Suppose that the edge $[k,c_1]$ is weightless. 
 Let us assume, without loss of generality, that
$k$ is connected to $a_2$. Then the edge $[k,a_2]$ is weighted 2
and the subdiagram with the vertices $\{k,a_2,b_1,b_2\}$ is of type $C_n^{(1)}$
(Fig.\ref{fig:Fhat11-12ii}).
\end{itemize}

\begin{figure}[ht]

\setlength{\unitlength}{1.8pt}

\begin{center}

\begin{picture}(70,40)(-10,-10)

\put(50,20){\line(1,-1){20}}
\put(50,20){\line(-1,0){10}}

\thicklines

\put(50,20){\circle*{2.0}}
\put(10,0){\circle*{2.0}$^{a_1}$}
\put(30,0){\circle*{2.0}$^{b_1}$}
\put(50,0){\circle*{2.0}$^{b_2}$}
\put(70,0){\circle*{2.0}$^{c_1}$}
\put(20,20){\circle*{2.0}$^{a_2}$}

\put(40,-3){\makebox(0,0){$2$}}
\put(60,23){\makebox(0,0){$k$}}
\put(62,10){\makebox(0,0){$2$}}

\put(10,0){\line(1,0){60}}
\put(20,20){\line(-1,-2){10}}
\put(20,20){\line(1,-2){10}}

\end{picture}

\end{center}

\caption{} 

\label{fig:Fhat11-12i}

\end{figure}

\begin{figure}[ht]

\setlength{\unitlength}{1.8pt}

\begin{center}

\begin{picture}(70,40)(-10,-10)

\put(50,20){\line(1,-1){20}}
\put(50,20){\line(-1,0){30}}

\thicklines

\put(50,20){\circle*{2.0}}
\put(10,0){\circle*{2.0}$^{a_1}$}
\put(30,0){\circle*{2.0}$^{b_1}$}
\put(50,0){\circle*{2.0}$^{b_2}$}
\put(70,0){\circle*{2.0}$^{c_1}$}
\put(20,20){\circle*{2.0}$^{a_2}$}

\put(40,-3){\makebox(0,0){$2$}}
\put(60,23){\makebox(0,0){$k$}}
\put(40,22){\makebox(0,0){$2$}}

\put(10,0){\line(1,0){60}}
\put(20,20){\line(-1,-2){10}}
\put(20,20){\line(1,-2){10}}

\end{picture}

\end{center}

\caption{} 

\label{fig:Fhat11-12ii}

\end{figure}

\credit{Case 2} \emph{$k$ is connected to both $b_1$ and $b_2$.}

\credit{Subcase 2.1}\emph{The edge $[k,b_1]$ is weighted 2.}

We first note that the edge $[k,b_2]$ is weightless. 

\begin{itemize}
\item[(i)] 
If $k$ is connected to $c_1$, then the subdiagram with the vertices
$\{k,b_1,b_2,c_1\}$ is of type $B_n^{(1)}(r)$ (Fig.\ref{fig:Fhat11-21i}).

\item[(ii)] If $k$ is connected to both $a_1$ and $a_2$, then
$k$ is contained in a non-oriented cycle

\item[(iii)] Suppose that neither (i) nor (ii) holds.
We may also assume, without loss of generality, that $k$ is
connected to $a_2$. 
Then, e.g., the subdiagram with the vertices
$\{k,b_1,a_1,a_2\}$ is of type $B_n^{(1)}(r)$ (Fig.\ref{fig:Fhat11-21iii}).
\end{itemize}

\begin{figure}[ht]

\setlength{\unitlength}{1.8pt}

\begin{center}

\begin{picture}(70,40)(-10,-10)

\put(50,20){\line(1,-1){20}}
\put(50,20){\line(-1,-1){20}}
\put(50,20){\line(0,-1){20}}
\put(50,20){\line(-1,0){10}}

\thicklines

\put(50,20){\circle*{2.0}}
\put(10,0){\circle*{2.0}}
\put(30,0){\circle*{2.0}}
\put(50,0){\circle*{2.0}}
\put(70,0){\circle*{2.0}}
\put(20,20){\circle*{2.0}}

\put(40,-3){\makebox(0,0){$2$}}
\put(50,23){\makebox(0,0){$k$}}
\put(40,12){\makebox(0,0){$2$}}
\put(30,-3){\makebox(0,0){$b_1$}}
\put(50,-3){\makebox(0,0){$b_2$}}
\put(70,-3){\makebox(0,0){$c_1$}}

\put(10,0){\line(1,0){60}}
\put(20,20){\line(-1,-2){10}}
\put(20,20){\line(1,-2){10}}

\end{picture}

\end{center}

\caption{} 

\label{fig:Fhat11-21i}

\end{figure}

\begin{figure}[ht]

\setlength{\unitlength}{1.8pt}

\begin{center}

\begin{picture}(70,40)(-10,-10)

\put(50,20){\line(-1,-1){20}}
\put(50,20){\line(0,-1){20}}
\put(50,20){\line(-1,0){30}}

\thicklines

\put(50,20){\circle*{2.0}}
\put(10,0){\circle*{2.0}}
\put(30,0){\circle*{2.0}}
\put(50,0){\circle*{2.0}}
\put(70,0){\circle*{2.0}}
\put(20,20){\circle*{2.0}}

\put(40,-3){\makebox(0,0){$2$}}
\put(50,23){\makebox(0,0){$k$}}
\put(30,22){\makebox(0,0){$2$}}
\put(40,12){\makebox(0,0){$2$}}
\put(30,-3){\makebox(0,0){$b_1$}}
\put(50,-3){\makebox(0,0){$b_2$}}
\put(70,-3){\makebox(0,0){$c_1$}}
\put(20,23){\makebox(0,0){$a_2$}}
\put(10,-3){\makebox(0,0){$a_1$}}

\put(10,0){\line(1,0){60}}
\put(20,20){\line(-1,-2){10}}
\put(20,20){\line(1,-2){10}}

\end{picture}

\end{center}

\caption{} 

\label{fig:Fhat11-21iii}

\end{figure}

\credit{Subcase 2.2}\emph{The edge $[k,b_2]$ is weighted 2.}

Let us consider the subcases.

\begin{itemize}
\item
If $k$ is connected to both $a_1$ and $a_2$, then
$k$ is contained in a non-oriented cycle

\item
Suppose that $k$ is connected to precisely one vertex, say $a_2$, in $\{a_1,a_2\}$.
Then the subdiagram with the vertices $\{k,a_2,b_1,b_2\}$
is of type $B_n^{(1)}(r)$ (Fig.~\ref{fig:Fhat11-22ii}).
\item
Suppose that $k$ is not connected to any vertex in $\{a_1,a_2\}$.
Let $D$ denote the subdiagram with the vertices $\{a_1,b_1,k,b_2,c_1\}$.
If $k$ is not connected to $c_1$, then $D$ is of type 
$F_4^{(1)}(3^{1};1;1^{2})$
, otherwise it is of type 
$F_4^{(1)}(3^{2};0;1^1)$. 
\end{itemize}

\begin{figure}[ht]

\setlength{\unitlength}{1.8pt}

\begin{center}

\begin{picture}(70,40)(-10,-10)

\put(50,20){\line(1,-1){10}}
\put(50,20){\line(-1,-1){20}}
\put(50,20){\line(0,-1){20}}
\put(50,20){\line(-1,0){30}}

\thicklines

\put(50,20){\circle*{2.0}}
\put(10,0){\circle*{2.0}}
\put(30,0){\circle*{2.0}}
\put(50,0){\circle*{2.0}}
\put(70,0){\circle*{2.0}}
\put(20,20){\circle*{2.0}}

\put(40,-3){\makebox(0,0){$2$}}
\put(50,23){\makebox(0,0){$k$}}
\put(52,10){\makebox(0,0){$2$}}
\put(30,-3){\makebox(0,0){$b_1$}}
\put(50,-3){\makebox(0,0){$b_2$}}
\put(70,-3){\makebox(0,0){$c_1$}}
\put(20,23){\makebox(0,0){$a_2$}}
\put(10,-3){\makebox(0,0){$a_1$}}

\put(10,0){\line(1,0){60}}
\put(20,20){\line(-1,-2){10}}
\put(20,20){\line(1,-2){10}}

\end{picture}

\end{center}

\caption{} 

\label{fig:Fhat11-22ii}

\end{figure}

\credit{Case 3} \emph{$k$ is connected to $b_1$ not $b_2$.}

Let us note that if the edge $[k,b_1]$ is weighted 2,
then the subdiagram with the vertices $\{k,b_1,b_2\}$ is of type $C_n^{(1)}$.
To consider the remaining subcases, we assume that the edge $[k,b_1]$
is weightless.

\credit{Subcase 3.1}\emph{k is not connected to $a_1$ or not connected to $a_2$.}

Let us assume, without loss of generality, that $k$ is not connected to $a_1$. 
Then the subdiagram with the vertices 
$\{k,b_1,a_1,b_2\}$ is of type $B_n^{(1)}$.

\begin{figure}[ht]

\setlength{\unitlength}{1.8pt}

\begin{center}

\begin{picture}(70,40)(-10,-10)

\put(50,20){\line(1,-1){10}}
\put(50,20){\line(-1,-1){20}}
\put(50,20){\line(-1,0){10}}

\thicklines

\put(50,20){\circle*{2.0}}
\put(10,0){\circle*{2.0}}
\put(30,0){\circle*{2.0}}
\put(50,0){\circle*{2.0}}
\put(70,0){\circle*{2.0}}
\put(20,20){\circle*{2.0}}

\put(40,-3){\makebox(0,0){$2$}}
\put(50,23){\makebox(0,0){$k$}}
\put(30,-3){\makebox(0,0){$b_1$}}
\put(50,-3){\makebox(0,0){$b_2$}}
\put(70,-3){\makebox(0,0){$c_1$}}
\put(20,23){\makebox(0,0){$a_2$}}
\put(10,-3){\makebox(0,0){$a_1$}}

\put(10,0){\line(1,0){60}}
\put(20,20){\line(-1,-2){10}}
\put(20,20){\line(1,-2){10}}

\end{picture}

\end{center}

\caption{} 

\label{fig:Fhat11-31}

\end{figure}

\credit{Subcase 3.2}\emph{$k$ is connected to both $a_1$ and $a_2$}

Then $k$ is contained in a non-oriented cycle.

\credit{Case 4} \emph{$k$ is connected to $b_2$ not $b_1$.}

We note that if the edge $[k,b_2]$ is weighted 2,
then the subdiagram with the vertices $\{k,b_1,b_2\}$ is of type $C_n^{(1)}$.
To consider the remaining subcases ,we assume that the edge $[k,b_2]$
is weightless.

\credit{Subcase 4.1}\emph{$k$ is connected to a vertex in $\{a_1,a_2\}$.}

\begin{itemize}
\item
Suppose that $k$ is connected to both $a_1$ and $a_2$.
Then the subdiagram with the vertices 
$\{k,a_1,a_2,b_1\}$ is of type $B_n^{(1)}(r)$. 
\item
Suppose that $k$ is connected to precisely one vertex, say $a_2$, in $\{a_1,a_2\}$.
Then the subdiagram with the vertices 
$\{a_2,a_1,b_1,b_2,k\}$ is of type $F_4^{(1)}(4^{1};3^1)$
(Fig.~\ref{fig:Fhat11-41}).
\end{itemize}

\begin{figure}[ht]

\setlength{\unitlength}{1.8pt}

\begin{center}

\begin{picture}(70,40)(-10,-10)

\put(50,20){\line(1,-1){10}}
\put(50,20){\line(0,-1){20}}
\put(50,20){\line(-1,0){30}}

\thicklines

\put(50,20){\circle*{2.0}}
\put(10,0){\circle*{2.0}}
\put(30,0){\circle*{2.0}}
\put(50,0){\circle*{2.0}}
\put(70,0){\circle*{2.0}}
\put(20,20){\circle*{2.0}}

\put(40,-3){\makebox(0,0){$2$}}
\put(50,23){\makebox(0,0){$k$}}
\put(30,22){\makebox(0,0){$2$}}
\put(30,-3){\makebox(0,0){$b_1$}}
\put(50,-3){\makebox(0,0){$b_2$}}
\put(70,-3){\makebox(0,0){$c_1$}}
\put(20,23){\makebox(0,0){$a_2$}}
\put(10,-3){\makebox(0,0){$a_1$}}

\put(10,0){\line(1,0){60}}
\put(20,20){\line(-1,-2){10}}
\put(20,20){\line(1,-2){10}}

\end{picture}

\end{center}

\caption{} 

\label{fig:Fhat11-41}

\end{figure}

\credit{Subcase 4.2}\emph{$k$ is not connected to any vertex in $\{a_1,a_2\}$}

Then the subdiagram with the vertices 
$\{a_1,a_2,b_1,b_2,k\}$ is of type $F_4^{(1)}(3^{1};1;2^{1})_1$ 
(Fig.~\ref{fig:Fhat11-42}).

\begin{figure}[ht]

\setlength{\unitlength}{1.8pt}

\begin{center}

\begin{picture}(70,40)(-10,-10)

\put(50,20){\line(1,-1){20}}
\put(50,20){\line(0,-1){20}}

\thicklines

\put(50,20){\circle*{2.0}}
\put(10,0){\circle*{2.0}}
\put(30,0){\circle*{2.0}}
\put(50,0){\circle*{2.0}}
\put(70,0){\circle*{2.0}}
\put(20,20){\circle*{2.0}}

\put(40,-3){\makebox(0,0){$2$}}
\put(50,23){\makebox(0,0){$k$}}
\put(30,-3){\makebox(0,0){$b_1$}}
\put(50,-3){\makebox(0,0){$b_2$}}
\put(70,-3){\makebox(0,0){$c_1$}}
\put(20,23){\makebox(0,0){$a_2$}}
\put(10,-3){\makebox(0,0){$a_1$}}

\put(10,0){\line(1,0){60}}
\put(20,20){\line(-1,-2){10}}
\put(20,20){\line(1,-2){10}}

\end{picture}

\end{center}

\caption{} 

\label{fig:Fhat11-42}

\end{figure}

\begin{lemma}
\label{lem:caseFhat12}

If $\Gamma$ is of type $F_4^{(1)}(3^{1};1;2^{1})_2$, then Lemma~\ref{lem:xxx} holds.

\end{lemma}

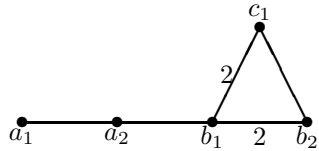
\begin{figure}[ht]

\setlength{\unitlength}{1.8pt}

\begin{center}

\begin{picture}(85,40)(-10,-10)

\thicklines

\put(10,0){\circle*{2.0}}
\put(30,0){\circle*{2.0}}
\put(50,0){\circle*{2.0}}
\put(70,0){\circle*{2.0}}
\put(60,20){\circle*{2.0}}

\put(60,-3){\makebox(0,0){$2$}}
\put(53,10){\makebox(0,0){$2$}}
\put(10,-3){\makebox(0,0){$a_1$}}
\put(30,-3){\makebox(0,0){$a_2$}}
\put(50,-3){\makebox(0,0){$b_1$}}
\put(70,-3){\makebox(0,0){$b_2$}}
\put(60,23){\makebox(0,0){$c_1$}}

\put(10,0){\line(1,0){60}}
\put(50,0){\line(1,2){10}}
\put(70,0){\line(-1,2){10}}

\end{picture}

\end{center}

\caption{The diagram $F_4^{(1)}(3^{1};1;2^{1})_2$}

\label{fig:Fhat12}

\end{figure}

\proof We assume that $\Gamma$ is indexed as in Fig.~\ref{fig:Fhat12}.

\credit{Case 1} \emph{$k$ is not connected to any vertex in $\{b_1,b_2\}$.}

\credit{Subcase 1.1}\emph{$k$ is not connected to $c_1$.}

In this case, the vertex $k$ is connected to both $a_1$ and $a_2$.
If the edge $[k,a_2]$ (hence the edge $[k,a_1]$) is weighted $2$
, then the subdiagram with the vertices $\{k,a_2,b_1,b_2\}$ is of type $C_n^{(1)}$ (Fig.~\ref{fig:Fhat12-11})
; otherwise the subdiagram with the vertices $\{k,a_2,b_1,b_2,c_1\}$ is of type $F_4^{(1)}(3^{1};1;2^{1})_2$. 

\begin{figure}[ht]

\setlength{\unitlength}{1.8pt}

\begin{center}

\begin{picture}(70,40)(-10,-10)

\put(10,0){\line(1,2){10}}
\put(30,0){\line(-1,2){10}}

\thicklines

\put(10,0){\circle*{2.0}}
\put(30,0){\circle*{2.0}}
\put(50,0){\circle*{2.0}}
\put(70,0){\circle*{2.0}}
\put(60,20){\circle*{2.0}}
\put(20,20){\circle*{2.0}}

\put(60,-3){\makebox(0,0){$2$}}
\put(53,10){\makebox(0,0){$2$}}
\put(10,-3){\makebox(0,0){$a_1$}}
\put(30,-3){\makebox(0,0){$a_2$}}
\put(50,-3){\makebox(0,0){$b_1$}}
\put(70,-3){\makebox(0,0){$b_2$}}
\put(60,23){\makebox(0,0){$c_1$}}
\put(20,23){\makebox(0,0){$k$}}
\put(13,10){\makebox(0,0){$2$}}
\put(27,10){\makebox(0,0){$2$}}

\put(10,0){\line(1,0){60}}
\put(50,0){\line(1,2){10}}
\put(70,0){\line(-1,2){10}}

\end{picture}

\end{center}

\caption{} 

\label{fig:Fhat12-11}

\end{figure}

\credit{Subcase 1.2}\emph{$k$ is connected to $c_1$.}

\begin{itemize}
\item
If the edge $[k,c_1]$ is weighted 2, then the subdiagram with the vertices 
$\{k,c_1,b_1\}$ is of type $C_n^{(1)}$
(Fig.~\ref{fig:Fhat12-12a}).

\item
Suppose that the edge $[k,c_1]$ is weightless. If $k$ is connected to $a_2$, then the subdiagram with the vertices $\{k,a_2,b_1,b_2\}$ is of type $C_n^{(1)}$ (Fig.~\ref{fig:Fhat12-12b});
otherwise the subdiagram with the vertices $\{k,a_1,a_2,b_1,b_2\}$ is of type $C_n^{(1)}$ 
(Fig.~\ref{fig:Fhat12-12c}).
\end{itemize}

\begin{figure}[ht]

\setlength{\unitlength}{1.8pt}

\begin{center}

\begin{picture}(70,40)(-10,-10)

\put(40,20){\line(1,0){20}}
\put(40,20){\line(-2,-1){10}}

\thicklines

\put(10,0){\circle*{2.0}}
\put(30,0){\circle*{2.0}}
\put(50,0){\circle*{2.0}}
\put(70,0){\circle*{2.0}}
\put(60,20){\circle*{2.0}}
\put(40,20){\circle*{2.0}}

\put(60,-3){\makebox(0,0){$2$}}
\put(53,10){\makebox(0,0){$2$}}
\put(10,-3){\makebox(0,0){$a_1$}}
\put(30,-3){\makebox(0,0){$a_2$}}
\put(50,-3){\makebox(0,0){$b_1$}}
\put(70,-3){\makebox(0,0){$b_2$}}
\put(60,23){\makebox(0,0){$c_1$}}
\put(40,23){\makebox(0,0){$k$}}
\put(50,23){\makebox(0,0){$2$}}

\put(10,0){\line(1,0){60}}
\put(50,0){\line(1,2){10}}
\put(70,0){\line(-1,2){10}}

\end{picture}

\end{center}

\caption{} 

\label{fig:Fhat12-12a}

\end{figure}

\begin{figure}[ht]

\setlength{\unitlength}{1.8pt}

\begin{center}

\begin{picture}(70,40)(-10,-10)

\put(40,20){\line(1,0){20}}
\put(40,20){\line(-1,-2){10}}
\put(40,20){\line(-3,-1){10}}

\thicklines

\put(10,0){\circle*{2.0}}
\put(30,0){\circle*{2.0}}
\put(50,0){\circle*{2.0}}
\put(70,0){\circle*{2.0}}
\put(60,20){\circle*{2.0}}
\put(40,20){\circle*{2.0}}

\put(60,-3){\makebox(0,0){$2$}}
\put(53,10){\makebox(0,0){$2$}}
\put(10,-3){\makebox(0,0){$a_1$}}
\put(30,-3){\makebox(0,0){$a_2$}}
\put(50,-3){\makebox(0,0){$b_1$}}
\put(70,-3){\makebox(0,0){$b_2$}}
\put(60,23){\makebox(0,0){$c_1$}}
\put(40,23){\makebox(0,0){$k$}}
\put(33,10){\makebox(0,0){$2$}}

\put(10,0){\line(1,0){60}}
\put(50,0){\line(1,2){10}}
\put(70,0){\line(-1,2){10}}

\end{picture}

\end{center}

\caption{} 

\label{fig:Fhat12-12b}

\end{figure}

\begin{figure}[ht]

\setlength{\unitlength}{1.8pt}

\begin{center}

\begin{picture}(70,40)(-10,-10)

\put(40,20){\line(1,0){20}}
\put(40,20){\line(-3,-2){30}}

\thicklines

\put(10,0){\circle*{2.0}}
\put(30,0){\circle*{2.0}}
\put(50,0){\circle*{2.0}}
\put(70,0){\circle*{2.0}}
\put(60,20){\circle*{2.0}}
\put(40,20){\circle*{2.0}}

\put(60,-3){\makebox(0,0){$2$}}
\put(53,10){\makebox(0,0){$2$}}
\put(10,-3){\makebox(0,0){$a_1$}}
\put(30,-3){\makebox(0,0){$a_2$}}
\put(50,-3){\makebox(0,0){$b_1$}}
\put(70,-3){\makebox(0,0){$b_2$}}
\put(60,23){\makebox(0,0){$c_1$}}
\put(40,23){\makebox(0,0){$k$}}
\put(25,13){\makebox(0,0){$2$}}

\put(10,0){\line(1,0){60}}
\put(50,0){\line(1,2){10}}
\put(70,0){\line(-1,2){10}}

\end{picture}

\end{center}

\caption{} 

\label{fig:Fhat12-12c}

\end{figure}

\credit{Case 2} \emph{$k$ is connected to both $b_1$ and $b_2$.}

\credit{Subcase 2.1}\emph{The edge $[k,b_1]$ is weighted 2.}

We first note that the edge $[k,b_2]$ is weightless.
\begin{itemize}
\item
Suppose that $k$ is not connected to $c_1$. Then the subdiagram with the vertices
$\{c_1,b_1,k\}$ is of type $C_n^{(1)}$ 
(Fig.~\ref{fig:Fhat12-21}).
\item
Suppose that $k$ is connected to $c_1$. Then
$k$ is contained in a non-oriented cycle.
\end{itemize}

\begin{figure}[ht]

\setlength{\unitlength}{1.8pt}

\begin{center}

\begin{picture}(70,50)(-10,-20)

\put(60,-20){\line(-1,2){10}}
\put(60,-20){\line(1,2){10}}
\put(60,-20){\line(-2,1){10}}

\thicklines

\put(10,0){\circle*{2.0}}
\put(30,0){\circle*{2.0}}
\put(50,0){\circle*{2.0}}
\put(70,0){\circle*{2.0}}
\put(60,20){\circle*{2.0}}
\put(60,-20){\circle*{2.0}}

\put(60,-3){\makebox(0,0){$2$}}
\put(53,10){\makebox(0,0){$2$}}
\put(9,-3){\makebox(0,0){$a_1$}}
\put(29,-3){\makebox(0,0){$a_2$}}
\put(49,-3){\makebox(0,0){$b_1$}}
\put(72,-3){\makebox(0,0){$b_2$}}
\put(60,23){\makebox(0,0){$c_1$}}
\put(60,-23){\makebox(0,0){$k$}}
\put(53,-10){\makebox(0,0){$2$}}

\put(10,0){\line(1,0){60}}
\put(50,0){\line(1,2){10}}
\put(70,0){\line(-1,2){10}}

\end{picture}

\end{center}

\caption{} 

\label{fig:Fhat12-21}

\end{figure}

\credit{Subcase 2.2}\emph{The edge $[k,b_1]$ is weightless.}

We note that the edge $[k,b_2]$ is weighted 2.
If $k$ is connected to $c_1$ then $k$ is contained in a non-oriented cycle.
Let us now assume that $k$ is not connected to $c_1$ and consider the subcases.
\begin{itemize}
\item
Suppose that $k$ is connected to $a_2$.
Then the subdiagram with the vertices $\{a_2,k,b_1,b_2\}$
is of type $B_n^{(1)}(r)$ 
(Fig.~\ref{fig:Fhat12-22i}).
 
\item
Suppose that $k$ is not connected to $a_2$.
Then the subdiagram with the vertices $\{c_1,b_1,a_2,k\}$ of type $B_n^{(1)}$ 
(Fig.~\ref{fig:Fhat12-22ii}).
\end{itemize}

\begin{figure}[ht]

\setlength{\unitlength}{1.8pt}

\begin{center}

\begin{picture}(70,50)(-10,-20)

\put(60,-20){\line(-1,2){10}}
\put(60,-20){\line(-3,1){10}}
\put(60,-20){\line(1,2){10}}
\put(60,-20){\line(-3,2){30}}

\thicklines

\put(10,0){\circle*{2.0}}
\put(30,0){\circle*{2.0}}
\put(50,0){\circle*{2.0}}
\put(70,0){\circle*{2.0}}
\put(60,20){\circle*{2.0}}
\put(60,-20){\circle*{2.0}}

\put(60,-3){\makebox(0,0){$2$}}
\put(53,10){\makebox(0,0){$2$}}
\put(9,-3){\makebox(0,0){$a_1$}}
\put(29,-3){\makebox(0,0){$a_2$}}
\put(49,-3){\makebox(0,0){$b_1$}}
\put(72,-3){\makebox(0,0){$b_2$}}
\put(60,23){\makebox(0,0){$c_1$}}
\put(60,-23){\makebox(0,0){$k$}}
\put(68,-10){\makebox(0,0){$2$}}

\put(10,0){\line(1,0){60}}
\put(50,0){\line(1,2){10}}
\put(70,0){\line(-1,2){10}}

\end{picture}

\end{center}

\caption{} 

\label{fig:Fhat12-22i}

\end{figure}

\begin{figure}[ht]

\setlength{\unitlength}{1.8pt}

\begin{center}

\begin{picture}(70,50)(-10,-20)

\put(60,-20){\line(-1,2){10}}
\put(60,-20){\line(-1,0){10}}
\put(60,-20){\line(1,2){10}}

\thicklines

\put(10,0){\circle*{2.0}}
\put(30,0){\circle*{2.0}}
\put(50,0){\circle*{2.0}}
\put(70,0){\circle*{2.0}}
\put(60,20){\circle*{2.0}}
\put(60,-20){\circle*{2.0}}

\put(60,-3){\makebox(0,0){$2$}}
\put(53,10){\makebox(0,0){$2$}}
\put(9,-3){\makebox(0,0){$a_1$}}
\put(29,-3){\makebox(0,0){$a_2$}}
\put(49,-3){\makebox(0,0){$b_1$}}
\put(72,-3){\makebox(0,0){$b_2$}}
\put(60,23){\makebox(0,0){$c_1$}}
\put(60,-23){\makebox(0,0){$k$}}
\put(68,-10){\makebox(0,0){$2$}}

\put(10,0){\line(1,0){60}}
\put(50,0){\line(1,2){10}}
\put(70,0){\line(-1,2){10}}

\end{picture}

\end{center}

\caption{} 

\label{fig:Fhat12-22ii}

\end{figure}

\credit{Case 3} \emph{$k$ is connected to $b_1$ not $b_2$.}

We immediately notice that if the edge $[k,b_1]$ is weighted 2,
then the subdiagram with the vertices $\{k,b_1,b_2\}$ is of type $C_n^{(1)}$.
To consider the remaining subcases, we assume that the edge $[k,b_1]$
is weightless.

\credit{Subcase 3.1}\emph{$k$ is not connected to $a_2$.}

Then the subdiagram with the vertices 
$\{k,b_1,a_2,b_2\}$ is of type $B_n^{(1)}$
(Fig.~\ref{fig:Fhat12-31}).

\begin{figure}[ht]

\setlength{\unitlength}{1.8pt}

\begin{center}

\begin{picture}(70,50)(-10,-20)

\put(40,20){\line(1,-2){10}}
\put(40,20){\line(1,0){10}}
\put(40,20){\line(-3,-1){10}}

\thicklines

\put(10,0){\circle*{2.0}}
\put(30,0){\circle*{2.0}}
\put(50,0){\circle*{2.0}}
\put(70,0){\circle*{2.0}}
\put(60,20){\circle*{2.0}}
\put(40,20){\circle*{2.0}}

\put(60,-3){\makebox(0,0){$2$}}
\put(53,10){\makebox(0,0){$2$}}
\put(9,-3){\makebox(0,0){$a_1$}}
\put(29,-3){\makebox(0,0){$a_2$}}
\put(49,-3){\makebox(0,0){$b_1$}}
\put(72,-3){\makebox(0,0){$b_2$}}
\put(60,23){\makebox(0,0){$c_1$}}
\put(40,23){\makebox(0,0){$k$}}

\put(10,0){\line(1,0){60}}
\put(50,0){\line(1,2){10}}
\put(70,0){\line(-1,2){10}}

\end{picture}

\end{center}

\caption{} 

\label{fig:Fhat12-31}

\end{figure}

\credit{Subcase 3.2}\emph{k is connected to $a_2$}

If $k$ is connected to $c_1$, then
the subdiagram with the vertices 
$\{k,b_1,c_1,a_2\}$ is of type $B_n^{(1)}(r)$
(Fig.~\ref{fig:Fhat12-32a}).
Let us now assume that $k$ is not connected to $c_1$.
Then the subdiagram with the vertices 
$\{k,a_2,b_1,b_2,c_1\}$ is of type $F_4^{(1)}(3^{2};0)$ 
(Fig.~\ref{fig:Fhat12-32b}).

\begin{figure}[ht]
\setlength{\unitlength}{1.8pt}

\begin{center}

\begin{picture}(70,50)(-10,-20)

\put(40,20){\line(1,-2){10}}
\put(40,20){\line(1,0){20}}
\put(40,20){\line(-3,-1){10}}
\put(40,20){\line(-1,-2){10}}

\thicklines

\put(10,0){\circle*{2.0}}
\put(30,0){\circle*{2.0}}
\put(50,0){\circle*{2.0}}
\put(70,0){\circle*{2.0}}
\put(60,20){\circle*{2.0}}
\put(40,20){\circle*{2.0}}

\put(60,-3){\makebox(0,0){$2$}}
\put(53,10){\makebox(0,0){$2$}}
\put(9,-3){\makebox(0,0){$a_1$}}
\put(29,-3){\makebox(0,0){$a_2$}}
\put(49,-3){\makebox(0,0){$b_1$}}
\put(72,-3){\makebox(0,0){$b_2$}}
\put(60,23){\makebox(0,0){$c_1$}}
\put(40,23){\makebox(0,0){$k$}}
\put(50,23){\makebox(0,0){$2$}}

\put(10,0){\line(1,0){60}}
\put(50,0){\line(1,2){10}}
\put(70,0){\line(-1,2){10}}

\end{picture}

\end{center}

\caption{} 

\label{fig:Fhat12-32a}

\end{figure}

\begin{figure}[ht]

\setlength{\unitlength}{1.8pt}

\begin{center}

\begin{picture}(70,50)(-10,-20)

\put(40,20){\line(1,-2){10}}
\put(40,20){\line(-3,-1){10}}
\put(40,20){\line(-1,-2){10}}

\thicklines

\put(10,0){\circle*{2.0}}
\put(30,0){\circle*{2.0}}
\put(50,0){\circle*{2.0}}
\put(70,0){\circle*{2.0}}
\put(60,20){\circle*{2.0}}
\put(40,20){\circle*{2.0}}

\put(60,-3){\makebox(0,0){$2$}}
\put(53,10){\makebox(0,0){$2$}}
\put(9,-3){\makebox(0,0){$a_1$}}
\put(29,-3){\makebox(0,0){$a_2$}}
\put(49,-3){\makebox(0,0){$b_1$}}
\put(72,-3){\makebox(0,0){$b_2$}}
\put(60,23){\makebox(0,0){$c_1$}}
\put(40,23){\makebox(0,0){$k$}}

\put(10,0){\line(1,0){60}}
\put(50,0){\line(1,2){10}}
\put(70,0){\line(-1,2){10}}

\end{picture}

\end{center}

\caption{} 

\label{fig:Fhat12-32b}

\end{figure}

\credit{Case 4} \emph{$k$ is connected to $b_2$ not $b_1$.}

We immediately notice that if the edge $[k,b_2]$ is weighted 2,
then the subdiagram with the vertices $\{k,b_1,b_2\}$ is of type $C_n^{(1)}$.
To consider the remaining subcases, we assume that the edge $[k,b_2]$
is weightless.

\credit{Subcase 4.1}\emph{$k$ is connected to $c_1$.}

Then the subdiagram with the vertices 
$\{b_1,b_2,c_1,k\}$ is of type $B_n^{(1)}(r)$
(Fig.~\ref{fig:Fhat12-41}).

\begin{figure}[ht]

\setlength{\unitlength}{1.8pt}

\begin{center}

\begin{picture}(70,50)(-10,-20)

\put(80,20){\line(-1,-2){10}}
\put(80,20){\line(-1,0){20}}

\thicklines

\put(10,0){\circle*{2.0}}
\put(30,0){\circle*{2.0}}
\put(50,0){\circle*{2.0}}
\put(70,0){\circle*{2.0}}
\put(60,20){\circle*{2.0}}
\put(80,20){\circle*{2.0}}

\put(60,-3){\makebox(0,0){$2$}}
\put(53,10){\makebox(0,0){$2$}}
\put(9,-3){\makebox(0,0){$a_1$}}
\put(29,-3){\makebox(0,0){$a_2$}}
\put(49,-3){\makebox(0,0){$b_1$}}
\put(72,-3){\makebox(0,0){$b_2$}}
\put(60,23){\makebox(0,0){$c_1$}}
\put(80,23){\makebox(0,0){$k$}}

\put(10,0){\line(1,0){60}}
\put(50,0){\line(1,2){10}}
\put(70,0){\line(-1,2){10}}

\end{picture}

\end{center}

\caption{} 

\label{fig:Fhat12-41}

\end{figure}

\credit{Subcase 4.2}\emph{$k$ is not connected to $c_1$}

\begin{itemize}
\item
Suppose that $k$ is connected to $a_2$.
Then the subdiagram with the vertices 
$\{k,a_2,b_1,c_1\}$ is of type $C_n^{(1)}$ 
(Fig.~\ref{fig:Fhat12-42}).
\item
Let us assume that $k$ is not connected to $a_2$.
 Then $k$ is
connected to $a_1$ and the subdiagram with the vertices 
$\{k,a_1,a_2,b_1,c_1\}$ is of type $C_n^{(1)}$. 
\end{itemize}

\begin{figure}[ht]

\setlength{\unitlength}{1.8pt}

\begin{center}

\begin{picture}(70,50)(-10,-20)

\put(60,-20){\line(-3,2){30}}
\put(60,-20){\line(1,2){10}}
\put(60,-20){\line(-5,2){10}}

\thicklines

\put(10,0){\circle*{2.0}}
\put(30,0){\circle*{2.0}}
\put(50,0){\circle*{2.0}}
\put(70,0){\circle*{2.0}}
\put(60,20){\circle*{2.0}}
\put(60,-20){\circle*{2.0}}

\put(60,-3){\makebox(0,0){$2$}}
\put(53,10){\makebox(0,0){$2$}}
\put(9,-3){\makebox(0,0){$a_1$}}
\put(29,-3){\makebox(0,0){$a_2$}}
\put(49,-3){\makebox(0,0){$b_1$}}
\put(72,-3){\makebox(0,0){$b_2$}}
\put(60,23){\makebox(0,0){$c_1$}}
\put(60,-23){\makebox(0,0){$k$}}
\put(43,-11){\makebox(0,0){$2$}}

\put(10,0){\line(1,0){60}}
\put(50,0){\line(1,2){10}}
\put(70,0){\line(-1,2){10}}

\end{picture}

\end{center}

\caption{} 

\label{fig:Fhat12-42}

\end{figure}

\endproof

\begin{lemma}
\label{lem:caseFhat22}

If $\Gamma$ is of type $F_4^{(1)}(3^{2};0)$, then Lemma~\ref{lem:xxx} holds.

\end{lemma}

\begin{figure}[ht]

\setlength{\unitlength}{1.8pt}

\begin{center}

\begin{picture}(70,40)(-10,-10)

\thicklines

\put(10,0){\circle*{2.0}}
\put(10,20){\circle*{2.0}}
\put(30,0){\circle*{2.0}}
\put(50,0){\circle*{2.0}}
\put(50,20){\circle*{2.0}}

\put(40,-3){\makebox(0,0){$2$}}
\put(39,12){\makebox(0,0){$2$}}
\put(10,-3){\makebox(0,0){$a_1$}}
\put(30,-3){\makebox(0,0){$b_1$}}
\put(50,-3){\makebox(0,0){$b_2$}}
\put(10,23){\makebox(0,0){$a_2$}}
\put(50,23){\makebox(0,0){$c_1$}}

\put(10,0){\line(1,0){40}}
\put(10,0){\line(0,1){20}}
\put(30,0){\line(1,1){20}}
\put(30,0){\line(-1,1){20}}
\put(50,0){\line(0,1){20}}

\end{picture}

\end{center}

\caption{The diagram $F_4^{(1)}(3^{2};0)$}

\label{fig:Fhat22}

\end{figure}
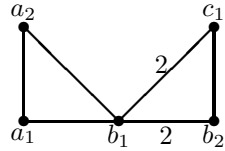

\proof We assume that $\Gamma$ is indexed as in Fig.~\ref{fig:Fhat22}.

\credit{Case 1} \emph{$k$ is not connected to any vertex in $\{b_1,b_2\}$.}

\credit{Subcase 1.1}\emph{$k$ is not connected to $c_1$.}

In this case the vertex $k$ is connected to both $a_1$ and $a_2$.
If the edge $[k,a_2]$ (hence the edge $[k,a_1]$) is weighted 2
, then the subdiagram with the vertices $\{k,a_1,a_2,b_1\}$ is of type $B_n^{(1)}(r)$ 
(Fig.~\ref{fig:Fhat22-11}).
; otherwise the subdiagram with the vertices $\{k,a_2,b_1,b_2,c_1\}$ is of type $F_4^{(1)}(3^{1};1;2^{1})_2$. 

\begin{figure}[ht]

\setlength{\unitlength}{1.8pt}

\begin{center}

\begin{picture}(70,40)(-10,-10)

\put(-10,20){\line(1,0){20}}
\put(-10,20){\line(1,-1){20}}
\put(0,23){\makebox(0,0){$2$}}
\put(0,8){\makebox(0,0){$2$}}

\thicklines

\put(10,0){\circle*{2.0}}
\put(10,20){\circle*{2.0}}
\put(30,0){\circle*{2.0}}
\put(50,0){\circle*{2.0}}
\put(50,20){\circle*{2.0}}
\put(-10,20){\circle*{2.0}}

\put(40,-3){\makebox(0,0){$2$}}
\put(39,12){\makebox(0,0){$2$}}
\put(10,-3){\makebox(0,0){$a_1$}}
\put(30,-3){\makebox(0,0){$b_1$}}
\put(50,-3){\makebox(0,0){$b_2$}}
\put(10,23){\makebox(0,0){$a_2$}}
\put(50,23){\makebox(0,0){$c_1$}}
\put(-13,20){\makebox(0,0){$k$}}

\put(10,0){\line(1,0){40}}
\put(10,0){\line(0,1){20}}
\put(30,0){\line(1,1){20}}
\put(30,0){\line(-1,1){20}}
\put(50,0){\line(0,1){20}}

\end{picture}

\end{center}

\caption{} 

\label{fig:Fhat22-11}

\end{figure}

\credit{Subcase 1.2}\emph{$k$ is connected to $c_1$.}

\begin{itemize}
\item
If the edge $[k,c_1]$ is weighted 2, then the subdiagram with the vertices 
$\{k,c_1,b_1\}$ is of type $C_n^{(1)}$.
\item
Suppose that the edge $[k,c_1]$ is weightless and assume, without loss of generality, that
$k$ is connected to $a_2$. Then
the subdiagram with the vertices $\{k,a_2,b_1,b_2\}$ is of type $C_n^{(1)}$
(Fig.~\ref{fig:Fhat22-12}).
\end{itemize}
\begin{figure}[ht]

\setlength{\unitlength}{1.8pt}

\begin{center}

\begin{picture}(70,60)(-10,-10)

\put(30,40){\line(1,-1){20}}
\put(30,40){\line(-1,-1){20}}
\put(20,33){\makebox(0,0){$2$}}
\put(30,40){\line(-1,0){10}}

\thicklines

\put(10,0){\circle*{2.0}}
\put(10,20){\circle*{2.0}}
\put(30,0){\circle*{2.0}}
\put(50,0){\circle*{2.0}}
\put(50,20){\circle*{2.0}}
\put(30,40){\circle*{2.0}}

\put(40,-3){\makebox(0,0){$2$}}
\put(39,12){\makebox(0,0){$2$}}
\put(10,-3){\makebox(0,0){$a_1$}}
\put(30,-3){\makebox(0,0){$b_1$}}
\put(50,-3){\makebox(0,0){$b_2$}}
\put(10,23){\makebox(0,0){$a_2$}}
\put(50,23){\makebox(0,0){$c_1$}}
\put(32,43){\makebox(0,0){$k$}}

\put(10,0){\line(1,0){40}}
\put(10,0){\line(0,1){20}}
\put(30,0){\line(1,1){20}}
\put(30,0){\line(-1,1){20}}
\put(50,0){\line(0,1){20}}

\end{picture}

\end{center}

\caption{}

\label{fig:Fhat22-12}

\end{figure}

\credit{Case 2} \emph{$k$ is connected to both $b_1$ and $b_2$.}

\credit{Subcase 2.1}\emph{The edge $[k,b_1]$ is weighted 2.}

We first note that the edge $[k,b_2]$ is weightless. Let us now consider the subcases.
\begin{itemize}
\item
Suppose that $k$ is not connected to $c_1$. Then the subdiagram with the vertices
$\{c_1,b_1,k\}$ is of type $C_n^{(1)}$ 
(Fig.~\ref{fig:Fhat22-21}).
\item 
Suppose that $k$ is connected to $c_1$. Then
$k$ is contained in a non-oriented cycle.
\end{itemize}

\begin{figure}[ht]

\setlength{\unitlength}{1.8pt}

\begin{center}

\begin{picture}(70,50)(-10,-20)

\put(40,-20){\line(1,2){10}}
\put(40,-20){\line(-1,2){10}}
\put(40,-20){\line(-3,1){10}}
\put(33,-11){\makebox(0,0){$2$}}

\thicklines

\put(10,0){\circle*{2.0}}
\put(10,20){\circle*{2.0}}
\put(30,0){\circle*{2.0}}
\put(50,0){\circle*{2.0}}
\put(50,20){\circle*{2.0}}
\put(40,-20){\circle*{2.0}}

\put(40,-2){\makebox(0,0){$2$}}
\put(39,12){\makebox(0,0){$2$}}
\put(9,-3){\makebox(0,0){$a_1$}}
\put(29,-3){\makebox(0,0){$b_1$}}
\put(52,-3){\makebox(0,0){$b_2$}}
\put(10,23){\makebox(0,0){$a_2$}}
\put(50,23){\makebox(0,0){$c_1$}}
\put(40,-23){\makebox(0,0){$k$}}

\put(10,0){\line(1,0){40}}
\put(10,0){\line(0,1){20}}
\put(30,0){\line(1,1){20}}
\put(30,0){\line(-1,1){20}}
\put(50,0){\line(0,1){20}}

\end{picture}

\end{center}

\caption{} 

\label{fig:Fhat22-21}

\end{figure}

\credit{Subcase 2.2}\emph{The edge $[k,b_1]$ is weightless.}
We note that the edge $[k,b_2]$ is weighted 2.
If $k$ is connected to $c_1$ then $k$ is contained in a non-oriented cycle.
Let us now assume that $k$ is not connected to $c_1$ and consider the subcases.
\begin{itemize}
\item
Suppose that $k$ is connected to a vertex in $\{a_1,a_2\}$.
If $k$ is connected to both $a_1$ and $a_2$, then
it is contained in a non-oriented cycle.
Let us now assume, without loss of generality that, the vertex $k$
is connected to $a_1$ but not connected to $a_2$.
Then the subdiagram with the vertices $\{a_1,k,b_1,b_2\}$
is of type $B_n^{(1)}(r)$ 
(Fig.~\ref{fig:Fhat22-22i}).
\item 
Suppose that $k$ is not connected to any vertex in $\{a_1,a_2\}$.
Then the subdiagram with the vertices $\{c_1,b_1,a_2,k\}$ of type $B_n^{(1)}$ 
(Fig.~\ref{fig:Fhat22-22ii}).
\end{itemize}

\begin{figure}[ht]
\setlength{\unitlength}{1.8pt}

\begin{center}

\begin{picture}(70,50)(-10,-20)

\put(40,-20){\line(1,2){10}}
\put(40,-20){\line(-1,2){10}}
\put(40,-20){\line(-3,2){30}}
\put(47,-11){\makebox(0,0){$2$}}

\thicklines

\put(10,0){\circle*{2.0}}
\put(10,20){\circle*{2.0}}
\put(30,0){\circle*{2.0}}
\put(50,0){\circle*{2.0}}
\put(50,20){\circle*{2.0}}
\put(40,-20){\circle*{2.0}}

\put(40,-2){\makebox(0,0){$2$}}
\put(39,12){\makebox(0,0){$2$}}
\put(9,-3){\makebox(0,0){$a_1$}}
\put(29,-3){\makebox(0,0){$b_1$}}
\put(52,-3){\makebox(0,0){$b_2$}}
\put(10,23){\makebox(0,0){$a_2$}}
\put(50,23){\makebox(0,0){$c_1$}}
\put(40,-23){\makebox(0,0){$k$}}

\put(10,0){\line(1,0){40}}
\put(10,0){\line(0,1){20}}
\put(30,0){\line(1,1){20}}
\put(30,0){\line(-1,1){20}}
\put(50,0){\line(0,1){20}}

\end{picture}

\end{center}

\caption{} 

\label{fig:Fhat22-22i}

\end{figure}

\begin{figure}[ht]

\setlength{\unitlength}{1.8pt}

\begin{center}

\begin{picture}(70,50)(-10,-20)

\put(40,-20){\line(1,2){10}}
\put(40,-20){\line(-1,2){10}}
\put(47,-11){\makebox(0,0){$2$}}

\thicklines

\put(10,0){\circle*{2.0}}
\put(10,20){\circle*{2.0}}
\put(30,0){\circle*{2.0}}
\put(50,0){\circle*{2.0}}
\put(50,20){\circle*{2.0}}
\put(40,-20){\circle*{2.0}}

\put(40,-2){\makebox(0,0){$2$}}
\put(39,12){\makebox(0,0){$2$}}
\put(9,-3){\makebox(0,0){$a_1$}}
\put(29,-3){\makebox(0,0){$b_1$}}
\put(52,-3){\makebox(0,0){$b_2$}}
\put(10,23){\makebox(0,0){$a_2$}}
\put(50,23){\makebox(0,0){$c_1$}}
\put(40,-23){\makebox(0,0){$k$}}

\put(10,0){\line(1,0){40}}
\put(10,0){\line(0,1){20}}
\put(30,0){\line(1,1){20}}
\put(30,0){\line(-1,1){20}}
\put(50,0){\line(0,1){20}}

\end{picture}

\end{center}

\caption{} 

\label{fig:Fhat22-22ii}

\end{figure}

\credit{Case 3} \emph{$k$ is connected to $b_1$ and not connected to $b_2$.}

We immediately note that if the edge $[k,b_1]$ is weighted 2,
then the subdiagram with the vertices $\{k,b_1,b_2\}$ is of type $C_n^{(1)}$.
To consider the remaining subcases we assume that the edge $[k,b_1]$
is weightless.

\credit{Subcase 3.1}\emph{$k$ is not connected to $a_1$ or not connected to $a_2$.}

Let us assume without loss of generality that $k$ is not connected to $a_2$.
Then the subdiagram with the vertices 
$\{k,b_1,a_2,b_2\}$ is of type $B_n^{(1)}$
(Fig.~\ref{fig:Fhat22-31}).

\begin{figure}[ht]

\setlength{\unitlength}{1.8pt}

\begin{center}

\begin{picture}(70,50)(-10,-20)

\put(30,-20){\line(0,1){20}}
\put(30,-20){\line(-1,1){20}}
\put(30,-20){\line(1,2){5}}

\thicklines

\put(10,0){\circle*{2.0}}
\put(10,20){\circle*{2.0}}
\put(30,0){\circle*{2.0}}
\put(50,0){\circle*{2.0}}
\put(50,20){\circle*{2.0}}
\put(30,-20){\circle*{2.0}}

\put(40,-2){\makebox(0,0){$2$}}
\put(39,12){\makebox(0,0){$2$}}
\put(9,-3){\makebox(0,0){$a_1$}}
\put(29,-3){\makebox(0,0){$b_1$}}
\put(52,-3){\makebox(0,0){$b_2$}}
\put(10,23){\makebox(0,0){$a_2$}}
\put(50,23){\makebox(0,0){$c_1$}}
\put(30,-23){\makebox(0,0){$k$}}

\put(10,0){\line(1,0){40}}
\put(10,0){\line(0,1){20}}
\put(30,0){\line(1,1){20}}
\put(30,0){\line(-1,1){20}}
\put(50,0){\line(0,1){20}}

\end{picture}

\end{center}

\caption{} 

\label{fig:Fhat22-31}

\end{figure}

\credit{Subcase 3.2}\emph{$k$ is connected to $a_1$ and $a_2$}

Then $k$ is contained in a non-oriented cycle.

\credit{Case 4} \emph{$k$ is connected to $b_2$ not $b_1$.}

We immediately note that if the edge $[k,b_2]$ is weighted 2,
then the subdiagram with the vertices $\{k,b_1,b_2\}$ is of type $C_n^{(1)}$.
To consider the remaining subcases we assume that the edge $[k,b_2]$
is weightless.

\credit{Subcase 4.1}\emph{$k$ is connected to $c_1$.}

Then the subdiagram with the vertices 
$\{b_1,b_2,c_1,k\}$ is of type $B_n^{(1)}(r)$ 
(Fig.~\ref{fig:Fhat22-41}).

\begin{figure}[ht]

\setlength{\unitlength}{1.8pt}

\begin{center}

\begin{picture}(70,50)(-10,-20)

\put(70,0){\line(-1,0){20}}
\put(70,0){\line(-1,1){20}}
\put(70,0){\line(0,1){5}}

\thicklines

\put(10,0){\circle*{2.0}}
\put(10,20){\circle*{2.0}}
\put(30,0){\circle*{2.0}}
\put(50,0){\circle*{2.0}}
\put(50,20){\circle*{2.0}}
\put(70,0){\circle*{2.0}}

\put(40,-2){\makebox(0,0){$2$}}
\put(39,12){\makebox(0,0){$2$}}
\put(9,-3){\makebox(0,0){$a_1$}}
\put(29,-3){\makebox(0,0){$b_1$}}
\put(52,-3){\makebox(0,0){$b_2$}}
\put(10,23){\makebox(0,0){$a_2$}}
\put(50,23){\makebox(0,0){$c_1$}}
\put(73,0){\makebox(0,0){$k$}}

\put(10,0){\line(1,0){40}}
\put(10,0){\line(0,1){20}}
\put(30,0){\line(1,1){20}}
\put(30,0){\line(-1,1){20}}
\put(50,0){\line(0,1){20}}

\end{picture}

\end{center}

\caption{} 

\label{fig:Fhat22-41}

\end{figure}

\credit{Subcase 4.2}\emph{$k$ is not connected to $c_1$.}

We may assume without loss of generality that $k$ is connected to $a_1$.
Then the subdiagram with the vertices 
$\{k,a_1,b_1,c_1\}$ is of type $C_n^{(1)}$ 
(Fig.~\ref{fig:Fhat22-42}).
\endproof

\begin{figure}[ht]

\setlength{\unitlength}{1.8pt}

\begin{center}

\begin{picture}(70,50)(-10,-20)

\put(30,-20){\line(1,1){20}}
\put(30,-20){\line(-1,1){20}}
\put(10,20){\line(1,-2){4}}
\put(19,-12){\makebox(0,0){$2$}}

\thicklines

\put(10,0){\circle*{2.0}}
\put(10,20){\circle*{2.0}}
\put(30,0){\circle*{2.0}}
\put(50,0){\circle*{2.0}}
\put(50,20){\circle*{2.0}}
\put(30,-20){\circle*{2.0}}

\put(40,-2){\makebox(0,0){$2$}}
\put(39,12){\makebox(0,0){$2$}}
\put(9,-3){\makebox(0,0){$a_1$}}
\put(29,-3){\makebox(0,0){$b_1$}}
\put(52,-3){\makebox(0,0){$b_2$}}
\put(10,23){\makebox(0,0){$a_2$}}
\put(50,23){\makebox(0,0){$c_1$}}
\put(30,-23){\makebox(0,0){$k$}}

\put(10,0){\line(1,0){40}}
\put(10,0){\line(0,1){20}}
\put(30,0){\line(1,1){20}}
\put(30,0){\line(-1,1){20}}
\put(50,0){\line(0,1){20}}

\end{picture}

\end{center}

\caption{} 

\label{fig:Fhat22-42}

\end{figure}

\begin{lemma}
\label{lem:caseFhat12b}

If $\Gamma$ is of type $F_4^{(1)}(3^{1};1;1^2)$, then Lemma~\ref{lem:xxx} holds.

\end{lemma}

\begin{figure}[ht]

\setlength{\unitlength}{1.8pt}

\begin{center}

\begin{picture}(82,40)(-10,-10)

\thicklines

\put(10,0){\circle*{2.0}}
\put(30,0){\circle*{2.0}}
\put(50,0){\circle*{2.0}}
\put(70,0){\circle*{2.0}}
\put(40,20){\circle*{2.0}}

\put(40,-2){\makebox(0,0){$2$}}
\put(47,10){\makebox(0,0){$2$}}
\put(10,-3){\makebox(0,0){$a_1$}}
\put(30,-3){\makebox(0,0){$b_1$}}
\put(50,-3){\makebox(0,0){$b_2$}}
\put(70,-3){\makebox(0,0){$a_2$}}
\put(40,23){\makebox(0,0){$c_1$}}

\put(10,0){\line(1,0){60}}
\put(30,0){\line(1,2){10}}
\put(50,0){\line(-1,2){10}}

\end{picture}

\end{center}

\caption{The diagram $F_4^{(1)}(3^{1};1;1^{2})$}

\label{fig:Fhat12b}

\end{figure}
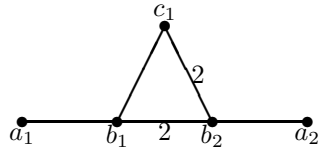

\proof We assume that $\Gamma$ is indexed as in Fig.~\ref{fig:Fhat12b}.

\credit{Case 1} \emph{$k$ is not connected to any vertex in $\{b_1,b_2\}$.}

\credit{Subcase 1.1}\emph{$k$ is not connected to $c_1$.}

We note that the vertex $k$ is connected to both $a_1$ and $a_2$.
If the edge $[k,a_2]$ is weighted 2
, then the subdiagram with the vertices $\{k,a_2,b_2,b_1\}$ is of type $C_n^{(1)}$
; otherwise the subdiagram with the vertices $\{k,a_1,b_1,b_2\}$ is of type $C_n^{(1)}$ 
(Fig.~\ref{fig:Fhat12b-11}).

\begin{figure}[ht]

\setlength{\unitlength}{1.8pt}

\begin{center}

\begin{picture}(70,40)(-10,-20)

\put(40,-20){\line(-3,2){30}}
\put(40,-20){\line(3,2){30}}
\put(26,-13){\makebox(0,0){$2$}}
\put(40,-20){\circle*{2.0}}
\put(40,-24){\makebox(0,0){$k$}}

\thicklines

\put(10,0){\circle*{2.0}}
\put(30,0){\circle*{2.0}}
\put(50,0){\circle*{2.0}}
\put(70,0){\circle*{2.0}}
\put(40,20){\circle*{2.0}}

\put(40,-2){\makebox(0,0){$2$}}
\put(47,10){\makebox(0,0){$2$}}
\put(9,-3){\makebox(0,0){$a_1$}}
\put(30,-3){\makebox(0,0){$b_1$}}
\put(50,-3){\makebox(0,0){$b_2$}}
\put(71,-3){\makebox(0,0){$a_2$}}
\put(40,23){\makebox(0,0){$c_1$}}

\put(10,0){\line(1,0){60}}
\put(30,0){\line(1,2){10}}
\put(50,0){\line(-1,2){10}}

\end{picture}

\end{center}

\caption{} 

\label{fig:Fhat12b-11}

\end{figure}

\credit{Subcase 1.2}\emph{$k$ is connected to $c_1$.}

\begin{itemize}
\item
If the edge $[k,c_1]$ is weighted 2, then the subdiagram with the vertices 
$\{k,c_1,b_2\}$ is of type $C_n^{(1)}$.
\item
Suppose that edge $[k,c_1]$ is weightless. If $k$ is connected to $a_2$, then the subdiagram with the vertices $\{k,a_2,b_2,b_1\}$ is of type $C_n^{(1)}$
(Fig.~\ref{fig:Fhat12b-12a}); 
otherwise the subdiagram with the vertices $\{k,a_1,b_1,c_1,b_2\}$ is of type $B_n^{(1)}(r)$ (Fig.~\ref{fig:Fhat12b-12b}). 
\end{itemize}

\begin{figure}[ht]

\setlength{\unitlength}{1.8pt}

\begin{center}

\begin{picture}(70,40)(-10,-10)

\put(70,20){\line(0,-1){20}}
\put(70,20){\line(-1,0){30}}
\put(72,10){\makebox(0,0){$2$}}
\put(70,20){\circle*{2.0}}
\put(70,24){\makebox(0,0){$k$}}
\put(10,0){\line(0,1){10}}

\thicklines

\put(10,0){\circle*{2.0}}
\put(30,0){\circle*{2.0}}
\put(50,0){\circle*{2.0}}
\put(70,0){\circle*{2.0}}
\put(40,20){\circle*{2.0}}

\put(40,-2){\makebox(0,0){$2$}}
\put(47,10){\makebox(0,0){$2$}}
\put(9,-3){\makebox(0,0){$a_1$}}
\put(30,-3){\makebox(0,0){$b_1$}}
\put(50,-3){\makebox(0,0){$b_2$}}
\put(71,-3){\makebox(0,0){$a_2$}}
\put(40,23){\makebox(0,0){$c_1$}}

\put(10,0){\line(1,0){60}}
\put(30,0){\line(1,2){10}}
\put(50,0){\line(-1,2){10}}

\end{picture}

\end{center}

\caption{} 

\label{fig:Fhat12b-12a}

\end{figure}

\begin{figure}[ht]

\setlength{\unitlength}{1.8pt}

\begin{center}

\begin{picture}(70,40)(-10,-10)

\put(10,20){\line(0,-1){20}}
\put(10,20){\line(1,0){30}}
\put(10,20){\circle*{2.0}}
\put(10,24){\makebox(0,0){$k$}}

\thicklines

\put(10,0){\circle*{2.0}}
\put(30,0){\circle*{2.0}}
\put(50,0){\circle*{2.0}}
\put(70,0){\circle*{2.0}}
\put(40,20){\circle*{2.0}}

\put(40,-2){\makebox(0,0){$2$}}
\put(47,10){\makebox(0,0){$2$}}
\put(9,-3){\makebox(0,0){$a_1$}}
\put(30,-3){\makebox(0,0){$b_1$}}
\put(50,-3){\makebox(0,0){$b_2$}}
\put(71,-3){\makebox(0,0){$a_2$}}
\put(40,23){\makebox(0,0){$c_1$}}

\put(10,0){\line(1,0){60}}
\put(30,0){\line(1,2){10}}
\put(50,0){\line(-1,2){10}}

\end{picture}

\end{center}

\caption{} 

\label{fig:Fhat12b-12b}

\end{figure}

\credit{Case 2} \emph{$k$ is connected to both $b_1$ and $b_2$.}

If $k$ is connected to $c_1$ then $k$ is contained in a non-oriented cycle.
Let us now assume that $k$ is not connected to $c_1$ and consider the subcases.

\credit{Subcase 2.1}\emph{The edge $[k,b_1]$ is weighted 2.}

We first note that the edge $[k,b_2]$ is weightless. 
\begin{itemize}
\item
Suppose that $k$ be not connected to $a_1$ or not connected to $a_2$. 
Let us first assume that $k$ is not
connected to $a_1$. Then the subdiagram with the vertices
$\{c_1,b_1,a_1,k\}$ is of type $B_n^{(1)}$ 
(Fig.~\ref{fig:Fhat12b-21ia}).
Let us now assume that $k$ is not
connected to $a_2$.
Then the subdiagram with the vertices
$\{c_1,b_2,a_2,k\}$ is of type $B_n^{(1)}$ 
(Fig.~\ref{fig:Fhat12b-21ib}).

\item
Suppose that $k$ is connected to $a_1$ and $a_2$. Then
the subdiagram with the vertices $\{k,b_1,b_2,a_2\}$ is of type $B_n^{(1)}(r)$
(Fig.~\ref{fig:Fhat12b-21ii}).
\end{itemize}

\begin{figure}[ht]

\setlength{\unitlength}{1.8pt}

\begin{center}

\begin{picture}(70,40)(-10,-20)

\put(40,-20){\line(-1,2){10}}
\put(40,-20){\line(1,2){10}}
\put(40,-20){\line(3,2){10}}
\put(34,-11){\makebox(0,0){$2$}}
\put(40,-20){\circle*{2.0}}
\put(40,-24){\makebox(0,0){$k$}}

\thicklines

\put(10,0){\circle*{2.0}}
\put(30,0){\circle*{2.0}}
\put(50,0){\circle*{2.0}}
\put(70,0){\circle*{2.0}}
\put(40,20){\circle*{2.0}}

\put(40,-2){\makebox(0,0){$2$}}
\put(47,10){\makebox(0,0){$2$}}
\put(9,-3){\makebox(0,0){$a_1$}}
\put(29,-3){\makebox(0,0){$b_1$}}
\put(51,-3){\makebox(0,0){$b_2$}}
\put(71,-3){\makebox(0,0){$a_2$}}
\put(40,23){\makebox(0,0){$c_1$}}

\put(10,0){\line(1,0){60}}
\put(30,0){\line(1,2){10}}
\put(50,0){\line(-1,2){10}}

\end{picture}

\end{center}

\caption{} 

\label{fig:Fhat12b-21ia}

\end{figure}

\begin{figure}[ht]

\setlength{\unitlength}{1.8pt}

\begin{center}

\begin{picture}(70,40)(-10,-20)

\put(40,-20){\line(-1,2){10}}
\put(40,-20){\line(1,2){10}}
\put(40,-20){\line(-3,2){10}}

\put(34,-11){\makebox(0,0){$2$}}
\put(40,-20){\circle*{2.0}}
\put(40,-24){\makebox(0,0){$k$}}

\thicklines

\put(10,0){\circle*{2.0}}
\put(30,0){\circle*{2.0}}
\put(50,0){\circle*{2.0}}
\put(70,0){\circle*{2.0}}
\put(40,20){\circle*{2.0}}

\put(40,-2){\makebox(0,0){$2$}}
\put(47,10){\makebox(0,0){$2$}}
\put(9,-3){\makebox(0,0){$a_1$}}
\put(29,-3){\makebox(0,0){$b_1$}}
\put(51,-3){\makebox(0,0){$b_2$}}
\put(71,-3){\makebox(0,0){$a_2$}}
\put(40,23){\makebox(0,0){$c_1$}}

\put(10,0){\line(1,0){60}}
\put(30,0){\line(1,2){10}}
\put(50,0){\line(-1,2){10}}

\end{picture}

\end{center}

\caption{} 

\label{fig:Fhat12b-21ib}

\end{figure}

\begin{figure}[ht]

\setlength{\unitlength}{1.8pt}

\begin{center}

\begin{picture}(70,40)(-10,-20)

\put(40,-20){\line(-1,2){10}}
\put(40,-20){\line(1,2){10}}
\put(40,-20){\line(3,2){30}}
\put(40,-20){\line(-3,2){30}}

\put(34,-11){\makebox(0,0){$2$}}
\put(21,-11){\makebox(0,0){$2$}}
\put(40,-20){\circle*{2.0}}
\put(40,-24){\makebox(0,0){$k$}}

\thicklines

\put(10,0){\circle*{2.0}}
\put(30,0){\circle*{2.0}}
\put(50,0){\circle*{2.0}}
\put(70,0){\circle*{2.0}}
\put(40,20){\circle*{2.0}}

\put(40,-2){\makebox(0,0){$2$}}
\put(47,10){\makebox(0,0){$2$}}
\put(9,-3){\makebox(0,0){$a_1$}}
\put(29,-3){\makebox(0,0){$b_1$}}
\put(51,-3){\makebox(0,0){$b_2$}}
\put(71,-3){\makebox(0,0){$a_2$}}
\put(40,23){\makebox(0,0){$c_1$}}

\put(10,0){\line(1,0){60}}
\put(30,0){\line(1,2){10}}
\put(50,0){\line(-1,2){10}}

\end{picture}

\end{center}

\caption{} 

\label{fig:Fhat12b-21ii}

\end{figure}

\credit{Subcase 2.2}\emph{The edge $[k,b_2]$ is weighted 2.}

Then the subdiagram with the vertices $\{k,b_2,c_1\}$
is of type $C_n^{(1)}$.

\credit{Case 3} \emph{$k$ is connected to $b_1$ not $b_2$.}

We immediately note that if the edge $[k,b_1]$ is weighted 2,
then the subdiagram with the vertices $\{k,b_1,b_2\}$ is of type $C_n^{(1)}$.
To consider the remaining subcases we assume that the edge $[k,b_1]$
is weightless.

\credit{Subcase 3.1}\emph{$k$ is connected to $c_1$.}

Then the subdiagram with the vertices 
$\{k,b_1,c_1,b_2\}$ is of type $B_n^{(1)}(r)$ 
(Fig.~\ref{fig:Fhat12b-31}).

\begin{figure}[ht]

\setlength{\unitlength}{1.8pt}

\begin{center}

\begin{picture}(70,40)(-10,-10)

\put(20,20){\line(1,-2){10}}
\put(20,20){\line(1,0){20}}
\put(20,20){\line(1,1){10}}

\put(20,20){\circle*{2.0}}
\put(20,24){\makebox(0,0){$k$}}

\thicklines

\put(10,0){\circle*{2.0}}
\put(30,0){\circle*{2.0}}
\put(50,0){\circle*{2.0}}
\put(70,0){\circle*{2.0}}
\put(40,20){\circle*{2.0}}

\put(40,-2){\makebox(0,0){$2$}}
\put(47,10){\makebox(0,0){$2$}}
\put(9,-3){\makebox(0,0){$a_1$}}
\put(29,-3){\makebox(0,0){$b_1$}}
\put(51,-3){\makebox(0,0){$b_2$}}
\put(71,-3){\makebox(0,0){$a_2$}}
\put(40,23){\makebox(0,0){$c_1$}}

\put(10,0){\line(1,0){60}}
\put(30,0){\line(1,2){10}}
\put(50,0){\line(-1,2){10}}

\end{picture}

\end{center}

\caption{} 

\label{fig:Fhat12b-31}

\end{figure}

\credit{Subcase 3.2}\emph{$k$ is not connected to $c_1$.}

If $k$ is connected to $a_2$, then
the subdiagram with the vertices 
$\{k,a_2,b_2,c_1\}$ is of type $C_n^{(1)}$ 
(Fig.~\ref{fig:Fhat12b-32a});
otherwise the subdiagram with the vertices 
$\{k,b_1,b_2,c_1,a_1\}$ is of type $F_4^{(1)}(3^{1};1;1^2)$ 
(Fig.~\ref{fig:Fhat12b-32b}).

\begin{figure}[ht]

\setlength{\unitlength}{1.8pt}

\begin{center}

\begin{picture}(70,40)(-10,-20)

\put(40,-20){\line(-1,2){10}}
\put(40,-20){\line(3,2){30}}
\put(40,-20){\line(-3,2){10}}

\put(60,-10){\makebox(0,0){$2$}}
\put(40,-20){\circle*{2.0}}
\put(40,-24){\makebox(0,0){$k$}}

\thicklines

\put(10,0){\circle*{2.0}}
\put(30,0){\circle*{2.0}}
\put(50,0){\circle*{2.0}}
\put(70,0){\circle*{2.0}}
\put(40,20){\circle*{2.0}}

\put(40,-2){\makebox(0,0){$2$}}
\put(47,10){\makebox(0,0){$2$}}
\put(9,-3){\makebox(0,0){$a_1$}}
\put(29,-3){\makebox(0,0){$b_1$}}
\put(51,-3){\makebox(0,0){$b_2$}}
\put(71,-3){\makebox(0,0){$a_2$}}
\put(40,23){\makebox(0,0){$c_1$}}

\put(10,0){\line(1,0){60}}
\put(30,0){\line(1,2){10}}
\put(50,0){\line(-1,2){10}}

\end{picture}

\end{center}

\caption{} 

\label{fig:Fhat12b-32a}

\end{figure}

\begin{figure}[ht]

\setlength{\unitlength}{1.8pt}

\begin{center}

\begin{picture}(70,40)(-10,-20)

\put(30,-20){\line(0,1){20}}
\put(30,-20){\line(-1,1){20}}

\put(30,-20){\circle*{2.0}}
\put(30,-24){\makebox(0,0){$k$}}

\thicklines

\put(10,0){\circle*{2.0}}
\put(30,0){\circle*{2.0}}
\put(50,0){\circle*{2.0}}
\put(70,0){\circle*{2.0}}
\put(40,20){\circle*{2.0}}

\put(40,-2){\makebox(0,0){$2$}}
\put(47,10){\makebox(0,0){$2$}}
\put(9,-3){\makebox(0,0){$a_1$}}
\put(29,-3){\makebox(0,0){$b_1$}}
\put(51,-3){\makebox(0,0){$b_2$}}
\put(71,-3){\makebox(0,0){$a_2$}}
\put(40,23){\makebox(0,0){$c_1$}}

\put(10,0){\line(1,0){60}}
\put(30,0){\line(1,2){10}}
\put(50,0){\line(-1,2){10}}

\end{picture}

\end{center}

\caption{} 

\label{fig:Fhat12b-32b}

\end{figure}

\credit{Case 4} \emph{$k$ is connected to $b_2$ not $b_1$.}

We note that if the edge $[k,b_2]$ is weighted 2,
then the subdiagram with the vertices $\{k,b_1,b_2\}$ is of type $C_n^{(1)}$.
To consider the remaining subcases we assume that the edge $[k,b_2]$
is weightless.

\credit{Subcase 4.1}\emph{$k$ is not connected to $a_2$.}

Then the subdiagram with the vertices 
$\{c_1,b_2,k,a_2\}$ is of type $B_n^{(1)}$ 
(Fig.~\ref{fig:Fhat12b-41}).

\begin{figure}[ht]

\setlength{\unitlength}{1.8pt}

\begin{center}

\begin{picture}(70,40)(-10,-20)

\put(30,-20){\line(1,1){20}}
\put(30,-20){\line(-1,1){20}}

\put(18,-11){\makebox(0,0){$2$}}
\put(30,-20){\circle*{2.0}}
\put(30,-24){\makebox(0,0){$k$}}

\thicklines

\put(10,0){\circle*{2.0}}
\put(30,0){\circle*{2.0}}
\put(50,0){\circle*{2.0}}
\put(70,0){\circle*{2.0}}
\put(40,20){\circle*{2.0}}

\put(40,-2){\makebox(0,0){$2$}}
\put(47,10){\makebox(0,0){$2$}}
\put(9,-3){\makebox(0,0){$a_1$}}
\put(29,-3){\makebox(0,0){$b_1$}}
\put(51,-3){\makebox(0,0){$b_2$}}
\put(71,-3){\makebox(0,0){$a_2$}}
\put(40,23){\makebox(0,0){$c_1$}}

\put(10,0){\line(1,0){60}}
\put(30,0){\line(1,2){10}}
\put(50,0){\line(-1,2){10}}

\end{picture}

\end{center}

\caption{} 

\label{fig:Fhat12b-41}

\end{figure}

\credit{Subcase 4.2}\emph{$k$ is connected to $a_2$}

Then the subdiagram with the vertices 
$\{k,a_2,b_2,b_1,c_1\}$ is of type $F_4^{(1)}(3^{2};0)$ 
(Fig.~\ref{fig:Fhat12b-42}).
\endproof

\begin{figure}[ht]

\setlength{\unitlength}{1.8pt}

\begin{center}

\begin{picture}(70,40)(-10,-20)

\put(50,-20){\line(1,1){20}}
\put(50,-20){\line(0,1){20}}
\put(50,-20){\line(-1,0){5}}

\put(50,-20){\circle*{2.0}}
\put(50,-24){\makebox(0,0){$k$}}

\thicklines

\put(10,0){\circle*{2.0}}
\put(30,0){\circle*{2.0}}
\put(50,0){\circle*{2.0}}
\put(70,0){\circle*{2.0}}
\put(40,20){\circle*{2.0}}

\put(40,-2){\makebox(0,0){$2$}}
\put(47,10){\makebox(0,0){$2$}}
\put(9,-3){\makebox(0,0){$a_1$}}
\put(29,-3){\makebox(0,0){$b_1$}}
\put(51,-3){\makebox(0,0){$b_2$}}
\put(71,-3){\makebox(0,0){$a_2$}}
\put(40,23){\makebox(0,0){$c_1$}}

\put(10,0){\line(1,0){60}}
\put(30,0){\line(1,2){10}}
\put(50,0){\line(-1,2){10}}

\end{picture}

\end{center}

\caption{} 

\label{fig:Fhat12b-42}

\end{figure}

\begin{lemma}
\label{lem:caseFhat23}

If $\Gamma$ is of type $F_4^{(1)}(3^{2};0;1^1)$, then Lemma~\ref{lem:xxx} holds.

\end{lemma}

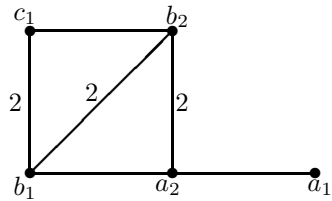
\begin{figure}[ht]

\setlength{\unitlength}{1.8pt}

\begin{center}

\begin{picture}(85,50)(-10,-10)

\thicklines

\put(70,0){\circle*{2.0}}
\put(40,0){\circle*{2.0}}
\put(10,0){\circle*{2.0}}
\put(40,30){\circle*{2.0}}
\put(10,30){\circle*{2.0}}

\put(7,15){\makebox(0,0){$2$}}
\put(23,17){\makebox(0,0){$2$}}
\put(42,15){\makebox(0,0){$2$}}
\put(9,-3){\makebox(0,0){$b_1$}}
\put(39,-3){\makebox(0,0){$a_2$}}
\put(71,-3){\makebox(0,0){$a_1$}}
\put(9,33){\makebox(0,0){$c_1$}}
\put(41,33){\makebox(0,0){$b_2$}}

\put(10,0){\line(1,0){60}}
\put(10,0){\line(0,1){30}}
\put(10,0){\line(1,1){30}}
\put(40,0){\line(0,1){30}}
\put(40,30){\line(-1,0){30}}

\end{picture}

\end{center}

\caption{The diagram $F_4^{(1)}(3^{2};0;1^{1})$}

\label{fig:Fhat23}

\end{figure}

\proof We assume that $\Gamma$ is indexed as in Fig.~\ref{fig:Fhat23}.

\credit{Case 1} \emph{$k$ is not connected to any vertex in $\{b_1,b_2\}$.}

\credit{Subcase 1.1}\emph{$k$ is not connected to $c_1$.}

In this case the vertex $k$ is connected to both $a_1$ and $a_2$.
If the edge $[k,a_2]$ (hence the edge $[k,a_1]$ is weighted 2
, then the subdiagram with the vertices $\{k,a_2,b_2\}$ is of type $C_n^{(1)}$
; otherwise the subdiagram with the vertices $\{k,a_2,b_1,b_2,c_1\}$ is of type $F_4^{(1)}(3^{2};0;1^1)$ 
(Fig.~\ref{fig:Fhat23-11}).

\begin{figure}[ht]

\setlength{\unitlength}{1.8pt}

\begin{center}

\begin{picture}(70,50)(-10,-10)

\put(70,30){\line(0,-1){30}}
\put(70,30){\line(-1,-1){30}}
\put(70,30){\circle*{2.0}}
\put(71,33){\makebox(0,0){$k$}}

\thicklines

\put(70,0){\circle*{2.0}}
\put(40,0){\circle*{2.0}}
\put(10,0){\circle*{2.0}}
\put(40,30){\circle*{2.0}}
\put(10,30){\circle*{2.0}}

\put(7,15){\makebox(0,0){$2$}}
\put(23,17){\makebox(0,0){$2$}}
\put(42,15){\makebox(0,0){$2$}}
\put(9,-3){\makebox(0,0){$b_1$}}
\put(39,-3){\makebox(0,0){$a_2$}}
\put(71,-3){\makebox(0,0){$a_1$}}
\put(9,33){\makebox(0,0){$c_1$}}
\put(41,33){\makebox(0,0){$b_2$}}

\put(10,0){\line(1,0){60}}
\put(10,0){\line(0,1){30}}
\put(10,0){\line(1,1){30}}
\put(40,0){\line(0,1){30}}
\put(40,30){\line(-1,0){30}}

\end{picture}

\end{center}

\caption{} 

\label{fig:Fhat23-11}

\end{figure}

\credit{Subcase 1.2}\emph{$k$ is connected to $c_1$.}

\begin{itemize}
\item
If the edge $[k,c_1]$ is weighted 2, then the subdiagram with the vertices 
$\{k,c_1,b_2\}$ is of type $C_n^{(1)}$
\item
If the edge $[k,c_1]$ is weightless, then $k$ is contained in a non-oriented cycle
(Fig.~\ref{fig:Fhat23-12}). 
\end{itemize}

\begin{figure}[ht]

\setlength{\unitlength}{1.8pt}

\begin{center}

\begin{picture}(70,80)(-10,-10)

\put(70,60){\line(0,-1){10}}
\put(70,60){\line(-2,-1){60}}
\put(70,60){\circle*{2.0}}
\put(73,60){\makebox(0,0){$k$}}
\put(56,28){\makebox(0,0){$2$}}
\put(70,60){\line(-1,-2){30}}

\thicklines

\put(70,0){\circle*{2.0}}
\put(40,0){\circle*{2.0}}
\put(10,0){\circle*{2.0}}
\put(40,30){\circle*{2.0}}
\put(10,30){\circle*{2.0}}

\put(7,15){\makebox(0,0){$2$}}
\put(23,17){\makebox(0,0){$2$}}
\put(42,15){\makebox(0,0){$2$}}
\put(9,-3){\makebox(0,0){$b_1$}}
\put(39,-3){\makebox(0,0){$a_2$}}
\put(71,-3){\makebox(0,0){$a_1$}}
\put(9,33){\makebox(0,0){$c_1$}}
\put(41,33){\makebox(0,0){$b_2$}}

\put(10,0){\line(1,0){60}}
\put(10,0){\line(0,1){30}}
\put(10,0){\line(1,1){30}}
\put(40,0){\line(0,1){30}}
\put(40,30){\line(-1,0){30}}

\end{picture}

\end{center}

\caption{} 

\label{fig:Fhat23-12}

\end{figure}

\credit{Case 2} \emph{$k$ is connected to both $b_1$ and $b_2$.}

If $k$ is connected to a vertex in $C=\{c_1,a_2,a_1\}$, 
then $k$ is contained in a non-oriented cycle.
Let us now assume that $k$ is not connected to 
any vertex in $C$ and consider the subcases.

\credit{Subcase 2.1}\emph{The edge $[k,b_1]$ is weighted 2.}

Then the subdiagram with the vertices
$\{k,b_1,c_1\}$ is of type $C_n^{(1)}$. 

\credit{Subcase 2.2}\emph{The edge $[k,b_1]$ is weightless.}

Then the edge $[k,b_2]$ is weighted 2 and the 
subdiagram with the vertices $\{k,b_2,a_2\}$
is of type $C_n^{(1)}$.

\credit{Case 3} \emph{$k$ is connected to $b_1$ not $b_2$.}

We note that if the edge $[k,b_1]$ is weighted 2,
then the subdiagram with the vertices $\{k,b_1,b_2\}$ is of type $C_n^{(1)}$.
To consider the remaining subcases we assume that the edge $[k,b_1]$
is weightless.

\credit{Subcase 3.1}\emph{$k$ is connected to $a_1$ or $a_2$.}

If $k$ is connected to $a_2$,
then 
the subdiagram with the vertices 
$\{k,b_1,b_2,a_2\}$ is of type $B_n^{(1)}(r)$, 
otherwise the subdiagram with the vertices 
$\{k,b_1,b_2,a_2,a_1\}$ is of type $B_n^{(1)}(r)$ 
(Fig.~\ref{fig:Fhat23-31}).

\begin{figure}[ht]

\setlength{\unitlength}{1.8pt}

\begin{center}

\begin{picture}(70,60)(-10,-30)

\put(40,-30){\line(-1,1){30}}
\put(40,-30){\line(1,1){30}}
\put(40,-30){\circle*{2.0}}
\put(40,-33){\makebox(0,0){$k$}}

\thicklines

\put(70,0){\circle*{2.0}}
\put(40,0){\circle*{2.0}}
\put(10,0){\circle*{2.0}}
\put(40,30){\circle*{2.0}}
\put(10,30){\circle*{2.0}}

\put(7,15){\makebox(0,0){$2$}}
\put(23,17){\makebox(0,0){$2$}}
\put(42,15){\makebox(0,0){$2$}}
\put(9,-3){\makebox(0,0){$b_1$}}
\put(39,-3){\makebox(0,0){$a_2$}}
\put(71,-3){\makebox(0,0){$a_1$}}
\put(9,33){\makebox(0,0){$c_1$}}
\put(41,33){\makebox(0,0){$b_2$}}

\put(10,0){\line(1,0){60}}
\put(10,0){\line(0,1){30}}
\put(10,0){\line(1,1){30}}
\put(40,0){\line(0,1){30}}
\put(40,30){\line(-1,0){30}}

\end{picture}

\end{center}

\caption{} 

\label{fig:Fhat23-31}

\end{figure}

\credit{Subcase 3.2}\emph{k is not connected $a_1$ and not connected $a_2$.}

Then $k$ is connected to $c_1$ and the edge $[k,c_1]$ is weighted 2.
We note that the subdiagram with the vertices 
$\{k,c_1,b_2,a_2\}$ is of type $C_n^{(1)}$ 
(Fig.~\ref{fig:Fhat23-32}).

\begin{figure}[ht]

\setlength{\unitlength}{1.8pt}

\begin{center}

\begin{picture}(70,60)(-10,-10)

\put(-20,0){\line(1,0){30}}
\put(-20,0){\line(1,1){30}}
\put(-20,0){\circle*{2.0}}
\put(-23,0){\makebox(0,0){$k$}}
\put(-5,17){\makebox(0,0){$2$}}

\thicklines

\put(70,0){\circle*{2.0}}
\put(40,0){\circle*{2.0}}
\put(10,0){\circle*{2.0}}
\put(40,30){\circle*{2.0}}
\put(10,30){\circle*{2.0}}

\put(7,15){\makebox(0,0){$2$}}
\put(23,17){\makebox(0,0){$2$}}
\put(42,15){\makebox(0,0){$2$}}
\put(9,-3){\makebox(0,0){$b_1$}}
\put(39,-3){\makebox(0,0){$a_2$}}
\put(71,-3){\makebox(0,0){$a_1$}}
\put(9,33){\makebox(0,0){$c_1$}}
\put(41,33){\makebox(0,0){$b_2$}}

\put(10,0){\line(1,0){60}}
\put(10,0){\line(0,1){30}}
\put(10,0){\line(1,1){30}}
\put(40,0){\line(0,1){30}}
\put(40,30){\line(-1,0){30}}

\end{picture}

\end{center}

\caption{} 

\label{fig:Fhat23-32}

\end{figure}

\credit{Case 4} \emph{$k$ is connected to $b_2$ not $b_1$.}

We note that if the edge $[k,b_2]$ is weighted 2,
then the subdiagram with the vertices $\{k,b_2,b_1\}$ is of type $C_n^{(1)}$.
To consider the remaining subcases we assume that the edge $[k,b_2]$
is weightless.

\credit{Subcase 4.1}\emph{$k$ is connected to $c_1$.}

Then the subdiagram with the vertices 
$\{b_1,c_1,b_2,k\}$ is of type $B_n^{(1)}(r)$ 
(Fig.~\ref{fig:Fhat23-41}).

\begin{figure}[ht]
\setlength{\unitlength}{1.8pt}

\begin{center}

\begin{picture}(70,80)(-10,-10)

\put(40,60){\line(0,-1){30}}
\put(40,60){\line(-1,-1){30}}
\put(40,60){\circle*{2.0}}
\put(40,63){\makebox(0,0){$k$}}
\put(40,60){\line(2,-1){10}}

\thicklines

\put(70,0){\circle*{2.0}}
\put(40,0){\circle*{2.0}}
\put(10,0){\circle*{2.0}}
\put(40,30){\circle*{2.0}}
\put(10,30){\circle*{2.0}}

\put(7,15){\makebox(0,0){$2$}}
\put(23,17){\makebox(0,0){$2$}}
\put(42,15){\makebox(0,0){$2$}}
\put(9,-3){\makebox(0,0){$b_1$}}
\put(39,-3){\makebox(0,0){$a_2$}}
\put(71,-3){\makebox(0,0){$a_1$}}
\put(9,33){\makebox(0,0){$c_1$}}
\put(42,33){\makebox(0,0){$b_2$}}

\put(10,0){\line(1,0){60}}
\put(10,0){\line(0,1){30}}
\put(10,0){\line(1,1){30}}
\put(40,0){\line(0,1){30}}
\put(40,30){\line(-1,0){30}}

\end{picture}

\end{center}

\caption{} 

\label{fig:Fhat23-41}

\end{figure}

\credit{Subcase 4.2}\emph{k is not connected to $c_1$}

If $k$ is connected to $a_1$, 
then the subdiagram with the vertices 
$\{b_1,b_2,k,a_1\}$ is of type $C_n^{(1)}$ 
(Fig.~\ref{fig:Fhat23-42a});
otherwise 
the subdiagram with the vertices 
$\{c_1,b_1,a_2,k\}$ is of type $C_n^{(1)}$ 
(Fig.~\ref{fig:Fhat23-42b}).
\endproof

\begin{figure}[ht]

\setlength{\unitlength}{1.8pt}

\begin{center}

\begin{picture}(70,60)(-10,-10)

\put(70,30){\line(0,-1){30}}
\put(70,30){\line(-1,0){30}}
\put(70,30){\circle*{2.0}}
\put(70,33){\makebox(0,0){$k$}}
\put(73,15){\makebox(0,0){$2$}}
\put(70,30){\line(-1,-1){10}}

\thicklines

\put(70,0){\circle*{2.0}}
\put(40,0){\circle*{2.0}}
\put(10,0){\circle*{2.0}}
\put(40,30){\circle*{2.0}}
\put(10,30){\circle*{2.0}}

\put(7,15){\makebox(0,0){$2$}}
\put(23,17){\makebox(0,0){$2$}}
\put(42,15){\makebox(0,0){$2$}}
\put(9,-3){\makebox(0,0){$b_1$}}
\put(39,-3){\makebox(0,0){$a_2$}}
\put(71,-3){\makebox(0,0){$a_1$}}
\put(9,33){\makebox(0,0){$c_1$}}
\put(42,33){\makebox(0,0){$b_2$}}

\put(10,0){\line(1,0){60}}
\put(10,0){\line(0,1){30}}
\put(10,0){\line(1,1){30}}
\put(40,0){\line(0,1){30}}
\put(40,30){\line(-1,0){30}}

\end{picture}

\end{center}

\caption{} 

\label{fig:Fhat23-42a}

\end{figure}

\begin{figure}[ht]

\setlength{\unitlength}{1.8pt}

\begin{center}

\begin{picture}(70,60)(-10,-10)

\put(70,30){\line(-1,0){30}}
\put(70,30){\circle*{2.0}}
\put(70,33){\makebox(0,0){$k$}}
\put(58,15){\makebox(0,0){$2$}}
\put(70,30){\line(-1,-1){30}}

\thicklines

\put(70,0){\circle*{2.0}}
\put(40,0){\circle*{2.0}}
\put(10,0){\circle*{2.0}}
\put(40,30){\circle*{2.0}}
\put(10,30){\circle*{2.0}}

\put(7,15){\makebox(0,0){$2$}}
\put(23,17){\makebox(0,0){$2$}}
\put(42,15){\makebox(0,0){$2$}}
\put(9,-3){\makebox(0,0){$b_1$}}
\put(39,-3){\makebox(0,0){$a_2$}}
\put(71,-3){\makebox(0,0){$a_1$}}
\put(9,33){\makebox(0,0){$c_1$}}
\put(42,33){\makebox(0,0){$b_2$}}

\put(10,0){\line(1,0){60}}
\put(10,0){\line(0,1){30}}
\put(10,0){\line(1,1){30}}
\put(40,0){\line(0,1){30}}
\put(40,30){\line(-1,0){30}}

\end{picture}

\end{center}

\caption{} 

\label{fig:Fhat23-42b}

\end{figure}

\begin{lemma}
\label{lem:caseFhat224}

If $\Gamma$ is of type $F_4^{(1)}(4^{1};3^1)$, then Lemma~\ref{lem:xxx} holds.

\end{lemma}

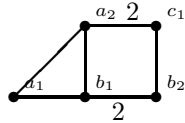
\begin{figure}[ht]

\setlength{\unitlength}{1.8pt}

\begin{center}

\begin{picture}(70,30)(-10,-10)

\thicklines

\put(10,0){\circle*{2.0}$^{a_1}$}
\put(25,15){\circle*{2.0}$^{a_2}$}
\put(25,0){\circle*{2.0}$^{b_1}$}
\put(40,0){\circle*{2.0}$^{b_2}$}
\put(40,15){\circle*{2.0}$^{c_1}$}

\put(32,-3){\makebox(0,0){$2$}}
\put(35,18){\makebox(0,0){$2$}}

\put(10,0){\line(1,0){30}}
\put(10,0){\line(1,1){15}}
\put(25,0){\line(0,1){15}}
\put(25,15){\line(1,0){15}}
\put(40,0){\line(0,1){15}}

\end{picture}

\end{center}

\caption{The diagram $F_4^{(1)}(4^{1};3^1) \,$.} 

\label{fig:Fhat224}

\end{figure}

\proof We assume that $\Gamma$ is indexed as in Fig.~\ref{fig:Fhat224}.

\credit{Case 1} \emph{$k$ is not connected to any vertex in $\{b_1,b_2\}$.}

\credit{Subcase 1.1}\emph{$k$ is not connected to $c_1$.}

In this case the vertex $k$ is connected to both $a_1$ and $a_2$.
If the edge $[k,a_2]$ (hence the edge $[k,a_1]$) is weighted 2
, then the subdiagram with the vertices $\{k,a_2,c_1\}$ is of type $B_n^{(1)}(r)$ 
; otherwise the subdiagram with the vertices $\{k,a_2,b_1,c_1\}$ is of type $B_n^{(1)}$ .

\credit{Subcase 1.2}\emph{$k$ is connected to $c_1$.}
\begin{itemize}
\item
If the edge $[k,c_1]$ is weighted 2, then the subdiagram with the vertices 
$\{k,c_1,b_2,b_1\}$ is of type $B_n^{(1)}$.
\item
Suppose that the edge $[k,c_1]$ is weightless.
If $k$ is connected to $a_2$, then
the subdiagram with the vertices $\{k,a_2,b_1,b_2\}$ is of type $C_n^{(1)}$;
otherwise $k$ is contained in a non-oriented cycle. 
\end{itemize}

\credit{Case 2} \emph{$k$ is connected to both $b_1$ and $b_2$.}

\credit{Subcase 2.1}\emph{The edge $[k,b_1]$ is weighted 2.}

We first note that the edge $[k,b_2]$ is weightless. 
\begin{itemize}
\item[(i)]
Suppose that $k$ is connected to $c_1$ or $a_2$. 
Then $k$ is contained in a non-oriented cycle.
\item[(ii)]
Suppose that (i) does not hold.
Then the subdiagram with the vertices $\{k,b_1,a_2,c_1\}$
is of type $C_n^{(1)}$. 
\end{itemize}

\credit{Subcase 2.2}\emph{The edge $[k,b_2]$ is weighted 2.}

If $k$ is connected to $c_1$, then $k$ is contained in a non-oriented cycle.
Let us now assume that $k$ is not connected to $c_1$ and consider the subcases.
\begin{itemize}
\item[(i)] 
Suppose that $k$ is not connected to $c_1$ or $a_2$. 
Then $k$ is contained in a non-oriented cycle.
\item[(ii)]
Suppose that (i) does not hold.
Then the subdiagram with the vertices $\{k,b_2,c_1,a_2\}$
is of type $C_n^{(1)}$. 
\end{itemize}

\credit{Case 3} \emph{$k$ is connected to $b_1$ not $b_2$.}

We immediately notice that if the edge $[k,b_1]$ is weighted 2,
then the subdiagram with the vertices $\{k,b_1,b_2\}$ is of type $C_n^{(1)}$.
To consider the remaining subcases we assume that the edge $[k,b_1]$
is weightless.

\credit{Subcase 3.1}\emph{k is not connected to $a_1$.}

Then the subdiagram with the vertices $\{b_2,b_1,k,a_1\}$ is of type
$B_n^{(1)}$. 

\credit{Subcase 3.2}\emph{k is connected to $a_1$}

(i) If $k$ is connected to $c_1$ or $a_2$, then it is contained
in a non-oriented cycle. 

(ii) Let us assume that (i) does not hold. Then
the subdiagram with the vertices $\{b_2,b_1,k,a_2\}$ is of type
$B_n^{(1)}$. 

\credit{Case 4} \emph{$k$ is connected to $b_2$ not $b_1$.}

We immediately notice that if the edge $(k,b_2)$ is weighted 2,
then the subdiagram with the vertices $\{k,b_1,b_2\}$ is of type $C_n^{(1)}$.
To consider the remaining subcases we assume that the edge $[k,b_2]$
is weightless.

\credit{Subcase 4.1}\emph{k is connected to $c_1$.}

If $k$ us connected to $a_1$ or $a_2$ then $k$ is contained 
in a non-oriented cycle; otherwise
the subdiagram with the vertices 
$\{b_1,b_2,k,c_1,a_2\}$ is of type $F_4^{(1)}(4^{1};3^1)$. 

\credit{Subcase 4.2}\emph{k is not connected to $c_1$}

Then the subdiagram with the vertices 
$\{b_1,b_2,k,c_1\}$ is of type $C_n^{(1)}$. 
\endproof

\subsection{Proof of Theorem~\ref{cor:main}}
\label{subsec:cor-main}
In view of Theorem~\ref{th:main}, the theorem is the same as 
Lemma~\ref{lem:min-eDynkin}, which we already proved.

\begin{remark}

\label{th:ps-main}
\rm{
By enlarging the set of representatives for mutation classes of minimal 2-infinite diagrams, we obtain the following recognition criterion for 2-finite diagrams:}

Using at most 9 mutations, any 2-infinite diagram can be transformed into a diagram which contains a either a diagram from Table 1 
or one of the following extended Dynkin diagrams: $E_6^{(1)}$, $E_7^{(1)}$, $E_8^{(1)}$, $F_4^{(1)}$. 

\rm{
To prove this statement, one may note that any (minimal 2-infinite) diagram which is not in Table 1 has at most 9 vertices. We also note that it is enough to prove the theorem for any diagram $\Gamma$ (which is not in Table 1) from our list, because any 2-infinite diagram contains a subdiagram from our list (Theorem~\ref{th:main}). This is manageable despite the relatively large number of diagrams in Tables 2-6. Although we checked this by hand we do not give the proof here.
}

\end{remark}

\pagebreak

\section{The list of minimal 2-infinite diagrams}
\label{subsec:list}
\emph{Unless otherwise stated, any diagram in this list is assumed to have an arbitrary orientation which does not contain any non-oriented cycle.}
\[ 
\begin{array}{ccl} 
A_n^{(1)}
&& 
\setlength{\unitlength}{1.0pt} 
\begin{picture}(140,60)(0,-2) 
\put(60,0){\circle*{2.0}} 
\put(60,20){\circle*{2.0}}
\put(40,40){\circle*{2.0}}
\put(20,40){\circle*{2.0}}
\put(0,20){\circle*{2.0}}
\put(20,-20){\circle*{2.0}}
\put(0,0){\circle*{2.0}}
\put(40,-20){\circle*{2.0}}

\put(60,0){\line(-1,-1){20}}
\put(60,0){\line(0,1){20}}
\put(60,20){\line(-1,1){20}}
\put(40,40){\line(-1,0){20}}
\put(20,40){\line(-1,-1){20}}
\put(0,20){\line(0,-1){20}}
\put(0,0){\line(1,-1){20}}
\put(20,-20){\line(1,0){20}}

\put(104,0){\makebox(0,0){non-oriented}}
\end{picture} 
\\[.3in] 
B_n^{(1)} 
&& 
\setlength{\unitlength}{1.0pt} 
\begin{picture}(140,20)(0,-2) 
\put(20,0){\line(1,0){120}} 
\put(0,10){\line(2,-1){20}} 
\put(0,-10){\line(2,1){20}} 
\multiput(20,0)(20,0){7}{\circle*{2}} 
\put(0,10){\circle*{2}} 
\put(0,-10){\circle*{2}} 
\put(130,4){\makebox(0,0){$2$}} 
\end{picture} 
\\[.1in]
B_n^{(1)}(m,r)
&& 
\setlength{\unitlength}{1.0pt} 
\begin{picture}(140,60)(20,-2) 
\put(100,0){\circle*{2.0}}

\put(80,20){\circle*{2.0}}
\put(80,0){\circle*{2.0}} 
\put(120,0){\circle*{2.0}}
\put(140,0){\circle*{2.0}}
\put(160,0){\circle*{2.0}}
\put(180,0){\circle*{2.0}}
\put(200,0){\circle*{2.0}}
\put(80,20){\circle*{2.0}}
\put(60,40){\circle*{2.0}}
\put(40,40){\circle*{2.0}}
\put(20,20){\circle*{2.0}}
\put(40,-20){\circle*{2.0}}
\put(20,0){\circle*{2.0}}
\put(60,-20){\circle*{2.0}}

\put(80,0){\line(1,0){120}}
\put(80,0){\line(-1,-1){20}}
\put(80,0){\line(0,1){20}}
\put(100,0){\line(-1,1){20}}
\put(80,20){\line(-1,1){20}}
\put(60,40){\line(-1,0){20}}
\put(40,40){\line(-1,-1){20}}
\put(20,20){\line(0,-1){20}}
\put(20,0){\line(1,-1){20}}
\put(40,-20){\line(1,0){20}}

\put(190,-4){\makebox(0,0){2}}
\put(81,-4){\makebox(0,0){$b_1$}}
\put(60,43){\makebox(0,0){$b_3$}}
\put(17,0){\makebox(0,0){$b_i$}}
\put(60,-23){\makebox(0,0){$b_r$}}
\put(81,23){\makebox(0,0){$b_2$}}
\put(-10,-10){\makebox(0,0){{\small $m\geq 1,~r\geq 3$}}}

\put(103,-4){\makebox(0,0){$a_1$}}
\put(123,-4){\makebox(0,0){$a_2$}}
\put(177,-4){\makebox(0,0){$a_m$}}

\end{picture} 
\\[.4in] 
B_n^{(1)}(r)
&& 
\setlength{\unitlength}{1.0pt} 
\begin{picture}(140,60)(20,-2) 
\put(80,0){\circle*{2.0}}

\put(80,0){\circle*{2.0}} 
\put(80,20){\circle*{2.0}}
\put(60,40){\circle*{2.0}}
\put(40,40){\circle*{2.0}}
\put(20,20){\circle*{2.0}}
\put(40,-20){\circle*{2.0}}
\put(20,0){\circle*{2.0}}
\put(60,-20){\circle*{2.0}}

\put(80,0){\line(1,0){20}}
\put(80,0){\line(-1,-1){20}}
\put(80,0){\line(0,1){20}}
\put(100,0){\line(-1,1){20}}
\put(80,20){\line(-1,1){20}}
\put(60,40){\line(-1,0){20}}
\put(40,40){\line(-1,-1){20}}
\put(20,20){\line(0,-1){20}}
\put(20,0){\line(1,-1){20}}
\put(40,-20){\line(1,0){20}}

\put(90,-4){\makebox(0,0){$2$}}
\put(93,13){\makebox(0,0){$2$}}
\put(81,-4){\makebox(0,0){$b_1$}}
\put(60,43){\makebox(0,0){$b_3$}}
\put(17,0){\makebox(0,0){$b_i$}}
\put(60,-23){\makebox(0,0){$b_r$}}
\put(81,23){\makebox(0,0){$b_2$}}
\put(104,0){\makebox(0,0){$c_1$}}
\put(-10,-10){\makebox(0,0){{$\small r\geq 3$}}}

\end{picture} 
\\[.4in] 
C_n^{(1)}
&& 
\setlength{\unitlength}{1.0pt} 
\begin{picture}(140,17)(0,-2) 
\put(0,0){\line(1,0){140}} 
\multiput(0,0)(20,0){8}{\circle*{2}} 
\put(10,4){\makebox(0,0){$2$}} 
\put(130,4){\makebox(0,0){$2$}} 
\end{picture} 
\\[.2in] 
D_n^{(1)} 
&& 
\setlength{\unitlength}{1.0pt} 
\begin{picture}(140,17)(0,-2) 
\put(20,0){\line(1,0){100}} 
\put(0,10){\line(2,-1){20}} 
\put(0,-10){\line(2,1){20}} 
\put(120,0){\line(2,-1){20}} 
\put(120,0){\line(2,1){20}} 
\multiput(20,0)(20,0){6}{\circle*{2}} 
\put(0,10){\circle*{2}} 
\put(0,-10){\circle*{2}} 
\put(140,10){\circle*{2}} 
\put(140,-10){\circle*{2}} 
\end{picture} 
\\
&& 
\setlength{\unitlength}{1.0pt} 
\begin{picture}(140,17)(0,-2) 
\end{picture} 
\\

D_n^{(1)}(m,r):
&& 
\setlength{\unitlength}{1.0pt}
\begin{picture}(180,63)(20,-2)
\put(40,0){\circle*{2.0}}
\put(20,10){\circle*{2.0}}  \put(24,13){\makebox(0,0){$b_2$}}
\put(20,-10){\circle*{2.0}} \put(24,-13){\makebox(0,0){$b_1$}}
\put(60,0){\circle*{2.0}}
\put(80,0){\circle*{2.0}}
\put(100,0){\circle*{2.0}}
\put(120,0){\circle*{2.0}}
\put(140,0){\circle*{2.0}}
\put(160,20){\circle*{2.0}}
\put(160,0){\circle*{2.0}}
\put(180,40){\circle*{2.0}}
\put(200,40){\circle*{2.0}}
\put(220,20){\circle*{2.0}}
\put(220,0){\circle*{2.0}}
\put(180,-20){\circle*{2.0}}
\put(200,-20){\circle*{2.0}}

\put(158,23){\makebox(0,0){$c_2$}}
\put(158,-3){\makebox(0,0){$c_1$}}
\put(178,43){\makebox(0,0){$c_3$}}
\put(178,-23){\makebox(0,0){$c_r$}}
\put(-10,-10){\makebox(0,0){{\small $m\geq 1,~r\geq 3$}}}

\put(40,0){\line(1,0){120}}
\put(40,0){\line(-2,-1){20}}
\put(40,0){\line(-2,1){20}}
\put(140,0){\line(1,1){20}}
\put(160,0){\line(0,1){20}}
\put(160,0){\line(1,-1){20}}
\put(160,20){\line(1,1){20}}
\put(180,40){\line(1,0){20}}
\put(200,40){\line(1,-1){20}}
\put(220,20){\line(0,-1){20}}
\put(220,0){\line(-1,-1){20}}
\put(200,-20){\line(-1,0){20}}

\put(138,-3){\makebox(0,0){$a_m$}}
\put(43,-3){\makebox(0,0){$a_1$}}

\end{picture}
\\[.3in] 
{\small D_n^{(1)}(m,r,s):}
&& 
\setlength{\unitlength}{1.0pt} 
\begin{picture}(224,64)(20,-2) 

\put(-10,-10){\makebox(0,0){{\small $m\geq 1,~r,s\geq 3$}}}

\put(82,-4){\makebox(0,0){$b_1$}} 
\put(82,23){\makebox(0,0){$b_2$}} 
\put(17,23){\makebox(0,0){$b_i$}} 
\put(63,-23){\makebox(0,0){$b_r$}}

\put(158,-4){\makebox(0,0){$c_1$}} 
\put(158,23){\makebox(0,0){$c_2$}} 
\put(224,20){\makebox(0,0){$c_i$}} 
\put(178,-23){\makebox(0,0){$c_s$}} 

\put(103,-4){\makebox(0,0){$a_1$}} 
\put(140,-4){\makebox(0,0){$a_m$}}

\put(40,-20){\line(1,0){20}} 
\put(40,40){\line(1,0){20}} 
\put(20,0){\line(0,1){20}} 
\put(80,0){\line(0,1){20}} 
\put(20,0){\line(1,-1){20}} 
\put(60,40){\line(1,-1){20}} 
\put(20,20){\line(1,1){20}} 
\put(60,-20){\line(1,1){20}} 
\put(80,20){\line(1,-1){20}}
\put(140,0){\line(1,1){20}}
\put(140,0){\line(1,0){20}}
\put(160,0){\line(0,1){20}} 
\put(160,20){\line(1,1){20}} 
\put(180,40){\line(1,0){20}}
\put(200,40){\line(1,-1){20}}
\put(220,20){\line(0,-1){20}}
\put(220,0){\line(-1,-1){20}}
\put(200,-20){\line(-1,0){20}}
\put(180,-20){\line(-1,1){20}}

\put(40,-20){\line(1,0){12}} 
\put(60,40){\line(-1,0){12}} 
\put(20,20){\line(0,-1){12}} 
\put(80,0){\line(0,1){12}} 
\put(20,0){\line(1,-1){12}} 
\put(80,20){\line(-1,1){12}} 
\put(40,40){\line(-1,-1){12}} 
\put(60,-20){\line(1,1){12}}

\multiput(40,-20)(20,0){2}{\circle*{2}} 
\multiput(40,40)(20,0){2}{\circle*{2}} 
\multiput(20,0)(60,0){2}{\circle*{2}} 
\multiput(20,20)(60,0){2}{\circle*{2}} 


\put(110,0){\circle*{2}} 
\put(130,0){\circle*{2}} 
\put(160,20){\circle*{2}}
\put(160,0){\circle*{2}}
\put(180,-20){\circle*{2}}
\put(180,40){\circle*{2}} 
\put(200,-20){\circle*{2}}
\put(200,40){\circle*{2}} 
\put(220,0){\circle*{2}}
\put(220,20){\circle*{2}} 

\multiput(100,0)(20,0){3}{\circle*{2}} 
\put(80,0){\line(1,0){60}} 
\end{picture} 
\\[.25in] 
D_n^{(1)}(r): r\geq 3 
&& 
\setlength{\unitlength}{1.0pt} 
\begin{picture}(140,60)(20,-2) 
\put(-10,-10){\makebox(0,0){{\small $r\geq 3$}}}

\put(40,-30){\line(1,0){20}} 
\put(40,30){\line(1,0){20}} 
\put(20,-10){\line(0,1){20}} 
\put(80,-10){\line(0,1){20}} 
\put(20,-10){\line(1,-1){20}} 
\put(60,30){\line(1,-1){20}} 
\put(20,10){\line(1,1){20}} 
\put(60,-30){\line(1,1){20}} 
\put(80,10){\line(1,-1){20}}

\put(20,10){\circle*{2}} 
\put(20,-10){\circle*{2}} 
\put(40,-30){\circle*{2}}
\put(40,30){\circle*{2}} 
\put(60,-30){\circle*{2}} 
\put(60,-10){\circle*{2}} 
\put(60,30){\circle*{2}} 
\put(80,-10){\circle*{2}} 
\put(80,10){\circle*{2}} 
\put(100,-10){\circle*{2}} 

\put(58,-14){\makebox(0,0){$a_1$}} 
\put(104,-14){\makebox(0,0){$c_1$}} 
\put(82,-15){\makebox(0,0){$b_1$}} 
\put(82,15){\makebox(0,0){$b_2$}} 
\put(17,15){\makebox(0,0){$b_i$}} 
\put(63,-34){\makebox(0,0){$b_r$}}

\put(80,-10){\line(1,0){20}} 
\put(80,-10){\line(-1,0){20}}  
\put(60,-10){\line(1,1){20}} 
\end{picture}
\\[.4in]
I_2(a)
&& 
\setlength{\unitlength}{1.0pt} 
\begin{picture}(65,17)(20,-2) 
\put(20,0){\line(1,0){30}} 
\put(20,0){\circle*{2}} 
\put(50,0){\circle*{2}} 
\put(35,4){\makebox(0,0){$a$}} 
\put(80,0){\makebox(0,0){$a\geq 4$}} 
\end{picture}\\[.2in]
&& 
\setlength{\unitlength}{1.0pt} 
\begin{picture}(140,17)(0,-2) 
\put(80,0){\makebox(0,0){Table 1 (continued)}}
\end{picture}

\end{array} 
\]

\[ 
\begin{array}{ccl} 
&& 
\setlength{\unitlength}{1.5pt}
\begin{picture}(120,43)(0,-2)

\put(10,0){\circle*{2.0}}
\put(10,20){\circle*{2.0}}
\put(30,20){\circle*{2.0}}
\put(30,0){\circle*{2.0}}

\put(50,0){\circle*{2.0}}
\put(70,0){\circle*{2.0}}
\put(60,20){\circle*{2.0}}
\put(90,0){\circle*{2.0}}
\put(110,0){\circle*{2.0}}
\put(100,20){\circle*{2.0}}

\put(107,11){\makebox(0,0){$3$}}

\put(20,-3){\makebox(0,0){$2$}}
\put(20,23){\makebox(0,0){$2$}}

\put(60,-3){\makebox(0,0){$2$}}
\put(68,10){\makebox(0,0){$2$}}
\put(60,-15){\makebox(0,0){(each cycle is non-oriented)}}
\put(100,-3){\makebox(0,0){$3$}}

\put(10,0){\line(1,0){20}}
\put(10,0){\line(0,1){20}}
\put(10,20){\line(1,0){20}}
\put(30,0){\line(0,1){20}}

\put(50,0){\line(1,0){20}}
\put(60,20){\line(1,-2){10}}
\put(60,20){\line(-1,-2){10}}

\put(90,0){\line(1,0){20}}
\put(90,0){\line(1,2){10}}
\put(100,20){\line(1,-2){10}}



\end{picture}\\[.4in]
&& 
\setlength{\unitlength}{1.5pt} 
\begin{picture}(140,17)(0,-2) 
\put(80,0){\makebox(0,0){Table 2 }}
\end{picture}

\end{array} 
\] 

\vspace {1 in}
In the following tables, we indexed the diagrams as follows. For any exceptional extended Dynkin diagram $X$, we denote by  $$X(c_i^{m_i};t;b_i^{l_i};d_i^{n_i})$$ a minimal 2-infinite diagram which is mutation equivalent to $X$ and which has
\begin{itemize}
\item
$m_i$ cycles of length $c_i$,
\item
$t$ triangles which are not adjacent to any cycle,
\item
$l_i$ branhes of length $l_i$, here a branch is a maximal chain which is not adjacent to any cycle,
\item
$n_i$ vertices which are connected to precisely $d_i$ vertices.
\end{itemize}

\pagebreak
\begin{center}
\psfrag{Table3}{Table 3: Minimal 2-infinite diagrams which are mutation equivalent to $F_4^{(1)}$ or $G_2^{(1)}$.}
\psfrag{0}{$F_4^{(1)}$}
\psfrag{11}{$F_4^{(1)}(3^{1};1;2^{1})_{1}$}
\psfrag{12}{$F_4^{(1)}(3^{1};1;2^{1})_{2}$}
\psfrag{22}{$F_4^{(1)}(3^{2};0)$}
\psfrag{1,2b}{$F_4^{(1)}(3^{1};1;1^{2})$}
\psfrag{2,3}{$F_4^{(1)}(3^{2};0;1^{1})$}
\psfrag{2,2,4}{$F_4^{(1)}(4^{1};3^{1})$}
\psfrag{3-1}{$G_2^{(1)}(1)$}
\psfrag{3-2}{$G_2^{(1)}(2)$}
\psfrag{3-3}{$G_2^{(1)}(3)$}
\psfrag{3-3-1}{$G_2^{(1)}(3^{1})$}
\includegraphics[width=4.6in]{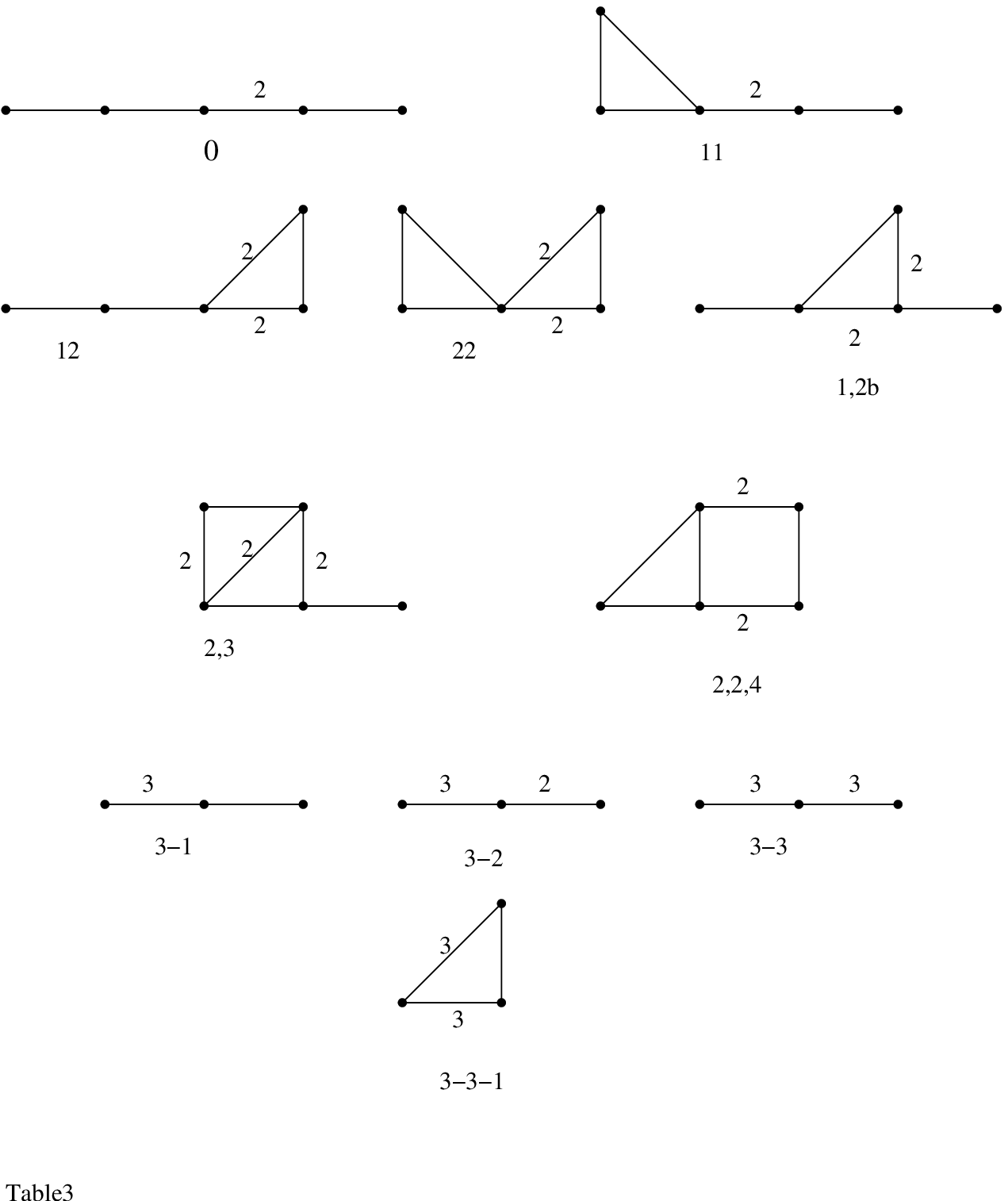}
\end{center}

\pagebreak
\begin{center}
\psfrag{p1}{$E_6^{(1)}$}
\psfrag{Table:71}{Table 4: Minimal 2-infinite diagrams which are mutation equivalent to $E_6^{(1)}$.}

\psfrag{p2}{$E_6^{(1)}(3^{1})$}
\psfrag{p3}{$E_6^{(1)}(3^{2};0)$}
\psfrag{p4}{$E_6^{(1)}(3^{2};2)$}
\psfrag{p5}{$E_6^{(1)}(3^{3})$}
\psfrag{p6}{$E_6^{(1)}(3^{4})$}
\psfrag{p7}{$E_6^{(1)}(4^{1};3^{1})$}
\psfrag{p8}{$E_6^{(1)}(4^{2};3^{1})$}
\includegraphics[width=4.6in]{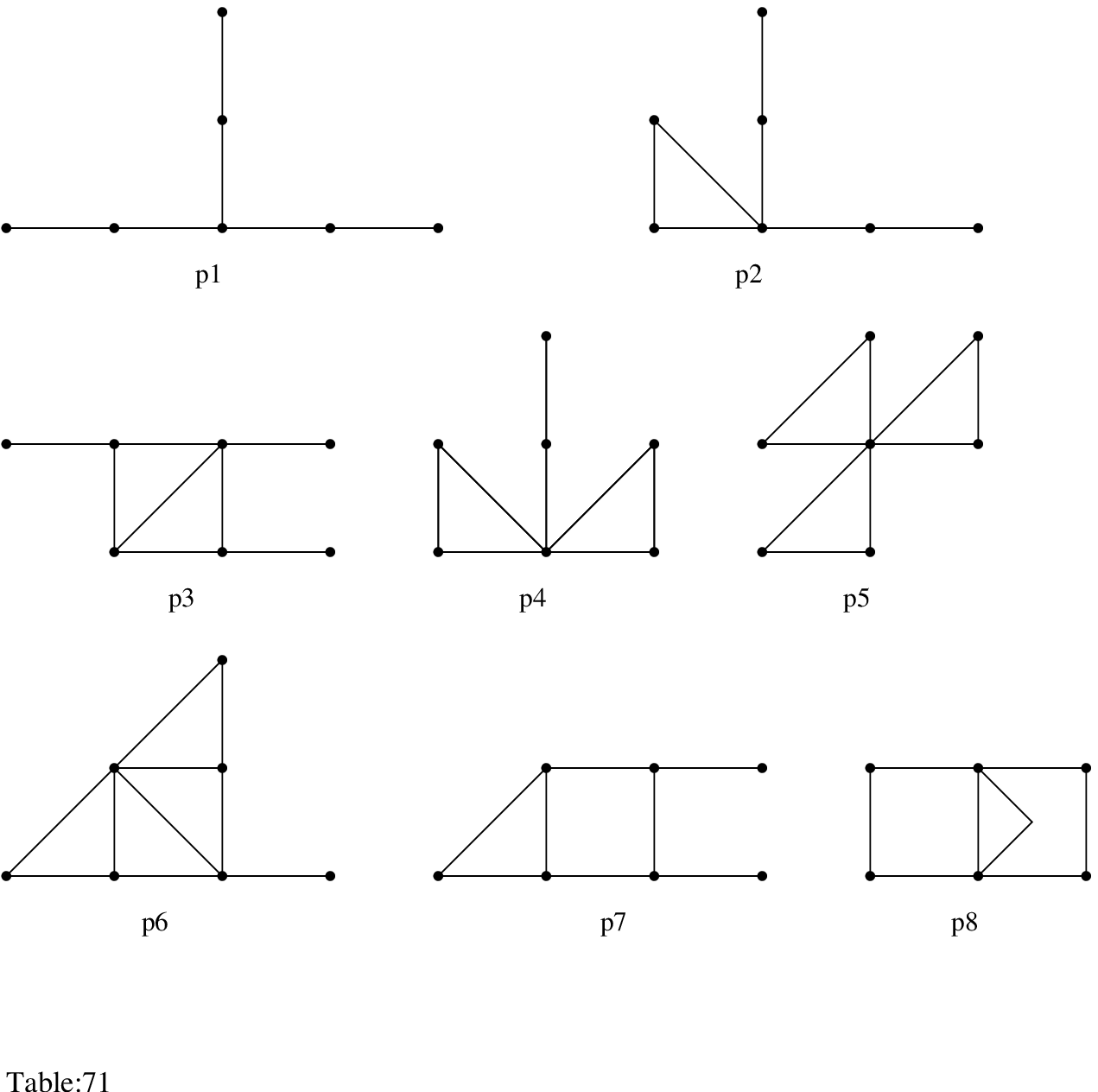}
\end{center}

\pagebreak
\begin{center}
\psfrag{Table: 81}{Table 5: Minimal 2-infinite diagrams which are mutation equivalent to $E_7^{(1)}$.}

\psfrag{p1}{$E_7^{(1)}$}
\psfrag{p2}{$E_7^{(1)}(3^{1};1;3^{1},1^{1})$}
\psfrag{p3}{$E_7^{(1)}(3^{1};1;3^{1},1^{2})$}
\psfrag{p4}{$E_7^{(1)}(3^{2};2;1^{1})$}
\psfrag{p5}{$E_7^{(1)}(3^{2};2;1^{2})$}
\psfrag{p6}{$E_7^{(1)}(3^{2};2;1^{3})$}
\psfrag{p7}{$E_7^{(1)}(3^{2};0;2^{1},1^{2})$}
\psfrag{p3-1}{$E_7^{(1)}(3^{3};1;1^{2})$}
\psfrag{p3-2}{$E_7^{(1)}(3^{3};0;2^{1},1^{1})$}
\psfrag{p3-3}{$E_7^{(1)}(3^{3};1;2^{1};4^{1})$}
\psfrag{p3-4}{$E_7^{(1)}(3^{3};1;2^{1};4^{2})$}
\psfrag{p3-5}{$E_7^{(1)}(3^{3};1;2^{1};5^{1})$}
\psfrag{P8}{$E_7^{(1)}(3^{2};0;2^{2};3^{2})$}
\psfrag{P9}{$E_7^{(1)}(3^{2};0;2^{2};4^{1})$}
\includegraphics[width=4.6in]{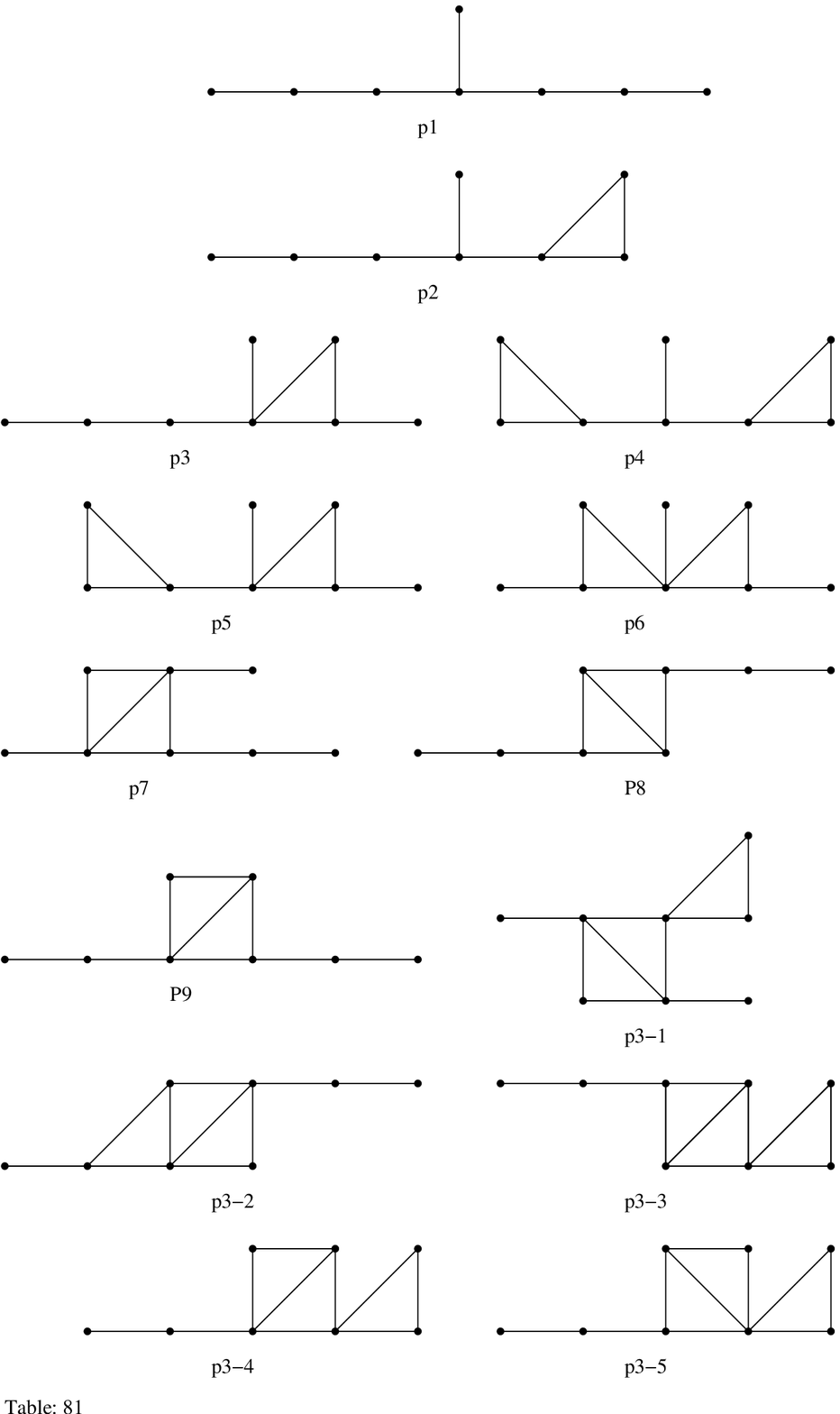}

\end{center}

\pagebreak
\begin{center}
\psfrag{Table: 82}{Table 5 (continued)}

\psfrag{T4-1}{$E_7^{(1)}(3^{4};2;0;4^{2})$}
\psfrag{T4-2}{$E_7^{(1)}(3^{4};2;0;5^{1})$}
\psfrag{T4-3}{$E_7^{(1)}(3^{4};1;1^{1})$}
\psfrag{T4-4}{$E_7^{(1)}(3^{4};0;1^{2};5^{1})_{1}$}
\psfrag{T4-5}{$E_7^{(1)}(3^{4};0;1^{2};4^{2})$}
\psfrag{T4-6}{$E_7^{(1)}(3^{4};0;1^{2};5^{1})_{2}$}
\psfrag{T4-7}{$E_7^{(1)}(3^{4};0;2^{1})$}
\psfrag{T5-1}{$E_7^{(1)}(3^{6};0)$}
\psfrag{T5-2}{$E_7^{(1)}(3^{5};1)$}
\includegraphics[width=4.6in]{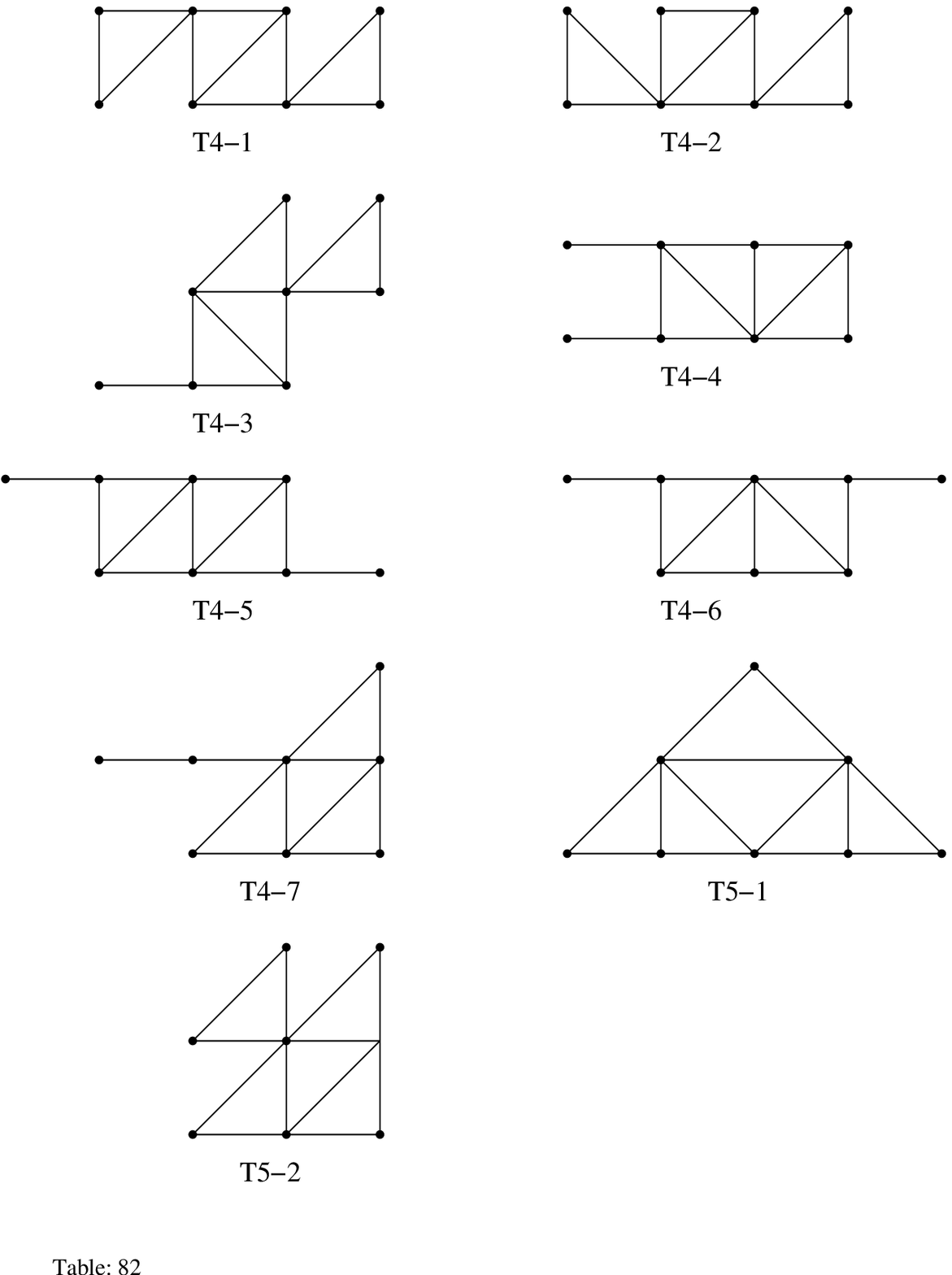}

\end{center}

\pagebreak
\begin{center}
\psfrag{Table: 83}{Table 5 (continued)}

\psfrag{S1}{$E_7^{(1)}(4^{1})$}
\psfrag{ST-1}{$E_7^{(1)}(4^{1},3^{1};0;1^{3})$}
\psfrag{ST-2}{$E_7^{(1)}(4^{1},3^{1};1;2^{1})$}
\psfrag{ST-3}{$E_7^{(1)}(4^{1},3^{1};0;2^{1})$}
\psfrag{ST2-0}{$E_7^{(1)}(4^{1},3^{2};2)$}
\psfrag{ST2-1}{$E_7^{(1)}(4^{1},3^{2};1)$}
\psfrag{ST2-2}{$E_7^{(1)}(4^{1},3^{2};0)$}
\psfrag{ST3-1}{$E_7^{(1)}(4^{2},3^{2};1)$}
\psfrag{ST3-2}{$E_7^{(1)}(4^{1},3^{3};0;1^{1};4^{1})$}
\psfrag{ST3-3}{$E_7^{(1)}(4^{1},3^{3};0;1^{1};5^{1},3^{4})$}
\psfrag{ST3-4}{$E_7^{(1)}(4^{1},3^{3};0;1^{1};5^{1},4^{1})$}

\includegraphics[width=4.6in]{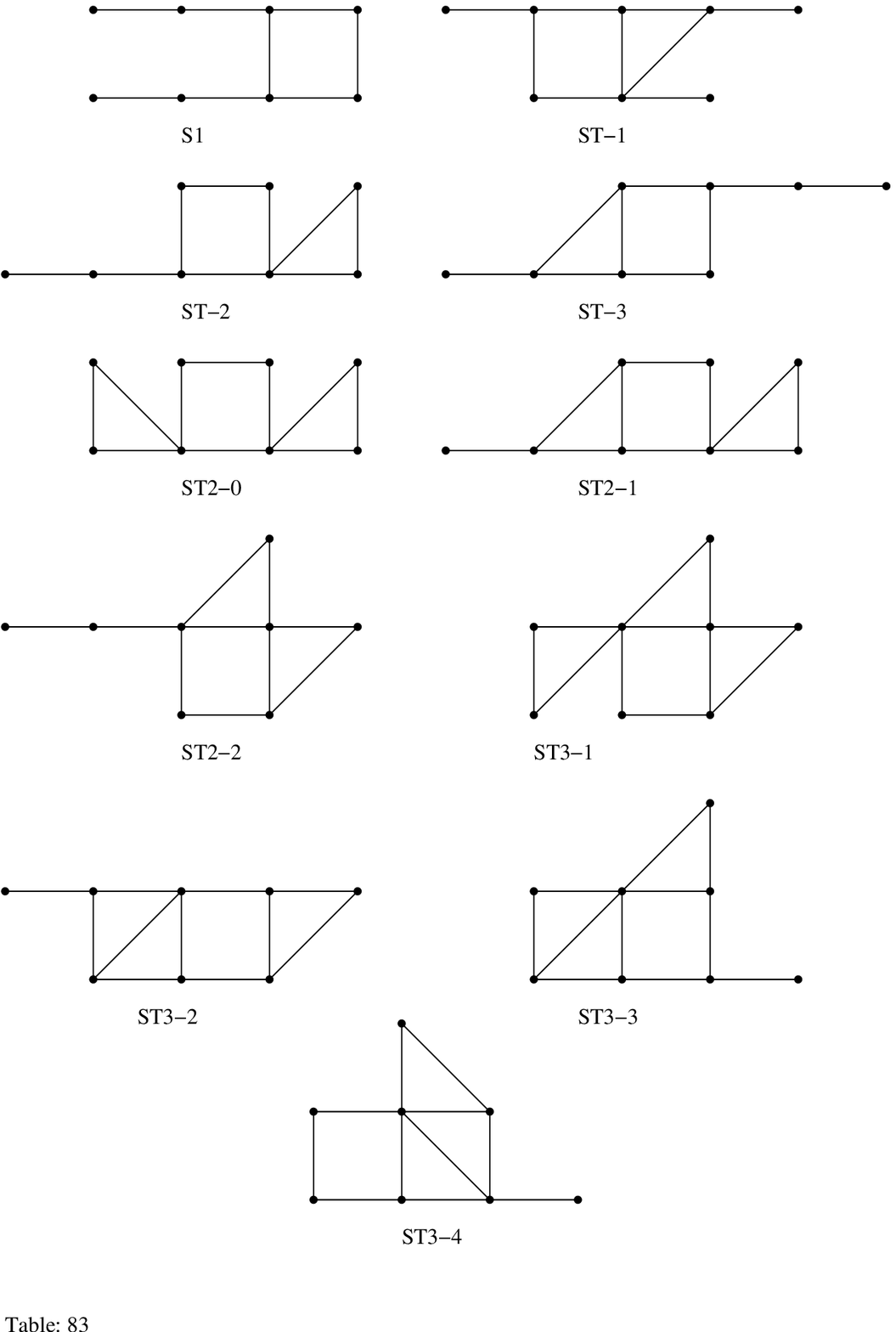}

\end{center}

\pagebreak
\begin{center}
\psfrag{Table:84}{Table 5 (continued)}

\psfrag{S2}{$E_7^{(1)}(4^{2},3^{0})$}
\psfrag{S2T2}{$E_7^{(1)}(4^{2},3^{2};0)$}
\psfrag{PT1}{$E_7^{(1)}(5^{1},4^{0},3^{1};0;1^{2})$}
\psfrag{PT2}{$E_7^{(1)}(5^{1},4^{0},3^{1};0;2^{1})$}
\psfrag{PT3}{$E_7^{(1)}(5^{1},4^{0},3^{2})$}
\psfrag{PST}{$E_7^{(1)}(5^{1},4^{1},3^{1})$}

\includegraphics[width=4.6in]{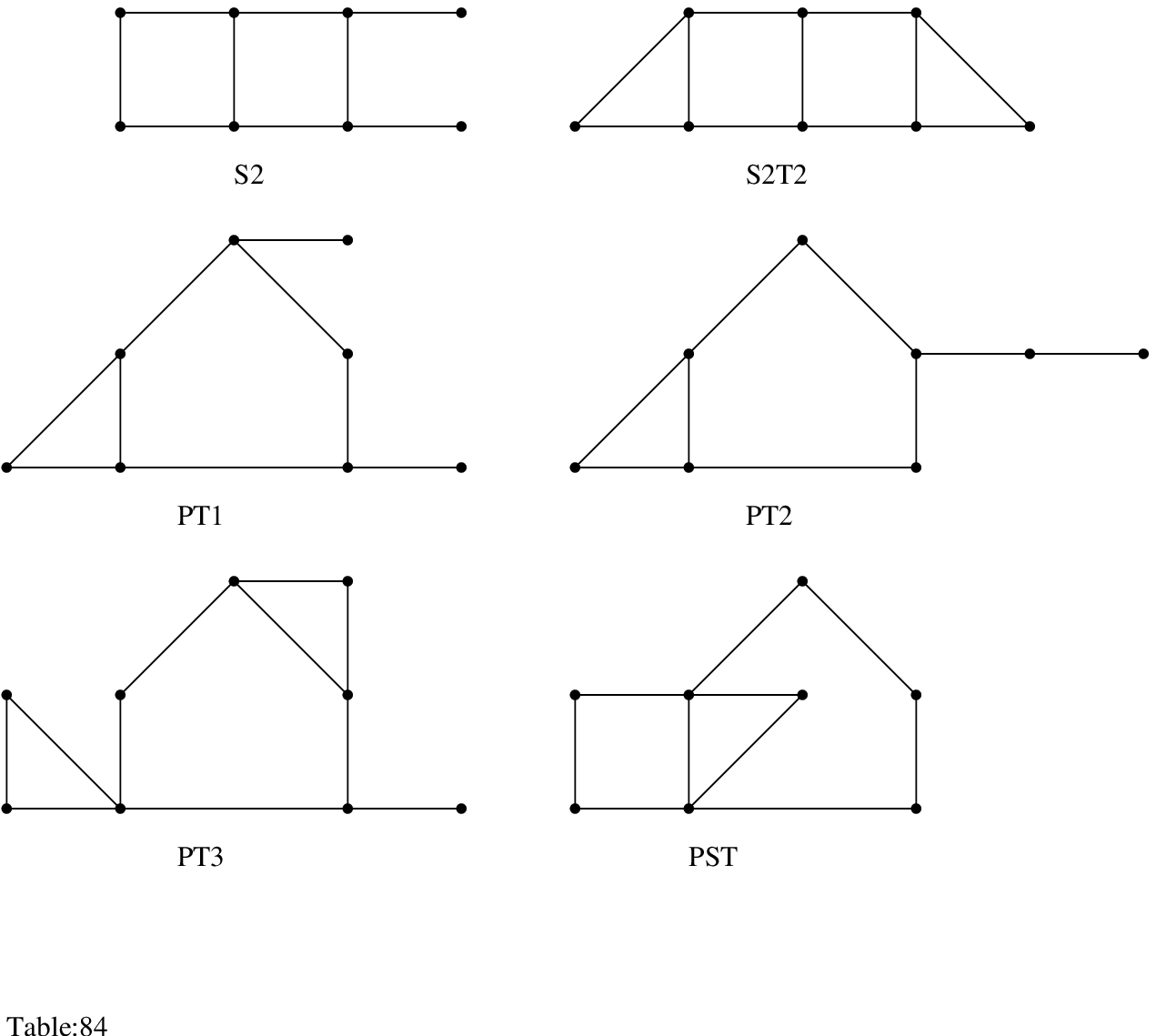}
\end{center}

\pagebreak
\begin{center}
\psfrag{Table:91}{Table 6: Minimal 2-infinite diagrams which are mutation equivalent to $E_8^{(1)}$.}

\psfrag{T1}{$E_8^{(1)}$}
\psfrag{T2}{$E_8^{(1)}(3^{1};1;2^{2},1^{2})$}
\psfrag{T3}{$E_8^{(1)}(3^{1};1;3^{1},2^{1},1^{1})$}
\psfrag{T4}{$E_8^{(1)}(3^{1};1;3^{1},2^{1},1^{2})$}
\psfrag{T5}{$E_8^{(1)}(3^{1};1;3^{1},2^{2})$}
\psfrag{T6}{$E_8^{(1)}(3^{1};1;4^{1},3^{1},1^{1})$}
\psfrag{T7}{$E_8^{(1)}(3^{1};1;5^{1},1^{1})$}
\psfrag{T8}{$E_8^{(1)}(3^{1};1;5^{1},4^{1})$}

\includegraphics[width=4.6in]{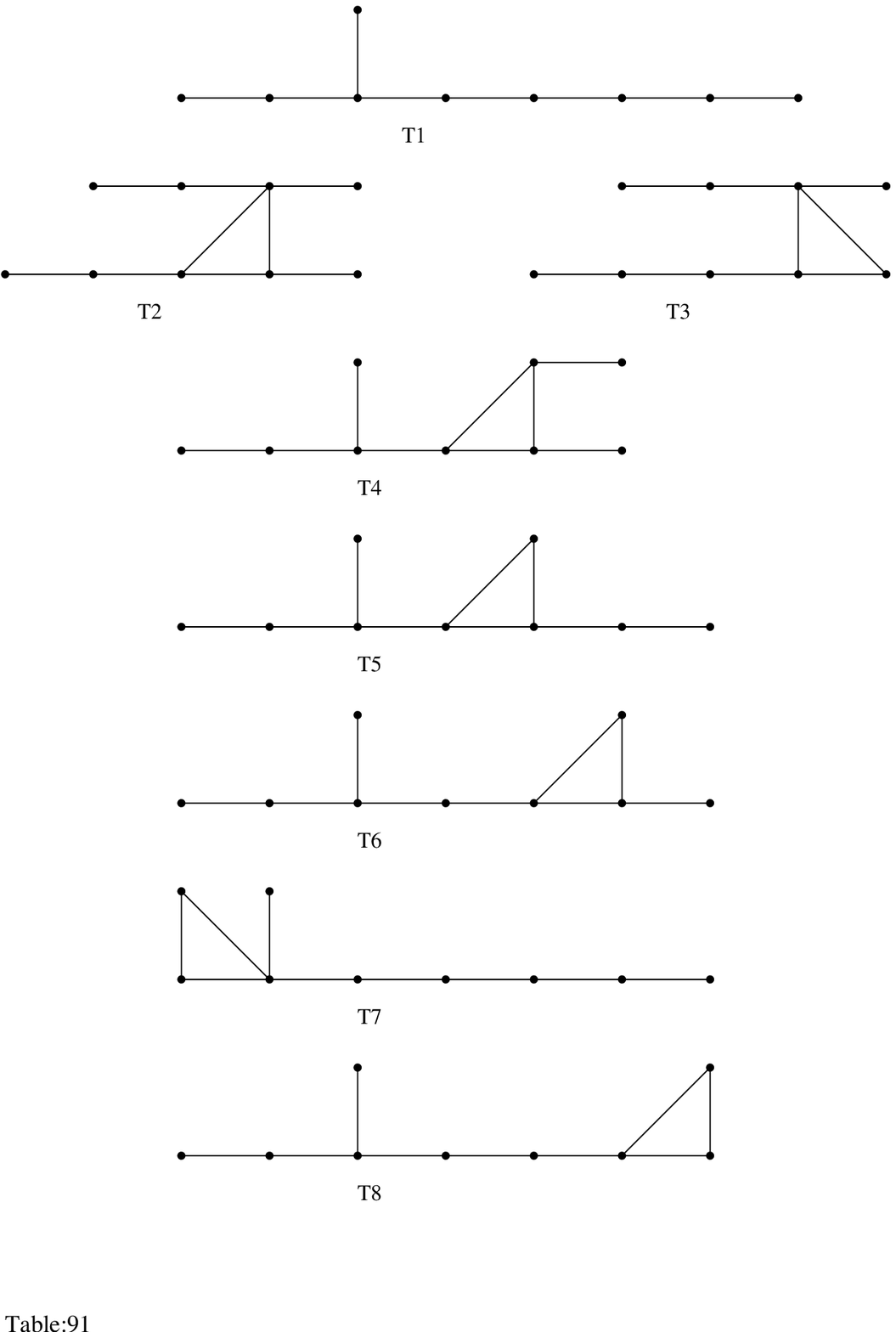}

\end{center}

\pagebreak
\begin{center}
\psfrag{Table: 92}{Table 6 (continued)}

\psfrag{T2-1}{$E_8^{(1)}(3^{2};2;1^{1})$}
\psfrag{T2-2}{$E_8^{(1)}(3^{2};2;1^{2})$}
\psfrag{T2-3}{$E_8^{(1)}(3^{2};2;1^{3})$}
\psfrag{T2-4}{$E_8^{(1)}(3^{2};2;2^{1},1^{1};4^{1},3^{2})_1$}
\psfrag{T2-5}{$E_8^{(1)}(3^{2};2;2^{1},1^{1};4^{1},3^{2})_2$}
\psfrag{T2-6}{$E_8^{(1)}(3^{2};2;2^{1},1^{1};5^{1})$}
\psfrag{T2-7}{$E_8^{(1)}(3^{2};2;2^{1},1^{2})$}
\psfrag{T2-7-2}{$E_8^{(1)}(3^{2};2;2^{1},1^{3})$}
\psfrag{T2-8}{$E_8^{(1)}(3^{2};0;2^{2})$}

\includegraphics[width=4.6in]{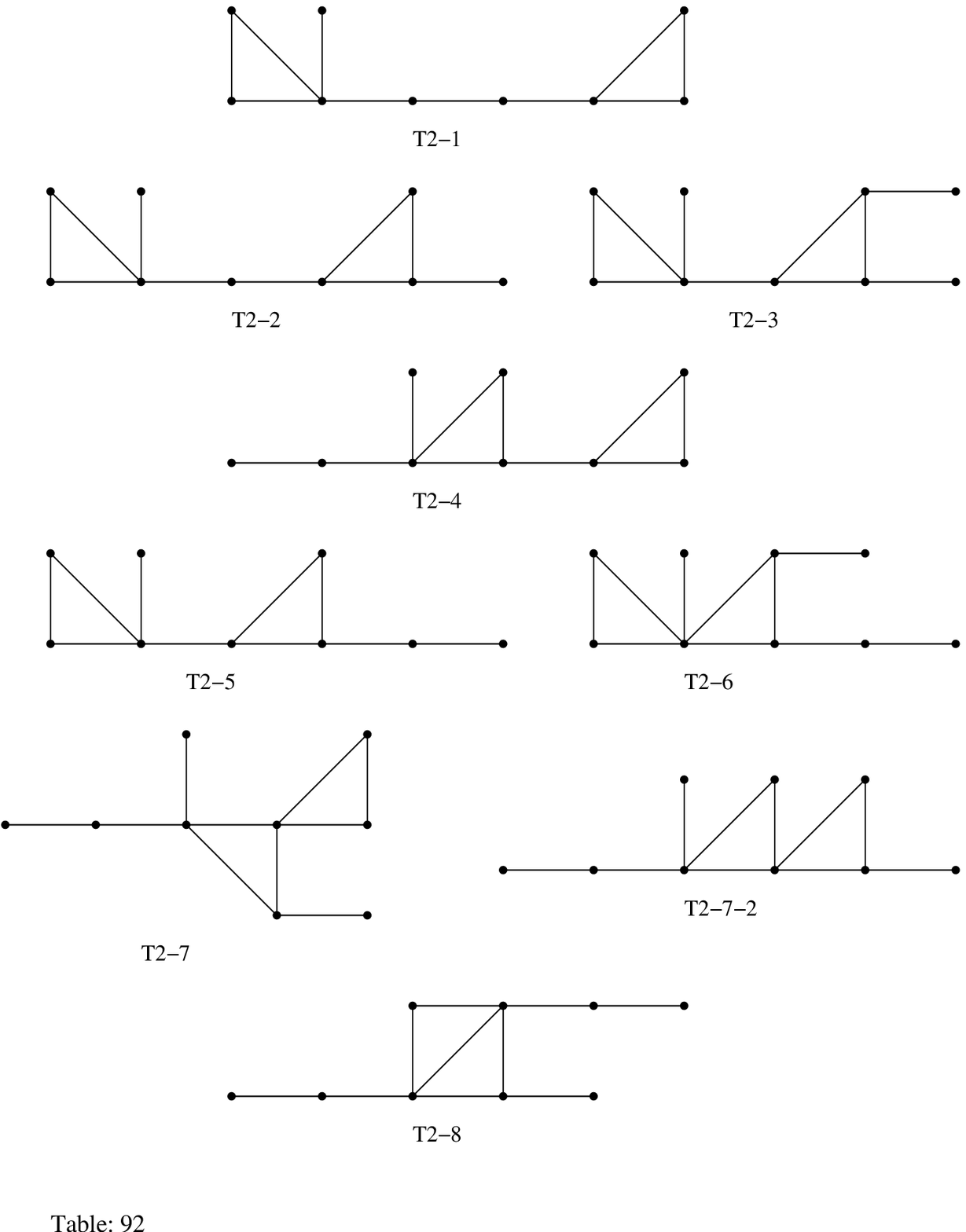}

\end{center}

\pagebreak
\begin{center}
\psfrag{Table: 93}{Table 6 (continued)}

\psfrag{T2-0-1}{$E_8^{(1)}(3^{2};0;3^{1},1^{2})$}
\psfrag{T2-0-2}{$E_8^{(1)}(3^{2};0;4^{1},1^{1};3^{4})$}
\psfrag{T2-0-3}{$E_8^{(1)}(3^{2};0;4^{1},1^{1};4^{1})_1$}
\psfrag{T2-0-4}{$E_8^{(1)}(3^{2};0;4^{1},1^{1};4^{1})_2$}
\psfrag{T2-2-1}{$E_8^{(1)}(3^{2};2;3^{1},1^{1})$}
\psfrag{T2-2-2}{$E_8^{(1)}(3^{2};2;3^{1},2^{1})$}

\includegraphics[width=4.6in]{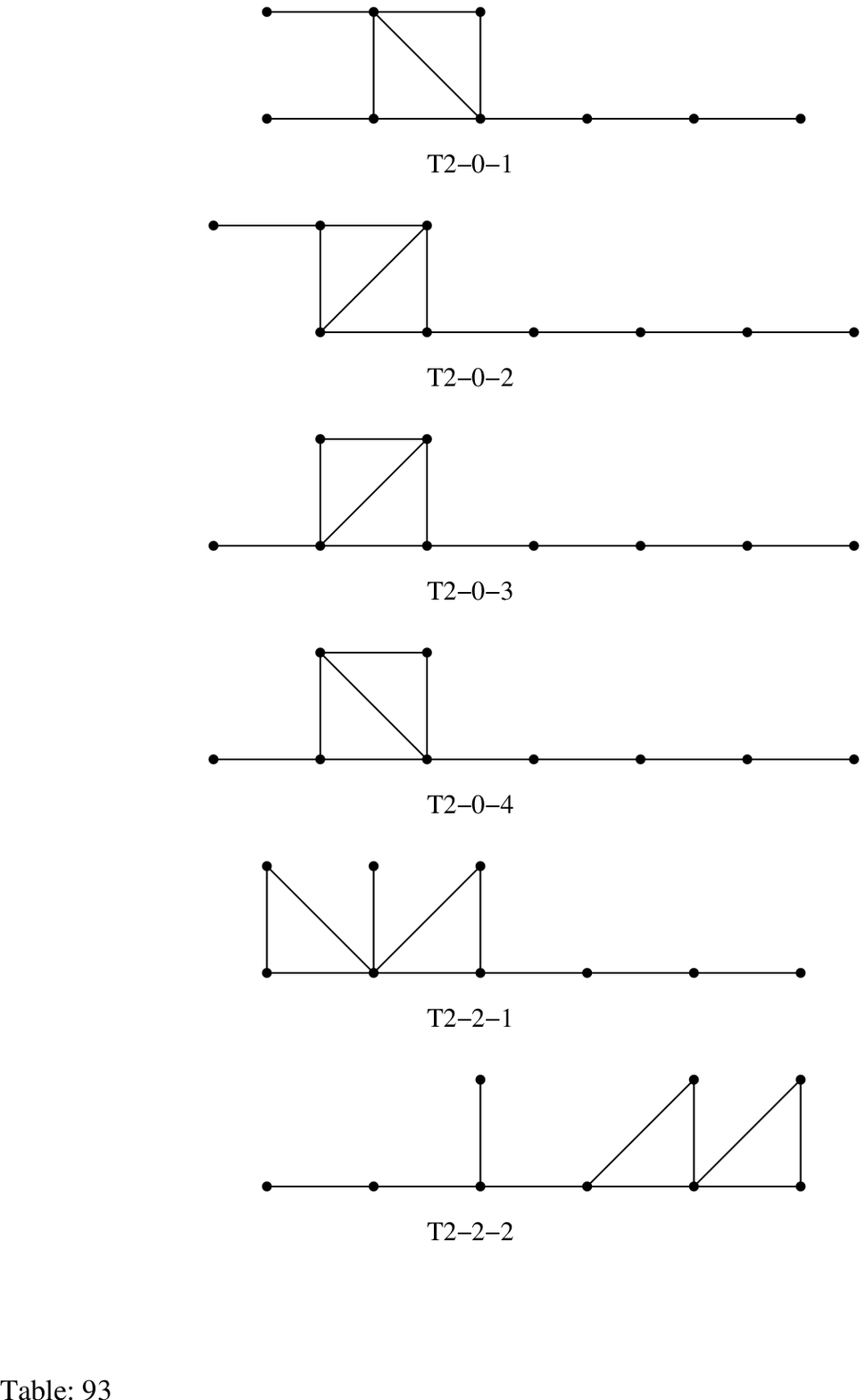}

\end{center}

\pagebreak
\begin{center}
\psfrag{Table: 94}{Table 6 (continued)}

\psfrag{T3-0-1}{$E_8^{(1)}(3^{3};0;2^{1},1^{2})$}
\psfrag{T3-0-2}{$E_8^{(1)}(3^{3};0;3^{1},1^{1};4^{1})$}
\psfrag{T3-0-3}{$E_8^{(1)}(3^{3};0;3^{1},1^{1};5^{1};3^{3})$}
\psfrag{T3-0-4}{$E_8^{(1)}(3^{3};0;3^{1},1^{1};5^{1},4^{1})$}

\psfrag{T3-0-5}{$E_8^{(1)}(3^{3};0;3^{1},1^{1};4^{2},3^{2})$}
\psfrag{T3-0-6}{$E_8^{(1)}(3^{3};0;3^{1},1^{1};4^{1},3^{2})$}
\psfrag{T3-0-7}{$E_8^{(1)}(3^{3};0;4^{1};4^{2})$}
\psfrag{T3-0-8}{$E_8^{(1)}(3^{3};0;4^{1};5^{1})$}

\includegraphics[width=4.6in]{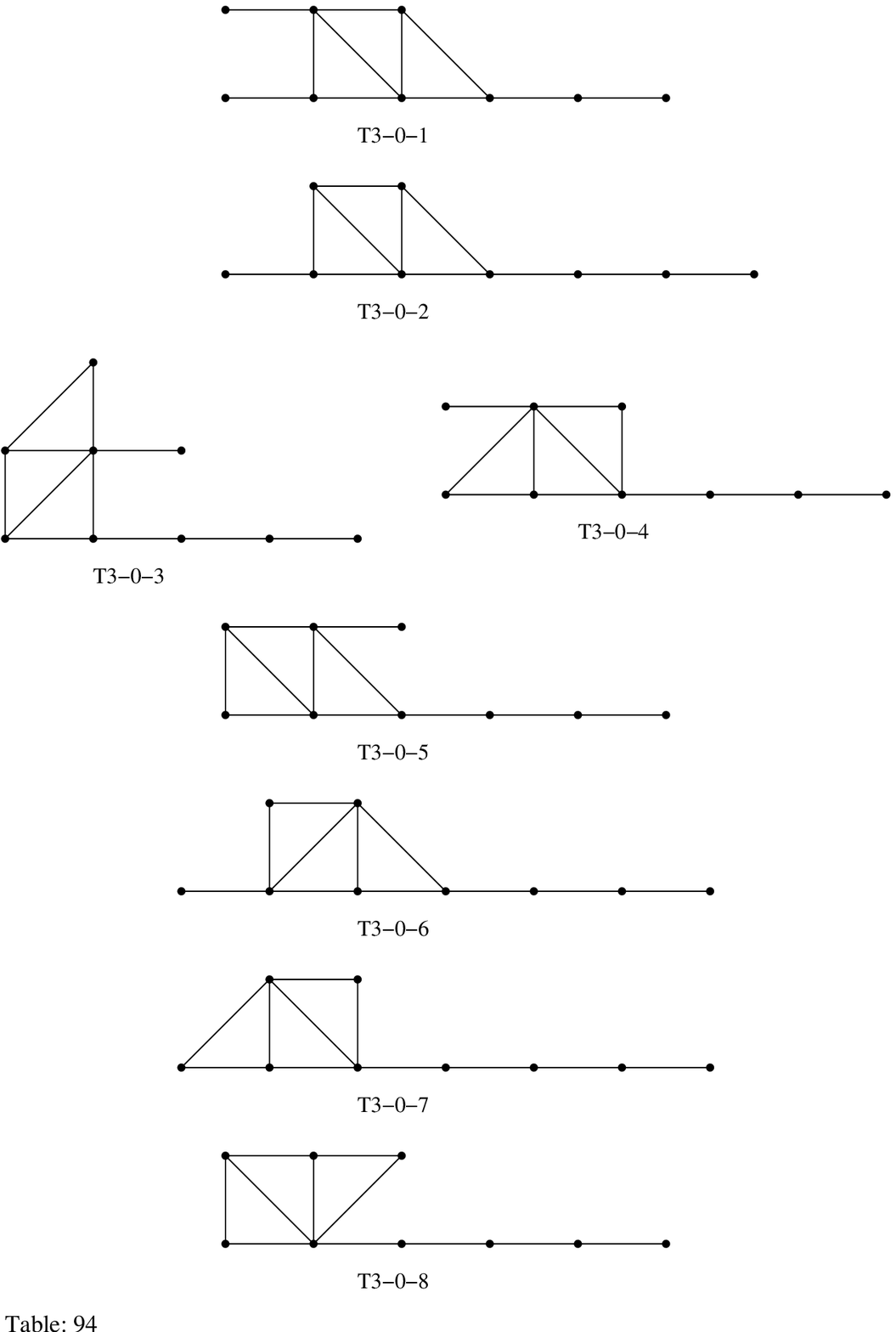}

\end{center}

\pagebreak
\begin{center}
\psfrag{Table:952}{Table 6 (continued)}
\psfrag{1-1}{$E_8^{(1)}(3^{3};1;1^{1};3^{4})$}
\psfrag{1-2}{$E_8^{(1)}(3^{3};1;1^{1};4^{1},3^{3})_{1}$}
\psfrag{1-3}{$E_8^{(1)}(3^{3};1;1^{1};4^{1},3^{3})_{2}$}
\psfrag{1-4}{$E_8^{(1)}(3^{3};1;1^{2};4^{2},3^{2})$}
\psfrag{1-5}{$E_8^{(1)}(3^{3};1;1^{2};4^{1},3^{4})_{1}$}
\psfrag{1-6}{$E_8^{(1)}(3^{3};1;1^{2};4^{1},3^{4})_{2}$}
\psfrag{1-6-2}{$E_8^{(1)}(3^{3};1;1^{2};4^{2},3^{2})$}
\psfrag{1-7}{$E_8^{(1)}(3^{3};1;1^{3};4^{1},3^{5})$}
\psfrag{1-8}{$E_8^{(1)}(3^{3};1;1^{3};4^{1},3^{3})$}

\includegraphics[width=4.6in]{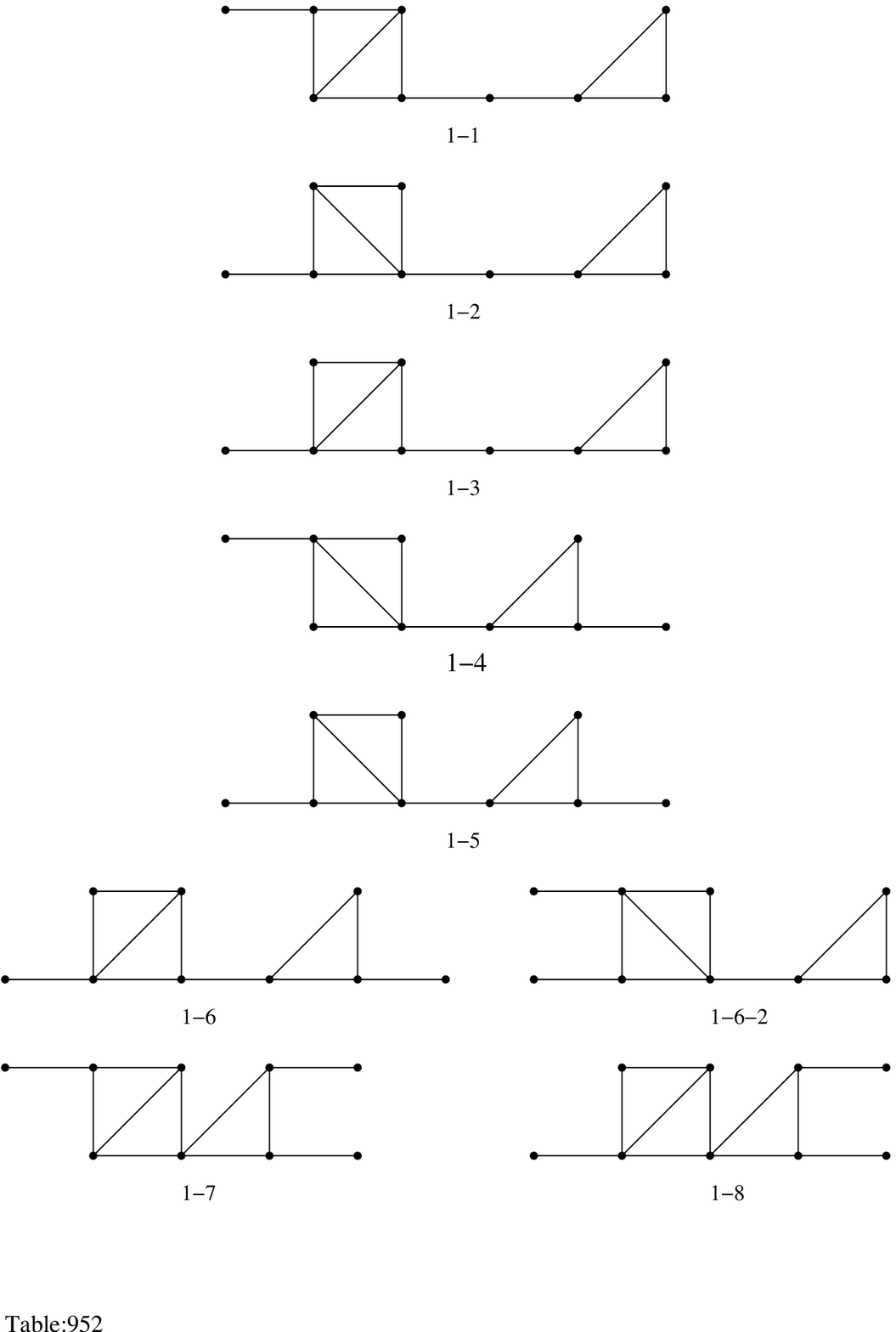}
\end{center}

\pagebreak
\begin{center}
\psfrag{Table:953}{Table 6 (continued)}
\psfrag{1-9}{$E_8^{(1)}(3^{3};1;1^{3};5^{1},3^{4})$}
\psfrag{1-9-2}{$E_8^{(1)}(3^{3};1;1^{3};5^{1},4^{1},3^{2})$}
\psfrag{1-10}{$E_8^{(1)}(3^{3};1;2^{1},1^{1};4^{2},3^{2})$}
\psfrag{1-11}{$E_8^{(1)}(3^{3};1;2^{1},1^{1};4^{1},3^{4})$}
\psfrag{1-12}{$E_8^{(1)}(3^{3};1;2^{1},1^{1};4^{2},3^{2})$}
\psfrag{1-13}{$E_8^{(1)}(3^{3};1;2^{1},1^{1};5^{1},3^{2})$}
\psfrag{3-1}{$E_8^{(1)}(3^{3};3;1^{1};4^{2},3^{1})$}
\psfrag{3-2}{$E_8^{(1)}(3^{3};3;1^{1};5^{1},3^{2})$}
\psfrag{3-3}{$E_8^{(1)}(3^{3};3;1^{2};5^{1},4^{1},3^{1})_{1}$}
\psfrag{3-4}{$E_8^{(1)}(3^{3};3;1^{2};5^{1},4^{1},3^{1})_{2}$}

\includegraphics[width=4.6in]{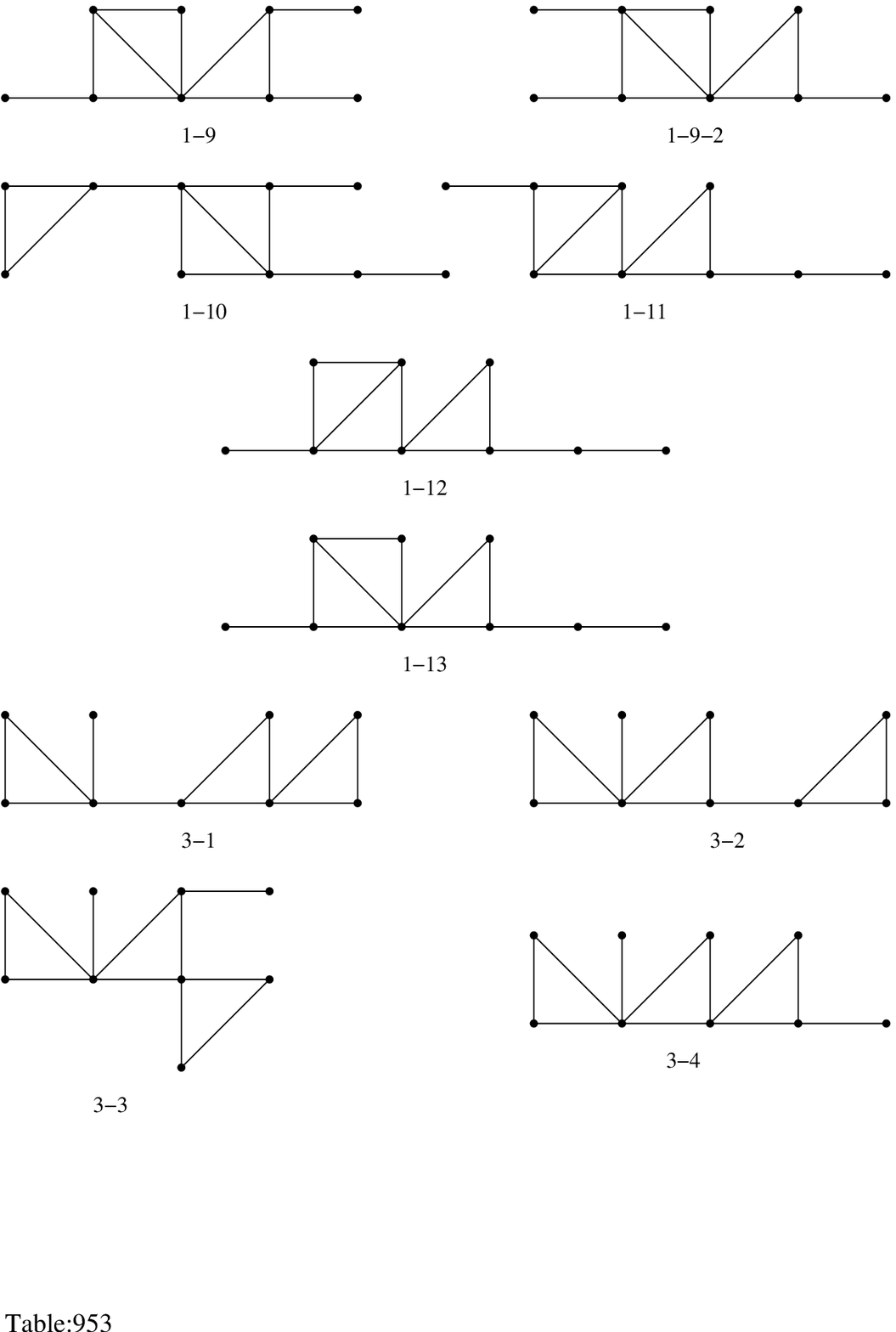}
\end{center}

\pagebreak
\begin{center}
\psfrag{Table: 95}{Table 6 (continued)}

\psfrag{T4-0-1}{$E_8^{(1)}(3^{4};1;1^{2};5^{1})$}
\psfrag{T4-0-2}{$E_8^{(1)}(3^{4};0;1^{3})_{1}$}
\psfrag{T4-0-3}{$E_8^{(1)}(3^{4};0;1^{3})_{2}$}
\psfrag{T4-0-4}{$E_8^{(1)}(3^{4};0;2^{1};4^{3})$}
\psfrag{T4-0-5}{$E_8^{(1)}(3^{4};0;2^{1},1^{1};4^{2})$}
\psfrag{T4-0-6}{$E_8^{(1)}(3^{4};0;2^{1},1^{1};5^{1})_{1}$}
\psfrag{T4-0-7}{$E_8^{(1)}(3^{4};0;2^{1},1^{1};5^{1})_{2}$}
\psfrag{T4-0-8}{$E_8^{(1)}(3^{4};0;3^{1})$}

\includegraphics[width=4.6in]{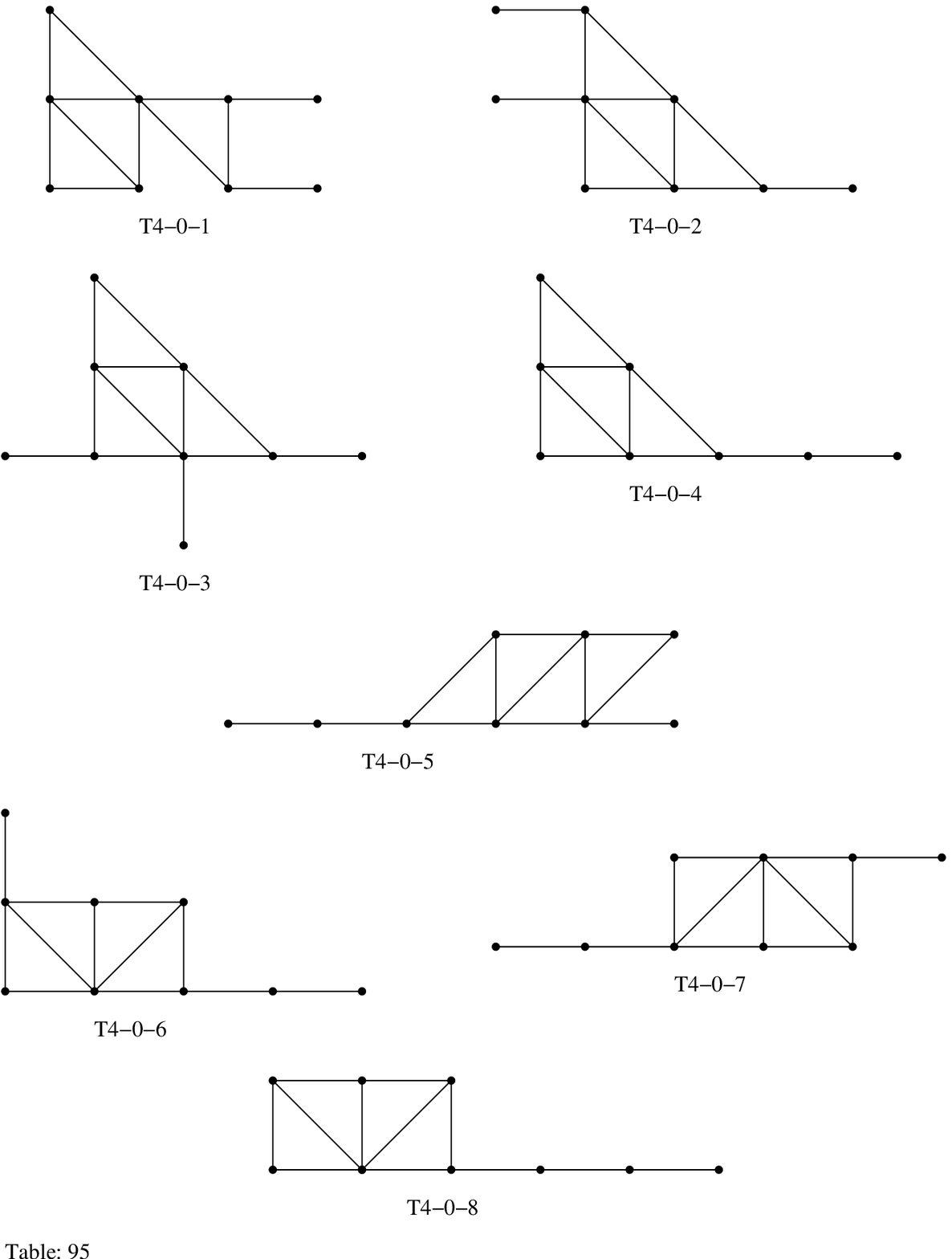}
\end{center}

\pagebreak
\begin{center}
\psfrag{Tablr: 96}{Table 6 (continued)}

\psfrag{T4-1-1}{$E_8^{(1)}(3^{4};1;0;4^{2})$}
\psfrag{T4-1-2}{$E_8^{(1)}(3^{4};1;0;5^{1})$}
\psfrag{T4-1-3}{$E_8^{(1)}(3^{4};1;1^{1};4^{1})$}
\psfrag{T4-1-4}{$E_8^{(1)}(3^{4};1;1^{2};4^{2},3^{3})_{1}$}
\psfrag{T4-1-5}{$E_8^{(1)}(3^{4};1;1^{2};4^{2},3^{3})_{2}$}
\psfrag{T4-1-6}{$E_8^{(1)}(3^{4};1;1^{1};5^{1},3^{3})$}
\psfrag{T4-1-7}{$E_8^{(1)}(3^{4};1;1^{1};5^{1},3^{4})$}
\psfrag{T4-1-8}{$E_8^{(1)}(3^{4};1;1^{1};5^{1},4^{1})$}
\psfrag{T4-1-9}{$E_8^{(1)}(3^{4};1;1^{2};4^{3})_{1}$}
\psfrag{T4-1-10}{$E_8^{(1)}(3^{4};1;1^{2};4^{3})_{2}$}

\includegraphics[width=4.6in]{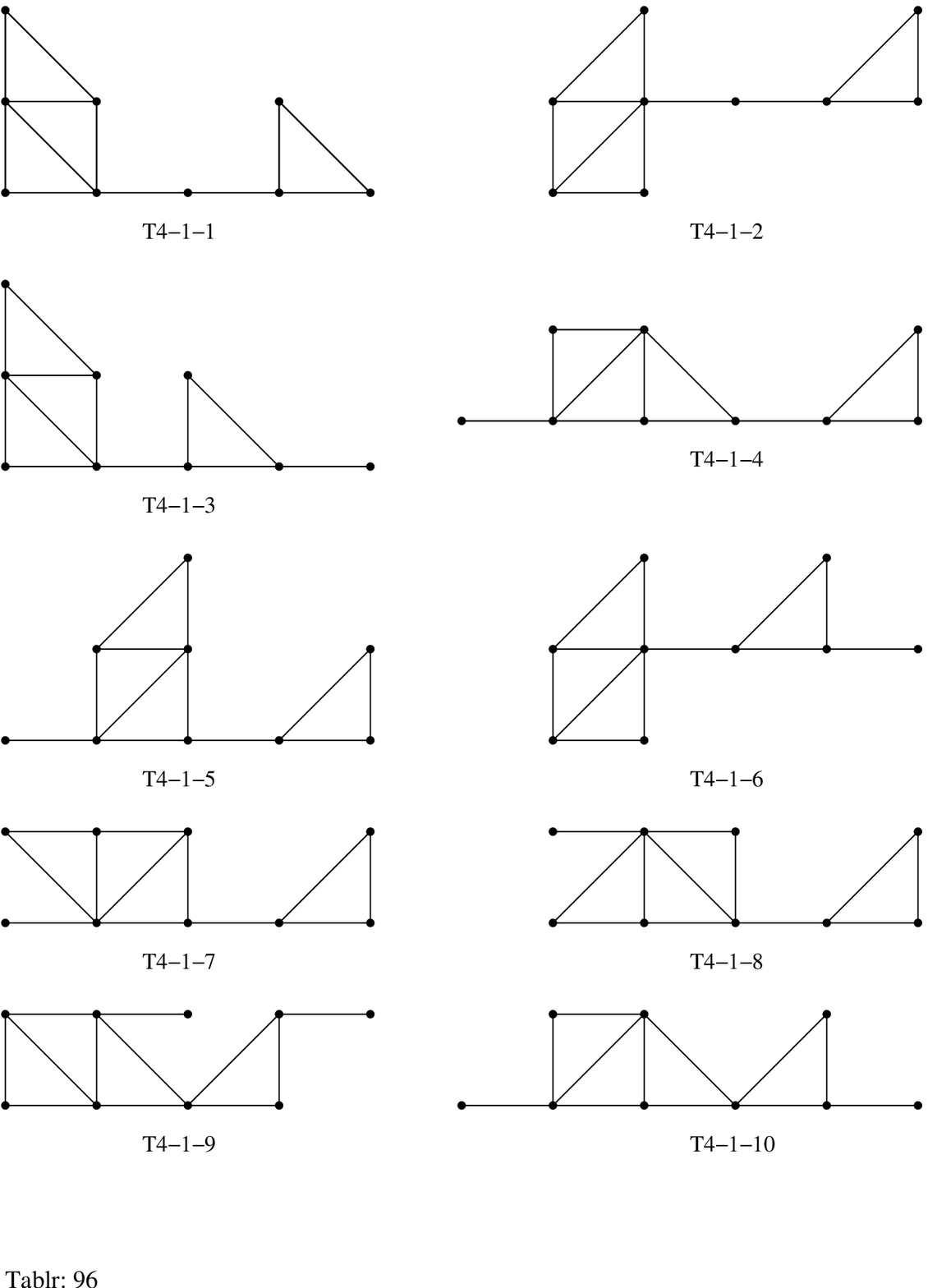}
\end{center}

\pagebreak
\begin{center}
\psfrag{Table:97}{Table 6 (continued)}

\psfrag{T4-1-11}{$E_8^{(1)}(3^{4};1;1^{2};4^{3})_{3}$}
\psfrag{T4-1-11-3}{$E_8^{(1)}(3^{4};1;1^{2};4^{3})_{4}$}
\psfrag{T4-1-11-1}{$E_8^{(1)}(3^{4};1;1^{2};5^{1},4^{1})$}
\psfrag{T4-1-12}{$E_8^{(1)}(3^{4};1;1^{2};6^{1})$}
\psfrag{T4-1-13}{$E_8^{(1)}(3^{4};1;2^{1};5^{1},4^{1})$}
\psfrag{T4-1-14}{$E_8^{(1)}(3^{4};1;2^{1};6^{1})$}
\psfrag{T4-2-1}{$E_8^{(1)}(3^{4};2;1^{1};4^{2},3^{3})$}
\psfrag{T4-2-2}{$E_8^{(1)}(3^{4};2;1^{1};4^{3})$}
\psfrag{T4-2-3}{$E_8^{(1)}(3^{4};2;1^{1};5^{1},4^{1},3^{1})$}
\psfrag{T4-2-4}{$E_8^{(1)}(3^{4};2;1^{1};5^{2})$}

\includegraphics[width=4.6in]{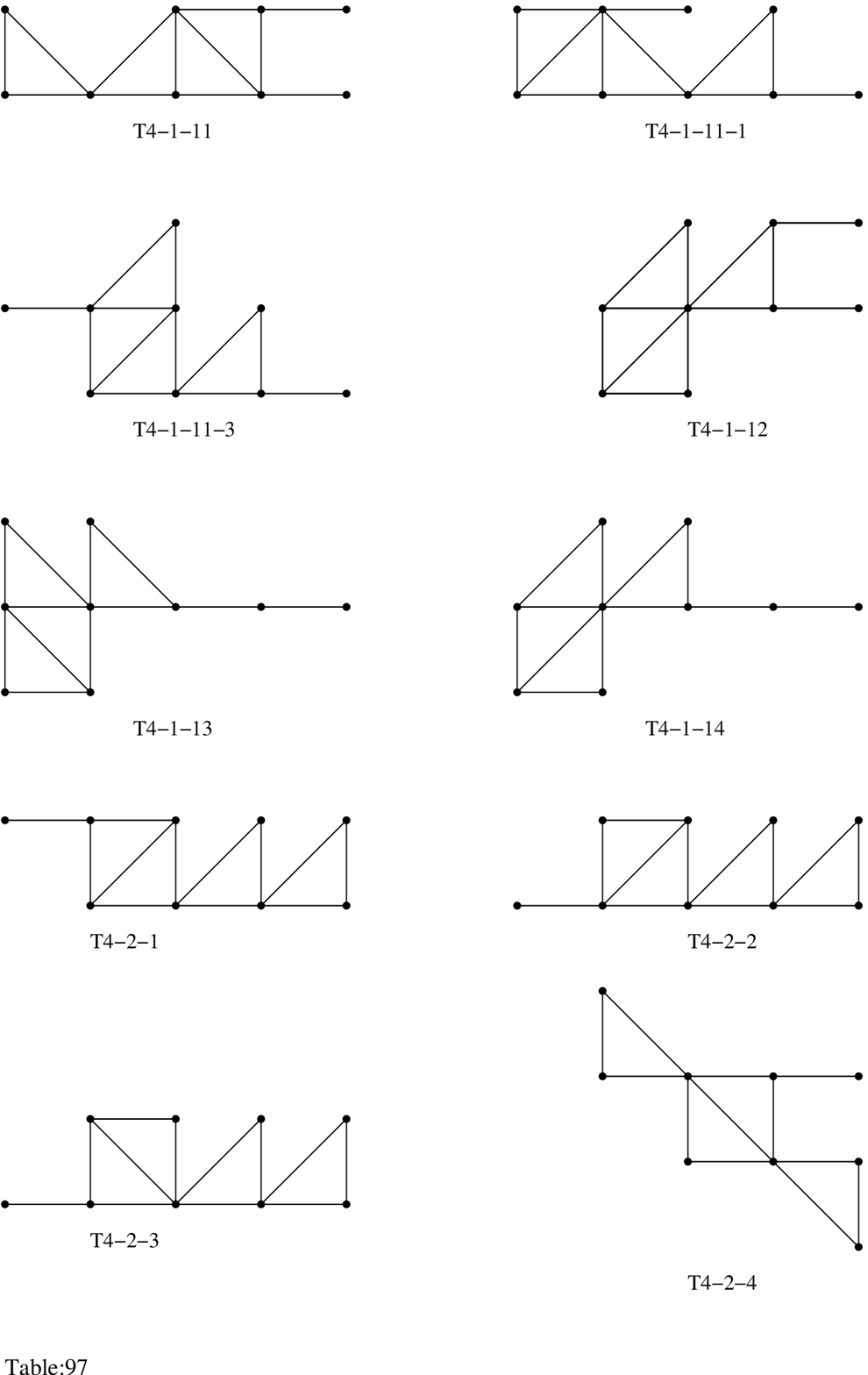}
\end{center}

\pagebreak
\begin{center}
\psfrag{Table: 98}{Table 6 (continued)}

\psfrag{T5-0-1}{$E_8^{(1)}(3^{5};0;1^{2};4^{3})$}
\psfrag{T5-0-2}{$E_8^{(1)}(3^{5};0;1^{2};5^{1})$}
\psfrag{T5-0-3}{$E_8^{(1)}(3^{5};0;1^{2};5^{1},4^{2},3^{2})_{1}$}
\psfrag{T5-0-4}{$E_8^{(1)}(3^{5};0;1^{2};5^{1},4^{2},3^{2})_{2}$}
\psfrag{T5-0-5}{$E_8^{(1)}(3^{5};0;2^{1};5^{1},4^{3})$}
\psfrag{T5-0-6}{$E_8^{(1)}(3^{5};0;2^{1};6^{1})$}
\psfrag{T5-1-1}{$E_8^{(1)}(3^{5};1;0;4^{1})$}
\psfrag{T5-1-2}{$E_8^{(1)}(3^{5};1;1^{1};4^{3})$}
\psfrag{T5-1-3}{$E_8^{(1)}(3^{5};1;1^{1};4^{4})$}
\psfrag{T5-1-4}{$E_8^{(1)}(3^{5};1;1^{1};5^{1},4^{1})$}
\psfrag{T4-0-4-2}{$E_8^{(1)}(3^{5};1;0;5^{1})$}
\includegraphics[width=4.6in]{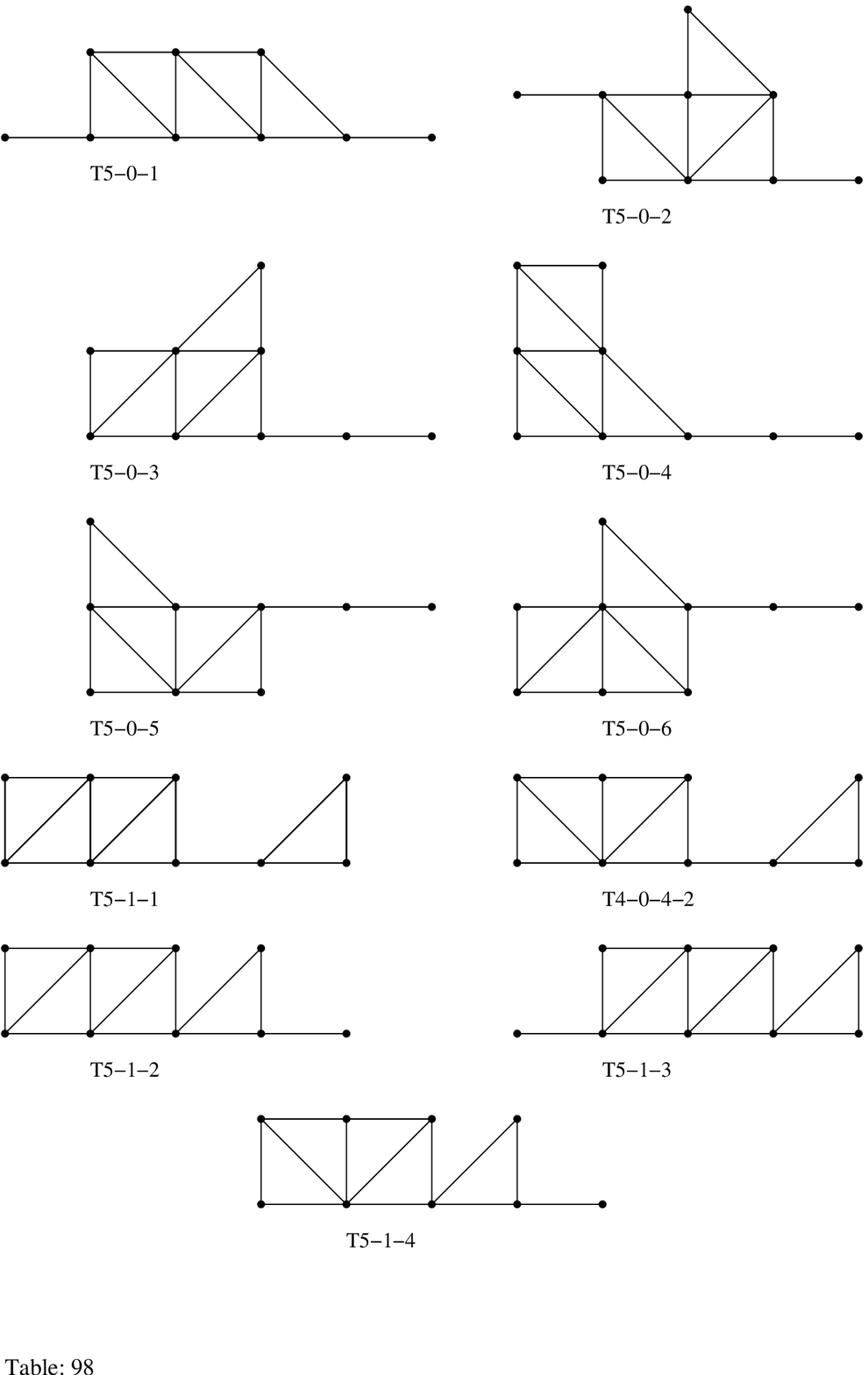}
\end{center}

\pagebreak
\begin{center}
\psfrag{Table: 99}{Table 6 (continued)}

\psfrag{T5-1-5}{$E_8^{(1)}(3^{5};1;1^{1};5^{1},4^{2})$}
\psfrag{T5-1-6}{$E_8^{(1)}(3^{5};1;1^{1};5^{1},4^{3})$}
\psfrag{T5-1-7}{$E_8^{(1)}(3^{5};1;1^{1};5^{2})$}
\psfrag{T5-2-1}{$E_8^{(1)}(3^{5};2;0;5^{1})$}
\psfrag{T5-2-2}{$E_8^{(1)}(3^{5};2;0;6^{1})$}
\psfrag{T6-5}{$E_8^{(1)}(3^{6};0;1^{1};4^{2})$}
\psfrag{T6-6}{$E_8^{(1)}(3^{6};0;1^{1};5^{1},4^{3})$}
\psfrag{T6-7}{$E_8^{(1)}(3^{6};0;1^{1};5^{2})$}
\psfrag{T6-1}{$E_8^{(1)}(3^{6};1;0;5^{1},4^{3})_{1}$}
\psfrag{T6-2}{$E_8^{(1)}(3^{6};1;0;5^{1},4^{3})_{2}$}
\psfrag{T6-3}{$E_8^{(1)}(3^{6};1;0;5^{2})$}
\psfrag{T6-4}{$E_8^{(1)}(3^{6};1;0;6^{1})$}
\includegraphics[width=4.6in]{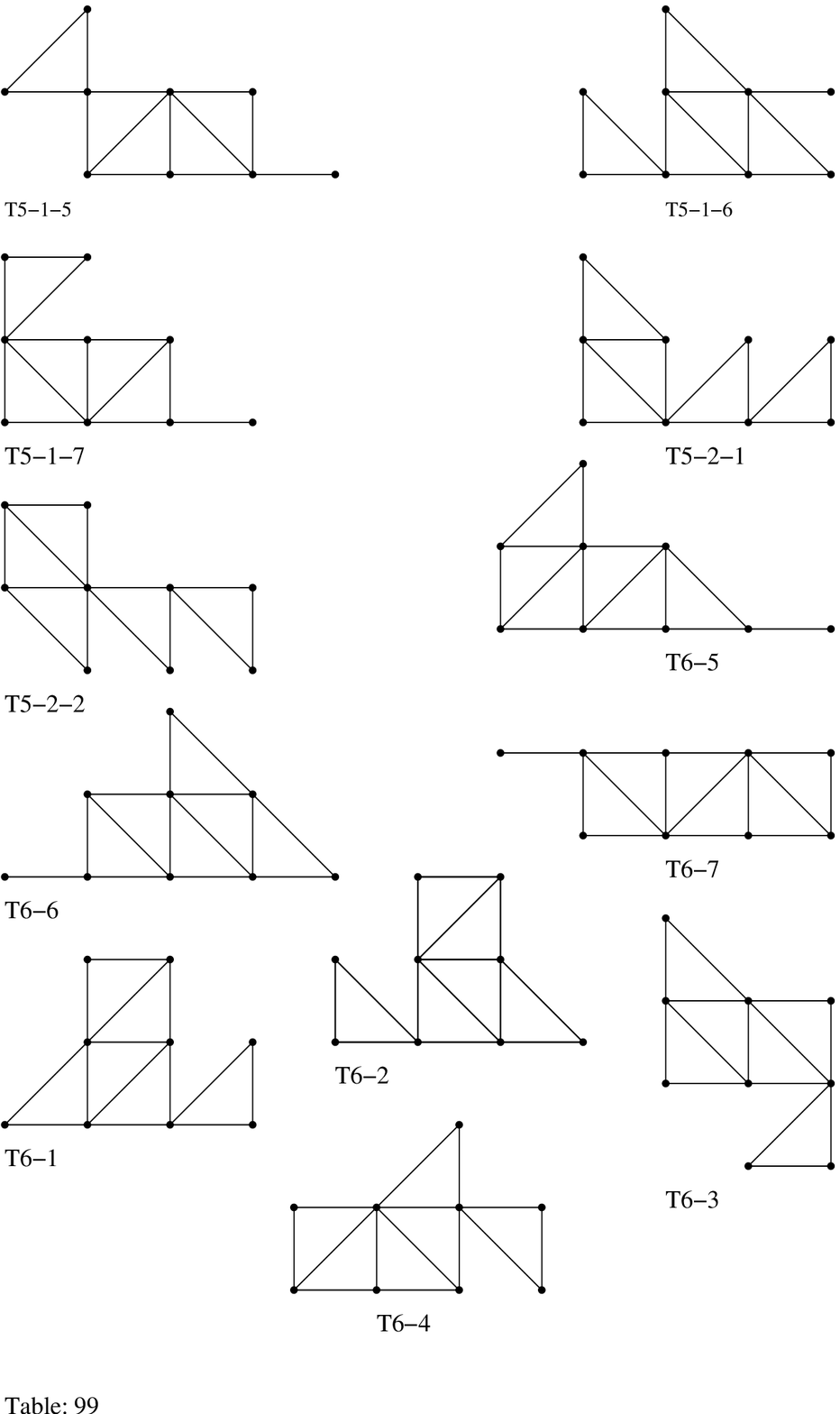}
\end{center}

\pagebreak
\begin{center}
\psfrag{Table: 9-10}{Table 6 (continued)}

\psfrag{T7-1}{$E_8^{(1)}(3^{7})$}
\psfrag{S-1}{$E_8^{(1)}(4^{1};0;4^{1},1^{1})$}
\psfrag{S-2}{$E_8^{(1)}(4^{1};0;3^{1},1^{2})$}
\psfrag{ST-1}{$E_8^{(1)}(4^{1},3^{1};0;2^{1},1^{2})$}
\psfrag{ST-2}{$E_8^{(1)}(4^{1},3^{1};0;2^{2})$}
\psfrag{ST3}{$E_8^{(1)}(4^{1},3^{1};0;3^{1},1^{1};3^{4})$}
\psfrag{ST-4}{$E_8^{(1)}(4^{1},3^{1};0;3^{1},1^{1};4^{1})_{1}$}
\psfrag{ST-5}{$E_8^{(1)}(4^{1},3^{1};0;3^{1},1^{1};4^{1})_{2}$}
\psfrag{ST-6}{$E_8^{(1)}(4^{1},3^{1};0;4^{1};3^{3})$}
\psfrag{ST-7}{$E_8^{(1)}(4^{1},3^{1};0;4^{1};4^{1})$}
\includegraphics[width=4.6in]{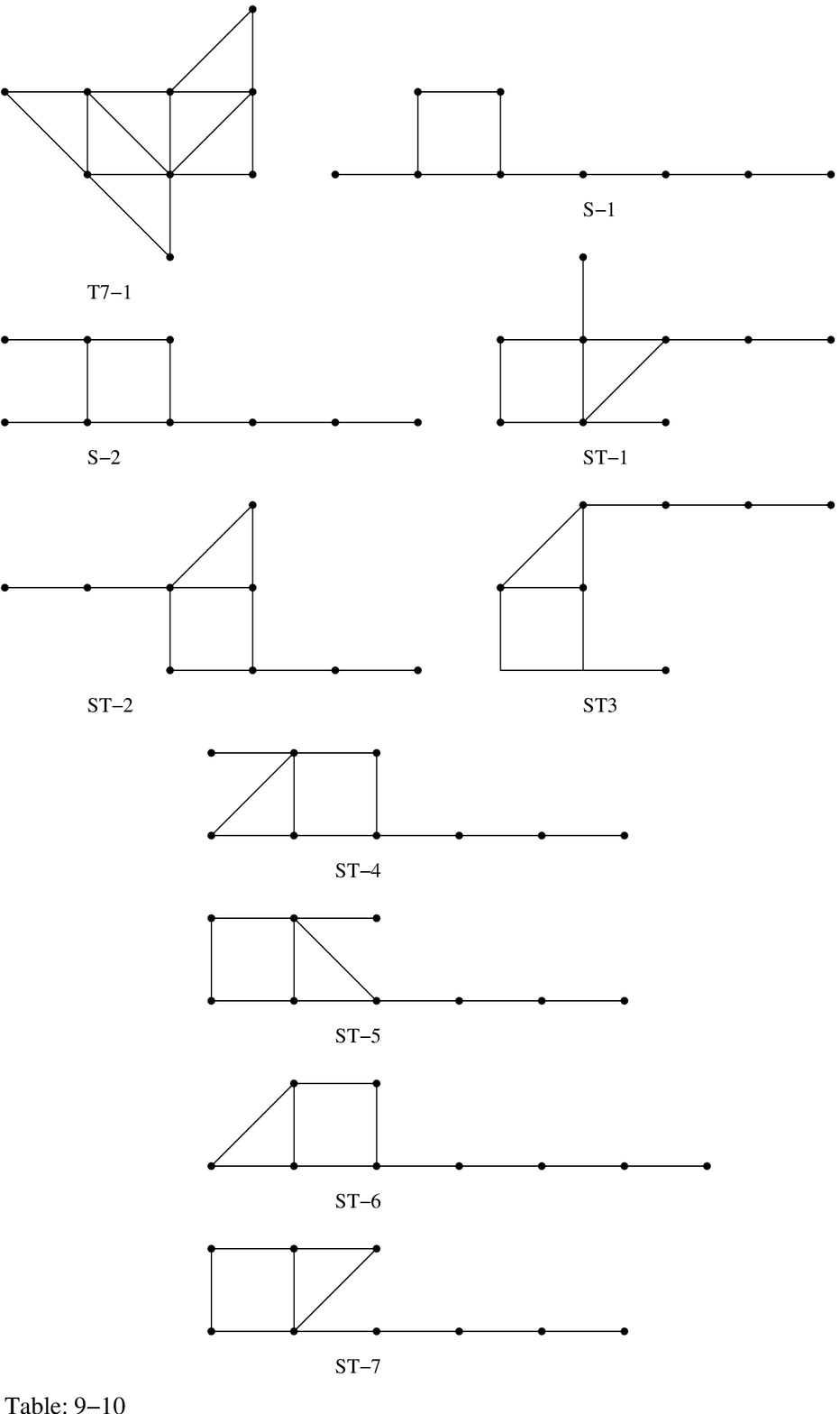}
\end{center}

\pagebreak
\begin{center}
\psfrag{Table: 9-11}{Table 6 (continued)}

\psfrag{ST-8}{$E_8^{(1)}(4^{1},3^{1};1;1^{1})$}
\psfrag{ST-9}{$E_8^{(1)}(4^{1},3^{1};1;1^{2};3^{4})_{1}$}
\psfrag{ST-10}{$E_8^{(1)}(4^{1},3^{1};1;1^{2};3^{4})_{2}$}
\psfrag{ST-11}{$E_8^{(1)}(4^{1},3^{1};1;1^{2};4^{1})$}
\psfrag{ST-12}{$E_8^{(1)}(4^{1},3^{1};1;1^{3})$}
\psfrag{ST-13}{$E_8^{(1)}(4^{1},3^{1};1;2^{1},1^{1};4^{1})$}
\psfrag{ST2-1}{$E_8^{(1)}(4^{1},3^{2};0;1^{3};4^{1})$}
\psfrag{ST2-2}{$E_8^{(1)}(4^{1},3^{2};0;2^{1},1^{1};4^{1},3^{4})_{1}$}
\psfrag{ST2-3}{$E_8^{(1)}(4^{1},3^{2};0;2^{1},1^{1};4^{1},3^{4})_{2}$}
\psfrag{ST2-4}{$E_8^{(1)}(4^{1},3^{2};0;2^{1},1^{1};4^{2})_{1}$}
\psfrag{ST2-5}{$E_8^{(1)}(4^{1},3^{2};0;2^{1},1^{1};4^{2})_{2}$}
\psfrag{ST2-6}{$E_8^{(1)}(4^{1},3^{2};0;2^{1},1^{2})$}
\psfrag{ST2-7}{$E_8^{(1)}(4^{1},3^{2};0;3^{1})$}
\includegraphics[width=4.6in]{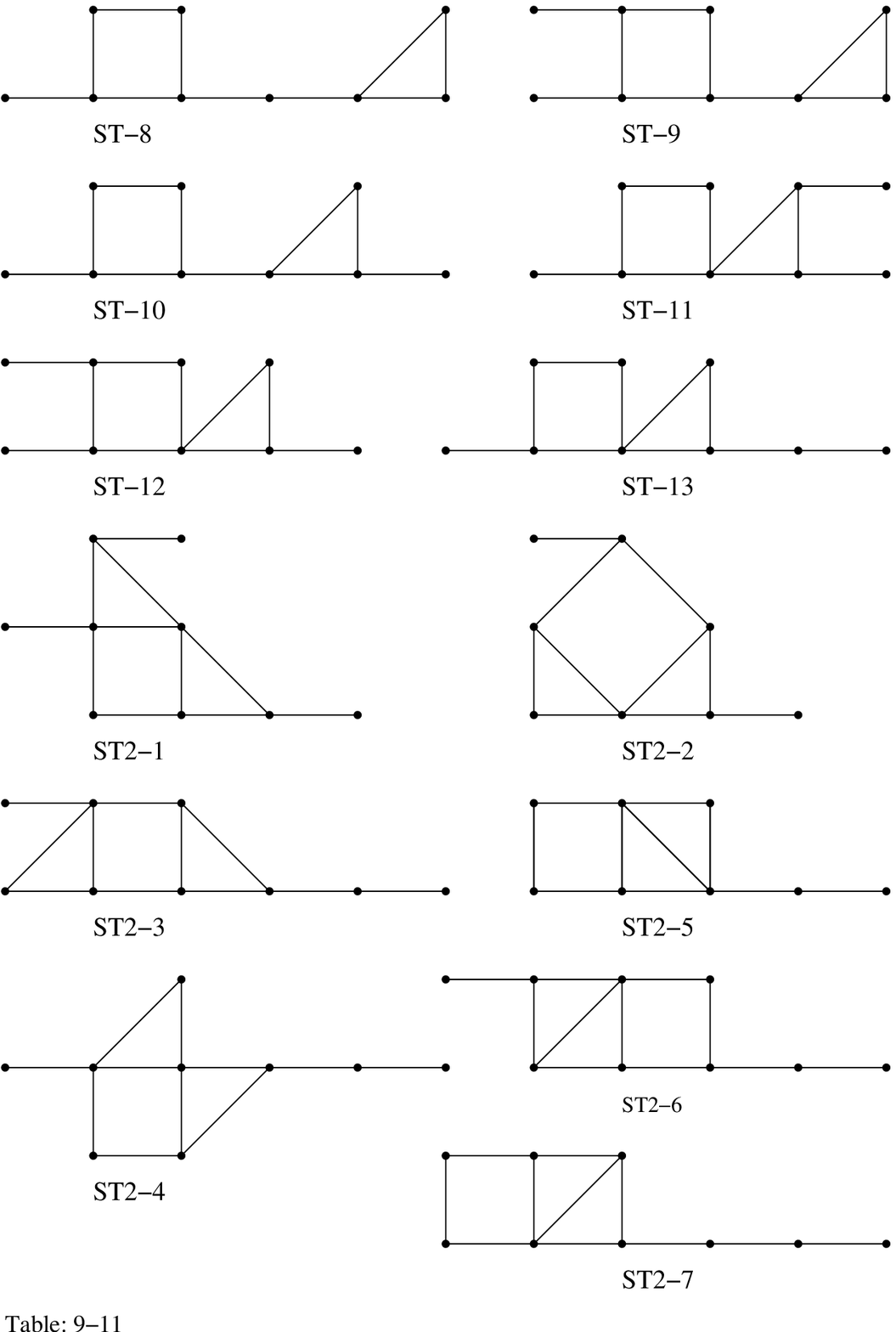}
\end{center}

\pagebreak

\begin{center}
\psfrag{Table:9-12}{Table 6 (continued)}

\psfrag{ST2-8}{$E_8^{(1)}(4^{1},3^{2};1;0;3^{4})$}
\psfrag{ST2-9}{$E_8^{(1)}(4^{1},3^{2};1;0;4^{1})$}
\psfrag{ST2-10}{$E_8^{(1)}(4^{1},3^{2};1;1^{1};3^{5})$}
\psfrag{ST2-11}{$E_8^{(1)}(4^{1},3^{2};1;1^{1};4^{1},3^{3})_{1}$}
\psfrag{ST2-12}{$E_8^{(1)}(4^{1},3^{2};1;1^{1};4^{1},3^{3})_{2}$}
\psfrag{ST2-13}{$E_8^{(1)}(4^{1},3^{2};1;1^{1};4^{1},3^{3})_{3}$}
\psfrag{ST2-14}{$E_8^{(1)}(4^{1},3^{2};1;1^{1};5^{1})$}
\psfrag{ST2-15}{$E_8^{(1)}(4^{1},3^{2};1;1^{2};4^{1},3^{4})_{1}$}
\psfrag{ST2-16}{$E_8^{(1)}(4^{1},3^{2};1;1^{2};4^{1},3^{4})_{2}$}
\psfrag{ST2-17}{$E_8^{(1)}(4^{1},3^{2};1;1^{2};4^{2},3^{2})_{1}$}
\psfrag{ST2-17-1}{$E_8^{(1)}(4^{1},3^{2};1;1^{2};4^{2},3^{2})_{2}$}
\psfrag{ST2-18}{$E_8^{(1)}(4^{1},3^{2};1;1^{2};4^{3})$}
\psfrag{ST2-19}{$E_8^{(1)}(4^{1},3^{2};1;1^{2};5^{1})$}
\psfrag{ST2-20}{$E_8^{(1)}(4^{1},3^{2};1;2^{1};4^{1})$}
\psfrag{ST2-21}{$E_8^{(1)}(4^{1},3^{2};1;2^{1};4^{2})$}
\psfrag{ST2-22}{$E_8^{(1)}(4^{1},3^{2};1;2^{1};5^{1})_{1}$}
\includegraphics[width=4.6in]{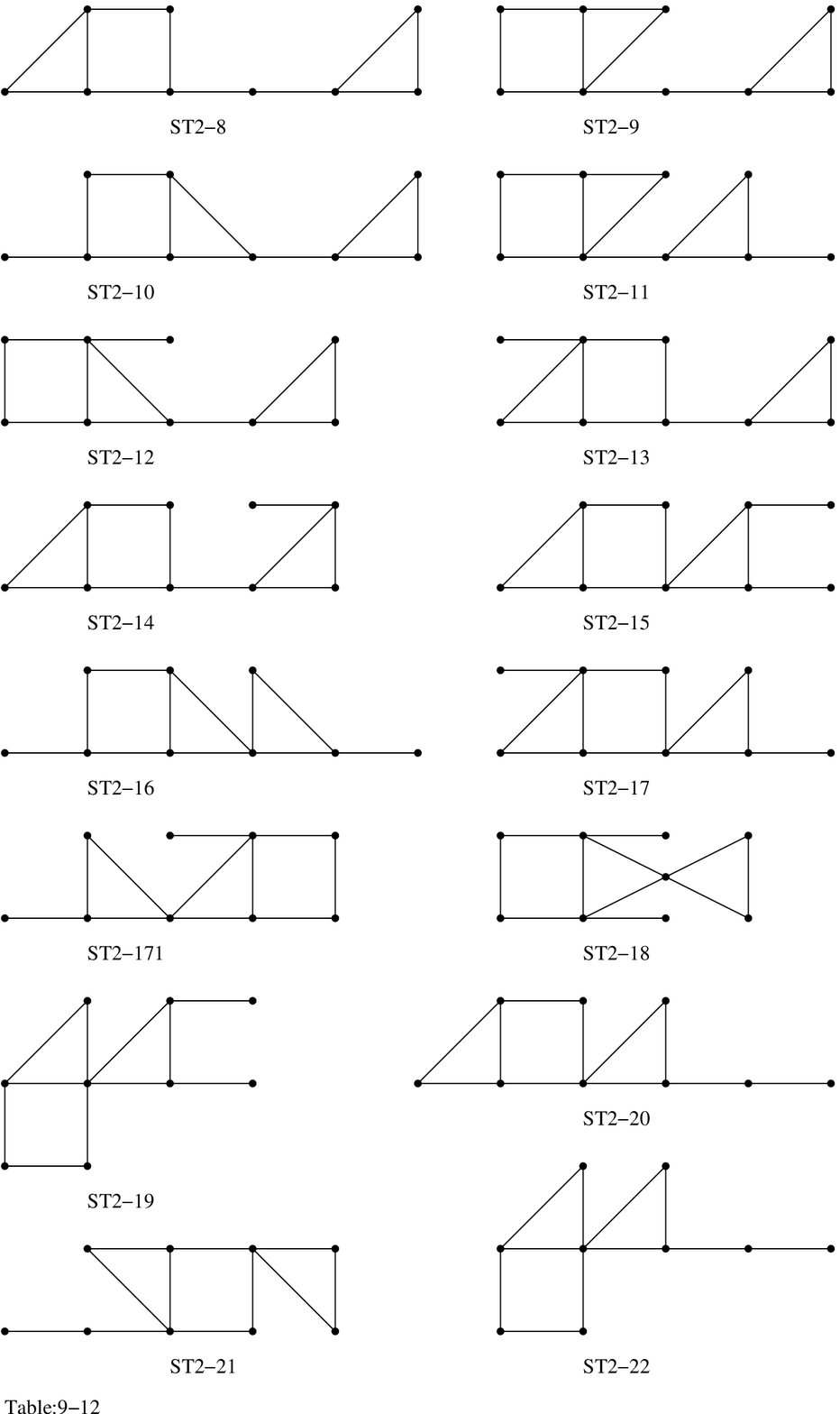}

\end{center}

\pagebreak
\begin{center}
\psfrag{Table:9-13}{Table 6 (continued)}

\psfrag{ST2-23}{$E_8^{(1)}(4^{1},3^{2};1;2^{1};5^{1})_{2}$}
\psfrag{ST2-24}{$E_8^{(1)}(4^{1},3^{2};2;1^{1})$}

\psfrag{ST3-2}{$E_8^{(1)}(4^{1},3^{3};0;1^{2};4^{2},3^{4})_{1}$}
\psfrag{ST3-3}{$E_8^{(1)}(4^{1},3^{3};0;1^{2};4^{2},3^{4})_{2}$}
\psfrag{ST3-4}{$E_8^{(1)}(4^{1},3^{3};0;1^{2};4^{3},3^{2})_{1}$}
\psfrag{ST3-5}{$E_8^{(1)}(4^{1},3^{3};0;1^{2};4^{3},3^{2})_{2}$}
\psfrag{ST3-6}{$E_8^{(1)}(4^{1},3^{3};0;1^{2};5^{1},3^{5})$}
\psfrag{ST3-7}{$E_8^{(1)}(4^{1},3^{3};0;1^{2};5^{1},4^{1},3^{3})_{1}$}
\psfrag{ST3-8}{$E_8^{(1)}(4^{1},3^{3};0;1^{2};5^{1},4^{1},3^{3})_{2}$}
\psfrag{ST3-9}{$E_8^{(1)}(4^{1},3^{3};0;1^{2};5^{1},4^{1},3^{3})_{3}$}
\includegraphics[width=4.6in]{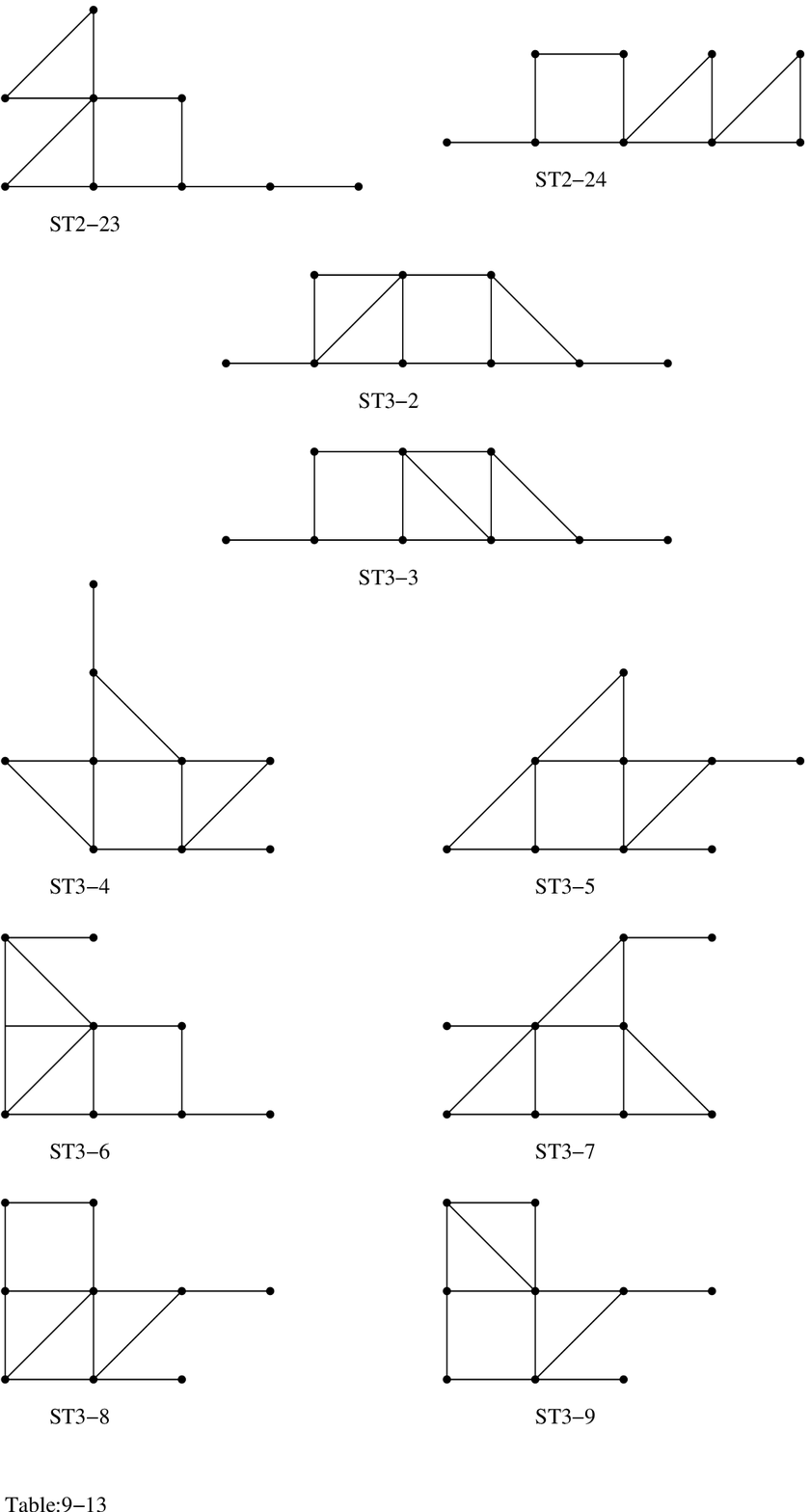}
\end{center}

\pagebreak
\begin{center}
\psfrag{Table:9-14}{Table 6 (continued)}
\psfrag{ST3-1}{$E_8^{(1)}(4^{1},3^{3};0;2^{1};5^{1},3^{4})$}
\psfrag{ST3-91}{$E_8^{(1)}(4^{1},3^{3};0;1^{2};5^{1},4^{1},3^{3})_{4}$}
\psfrag{ST3-92}{$E_8^{(1)}(4^{1},3^{3};0;1^{2};5^{1},4^{1},3^{3})_{5}$}
\psfrag{ST3-10}{$E_8^{(1)}(4^{1},3^{3};0;2^{1};4^{1})$}
\psfrag{ST3-11}{$E_8^{(1)}(4^{1},3^{3};0;2^{1};4^{2})$}
\psfrag{ST3-12}{$E_8^{(1)}(4^{1},3^{3};0;2^{1};4^{3})_{1}$}
\psfrag{ST3-13}{$E_8^{(1)}(4^{1},3^{3};0;2^{1};4^{3})_{2}$}
\psfrag{ST3-14}{$E_8^{(1)}(4^{1},3^{3};0;2^{1};5^{1},3^{3})$}
\psfrag{ST3-15}{$E_8^{(1)}(4^{1},3^{3};0;2^{1};5^{1},4^{1})_{1}$}
\psfrag{ST3-16}{$E_8^{(1)}(4^{1},3^{3};0;2^{1};5^{1},4^{1})_{2}$}
\psfrag{ST3-17}{$E_8^{(1)}(4^{1},3^{3};1;0)$}
\includegraphics[width=4.6in]{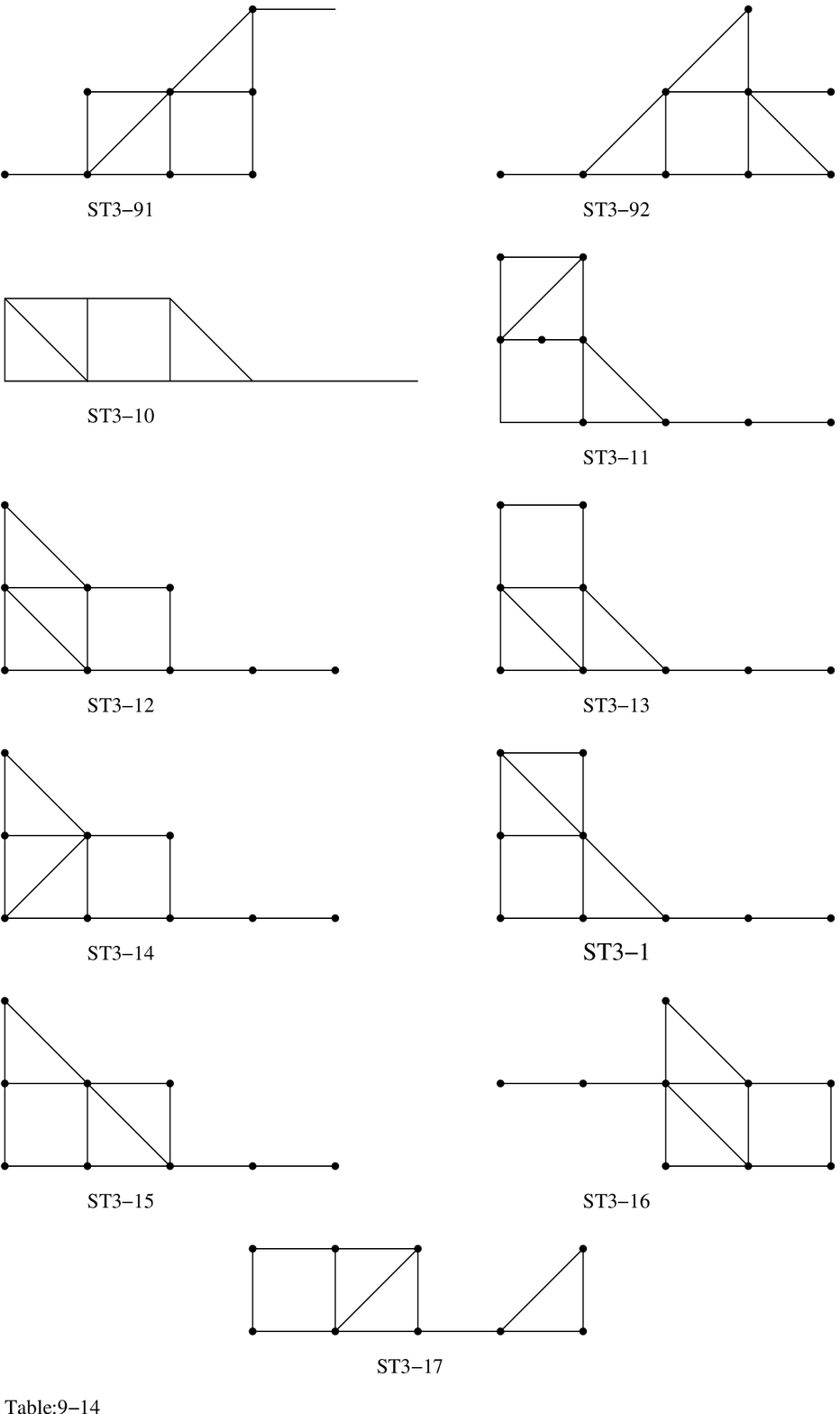}
\end{center}

\pagebreak
\begin{center}
\psfrag{Table:9-15}{Table 6 (continued)}
\psfrag{ST3-18}{$E_8^{(1)}(4^{1},3^{3};1;1^{1};4^{2},3^{2})$}
\psfrag{ST3-19}{$E_8^{(1)}(4^{1},3^{3};1;1^{1};4^{2},3^{3})_{1}$}
\psfrag{ST3-20}{$E_8^{(1)}(4^{1},3^{3};1;1^{1};4^{2},3^{3})_{2}$}
\psfrag{ST3-21}{$E_8^{(1)}(4^{1},3^{3};1;1^{1};4^{2},3^{3})_{3}$}
\psfrag{ST3-22}{$E_8^{(1)}(4^{1},3^{3};1;1^{1};4^{3})$}
\psfrag{ST3-23}{$E_8^{(1)}(4^{1},3^{3};1;1^{1};5^{1})$}
\psfrag{ST3-24}{$E_8^{(1)}(4^{1},3^{3};2;0;4^{2})$}
\psfrag{ST3-25}{$E_8^{(1)}(4^{1},3^{3};2;0;5^{1})_{1}$}
\psfrag{ST3-26}{$E_8^{(1)}(4^{1},3^{3};2;0;5^{1})_{2}$}
\psfrag{ST4-1}{$E_8^{(1)}(4^{1},3^{4};0;1^{1};4^{3},3^{3})_{1}$}
\psfrag{ST4-2}{$E_8^{(1)}(4^{1},3^{4};0;1^{1};4^{3},3^{3})_{2}$}
\psfrag{ST4-3}{$E_8^{(1)}(4^{1},3^{4};0;1^{1};4^{4})$}
\includegraphics[width=4.6in]{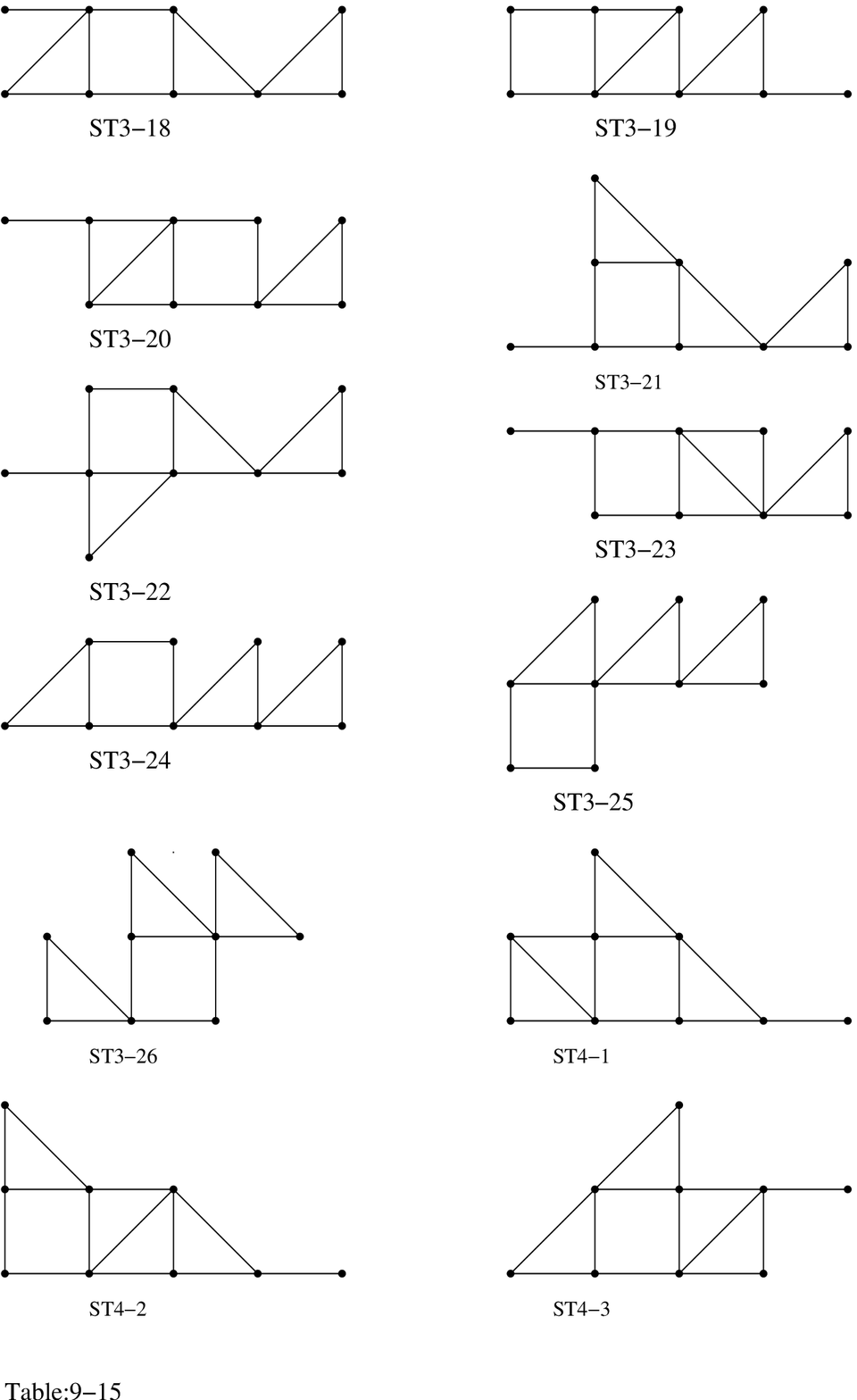}
\end{center}

\pagebreak
\begin{center}
\psfrag{Table:9-16}{Table 6 (continued)}

\psfrag{ST4-4}{$E_8^{(1)}(4^{1},3^{4};0;1^{1};5^{1})$}
\psfrag{ST4-5}{$E_8^{(1)}(4^{1},3^{4};1;0;4^{2})$}
\psfrag{ST4-6}{$E_8^{(1)}(4^{1},3^{4};1;0;4^{3})$}
\psfrag{ST4-7}{$E_8^{(1)}(4^{1},3^{4};1;0;5^{1},4^{1})_{1}$}
\psfrag{ST4-8}{$E_8^{(1)}(4^{1},3^{4};1;0;5^{1},4^{1})_{2}$}
\psfrag{ST4-9}{$E_8^{(1)}(4^{1},3^{4};1;0;4^{4})_{1}$}
\psfrag{ST4-10}{$E_8^{(1)}(4^{1},3^{4};1;0;4^{4})_{2}$}
\psfrag{ST4-11}{$E_8^{(1)}(4^{1},3^{4};1;0;5^{2})$}
\psfrag{ST4-12}{$E_8^{(1)}(4^{1},3^{4};1;0;6^{1})$}
\includegraphics[width=4.6in]{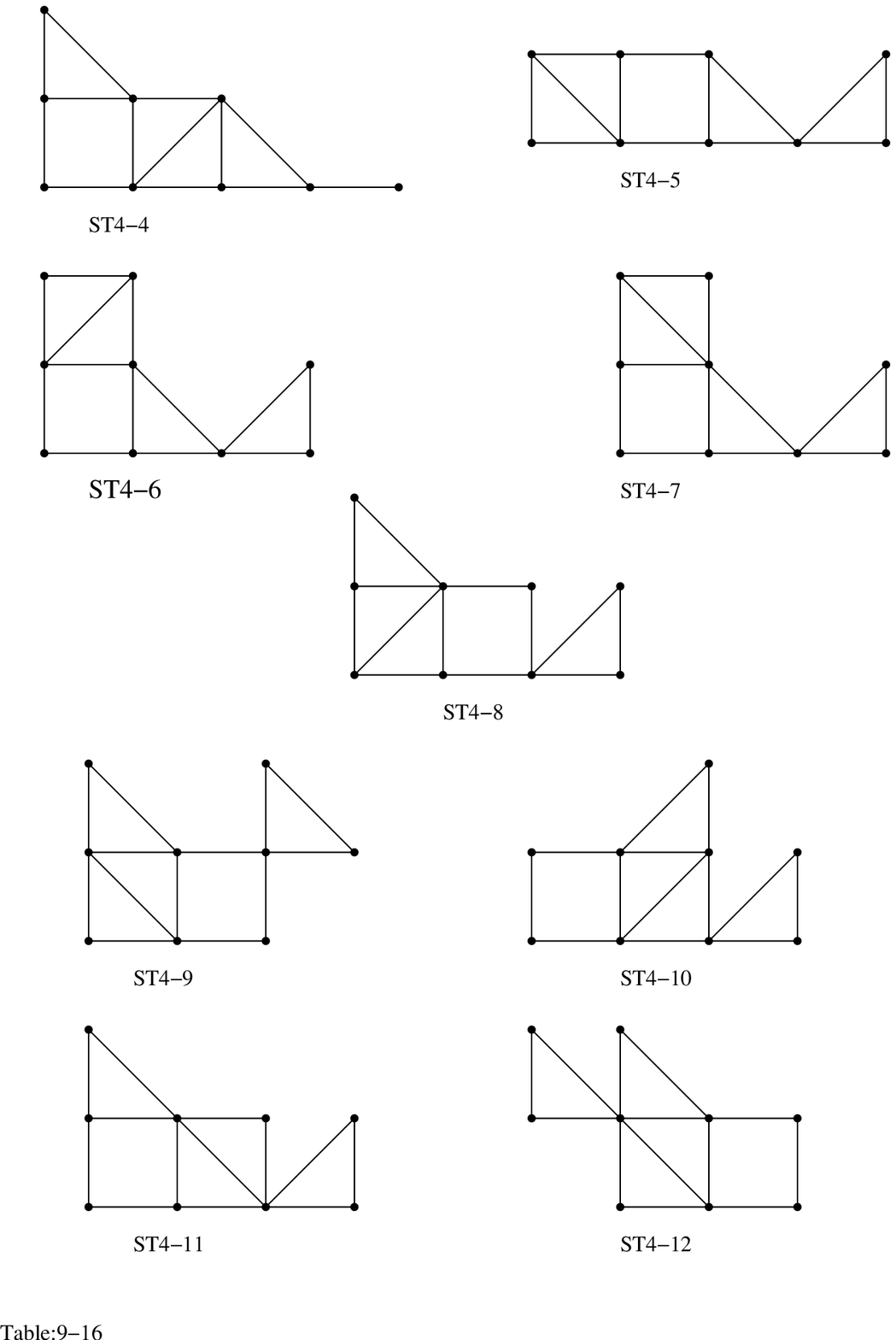}
\end{center}

\pagebreak
\begin{center}
\psfrag{Table: 9-17}{Table 6 (continued)}

\psfrag{ST5-1}{$E_8^{(1)}(4^{1},3^{5};0;0;5^{1},4^{1})$}
\psfrag{ST5-2}{$E_8^{(1)}(4^{1},3^{5};0;0;5^{1},4^{3})_{1}$}
\psfrag{ST5-21}{$E_8^{(1)}(4^{1},3^{5};0;0;5^{1},4^{3})_{2}$}
\psfrag{ST5-3}{$E_8^{(1)}(4^{1},3^{5};0;0;5^{2},3^{4})$}
\psfrag{ST5-4}{$E_8^{(1)}(4^{1},3^{5};0;0;5^{2},4^{1})$}
\psfrag{S2T-0}{$E_8^{(1)}(4^{2};0)$}
\psfrag{S2T-1}{$E_8^{(1)}(4^{1},3^{1};0;2^{1};3^{5})$}
\psfrag{S2T-2}{$E_8^{(1)}(4^{1},3^{1};0;2^{1};4^{1})_{1}$}
\psfrag{S2T-3}{$E_8^{(1)}(4^{1},3^{1};0;2^{1};4^{1})_{2}$}
\psfrag{S2T-4}{$E_8^{(1)}(4^{1},3^{1};0;2^{1};4^{2})$}
\psfrag{S2T-5}{$E_8^{(1)}(4^{1},3^{1};1;1^{1})$}
\includegraphics[width=4.6in]{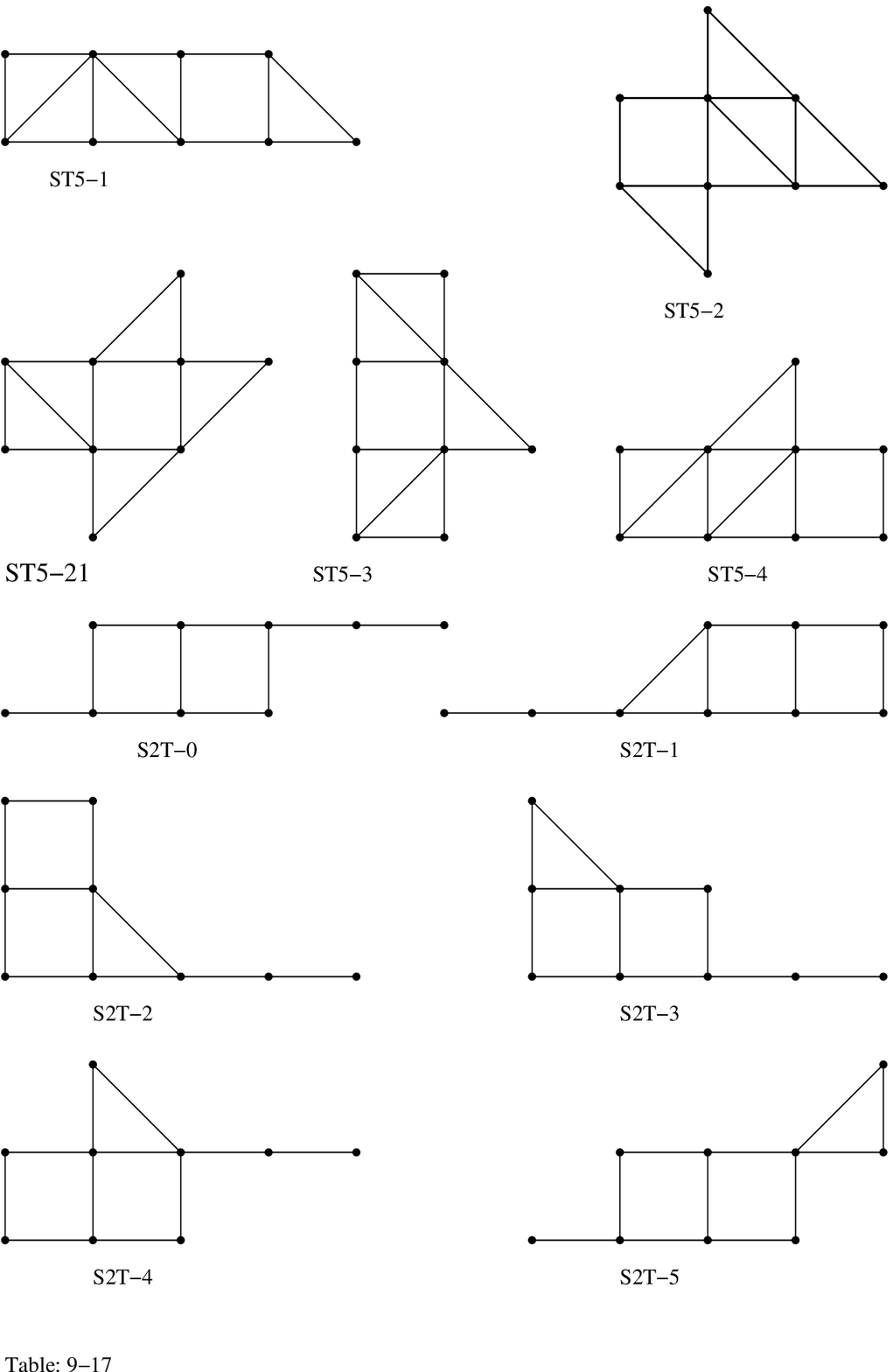}
\end{center}

\pagebreak
\begin{center}
\psfrag{Table: 9-18}{Table 6 (continued)}

\psfrag{S2T2-1}{$E_8^{(1)}(4^{2},3^{2};0;1^{1};4^{1})_{1}$}
\psfrag{S2T2-2}{$E_8^{(1)}(4^{2},3^{2};0;1^{1};4^{1})_{2}$}
\psfrag{S2T2-3}{$E_8^{(1)}(4^{2},3^{2};0;1^{1};4^{2})$}
\psfrag{S2T2-4}{$E_8^{(1)}(4^{2},3^{2};0;1^{1};5^{1})$}
\psfrag{S2T2-5}{$E_8^{(1)}(4^{2},3^{2};1;0;4^{1})$}
\psfrag{S2T2-6}{$E_8^{(1)}(4^{2},3^{2};1;0;4^{2})_{1}$}
\psfrag{S2T2-7}{$E_8^{(1)}(4^{2},3^{2};1;0;4^{2})_{2}$}
\psfrag{S2T2-8}{$E_8^{(1)}(4^{2},3^{2};1;0;5^{1})$}
\psfrag{S2T3-1}{$E_8^{(1)}(4^{2},3^{3};0;0;4^{3})$}
\psfrag{S2T3-2}{$E_8^{(1)}(4^{2},3^{3};0;0;4^{4})$}
\includegraphics[width=4.6in]{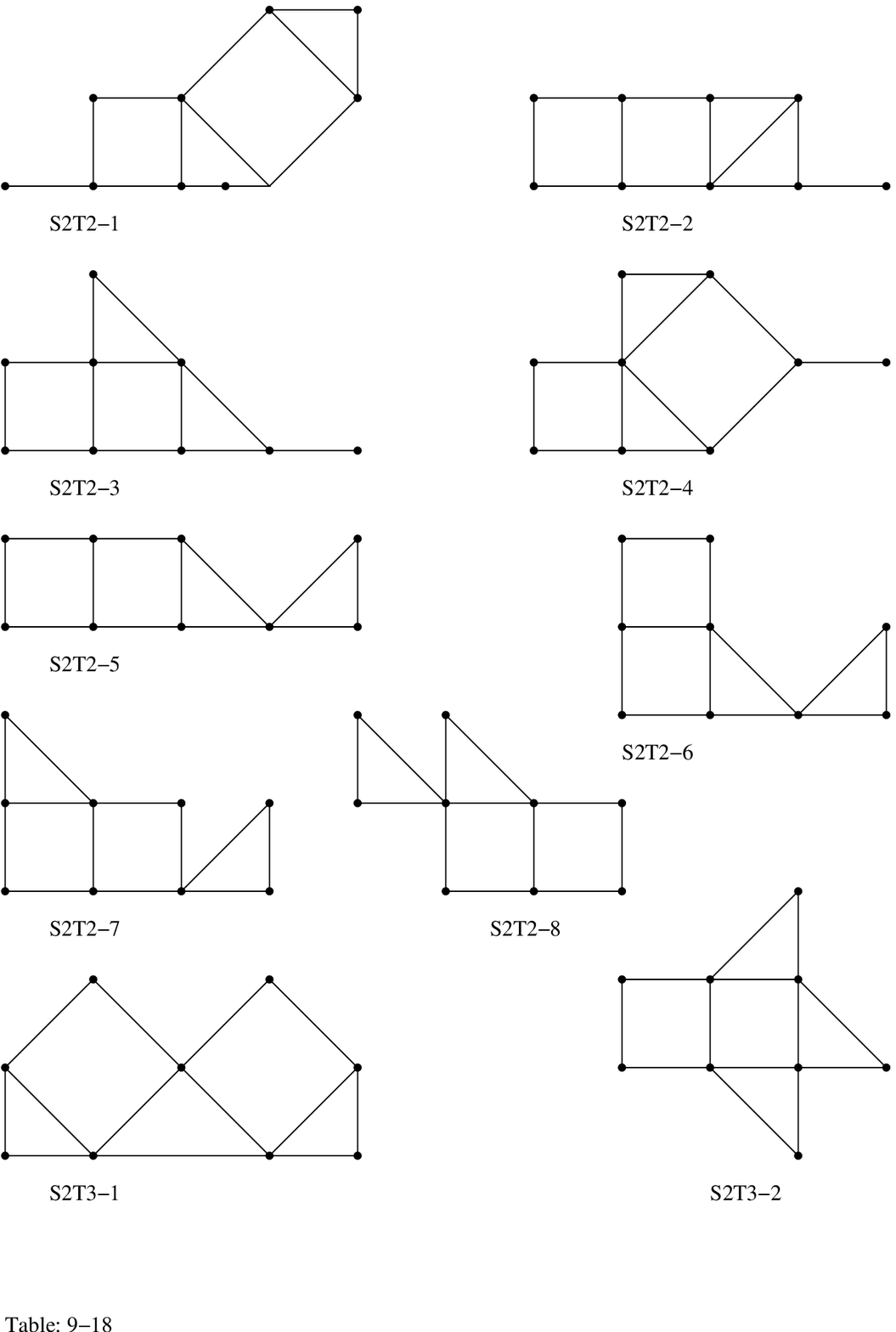}
\end{center}

\pagebreak
\begin{center}
\psfrag{Table: 9-19}{Table 6 (continued)}

\psfrag{S4T1}{$E_8^{(1)}(4^{3},3^{1})$}
\psfrag{P1}{$E_8^{(1)}(5^{1};0;2^{2})$}
\psfrag{P2}{$E_8^{(1)}(5^{1};0;4^{1})$}
\psfrag{P1T1-1}{$E_8^{(1)}(5^{1},3^{1};0;1^{3})$}
\psfrag{P1T1-2}{$E_8^{(1)}(5^{1},3^{1};0;2^{1},1^{1})_{1}$}
\psfrag{P1T1-3}{$E_8^{(1)}(5^{1},3^{1};0;2^{1},1^{1})_{2}$}
\psfrag{P1T1-4}{$E_8^{(1)}(5^{1},3^{1};1;0)$}
\psfrag{P1T1-5}{$E_8^{(1)}(5^{1},3^{1};1;1^{1})$}
\psfrag{P1T1-6}{$E_8^{(1)}(5^{1},3^{1};1;1^{2})$}
\psfrag{P1T1-7}{$E_8^{(1)}(5^{1},3^{1};1;2^{1};4^{1},3^{1})_{1}$}
\psfrag{P1T1-8}{$E_8^{(1)}(5^{1},3^{1};1;2^{1};4^{1},3^{1})_{2}$}
\includegraphics[width=4.6in]{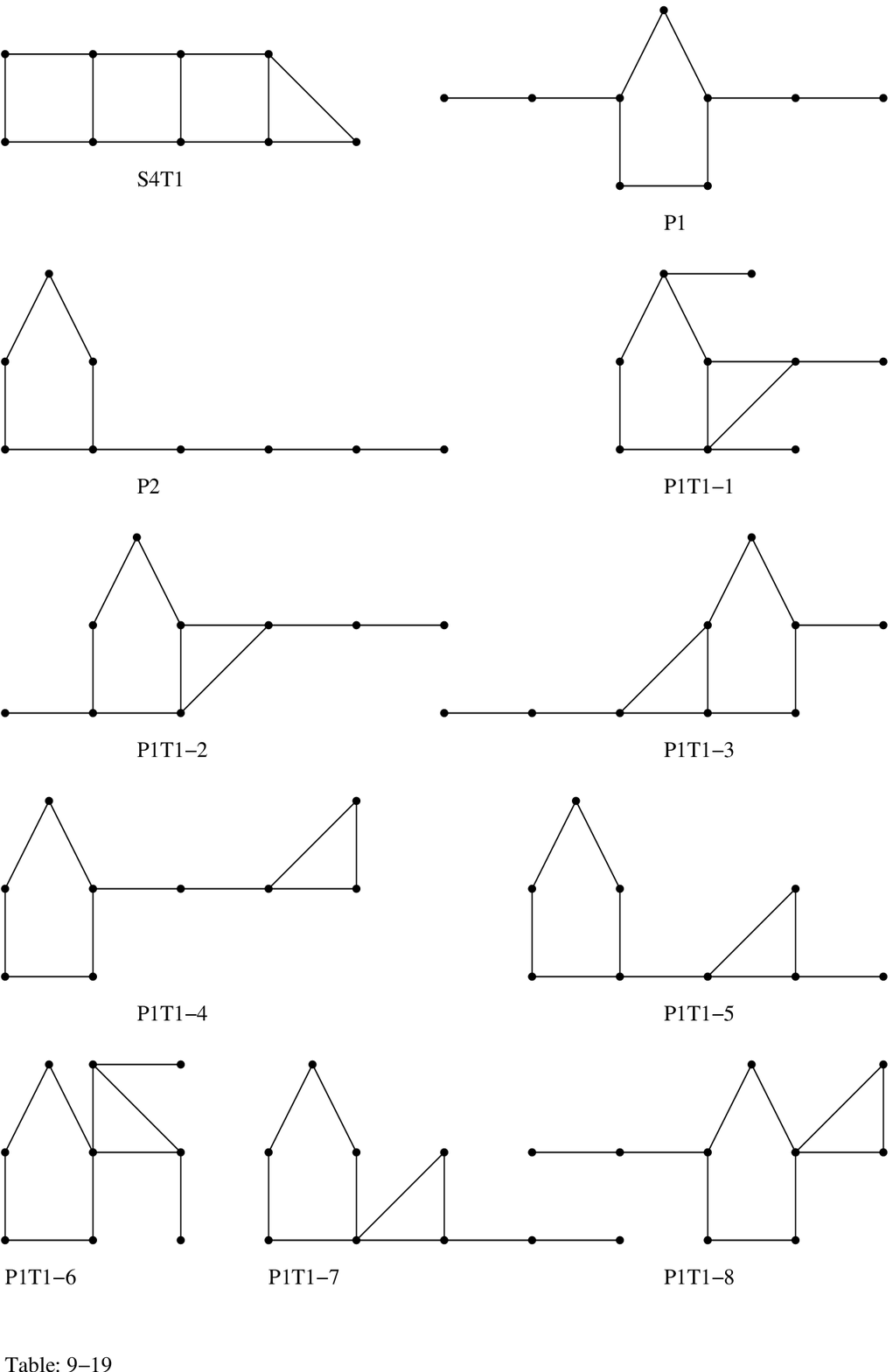}
\end{center}

\pagebreak
\begin{center}
\psfrag{Table:9-192}{Table 6 (continued)}

\psfrag{P1T2-1}{$E_8^{(1)}(5^{1},3^{2};0;1^{2};4^{1},3^{4})_{1}$}
\psfrag{P1T2-2}{$E_8^{(1)}(5^{1},3^{2};0;1^{2};4^{1},3^{4})_{2}$}
\psfrag{P1T2-3}{$E_8^{(1)}(5^{1},3^{2};0;1^{2};4^{2},3^{2})_{1}$}
\psfrag{P1T2-32}{$E_8^{(1)}(5^{1},3^{2};0;1^{2};4^{2},3^{2})_{2}$}
\psfrag{P1T2-4}{$E_8^{(1)}(5^{1},3^{2};1;1^{1};4^{1},3^{3})_{1}$}
\psfrag{P1T2-5}{$E_8^{(1)}(5^{1},3^{2};1;1^{1};4^{1},3^{3})_{2}$}
\psfrag{P1T2-6}{$E_8^{(1)}(5^{1},3^{2};2;0;4^{2})_{1}$}
\psfrag{P1T2-7}{$E_8^{(1)}(5^{1},3^{2};2;0;4^{2})_{2}$}

\includegraphics[width=4.6in]{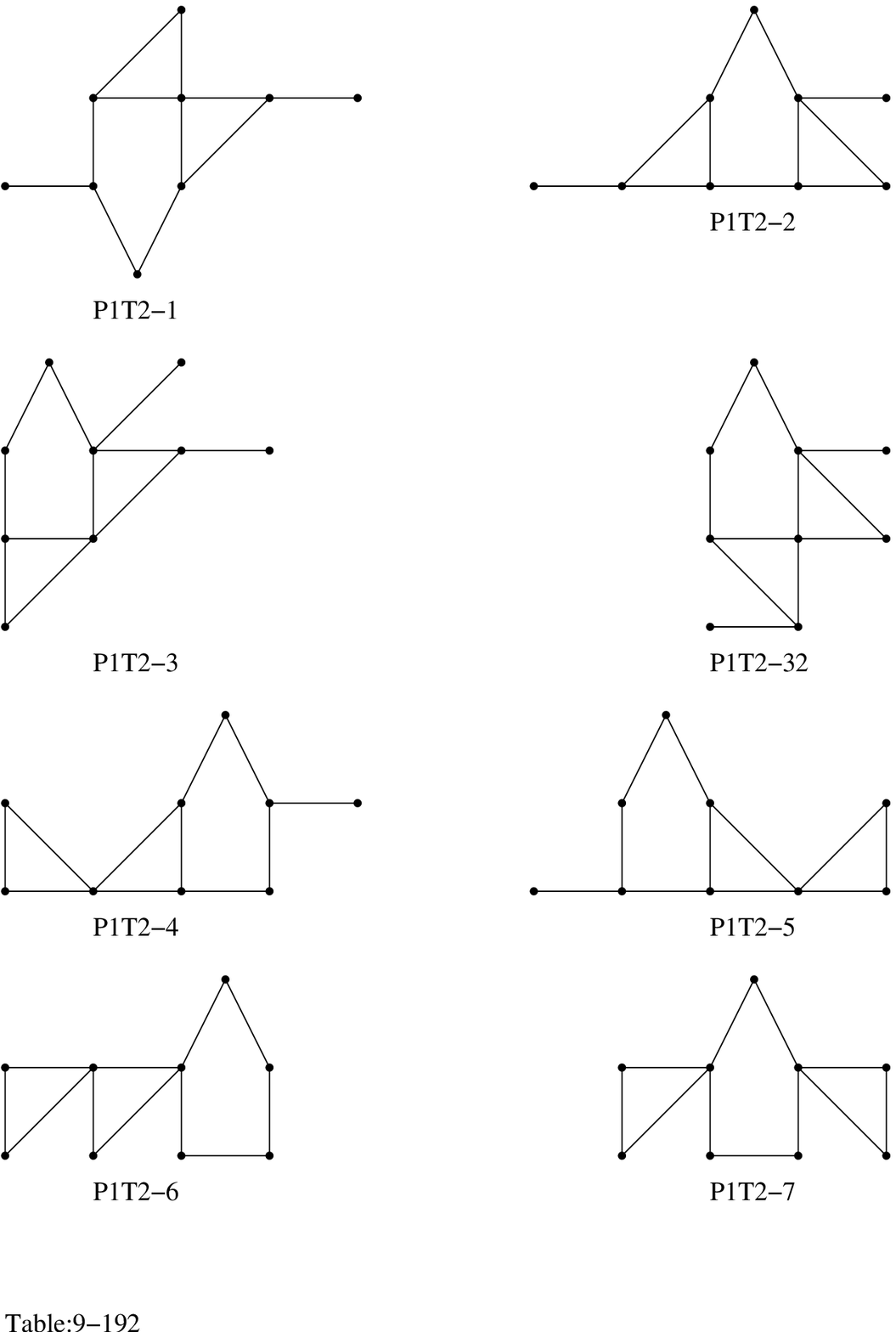}
\end{center}

\pagebreak
\begin{center}
\psfrag{Table: 9-20}{Table 6 (continued)}

\psfrag{P1T3-1}{$E_8^{(1)}(5^{1},3^{3};0;0)$}
\psfrag{P1T3-2}{$E_8^{(1)}(5^{1},3^{3};0;1^{1};4^{1})_{1}$}
\psfrag{P1T3-3}{$E_8^{(1)}(5^{1},3^{3};0;1^{1};4^{1})_{2}$}
\psfrag{P1T3-4}{$E_8^{(1)}(5^{1},3^{3};0;1^{1};4^{2})$}
\psfrag{P1T3-5}{$E_8^{(1)}(5^{1},3^{3};0;1^{1};5^{1})_{1}$}
\psfrag{P1T3-6}{$E_8^{(1)}(5^{1},3^{3};0;1^{1};5^{1})_{2}$}
\psfrag{P1T3-7}{$E_8^{(1)}(5^{1},3^{3};0;1^{1};5^{1})_{3}$}
\psfrag{P1T3-8}{$E_8^{(1)}(5^{1},3^{3};0;1^{1};5^{1},4^{1})$}
\includegraphics[width=4.6in]{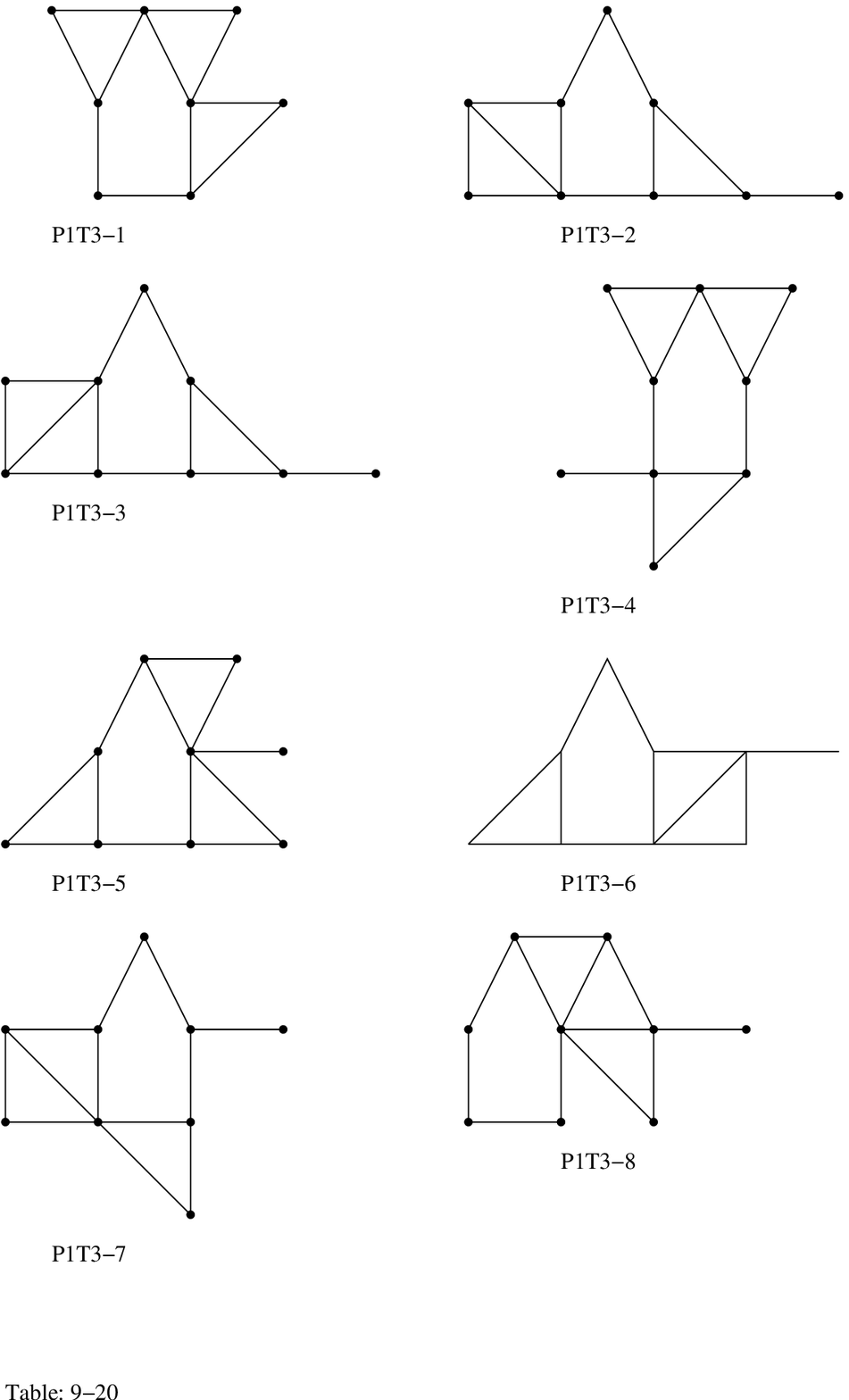}

\end{center}

\pagebreak
\begin{center}
\psfrag{Table:9-21}{Table 6 (continued)}

\psfrag{P1T4-1}{$E_8^{(1)}(5^{1},3^{4};0;0;4^{3})$}
\psfrag{P1T4-2}{$E_8^{(1)}(5^{1},3^{4};0;0;5^{1})$}
\psfrag{P1S1-1}{$E_8^{(1)}(5^{1},4^{1};0;1^{2})$}
\psfrag{P1S1-2}{$E_8^{(1)}(5^{1},4^{1};0;2^{1})$}
\psfrag{P1S1T1-1}{$E_8^{(1)}(5^{1},4^{1},3^{1};0;1^{1})$}
\psfrag{P1S1T1-2}{$E_8^{(1)}(5^{1},4^{1},3^{1};1;0)$}
\psfrag{P1S1T2-1}{$E_8^{(1)}(5^{1},4^{1},3^{2};0)$}
\includegraphics[width=4.6in]{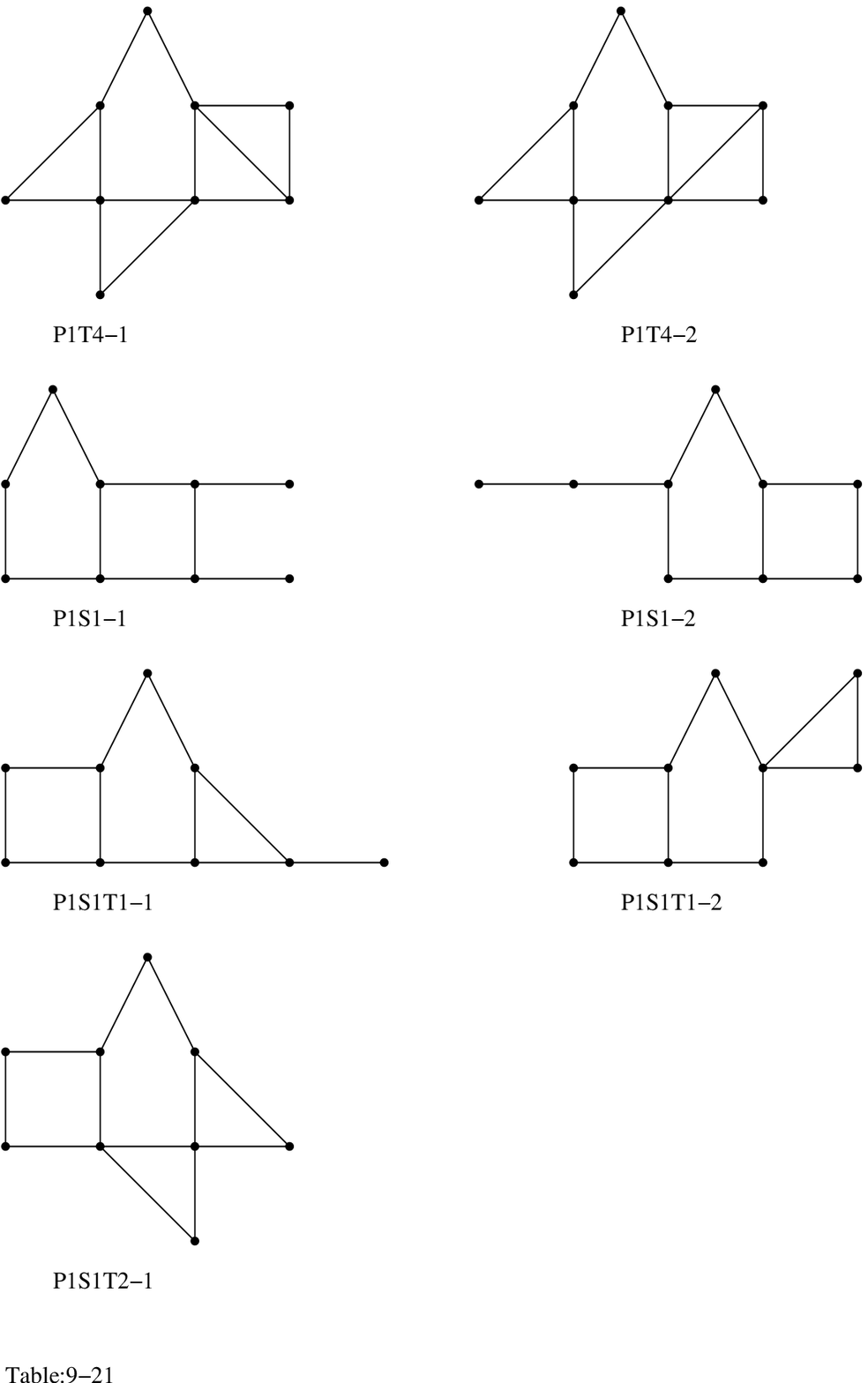}
\end{center}

\pagebreak
\begin{center}
\psfrag{Table:9-22}{Table 6 (continued)}

\psfrag{1}{$E_8^{(1)}(6^{1},3^{1};0;1^{2})_1$}
\psfrag{2}{$E_8^{(1)}(6^{1},3^{1};0;1^{2})_2$}
\psfrag{3}{$E_8^{(1)}(6^{1},3^{2};0;1^{1})_{1}$}
\psfrag{4}{$E_8^{(1)}(6^{1},3^{2};0;1^{1})_{2}$}
\psfrag{5}{$E_8^{(1)}(6^{1},3^{3};0)$}
\psfrag{6}{$E_8^{(1)}(6^{1},4^{1},3^{1};0;3^{4})$}
\psfrag{7}{$E_8^{(1)}(6^{1},4^{1},3^{1};0;4^{2})$}
\psfrag{8}{$E_8^{(1)}(7^{1},3^{1};0)$}
\includegraphics[width=4.6in]{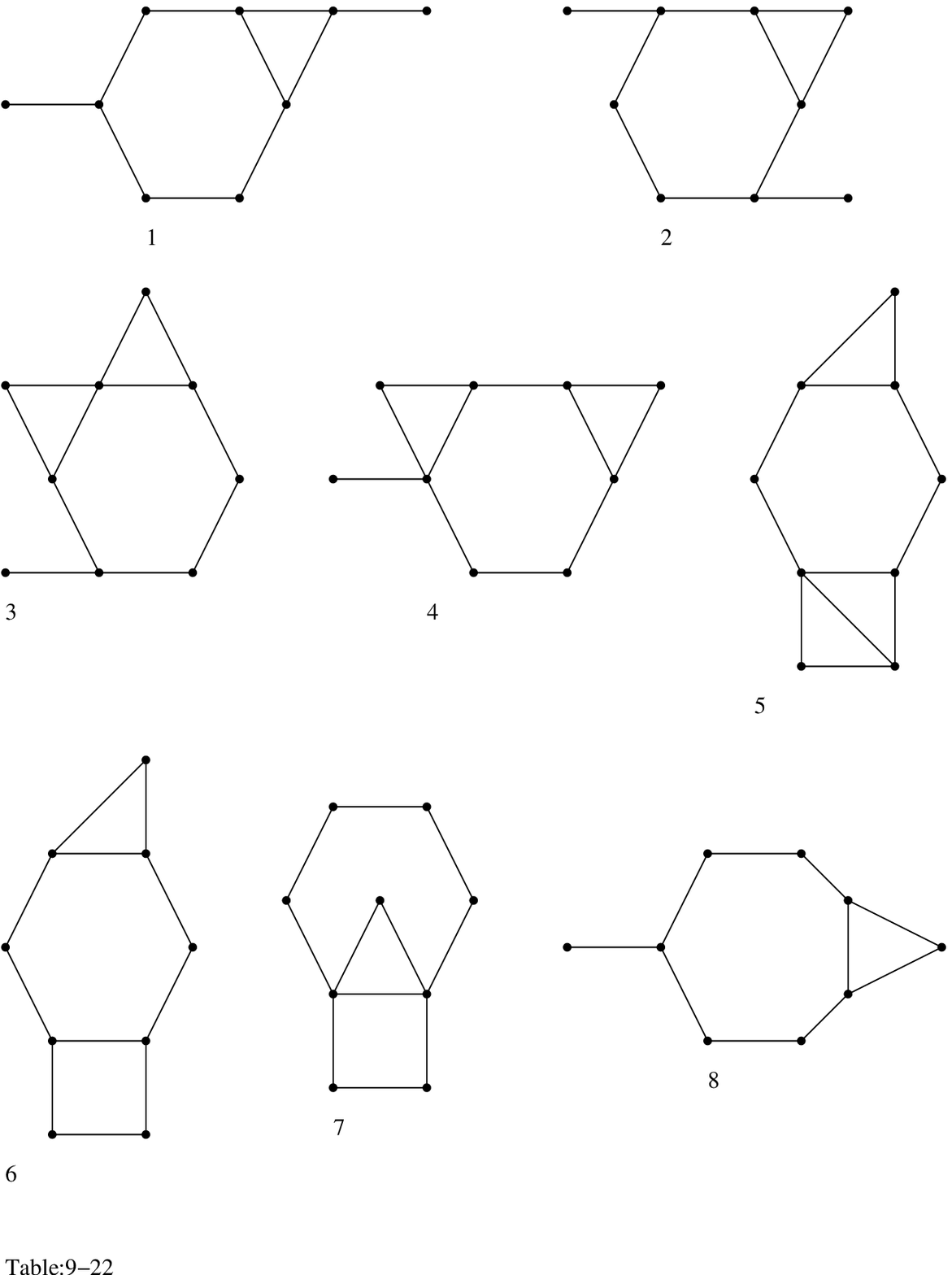}
\end{center}

\pagebreak

\section{Acknowledgements}
\label{sec:ack}

I am grateful to my graduate advisor A. Zelevinsky for his support and suggestions. In particular, it was he who suggested Problem~\ref{pr:pb}.

\end{document}